\documentclass[12pt]{article}
\usepackage[cp1251]{inputenc} %для WIN кодировки (под квантором символ)
\usepackage[english,russian]{babel}
\usepackage{amsfonts,amssymb}
\usepackage{latexsym}
\usepackage{srcltx}
\usepackage{cite}
\usepackage{bbding}
\usepackage{amsmath}
\usepackage{graphicx}

\newtheorem{theorem}{Теорема}

\newtheorem{proposition}{Предложение}
\newtheorem{definition}{Определение}

\newtheorem{zam}{Замечание}

\newtheorem{corollary}{Следствие}

\newtheorem{condition}{Условие}

\newcommand{\beq}{\begin{equation}\label}
\newcommand{\eeq}{\end{equation}}
\newcommand{\pp}{\mathcal{P}}
\newcommand{\pr}{\mathrm}
\newcommand{\ml}{\mathcal{L}}
\newcommand{\rr}{\mathbb{R}}

\newcommand{\tr}{\tau_{\mathbb{R}}}

\newcommand{\bem}{\begin{multline}\label}
\newcommand{\eem}{\end{multline}}
\newcommand{\zer}{\varnothing}
\newcommand{\bn}{\in]0,\infty[}

\begin{document}
\setcounter{page}{1}
 \thispagestyle{empty}

{\Large
\noindent А.\,G.~Chentsov

\vspace{1cm}

\noindent  Ju.\,V.~Shapar\\%[1cm]
}

\vspace{5cm}

  \begin{center}
 {\Large
\textbf{FINITELY ADDITIVE MEASURES IN}

\vspace{1cm}

\textbf{CONSTRUCTIONS OF EXTENSION FOR ABSTRACT}

\vspace{1cm}

\textbf{
 ATTAINABILITY PROBLEMS}}
 \end{center}

\newpage

\vspace{15cm}

{\Large
\noindent А.\,Г.~Ченцов

\vspace{1cm}

\noindent  Ю.\,В.~Шапарь\\%[1cm]
}

\vspace{5cm}

  \begin{center}
 {\Large
\textbf{КОНЕЧНО-АДДИТИВНЫЕ МЕРЫ В}

\vspace{1cm}

\textbf{КОНСТРУКЦИЯХ РАСШИРЕНИЙ АБСТРАКТНЫХ}

\vspace{1cm}

\textbf{ЗАДАЧ О ДОСТИЖИМОСТИ}}
 \end{center}

\newpage
\hspace{3 mm} \begin{center} {\LARGE {\bf Оглавление}}\end{center}
\vspace{0.3cm}

\large

 \noindent 1. \ Введение  \dotfill\hspace{2 mm} 5
\vspace{0.3cm}

\begin{center} \textbf{ Глава 1. \ \  Общие сведения} \end{center}
\vspace{0.3cm}

  \noindent 2. \ Элементы теории множеств: простейшие конструкции
    \dotfill \hspace{2 mm}11
\vspace{0.3cm}

 \noindent  3. \ Элементы топологии, 1 \dotfill \hspace{2 mm} 27
  \vspace{0.3cm}

 \noindent 4. \ Элементы топологии, 2  \dotfill\hspace{2 mm} 39
\vspace{0.3cm}

\noindent 5. \ Отношение  вписанности на пространстве разбиений  \dotfill\hspace{2 mm} 44
\vspace{0.3cm}

\noindent 6. \  Измеримые пространства с полуалгебрами и алгебрами множеств  \dotfill\hspace{2 mm} 54
\vspace{0.3cm}

\begin{center} \textbf{ Глава 2. \ \ Элементы конечно-аддитивной теории меры}\end{center}
\vspace{0.3cm}

 \noindent 7. \ Конечно-аддитивные меры: общие свойства  \dotfill\hspace{2 mm} 57
\vspace{0.3cm}

\noindent 8. \  Ступенчатые, ярусные и измеримые функции   \dotfill\hspace{2 mm} 65
\vspace{0.3cm}

\noindent 9. \ Интегрирование ярусных функций  по конечно-аддитивной мере \\ ограниченной вариации \dotfill\hspace{2 mm} 71
\vspace{0.3cm}

\noindent 10.  Интегральное представление линейных непрерывных функционалов; \\ оснащение $*$-слабой топологией   \dotfill\hspace{1 mm} 79
\vspace{0.3cm}

\noindent 11. \  Слабо абсолютно непрерывные конечно-аддитивные меры и их \\ приближение неопределенными интегралами, 1 \dotfill\hspace{2 mm} 85
\vspace{0.3cm}

\noindent 12. \  Слабо абсолютно непрерывные конечно-аддитивные меры и их \\ приближение неопределенными интегралами, 2 \dotfill\hspace{2 mm} 94
\vspace{0.3cm}

\begin{center}\textbf{ Глава 3. \ \ Абстрактная задача о достижимости с ограничениями асимптотического характера}\end{center}
\vspace{0.3cm}

 \noindent 13. \ Вариант задачи о достижимости при ограничениях моментного \\ характера \dotfill\hspace{2 mm} 110
\vspace{0.3cm}

\noindent 14. \ Абстрактная задача о достижимости и ее <<конечно-аддитивное>> \\ расширение   \dotfill\hspace{2 mm} 114
\vspace{0.3cm}

\noindent 15.  Множества притяжения и условия асимптотической  нечувствительности \\ при ослаблении части ограничений \dotfill\hspace{2 mm} 117
\vspace{0.3cm}

\begin{center} \textbf{ Глава 4. \ \ Линейные задачи управления с ограничениями импульсного характера}\end{center}
\vspace{0.3cm}

\noindent 16. \ Линейные управляемые системы с разрывностью в коэффициентах \\ при управляющих воздействиях (случай ограничений импульсного \\  характера)  \dotfill\hspace{2 mm} 134
\vspace{0.3cm}

\noindent 17. \ Управление материальной точкой  при  моментных  ограничениях \dotfill\hspace{2 mm} 139
\vspace{0.3cm}

\noindent 18. \ Одно достаточное  условие асимптотической нечувствительности при \\ ослаблении части ограничений   \dotfill\hspace{2 mm} 145
\vspace{0.3cm}

\noindent 19. \ Пример задачи с <<асимптотическими>> ограничениями    \dotfill\hspace{2 mm} 147
\vspace{0.3cm}

\noindent 20. \ Заключение   \dotfill\hspace{2 mm} 154
\vspace{0.5cm}

\noindent Список литературы \dotfill \hspace{2 mm} 156

\newpage
%\Large

\section{Введение}\setcounter{equation}{0} \setcounter{proposition}{0}
\ \ \ \ \ Предлагаемое  пособие посвящено в своем конечном итоге изучению задач, возникающих в теории управления и связанных с достижимостью тех или иных  состояний при наличии ограничений различного характера. Наиболее понятной в этом отношении является задача о построении области достижимости управляемой системы в заданный момент времени, имеющая серьезное практическое значение (данная  задача играет, в частности, важную роль в вопросах космической навигации) и, вместе с тем, представляющая интерес для теоретических исследований.

Однако нередко ограничения, используемые в постановке упомянутой задачи, не могут быть заданы точно; в частности, возможно <<малое>> ослабление этих ограничений, что объективно способствует улучшению достигаемого результата. При этом в ряде случаев данное улучшение имеет скачкообразный характер и, в частности, <<малым>> уже не является. В таком случае (неустойчивой задачи) практически интересными оказываются решения, достигаемые <<на грани  фола>> в смысле соблюдения ограничений: допускаются <<очень малые>> нарушения имеющихся ограничений.

По целому ряду причин во многих задачах такого типа не удается, однако, указать конкретную степень ослабления имеющихся условий. В этих случаях логичным представляется асмптотический подход, который на идейном уровне можно связать с приближенными решениями  \cite[гл.\,III]{44}. Данный подход получил широкое распространение в задачах оптимального управления; при незначительных в идейном отношении модификациях он может быть распространен и на задачи о достижимости, в которых уже не требуется оптимизировать какой-либо критерий качества.

Следует отметить, что подобные конструкции активно использовались в Свердловске-Екатеринбурге в работах Н.\,Н.~Красовского, А.\,И.~Субботина и их учеников. Особо отметим применение данных конструкций (по смыслу~--- расширений) в игровых задачах, включая задачи программного управления и позиционные дифференциальные игры. В теории дифференциальных игр исключительно важную роль сыграла основополагающая теорема об альтернативе Н.\,Н.~Красовского и А.\,И.~Субботина (см. \cite{20}), при доказательстве которой использовались, в частности, т.н. стабильные мосты, играющие роль фазовых ограничений, которые следовало соблюдать вплоть до встречи с целевым множеством; упомянутое соблюдение ограничений в   классе реализуемых пошаговых движений допускалось при этом всего лишь приближенным. В то же время, применяя в определении стабильности (в том или ином виде) обобщенные элементы, удается добиться (в сочетании с правилом экстремального сдвига) осуществление перемещений вблизи стабильных мостов, что оказывается достаточным для приемлемого решения соответствующей игровой задачи управления на аппроксимативном уровне. Используя упомянутую теорему об альтернативе, Н.\,Н.~Красовский и А.\,И.~Субботин установили затем важное положение о разрешимости нелинейной дифференциальной игры (с функцией платы) в смысле седловой точки, реализуемой в надлежащем классе позиционных стратегий.

Обсуждаемые выше  исследования относятся по большей части к задачам управления с геометрическими ограничениями на выбор управляющих воздействий, систематическое исследование которых было начато Л.\,С.~Понтрягиным  (см. \cite{21}). В связи с изучением вопросов, связанных с существованием оптимальных (обобщенных) управлений-мер, отметим работы Дж.~Варги \cite{44} и Р.\,В. Гамкрелидзе \cite{22}.

Более сложно (на наш взгляд) складывались аналогичные исследования в задачах управления с импульсными ограничениями. В этой связи особо отметим оригинальный подход, предложенный Н.\,Н.~Красовским и связанный с использованием аппарата обобщенных функций (см. \cite{23}). Этот подход послужил основой многих исследований в области импульсного управления; отметим, в частности, работы С.\,Т.~Завалищина и А.\,Н.~Сесекина (см. \cite{24}). Заметим, что импульсные ограничения играют важную роль в приложениях, поскольку именно посредством данных ограничений задаются условия на энергоресурс, т.е. на запас имеющегося топлива.

Одно из затруднений, возникающих в импульсных задачах, связано с эффектами, имеющими смысл произведения разрывной функции на обобщенную. В этой связи отметим направление, касающееся изучения линейных систем с разрывностью при управляющих воздействиях, которое связано с применением особых обобщенных элементов (обобщенных управлений)~--- скалярных или векторных конечно-аддитивных мер (см. \cite{25,46,47} и др.). Дело в том, что такие меры определяют представление линейных непрерывных функционалов на пространстве разрывных (точнее, ярусных) функций. При этом действие обобщенных элементов естественно связывать с упомянутыми функционалами, имея в виду <<хорошие>> условия компактности, определяемые теоремой Алаоглу (см. \cite{43}). Важно, конечно, позаботиться о том, чтобы создаваемые таким образом обобщенные элементы допускали приближение обычными управлениями в смысле соответствующей топологии, которая в данном случае является <<сужением>> $*$-слабой топологии пространства конечно-аддитивных мер ограниченной вариации. На этой основе в классе линейных систем удается формализовать эффекты, имеющие смысл произведения  разрывной функции на обобщенную.

В настоящем пособии предпринимается попытка дать систематическое и достаточно подробное изложение упомянутого подхода для одного типа ограничений импульсного характера: мы обсуждаем случай неотрицательных \\ управлений и требуем полного расходования имеющегося энергоресурса. Помимо этого предполагаем, что имеются и некоторые другие ограничения (моментного характера), которые должны соблюдаться с высокой, но все же конечной степенью точности. Точная формализация связывается при этом с введением асимптотических режимов, подобных приближенным решениям Дж. Варги (см. \cite[гл.\,III]{44}). Соответствующие области достижимости при данном подходе заменяются множествами притяжения, которые как раз и следует рассматривать как практически интересные варианты описания возможностей управляющих процедур при реализации моментных ограничений <<на грани фола>>.

Важно отметить, что при использовании обобщенных элементов, определяемых в виде тех или иных  конечно-аддитивных мер, соответствующие  множества притяжения совпадают с <<обычными>> множествами достижимости, но только в ситуации, когда эти обобщенные элементы используются в качестве своеобразных <<управлений>>. Итак,  мы  можем, избегая труднореализуемых предельных переходов, используемых при непосредственном построении множеств притяжения, пойти совсем другим и более простым в логическом отношении путем, а именно: создать <<хорошее>>, а точнее, компактное, пространство обобщенных управлений (далее будем использовать данный термин) и, оперируя (обычным образом) с элементами данного пространства как с управлениями, построить область достижимости в традиционном понимании. Здесь, конечно, возникают вопросы, связанные по сути дела с построением некой новой (по крайней мере формально) динамической системы, которая могла бы реагировать на конечно-аддитивные меры как на управления. Однако в рассматриваемом сейчас случае линейных управляемых систем эти вопросы легко решаются посредством соответствующего расширения формулы Коши (см. \cite{25}).

Полезно отметить, что вышеупомянутая схема исследования применима на самом деле в более общей постановке, уже не связанной обязательно с задачей управления. В этой связи обсудим сейчас одну общую схему, которая в настоящем пособии используется в качестве промежуточной, но может иметь и самостоятельное значение (например, в связи с задачами математического программирования). Итак, придерживаясь содержательного способа изложения, наметим данную схему, сохраняя при этом идею, связанную с исследованием вопросов о достижимости.

Пусть $X$\,--- непустое множество, даны функции $$g_1: X\rightarrow\rr,\ldots,g_k:X\rightarrow\rr, \ \  h_1: X\rightarrow\rr,\ldots,h_l:X\rightarrow\rr,$$ где $k,l$\,--- натуральные  числа, $\rr$\,--- вещественная прямая. В $k$-мерном арифметическом пространстве задано множество $\mathbf{Y},$ определяющее условие $$(g_1(x),\ldots,g_k(x))\in\mathbf{Y}$$ на выбор $x\in X.$  Требуется определить, каким при данном условии могут быть векторы $$(h_1(x),\ldots,h_l(x)),$$ реализующиеся в $l$-мерном арифметическом пространстве. Разумеется, ответ можно представить следующим образом: введем вектор-функционалы $\mathbf{G}$ и $\mathbf{H}$ в виде $$x\mapsto(g_1(x),\ldots,g_k(x)) \ \text{и} \ \ x\mapsto(h_1(x),\ldots,h_l(x))$$  соответственно. Тогда
\begin{multline}\label{0.1}
\mathbf{H}^1(\mathbf{G}^{-1}(\mathbf{Y}))=\{\mathbf{H}(x):x\in\mathbf{G}^{-1}(\mathbf{Y})\}=\{(h_1(x),\ldots,h_l(x)):  x\in X,\\ (g_1(x),\ldots,g_k(x))\in\mathbf{Y}\}
\end{multline}
($\mathbf{H}$-образ  $\mathbf{G}$-прообраза $\mathbf{Y}$) определяет множество всех интересующих нас <<целевых>> векторов. Представляет, однако, интерес следующий вопрос: если вместо $\mathbf{Y}$ будет использоваться <<большее>>  множество $\widetilde{\mathbf{Y}}$ (т.е. множество $\widetilde{\mathbf{Y}}$ в $k$-мерном пространстве со свойством $\mathbf{Y}\subset\widetilde{\mathbf{Y}}$), то как изменится интересующее нас множество (\ref{0.1}), а, точнее, каким (<<малым>> или нет) будет изменение
\beq{0.2}
\mathbf{H}^1(\mathbf{G}^{-1}(\mathbf{Y}))\rightarrow\mathbf{H}^1(\mathbf{G}^{-1}(\widetilde{\mathbf{Y}}))
\eeq
 для $\mathbf{Y}\thickapprox\widetilde{\mathbf{Y}};$ ясно, что $\mathbf{H}^1(\mathbf{G}^{-1}(\mathbf{Y}))\subset\mathbf{H}^1(\mathbf{G}^{-1}(\widetilde{\mathbf{Y}})).$ Итак, постулируя близость $\widetilde{\mathbf{Y}}$ к $\mathbf{Y},$ мы интересуемся преобразованием (\ref{0.2}), которое в целом ряде случаев является <<скачком>>: приближение $\widetilde{\mathbf{Y}}$ к $\mathbf{Y}$  извне (в рамках того или иного заданного правила) не приводит, вообще говоря, к приближению $\mathbf{H}^1(\mathbf{G}^{-1}(\widetilde{\mathbf{Y}}))$ к $\mathbf{H}^1(\mathbf{G}^{-1}(\mathbf{Y}));$ см. примеры в \cite{26} . В этих случаях мы интересуемся <<пределом>> множеств $\mathbf{H}^1(\mathbf{G}^{-1}(\widetilde{\mathbf{Y}})),$ используя для его точного определения множества притяжения (см. \cite{46,47,41}).
Тогда ранее упомянутый вариант задачи управления извлекается из схемы с (\ref{0.1}), (\ref{0.2})  в том случае, когда $X$ есть множество программных управлений, т.е. множество в функциональном пространстве. Однако, как уже отмечалось, схема исследования с использованием (\ref{0.1}), (\ref{0.2})  может быть применена и в других содержательных задачах; поэтому в данном пособии она, при некоторой конкретизации $\mathbf{H}$ и $\mathbf{G},$ выделена для отдельного рассмотрения.

В заключении данного введения приведем краткое описание излагаемого ниже материала по главам.

В главе 1 содержатся сведения общематематического  характера, относящиеся к теории множеств, топологии и к некоторым типам измеримых пространств. В краткой форме введены основные понятия, используемые в теории множеств, включая функции, отношения, вещественные числа. Подробно обсуждаются топологии и их построение посредством баз, окрестности и локальные базы (фундаментальные семейства окрестностей), операция замыкания, аксиомы отделимости, компактность. Особое внимание уделяется вопросам сходимости в топологическом пространстве (имеется в виду сходимость по Мору-Смиту). Рассматриваются алгебры и полуалгебры множеств, а также соотношения между этими двумя типами измеримых структур.

В главе 2 излагаются положения конечно-аддитивной теории меры. В частности, рассматриваются нужные в дальнейшем свойства конечно-аддитивных мер, а также выделяются наиболее важные классы таких мер. Вводятся некоторые топологии на пространствах конечно-аддитивных мер соотвествующего типа. Исследуются ступенчатые и ярусные функции (ярусные функции --- равномерные пределы ступенчатых). Рассматриваются линейные непрерывные функционалы на банаховом пространстве ярусных функций, которые затем отождествляются с конечно-аддитивными мерами ограниченной вариации (предварительно излагается простейшая схема интегрирования по упомянутым конечно-аддитивным мерам). Введена $*$-слабая топология; обсуждаются условия   $*$-слабой компактности. Для специального рассмотрения выделяются конечно-аддитивные меры со свойством слабой абсолютной непрерывности относительно фиксированной меры; подробно обсуждается связь с понятием неопределенного интеграла.

Глава 3 посвящена вопросам, связанным с построением и исследованием множеств достижимости в конечномерном арифметическом пространстве при точном и приближенном соблюдении ограничений моментного характера. При этом рассматриваются два различных, вообще говоря, варианта ослабления условий. Введено общее понятие множества притяжения, отвечающее идее достижимости в условиях приближенного соблюдения ограничений. Данное понятие детализировано для двух вышеупомянутых вариантов ослабления ограничений моментного характера. Установлены эффективно проверяемые условия асимптотической нечувствительности при ослаблении части ограничений.

В главе 4 общие конструкции главы 3 применяются для исследования вопросов, связанных с построением и исследованием областей достижимости линейных управляемых систем с разрывностью в коэффициентах при управляющих воздействиях и импульсных ограничениях на выбор самих управлений. Введены обобщенные управления, формализуемые в виде конечно-аддитивных мер со свойством слабой абсолютной неперывности относительно сужения меры Лебега. Более подробно рассматривается случай управления материальной точкой, для которого обсуждается естественный вариант асимптотической нечувствительности при ослаблении ограничений на скорость.

Заметим, что в настоящем пособии излагается  практически весь необходимый математический аппарат, включая сведения из теории множеств и общей топологии. Нам представляется, что знакомство  с этим материалом также будет полезным заинтересованному читателю.

\vspace{5mm}

\emph{Данное пособие подготовлено при финансовой поддержке РФФИ (проекты   № 16-01-00505, 16-01-00649).}
%\end{document}
%\begin{center} \section*{Глава 1} \end{center}

\newpage

\begin{center} \section*{Глава 1. Общие сведения } \end{center}
\section{Элементы теории множеств: простейшие конструкции}\setcounter{equation}{0} \setcounter{proposition}{0}

\ \ \ \ \ Основные конструкции в дальнейшем будут связаны с  вещественнозначными (в/з) функциями множеств, обладающими некоторыми специальными свойствами (конечной и счетной аддитивностью). В этой связи представляется необходимой краткая сводка понятий теории множеств. Эта сводка приводится в настоящем разделе (см. также \cite[гл.\,1]{30}, \cite{31,32}).

Рассматриваем объекты и свойства, которыми могут обладать те или иные объекты. Среди объектов особо будем выделять \emph{множества}. Множества могут содержать другие объекты в качестве своих элементов. То, что объект $a$ является элементом множества $A,$ обозначается посредством выражения $a\in A.$ Если же, напротив, $a$ не является элементом $A,$ используем обозначение $a\notin A.$

В дальнейшем следуем традиционным соглашениям, касающимся применения кванторов и связок для сокращенной записи словесных высказываний: $\forall\,\ldots$\,--- для всякого $\ldots,$ $\exists\ldots$\,--- существует $\ldots$, $\exists\,!\ldots$\,--- существует и единственно $\ldots,$ $\neg$\,--- не, $\Rightarrow$\,--- влечет, $\Leftrightarrow$\,--- эквивалентно, $\&$\,--- и, $\vee$\,--- или, $\pr{def}$ заменяет фразу <<по определению>>, $\triangleq$\,--- равенство по определению. Двоеточие в формулах заменяет фразу <<такой, что>> (<<такая, что>>, <<такое, что>>).

Если $A$  и $B$\,--- множества, то, как обычно,
$$(A\subset B)\stackrel{\pr{def}}{\Leftrightarrow}(a\in B \ \forall\,a\in A).$$ Если $P$ и $Q$\,--- множества, причем $P\subset Q,$ то называем $P$ \emph{подмножеством} (п/м) множества~$Q.$ Постулируем, что
\beq{2.1}
(P=Q)\stackrel{\pr{def}}{\Leftrightarrow}((P\subset Q)\ \& \ (Q\subset P)).
\eeq
Итак, два множества совпадают в том и только том случае, если они состоят из одних и тех же элементов. Соглашение (\ref{2.1}) является очень важным и неоднократно используется в последующих построениях.

Полагаем, что $\zer$ есть (единственное в силу (\ref{2.1})) множество, не содержащее ни одного элемента. Тогда в случае, если $X$\,--- какое-либо множество, запись $X\neq\zer$ означает, что $X$\,--- непустое множество, т.е. существует некоторый объект, являющийся элементом $X.$

Считаем известными основные теоретико-множественные операции (объединение, пересечение, разность и т.д.); см., в частности, \cite{31,32}. Поэтому лишь кратко напомним некоторые из них. Так, например, если $A$ и $B$\,--- множества, то  множество $A\cup B$ таково, что
\beq{2.1'}
(A\subset A\cup B)\ \& \ (B\subset A\cup B)\ \& \ (\forall\,x\in A\cup B \ (x\in A)\vee (x\in B)).
\eeq

Если $a$ и  $b$\,--- некоторые объекты, то $\{a;b\}$ есть (единственное в силу (\ref{2.1})) множество, содержащее в качестве своих элементов $a,$ $b$ и не содержащее никаких других элементов; $\{a;b\}$ есть, таким образом, неупорядоченная пара объектов $a$ и $b.$ Если же $x$\,--- объект, то  $$\{x\}\triangleq\{x;x\}$$ есть одноэлементное множество, содержащее (единственный) элемент $x.$

Множество, все элементы которого сами являются множествами, называем \emph{семейством}. Каждому семейству $\mathcal{A}$ сопоставляется единственное множество \beq{2.2}
\bigcup\limits_{A\in\mathcal{A}}A,
\eeq
для которого имеют место свойства
$$\biggl(B\subset\bigcup\limits_{A\in\mathcal{A}}A \ \forall\,B\in\mathcal{A}\biggl)\ \& \ \biggl(\forall\,x\in\bigcup\limits_{A\in\mathcal{A}}A \ \exists\,C\in\mathcal{A}:x\in C\biggl).$$
Заметим, что в (\ref{2.2}) вместо $A$ может использоваться произвольная буква; буква  $A$ играет роль своеобразного  <<бегающего индекса>>. При этом, конечно, для всякого множества $M$ в виде $\{M\}$ имеем (непустое) семейство и при этом
$$M=\bigcup\limits_{A\in\{M\}}A.$$
Аналогичным образом, для всяких двух множеств $U$ и $V$ в виде $\{U;V\}$ имеем непустое семейство и при этом (см. (\ref{2.1'}))
$$U\cup V=\bigcup\limits_{A\in\{U;V\}}A.$$

Если $T$\,--- множество, то через $\pp(T)$ обозначаем семейство всех п/м $T$ (см. \cite{31,32,33}; в частности, см. \cite[с.\,12]{33}), т.е. семейство всех множеств, содержащихся в $T.$  Итак, $\pp(T)$ есть единственное семейство, такое, что\\
1) \  для всякого множества $S$ истинна импликация $$(S\subset T)\Rightarrow(S\in\pp(T));$$
2) \  $M\subset T \ \ \forall\,M\in\pp(T).$\\
Ясно, что $\{\zer\}$ есть семейство, для которого $\zer\in\{\zer\}.$  Итак, $\{\zer\}$ есть непустое семейство.
Располагая множеством $X$ и некоторым свойством $\ldots,$   применимым к элементам $X,$ через
\beq{2.3} \{x\in X\mid\ldots\} \eeq  условимся обозначать множество всех элементов из $X,$ обладающих свойством $\ldots;$  ясно, что вместо $x$ в (\ref{2.3}) может использоваться произвольная буква. Отметим два простейших примера применения (\ref{2.3}):
если $A$ и $B$\,--- множества, то
\beq{2.4}
A\setminus B\triangleq\{x\in A\mid x\notin B\}
\eeq
  есть разность множеств $A$ и $B;$  для произвольных множеств $M$ и  $N$ полагаем,  учитывая  (\ref{2.4}), что $$M\cap N\triangleq M\setminus(M\setminus N),$$
получая пересечение этих двух множеств. Ясно, что $M\cap N$ есть единственное множество, для которого при всяком выборе объекта $x$
$$(x\in M\cap N)\Leftrightarrow((x\in M)\ \& \ (y\in N)).$$
Итак,  любым двум множествам $A$ и $B$ мы сопоставляем их объединение $A\cup B$ и пересечение $A\cap B.$ Если же $U,V$  и $W$\,---  три множества, то полагаем, что $$U\cup V\cup W\triangleq(U\cup V)\cup W$$ и, аналогичным образом,  $$U\cap V\cap W\triangleq(U\cap V)\cap W.$$
Легко видеть, что $U\cup V\cup W=U\cup( V\cup W)$ и $U\cap V\cap W= U\cap (V\cap W).$ Пересечение произвольного непустого семейства $\mathcal{A}$ есть множество
\beq{2.5}
\bigcap\limits_{A\in\mathcal{A}}A\triangleq\biggl\{x\in\bigcup\limits_{A\in\mathcal{A}}A\mid x\in\mathbb{A} \ \forall\,\mathbb{A}\in\mathcal{A}\biggl\},
\eeq
для которого имеют место свойства:\\
1)\  если $y\in\bigcap\limits_{A\in\mathcal{A}}A,$ то $y\in B \ \forall\,B\in\mathcal{A};$\\
2) \ для всякого объекта $x$
$$\left(x\in\widetilde{A} \ \forall\,\widetilde{A}\in\mathcal{A}\right)\Rightarrow\biggl(x\in\bigcap\limits_{A\in\mathcal{A}}A\biggl).$$
Для всякого множества $X$ введем  семейство всех непустых его подмножеств:  $$\pp'(X)\triangleq\pp(X)\setminus\{\zer\}.$$

Исключительно важным является следующее определение \emph{упорядоченной пары} (см. \cite{34}): если $x$ и $y$\,--- объекты, то определено непустое семейство
\beq{2.6}
(x,y)\triangleq\{\{x\};\{x;y\}\}
\eeq
(неупорядоченная пара специального вида). Определение (\ref{2.6}) позволяет установить следующее свойство: если $a,b,c$ и $d$\,--- объекты, то
\beq{2.6'}
((a,b)=(c,d))\Leftrightarrow((a=c)\ \& (b=d)).
\eeq
С понятием упорядоченной пары естественно связана конструкция декартова произведения: если $A$ и $B$\,---  множества, то
$$A\times B\triangleq\{z\in\pp(\pp(A\cup B))\mid\exists\,a\in A \ \exists\,b\in B:\ z=(a,b)\}$$
есть единственное множество, для которого:\\
1') \ \ $(a,b)\in A\times B \ \forall\,a\in A \ \forall\,b\in B;$\\
2') \ \ $\forall\,z\in A\times B \ \exists\,a\in A \ \exists\,b\in B: \ z=(a,b).$\\
С учетом 1') и 2') полезно иметь в виду часто используемое выражение $$A\times B=\{(a,b):a\in A, \ b\in B\},$$
толкование которого сводится к сочетанию свойств 1'), 2'). Учитывая эти свойства, легко проверяется, что для всяких четырех множеств $M,N,P$ и $Q$
$$(M\times N)\cap(P\times Q)=(M\cap P)\times(N\cap Q).$$

Если $x,y$ и $z$\,---  объекты, то через $(x,y,z)$ обозначается триплет, определяемый в виде $$(x,y,z)\triangleq((x,y),z),$$
т.е. как упорядоченная пара специального вида. Из общих свойств упорядоченных пар  следует, что для всяких объектов $x,y,z,u,v$ и $w$
$$((x,y,z)=(u,v,w))\Leftrightarrow \left((x=u)\ \& (y=v)\ \& (z=w)\right).$$
В этой связи отметим естественное соглашение, касающееся произведения трех произвольных множеств $A,B$ и $C:$
$$A\times B\times C\triangleq(A\times B)\times C.$$

\begin{center}
\textbf{Отношения и функции}
\end{center}

\emph{Отношением} называем п/м декартова произведения двух произвольных множеств. Отношения суть элементы семейств $\pp(A\times B),$ где $A$ и $B$\,--- какие-либо множества (множество-произведение $A\times B$ само является отношением).

Если $z$\,--- какая-либо упорядоченная пара, то через $\pr{pr}_1(z)$ и $\pr{pr}_2(z)$  обозначаем соответственно первый и второй элементы $z,$ однозначно (см. (\ref{2.6'})) определяемые условием
\beq{2.7}
z=(\pr{pr}_1(z),\pr{pr}_2(z))
\eeq
(более подробно см. в \cite{31,32,33}; в частности, см. \cite[c.\,23-25]{33}). В связи с (\ref{2.7}) отметим, что для всяких двух объектов $x$ и $y$
$$(\pr{pr}_1((x,y))=x)\ \& \ (\pr{pr}_2((x,y))=y).$$
Легко понять, что для любых двух множеств $A$ и $B,$ а также упорядоченной пары $z\in A\times B$
 $$(\pr{pr}_1(z)\in A)\ \& \ (\pr{pr}_2(z)\in B);$$
см. (\ref{2.6'}), а также свойства 1'), \ 2').

Напомним, что отношение $f$ называется \cite[гл.\,2]{34}  \emph{функцией}, если  $\forall\,u\in f \ \forall\,v\in f$
$$(\pr{pr}_1(u)=\pr{pr}_1(v))\Rightarrow(\pr{pr}_2(u)=\pr{pr}_2(v)).$$
Данное общее определение (см. подробнее \cite{31,32,33}) и, в частности, \cite[определение 2.2]{33} ниже практически не используется и мы воспользуемся более частным вариантом \cite[(1.1.9)]{30}: если $X$ и $Y$\,---  множества, то
\beq{2.8}
Y^X\triangleq\{f\in\pp(X\times Y)\mid \forall\,x\in X \ \exists\,!\, y\in Y: (x,y)\in f\}.
\eeq
Элементами данного множества являются отображения из $X$ в $Y$ и только они. При этом $\forall\,f\in Y^X\ \forall\,x\in X \ \exists\,!\, y\in Y: (x,y)\in f.$
С учетом этого полагаем, что $\forall\,f\in Y^X \ \forall\,x\in~X  \ \pr{def} \ f(x)\in Y:$ $$(x,f(x))\in f.$$

Тем самым определено значение  отображения $f$ в точке $x;$  более подробное обсуждение см. в \cite[с.\,30--32]{33}.
Итак, если $X$ и $Y$\,--- множества, $f\in Y^X$ (т.е. $f:X\rightarrow Y$) и $x\in X,$ то $f(x)\in Y$ есть значения отображения $f$ в точке $x.$
В связи с (\ref{2.8}) отметим,  что вместо выражения $f\in Y^X$  будет использоваться также традиционная запись $$f: X\rightarrow Y.$$
Отметим одно очевидное свойство \cite[(1.1.23)]{30}: для всяких множеств $X,Y$ и отображений $f\in Y^X$ и $g\in Y^X$
\beq{2.9}
(f=g)\Leftrightarrow(f(x)=g(x) \ \forall\,x\in X).
\eeq

\begin{center}
\textbf{Индексная форма записи отображений}
\end{center}

Пусть $A$ и $B$\,--- множества, $A\neq\zer,$ и, кроме того, $$b_x\in B \ \forall\,x\in A.$$
Тогда через $(b_x)_{x\in A}$ (вместо $x$ может использоваться произвольная буква) условимся обозначать единственное (в силу (\ref{2.9})) отображение $f\in B^A,$ для которого $f(a)=b_a \ \forall\,a\in A.$ В качестве дублирующего (по отношению к $(b_x)_{x\in A}$) обозначения используем также следующее:
\beq{2.10}
x\mapsto b_x: \ A\rightarrow B
\eeq
(разумеется, в (\ref{2.10}) вместо $x$ может использоваться произвольная буква).

С использованием индексной формы введем понятие суперпозиции отображений (подробнее см. в \cite[раздел 5]{33}): если $X,Y$ и $Z$\,--- множества, $g\in Y^X$ и $h\in Z^Y,$ то $$h\circ g=(h(g(x)))_{x\in X}\in Z^X;$$
иными словами  $h\circ g$ (суперпозиция отображений $g$ и $h$) есть $$x\mapsto h(g(x)): X\rightarrow Z.$$
\begin{center}\textbf{Сужение отображений}\end{center} Если $A$ и $B$\,--- множества, $f\in B^A$ и $C\in\pp(A),$ то \emph{сужение} $$(f|C)\in B^C$$  $f$ на $C$ (т.е. $(f|C): C\rightarrow B$) есть отображение, для которого $$(f|C)(x)\triangleq f(x) \ \ \forall\,x\in C.$$
Разумеется, можно было бы определить $(f|C)$ следующим образом (см.  \cite[с.\,62-63]{33}): $$(f|C)=f\cap(C\times B).$$
\newpage
\begin{center}
\textbf{Образы и прообразы множеств}
\end{center}

Если $X$ и $Y$\,--- множества, $f\in Y^X$ и $M\in\pp(X),$ то $$f^1(M)\triangleq\{f(x): x\in M\}\in\pp(Y)$$ есть единственное множество, обладающее свойствами:\\
a) \ \ $f(x)\in f^1(M) \ \forall\,x\in M;$\\
b) \ \ $\forall\,y\in f^1(M) \ \exists\,x\in M: \ y=f(x).$

Тем самым определен \emph{образ} множества $M$ при действии $f:$ равенство $f^1(M)=\{f(x):x\in M\}$ как раз и означает справедливость свойств a),b).  В качестве $M$ может использоваться $X.$  Если $f^1(X)=Y,$  то отображение $f$ называется \emph{сюръективным} (или отображением $X$ на $Y,$  сюръекцией $X$ на $Y$).
Легко понять (см. \cite[раздел~6]{33}), что при $M_1\in\pp(X)$ и $M_2\in\pp(X)$
\begin{multline*}
\left(f^1(M_1\cap M_2)\subset f^1(M_1)\cap f^1(M_2)\right) \ \& \\ \& \ \left(f^1(M_1\cup M_2)= f^1(M_1)\cup f^1(M_2)\right)\ \& \\ \& \ \left(f^1(M_1)\setminus f^1(M_2)\subset f^1(M_1\setminus M_2)\right)
\end{multline*}

Если $X$ и $Y$\,--- множества, $f\in Y^X$ и $L\in\pp(Y),$ то
 \beq{2.11}
 f^{-1}(L)\triangleq\{x\in X\mid f(x)\in L\}\in\pp(X);
 \eeq
 $f^{-1}(L)$ есть \emph{прообраз} $L$ при действии отображения $f.$  Напомним (см. \cite[раздел 7]{33}) простейшие свойства операции взятия прообраза: если $A\in\pp(Y)$ и $B\in\pp(Y),$ то
\begin{multline*}
\left(f^{-1}(A\cap B)= f^{-1}(A)\cap f^{-1}(B)\right) \ \& \\ \& \ \left(f^{-1}(A\cup B)= f^{-1}(A)\cup f^{-1}(B)\right)\ \& \\ \& \ \left(f^{-1}(A\setminus B)= f^{-1}(A)\setminus f^{-1}(B)\right).
\end{multline*}

Если $S\in\pp(Y),$ то $f^{-1}(Y\setminus S)=X\setminus f^{-1}(S).$  Ряд других (в том числе более общих) свойств см. в \cite[(7.33), (7.43)]{33}.
Отметим, используя  a),b) и (\ref{2.11}), что для всяких множеств $X,Y$ и отображения $f\in Y^X$
$$
\left(A\subset f^{-1}\left(f^1(A)\right) \ \forall\,A\in\pp(X)\right) \ \& \ \left(f^1\left(f^{-1}(B)\right)\subset B \  \forall\,B\in\pp(Y)\right).
$$
\begin{center}
\textbf{След семейства на заданное множество}
\end{center}

Если $\mathcal{X}$\,--- непустое семейство, а $Y$\,--- множество, то полагаем, что
\beq{2.219}
\mathcal{X}|_Y\triangleq\{X\cap Y:X\in\mathcal{\mathcal{X}}\}\in\pp'(\pp(Y)).
\eeq
В связи с корректностью данного определения заметим, что введенный таким образом объект есть на самом деле образ семейства $\mathcal{X}$ при действии отображения из множества $\pp(Y)^{\mathcal{X}}.$ В самом деле, используя индексную форму, рассмотрим отображение
\beq{2.219'}
f\triangleq(X\cap Y)_{X\in\mathcal{X}}\in\pp(Y)^{\mathcal{X}}.
\eeq
Тогда отображение $f:\mathcal{X}\rightarrow\pp(Y)$ таково, что $f(X)=X\cap Y \ \forall\,X\in\mathcal{X}.$  Иными словами, в нашем случае $f$ есть отображение
$$X\rightarrow X\cap Y:\mathcal{X}\rightarrow\pp(Y).$$
Тогда в рамках принятого ранее соглашения имеем в виде
\beq{2.200}
f^1(\mathcal{X})=\{X\cap Y: X\in\mathcal{X}\}
\eeq
<<обычный>> образ множества (точнее семейства) $\mathcal{X}$ при действии $f.$ Варианты (\ref{2.200}) будут в дальнейшем часто использоваться (имеется в виду употребление множеств, подобных множеству в правой части (\ref{2.200})). В частности, выражение, аналогичное применяемому в правой части (\ref{2.200}) потребуется при рассмотрении топологии, индуцированной из заданного пространства на п/м данного пространства. Сейчас же мы использовали (\ref{2.200}) для пояснения (\ref{2.219}): имеем $\mathcal{X}|_Y=f^1(\mathcal{X}),$ где $f$ есть отображение (\ref{2.219'}).

  Если $\mathbb{X}$\,--- множество, $\mathcal{X}\in\pp'(\pp(\mathbb{X}))$  и $Y\in\pp(\mathbb{X}),$  то, как легко видеть,  определено (непустое) семейство (\ref{2.219}), а, точнее, $\mathcal{X}|_Y\in\pp'(\pp(Y)).$

Отметим, что представления, подобные (\ref{2.200}) и реализуемые с помощью других операций над множествами, будут использоваться в дальнейшем неоднократно. Так, в частности, для произвольных множества $U$ и семейства $\mathcal{U}\in\pp'(\pp(U))$
\beq{2.201}
\mathbf{C}_U[\mathcal{U}]\triangleq\{U\setminus\widetilde{U}:\widetilde{U}\in\mathcal{U}\}\in\pp'(\pp(U))
\eeq
 есть (непустое) семейство п/м $U,$  двойственное к $\mathcal{U}:$
 $$(U\setminus V\in\mathbf{C}_U[\mathcal{U}] \ \forall\,V\in \mathcal{U}) \ \& \ (\forall\,A\in\mathbf{C}_U[\mathcal{U}] \ \exists\,B\in\mathcal{U}: \ A=U\setminus B)$$
 (семейство (\ref{2.201}) также может быть представлено в виде образа семейства $\mathcal{U}$ при действии отображения
 $$\widetilde{U}\mapsto U\setminus\widetilde{U}:\mathcal{U}\rightarrow\pp(U),$$
 но мы это представление сейчас подробно не рассматриваем).

\newpage
\begin{center}
\textbf{Аксиома выбора в форме Рассела (аксиома мультипликативности)}
\end{center}

Если $X$ и $Y$\,--- множества, а $\Phi\in\pp(Y)^X,$ то
\beq{2.12}
\prod\limits_{x\in X}\Phi(x)\triangleq\{f\in Y^X\mid f(x)\in \Phi(x) \ \ \forall\,x\in X\}
\eeq
есть декартово произведение всех множеств $\Phi(x), \ x\in X.$  Заметим, что отображение  $\Phi$ (мультифункция из $X$ в $Y$) часто записывается при использовании  (\ref{2.12}) в индексной форме. Тогда, как вариант (\ref{2.12}), имеем при заданных множествах $X,$ $X\neq\zer,$ и $Y,$ а также при
  $(A_x)_{x\in X}\in\pp(Y)^X,$ что
\beq{2.13}
\prod\limits_{x\in X}A_x=\{f\in Y^X\mid f(x)\in A_x \ \ \forall\,x\in X\}.
\eeq
Здесь (см. (\ref{2.12})) $\Phi=(A_x)_{x\in X},$ что означает справедливость системы равенств
\beq{2.14}
\Phi(y)=A_y \ \ \forall\,y\in X.
\eeq

В связи с (\ref{2.12}), (\ref{2.13}) особый интерес представляет случай, когда все множества в (\ref{2.14}) являются непустыми. Принимаем в этой связи следующую аксиому: если $X$ и $Y$\,--- множества, $X\neq\zer,$ и при этом $\Phi\in\pp'(Y)^X,$  то \cite{34}
\beq{2.15}\prod\limits_{x\in X}\Phi(x)\in\pp'(Y^X).\eeq
Свойство, связанное с (\ref{2.15}), можно рассматривать как аксиому выбора в форме Б.Рассела (аксиома мультипликативности). Итак, если $A_x\in\pp'(Y)$ при $x\in X,$ то
$$\prod\limits_{x\in X}A_x\in\pp'(Y^X);$$
разумеется (см. (\ref{2.13})), при этом  $f(y)\in A_y \ \ \forall\,f\in \prod\limits_{x\in X}A_x.$
\begin{center}
\textbf{Сюръективные, биективные и  инъективные отображения}\end{center}

В пределах настоящего пункта фиксируем непустые множества $X$ и $Y.$  Элементы множества
\beq{2.15'}
(\pr{su})[X;Y]\triangleq\{f\in Y^X\mid f^1(X)=Y\}\in\pp(Y^X)
\eeq
называют сюръективными отображениями из $X$ в $Y,$ или  \emph{сюръекциями} $X$ на $Y.$ Таким образом,  в (\ref{2.15'}) речь идет об отображениях из $Y^X,$ для  каждого из которых образ множества $X$ совпадает с $Y.$ Далее полагаем, что
\beq{2.15''}
(\pr{bi})[X;Y]\triangleq\{f\in(\pr{su})[X;Y]\mid\forall\,x_1\in X \ \forall\,x_2\in X \ (f(x_1)=f(x_2))\Rightarrow(x_1=x_2)\}
\eeq
получая  множество всех \emph{биективных} (взаимно однозначных) отображений $X$ на $Y.$ Если $X=Y,$ то элементы множества (\ref{2.15''}) называем перестановками множества $X$ (см. \cite[c.\,87]{35}).

Если $f\in Y^X$ таково, что $\forall\,x_1\in X \ \forall\,x_2\in X$ $$(f(x_1)=f(x_2))\Rightarrow(x_1=x_2),$$ то отображение $f$ называется \emph{инъективным}.

\begin{center}
\textbf{Вещественные числа}
\end{center}

Всюду в дальнейшем $\rr$\,--- вещественная прямая, элементы множества $\rr$ суть вещественные числа. Будем полагать, что элементы $\rr$ не являются множествами с тем, чтобы избежать двусмысленности в традиционных обозначениях; разумеется, $\rr$\,--- непустое множество, и мы считаем известными основные операции с вещественными числами (сложение, умножение, деление), а также свойства этих операций. Через $\leqslant$ обозначаем обычную упорядоченность $\rr$ (строго говоря, $\leqslant\in\pp(\rr\times\rr)$ и при $x\in\rr, \ y\in\rr$ выражение $x\leqslant y$ эквивалентно включению  $(x,y)\in\leqslant;$ данное отношение $\leqslant$ рефлексивно, транзитивно и антисимметрично; см. \cite[c.\,26]{30}).

Конечные и бесконечные промежутки в $\rr$ условимся обозначать только квадратными скобками (см. \cite[c.\,35,36]{30}).

Считаем известными свойства точных верхней и нижней граней непустых п/м $\rr,$ ограниченных соответственно сверху и снизу; подробнее см. в  \cite[c.\,37]{30}. В соответствующих обозначениях используем традиционные символы $\pr{sup}$ и $\pr{inf}.$ В частности, при $x\in\rr$  имеем $|x|=\pr{sup}(\{x;-x\})\in[0;\infty[,$ где $-x=(-1)x.$

Как обычно, $\mathbb{N}\in\pp'(\rr)$ есть натуральный ряд, $\mathbb{N}\triangleq\{1;2;\ldots\}$\,--- множество, для которого
\begin{multline*}
(1\in\mathbb{N})\ \& \ (k+1\in\mathbb{N} \ \forall\,k\in\mathbb{N}) \ \& \ (\,]n,n+1[\,\cap\,\mathbb{N}= \\ =\{\xi\in\rr\mid(n<\xi)\ \&\ (\xi<n+1)\}\cap\mathbb{N}=\zer\ \forall\,n\in\mathbb{N})\ \& \\ \& \ (\forall\,x\in\rr \  \exists\,N\in\mathbb{N}:x\leqslant N).
\end{multline*}

В дальнейшем будут использоваться следующие обозначения:
 \begin{multline*}
 \left(\,\overline{p,q}\triangleq\{i\in\mathbb{N}\mid(p\leqslant i) \ \& \ (i\leqslant q)\} \ \forall\,p\in\mathbb{N} \ \forall\,q\in\mathbb{N}\right) \ \& \\ \& \left(\overrightarrow{n,\infty}\triangleq\{s\in\mathbb{N}\mid n\leqslant s\} \ \forall\,n\in\mathbb{N}\right).
 \end{multline*}
 Заметим, что $\overline{1,m}=\{i\in\mathbb{N}\mid i\leqslant m\}$ при $m\in\mathbb{N};$ кроме того, при $p\in\mathbb{N}$ и $q\in\mathbb{N}$ допускается возможность $\overline{p,q}=\zer,$  реализующаяся при $q<p.$ В  этой связи полезно отметить, что  при $a\in\rr$ и $b\in\rr$ случай $[a,b]=\zer$ также допускается.
 Этот случай реализуется при $b<a;$  аналогичным образом заметим,  что $]a,b]=\zer,$ $[a,b[=\zer$ и $]a,b[=\zer$  при $b\leqslant a$.  Заметим также, что (см. \cite[(1.3.4)]{30}) $\forall\,E\in\pp'(\mathbb{N})$
 $$((1\in E)\ \& \ (k+1 \in E \ \forall\,k\in E))\Rightarrow(E=\mathbb{N}).$$
  Если $T$\,--- множество, то последовательностью в $T$ называется всякое отображение из $\mathbb{N}$ в $T;$  при этом $T^{\mathbb{N}}$\,---  множество всех последовательностей в $T.$ Для обозначения последовательностей традиционно используется индексная форма записи: $(t_i)_{i\in\mathbb{N}}\in T^{\mathbb{N}}$ означает, что $t_j\in T \ \forall\,j\in\mathbb{N}.$  Более подробно рассмотрим случай $T=\rr,$ имея в виду последовательности вещественных чисел. Как обычно \cite[(1.3.6)]{30}, при $(\xi_i)_{i\in\mathbb{N}}\in\rr^{\mathbb{N}}$ и $\xi\in\rr$
  \beq{2.16}
  \left((\xi_i)_{i\in\mathbb{N}}\rightarrow\xi\right)\stackrel{\pr{def}}\Leftrightarrow\left(\forall\,\varepsilon\in]0,\infty[ \ \exists\,m\in\mathbb{N}: \ |\xi_k-\xi|<\varepsilon \ \forall\,k\in\overrightarrow{m,\infty}\right).
  \eeq
  Тем самым определена \emph{секвенциальная сходимость} в $\rr.$ Отметим важное свойство полноты $\rr:$ $\forall\,(\xi_i)_{i\in\mathbb{N}}\in\rr^{\mathbb{N}}$
  \begin{multline*}\left(\forall\,\varepsilon\in]0,\infty[ \ \exists\,n\in\mathbb{N}: \ |\xi_p-\xi_q|<\varepsilon \ \forall\,p\in\overrightarrow{n,\infty} \ \forall\,q\in\overrightarrow{n,\infty}\right)\Leftrightarrow\\ \Leftrightarrow \left(\exists\,\xi\in\rr: \ (\xi_i)_{i\in\mathbb{N}}\rightarrow\xi\right).\end{multline*}

 Кроме того, в условиях, определяющих  (\ref{2.16}), число $\xi$\,--- предел последовательности $(\xi_i)_{i\in\mathbb{N}}$\,--- определяется единственным образом. Итак, $\forall\,(t_i)_{i\in\mathbb{N}}\in\rr^{\mathbb{N}} \ \forall\,t^*\in\rr \ \forall\,t_*\in\rr$
 $$\left(\left((t_i)_{i\in\mathbb{N}}\rightarrow t_*\right) \ \& \ \left((t_i)_{i\in\mathbb{N}}\rightarrow t^*\right)\right)\Rightarrow\left(t_*=t^*\right).$$
 Используем обычные обозначения для (конечных) сумм наборов вещественных чисел. В этой связи условимся об общем соглашении: если $T$\,--- множество и $m\in\mathbb{N},$ то, как обычно, вместо $T^{\overline{1,m}}$ будем использовать запись $T^m,$ полагая всюду в дальнейшем, что элементы $\mathbb{N}$ (натуральные числа) не являются множествами. Последнее соглашение принимается в целях исключения <<двойных>> обозначений (обозначений с двойным толкованием). Итак, элементы $T^m$ суть отображения (кортежи)
 \beq{2.17}
 (t_i)_{i\in\overline{1,m}}:\overline{1,m}\rightarrow T.
 \eeq
  В качестве $T$  может использоваться семейство и в этом случае элементы $T^m$ являются кортежами множеств. С другой стороны, для нас будет важен случай $T=\rr.$ Тогда при $m\in\mathbb{N}$ элементы $\rr^m$ рассматриваем в качестве $m$-мерных векторов; их толкование как кортежей~(\ref{2.17}),  где $T=\rr,$ будет более удобным по соображениям методического характера. Если $m\in\mathbb{N}$ и $(t_i)_{i\in\overline{1,m}}\in\rr^m,$ то обычным образом определяется сумма $$\sum\limits_{i=1}^mt_i\in\rr$$
 чисел $t_1,\ldots,t_m$ (данное понятие считаем известным; см. также \cite[с.\,43,44]{30}).  Отметим только, что при $n\in\mathbb{N},$ $(x_i)_{i\in\overline{1,n}}\in\rr^n$ и $l\in(\pr{bi})[\overline{1,n};\overline{1,n}]$
 \beq{2.17'}
 \sum\limits_{i=1}^nx_i=\sum\limits_{i=1}^nx_{l(i)}.
 \eeq
\begin{center}
\textbf{Конечные множества}\end{center}

Непустое множество $K$ называем \emph{конечным}, если
\beq{2.18}
\exists\,n\in\mathbb{N}: \ (\pr{bi})[\,\overline{1,n};K]\neq\zer.
\eeq
Полагаем также, следуя традиции, что пустое множество конечно. Итак, называем множество $K$ конечным, когда
$$(K=\zer)\vee(\exists\,n\in\mathbb{N}: \ (\pr{bi})[\,\overline{1,n};K]\neq\zer).$$

Отметим, что для всякого непустого множества $K$  эквивалентны свойства:  a) $K$ конечно; b) $\exists\,n\in\mathbb{N}: \ (\pr{su})[\,\overline{1,n};K]\neq\zer.$  Таким образам, требование инъективной нумерации элементов конечных множеств может быть, в принципе, ослаблено.
%\end{document}
Если $X$\,--- произвольное множество, то
\beq{2.19}
\pr{Fin}(X)\triangleq\left\{K\in\mathcal{P}'(X)\mid \exists\,n\in\mathbb{N}: \ (\pr{bi})[\,\overline{1,n};K]\neq\zer \right\}
\eeq
есть семейство всех непустых конечных п/м $X.$

Заметим, что \cite[\S\,1.4]{30}  на самом деле для каждого непустого конечного множества $\mathbb{K}$
\beq{2.20}
\exists\,!\, n\in\mathbb{N}: \ (\pr{bi})[\,\overline{1,n};\mathbb{K}]\neq\zer;
\eeq
с учетом этого полагаем, что число $|\mathbb{K}|\in\mathbb{N}$  (мощность $\mathbb{K}$) таково, что
\beq{2.20}
(\pr{bi})[\,\overline{1,|\mathbb{K}|};\mathbb{K}\,]\neq\zer.
\eeq
Напомним, что  \cite[\S\,1.4]{30} для всякого множества $A$
$$
\left(X\!\cup\! Y\!\!\in\!\pr{Fin}(A) \ \forall\,X\!\in\!\pr{Fin}(A) \ \forall\,Y\!\in\!\pr{Fin}(A)\right)\ \& \ \left(\pp'(K)\!\subset\!\pr{Fin}(A) \ \forall\,K\!\in\!\pr{Fin}(A)\right).
$$
Если же $X$ и $Y$\,--- множества, $P\in\pr{Fin}(X)$ и $Q\in\pr{Fin}(X),$ то $$P\times Q\in\pr{Fin}(X\times Y).$$
\begin{center}
\textbf{Конечные разбиения}\end{center}

Условимся о соглашении: если $A$\,--- множество, $\mathcal{A}\in\pp'(\pp(A))$   (непустое семейство п/м $A$), $B\in\pp(A)$ и $n\in\mathbb{N},$ то
\beq{2.21}
\triangle_n(B\!,\mathcal{A})\!\triangleq\!\biggl\{\!(A_i)_{i\in\overline{1,n}}\!\in\!\mathcal{A}^n\mid(B\!=\!\bigcup\limits_{i=1}^nA_i) \&  \left(A_p\!\cap\! A_q\!=\!\zer  \forall\,p\!\in\!\overline{1,n} \, \forall\,q\!\in\!\overline{1,\!n}\!\setminus\!\{p\}\!\right)\!\biggl\}
\eeq
(см. в этой связи общее определение в \cite[c.\,48]{30}) есть множество всех упорядоченных $\mathcal{A}$-разбиений множества $A,$ имеющих <<длину>> $n.$ Заметим, что (\ref{2.21}) соответствует \cite[(1.4.9)]{30}. Наряду с упорядоченными будем использовать неупорядоченные конечные разбиения, определяемые в виде непустых конечных подсемейств заданного семейства: если $S$\,--- множество, $\mathcal{S}\in\pp'(\pp(S))$ и $H\in\pp(S),$ то
\begin{multline}\label{2.22}
\mathbf{D}(H,\mathcal{S})\triangleq\biggl\{\mathcal{K}\in\pr{Fin}(\mathcal{S})\mid(H=\bigcup\limits_{L\in\mathcal{K}}L)\ \& \\ \& \bigl(\forall\,A\in\mathcal{K} \ \forall\,B\in\mathcal{K} \ (A\cap B\neq\zer)\Rightarrow (A=B)\bigl)\biggl\};
\end{multline}
см. \cite[c.\,58]{36}.  Конструкции (\ref{2.21}), (\ref{2.22}) взаимосвязаны, что проявляется в двух свойствах, приводимых ниже.

\noindent\emph{Свойство 1*)}   Если $X$\,--- множество,  $\mathcal{X}\in\pp'(\pp(X)),$ $Y\in\pp(X),$ $m\in\mathbb{N}$  и $(X_i)_{i\in\overline{1,m}}\in\triangle_m(Y,\mathcal{X}),$ то
\beq{2.23}
\{X_i:i\in\overline{1,m}\,\}\in\mathbf{D}(Y,\mathcal{X}).
\eeq
В самом деле, через $\mathfrak{X}$ обозначим семейство в левой части (\ref{2.23}).  Ясно, что $\mathfrak{X}\in\pr{Fin}(\mathcal{X}).$  Здесь мы учитываем, что
$(X_i)_{i\in\overline{1,m}}\in(\pr{su})[\,\overline{1,m};\mathfrak{X}].$  Кроме того,
\beq{2.23'}
Y=\bigcup\limits_{i=1}^mX_i=\bigcup\limits_{\mathbb{X}\in\mathfrak{X}}\mathbb{X}.
\eeq

Пусть теперь $\widetilde{A}\in\mathfrak{X}$ и $\widetilde{B}\in\mathfrak{X}$ таковы, что
\beq{2.24}
\widetilde{A}\cap\widetilde{B}\neq\zer.
\eeq
Подберем (см. (\ref{2.23})) $p\in\overline{1,m}$ и $q\in\overline{1,m}$ так, что
\beq{2.25}
\left(\widetilde{A}=X_p\right) \ \& \ \left(\widetilde{B}=X_q\right).
\eeq
Если $p\neq q,$ то согласно (\ref{2.21}) $X_p\cap X_q=\zer,$   а потому $\widetilde{A}\cap\widetilde{B}=\zer,$  что невозможно (см. (\ref{2.24})). Полученное (при условии $p\neq q$) противоречие означает, что $p=q,$ а потому согласно (\ref{2.25})
\beq{2.26}
\widetilde{A}=\widetilde{B}.
\eeq
Итак, (\ref{2.24})$\Rightarrow$(\ref{2.25}).   Поскольку выбор $\widetilde{A}, \widetilde{B}$  был произвольным, установлено, что $\forall\,A\in\mathfrak{X}$ $\forall\,B\in\mathfrak{X}$
\beq{2.27}
(A\cap B\neq\zer)\Rightarrow(A=B).
\eeq
Из  (\ref{2.22}),  (\ref{2.23'}) и  (\ref{2.27}) вытекает, что $\mathfrak{X}\in\mathbf{D}(Y,\mathcal{X}),$ что означает справедливость  (\ref{2.23}).

\noindent\emph{Свойство 2*)} Если $S$\,-- множество, $\mathcal{S}\!\!\in\!\pp'(\pp(S)),$ $H\!\!\in\!\pp(S)$ и $\mathcal{K}\!\in\!\mathbf{D}(H,\mathcal{S}),$ то
\beq{2.28}
(\pr{bi})[\overline{1,|\mathcal{K}|};\mathcal{K}]\subset\bigtriangleup_{|\mathcal{K}|}(H,\mathcal{S}).
\eeq
Действительно, пусть $n\triangleq|\mathcal{K}|$ (тогда $n\in\mathbb{N}$). Выберем произвольно
\beq{2.29}
(S_i)_{i\in\overline{1,n}}\in(\pr{bi})[\overline{1,n};\mathcal{K}].
\eeq
Тогда (см. (\ref{2.15''})) имеем, в частности, что
\beq{2.30}
(S_i)_{i\in\overline{1,n}}\in(\pr{su})[\overline{1,n};\mathcal{K}]
\eeq
и, в частности (см. (\ref{2.15'})),
\beq{2.31}
(S_i)_{i\in\overline{1,n}}:\overline{1,n}\rightarrow\mathcal{K}.
\eeq
При этом согласно (\ref{2.15''}) $\forall\,i_1\in\overline{1,n} \ \forall\,i_2\in\overline{1,n}$
\beq{2.32}
(S_{i_1}=S_{i_2})\Rightarrow(i_1=i_2).
\eeq
Полагая для краткости  $\varphi\triangleq(S_i)_{i\in\overline{1,n}},$ имеем свойство $\varphi\in\mathcal{K}^n,$ т.е. $\varphi:\overline{1,n}\rightarrow\mathcal{K}.$ Тогда $\varphi^1(\overline{1,n})\in\pp(\mathcal{K})$ обладает свойствами
\beq{2.33}
\left(S_i\in\varphi^1(\overline{1,n}) \  \forall\,i\in\overline{1,n}\right) \ \& \ \left(\forall\,\Phi\in\varphi^1(\overline{1,n}) \ \exists\,j\in\overline{1,n}:\Phi=S_j\right).
\eeq
Свойства (\ref{2.33}) кратко записываются в виде равенства
 $$\varphi^1(\overline{1,n})=\{S_t:t\in\overline{1,n}\}$$ (здесь $t$\, --- <<бегающий>> индекс). Из (\ref{2.15}) и (\ref{2.30}) следует, что $\varphi\in(\pr{su})[\overline{1,n};\mathcal{K}],$  а потому $\varphi^1(\overline{1,n})=\mathcal{K}.$ С учетом (\ref{2.33}) имеем
 \beq{2.34}
 \left(S_i\in\mathcal{K} \ \forall\,i\in\overline{1,n}\,\right)\ \& \ \left(\forall\,\Phi\in\mathcal{K} \ \exists\,j\in\overline{1,n}:\Phi=S_j\right).
 \eeq
Свойства (\ref{2.34}) в краткой форме можно записать следующим образом:
$$\mathcal{K}=\{S_t:t\in\overline{1,n}\}.$$
Покажем теперь, что $(S_i)_{i\in\overline{1,n}}\in\triangle_n(H,\mathcal{S}),$  где согласно (\ref{2.21})
\begin{multline}\label{2.35}
\triangle_n(H,\mathcal{S})=\biggl\{(A_i)_{i\in\overline{1,n}}\in\mathcal{S}^n\mid\bigl(H=\bigcup\limits_{i=1}^nA_i\bigl)\ \& \\ \left(A_p\cap A_q=\zer \ \ \forall\,p\in\overline{1,n} \ \ \forall\,q\in\overline{1,n}\setminus\{p\}\right)\biggl\}.
\end{multline}
По выбору $\mathcal{K}$ имеем, что (см. (\ref{2.22})) $\mathcal{K}\in\pr{Fin}(\mathcal{S})$ и согласно (\ref{2.19}) $\mathcal{K}\in\pp'(\mathcal{S}).$ Тогда $\mathcal{K}\neq\zer$ и $\mathcal{K}\subset\mathcal{S}.$ Согласно (\ref{2.31}) имеем теперь, что $$(S_i)_{i\in\overline{1,n}}:\overline{1,n}\rightarrow\mathcal{S},$$ т.е. $(S_i)_{i\in\overline{1,n}}\in\mathcal{S}^n.$ Далее в силу (\ref{2.22})
\beq{2.36}
\bigl(H=\bigcup\limits_{L\in\mathcal{K}}L\bigl)\ \& \ (\forall\,A\in\mathcal{K} \ \forall\,B\in\mathcal{K} \ \ (A\cap B\neq\zer)\Rightarrow(A=B)).
\eeq
Из (\ref{2.34}) следует, однако, равенство $$\bigcup\limits_{L\in\mathcal{K}}L=\bigcup\limits_{i=1}^nS_i;$$
см. \cite[(1.1.27),(1.1.30)]{30}. С учетом (\ref{2.36}) получаем, что
\beq{2.37}
H=\bigcup\limits_{i=1}^nS_i.
\eeq
Выберем произвольно $\mu\in\overline{1,n}$ и $\nu\in\overline{1,n}\setminus\{\mu\}.$ Тогда $\nu\in\overline{1,n}$ и при этом
\beq{2.38}
\mu\neq\nu.
\eeq
В силу (\ref{2.32}) $(S_\mu=S_\nu)\Rightarrow (\mu=\nu),$ а потому (см. (\ref{2.38}))
\beq{2.39}
S_\mu\neq S_\nu,
\eeq
где (согласно (\ref{2.34})) $S_\mu\in\mathcal{K}$ и $S_\nu\in\mathcal{K}.$ С учетом (\ref{2.36}) имеем импликацию $(S_\mu\cap S_\nu\neq~\zer)\Rightarrow(S_\mu=S_\nu)$
и, используя (\ref{2.39}), получаем, что $$S_\mu\cap S_\nu=\zer.$$

Поскольку выбор $\mu$ и $\nu$ был произвольным,
\beq{2.40}
S_p\cap S_q=\zer \ \forall\,p\in\overline{1,n} \ \forall\,q\in\overline{1,n}\setminus\{p\}.
\eeq
Таким образом (см. (\ref{2.37}), (\ref{2.40})) $(S_i)_{i\in\overline{1,n}}\in\mathcal{S}^n:$
$$\biggl(H=\bigcup\limits_{i=1}^nS_i\biggl)\ \& \ \left(S_p\cap S_q=\zer \ \forall\,p\in\overline{1,n} \ \forall\,q\in\overline{1,n}\setminus\{p\}\right).$$
С учетом (\ref{2.35}) получаем включение
$$(S_i)_{i\in\overline{1,n}}\in\triangle_n(H,\mathcal{S}).$$

Поскольку выбор $(S_i)_{i\in\overline{1,n}}$ (\ref{2.29}) был произвольным, установлено вложение
$$(\pr{bi})[\overline{1,n};\mathcal{K}]\subset\triangle_n(H,\mathcal{S}),$$
из которого (по определению $n$) извлекается требуемое свойство:
$$(\pr{bi})[\overline{1,|\mathcal{K}|};\mathcal{K}]\subset\triangle_{|\mathcal{K}|}(H,\mathcal{S}).$$
Итак, имеем вложение (\ref{2.28}), что и требовалось доказать.
Полученные свойства $1^*), 2^*)$ означают фактически, что упорядоченные и неупорядоченные разбиения по сути\, --- одно и то же.

Возвращаясь к (\ref{2.17'}), отметим (см. \cite[c.\,45]{30}), что для всякого непустого конечного множества $K$ и функции $(f_k)_{k\in K}\in\rr^{K}$
\beq{2.41}
\exists\,!\,c\in\rr: \ c=\sum\limits_{j=1}^{|K|}f_{l(j)} \ \forall\,l\in(\pr{bi})[\overline{1,|K|};K]
\eeq
(свойство (\ref{2.41}) легко извлекается из (\ref{2.17'})). С учетом (\ref{2.41}) полагаем, как обычно, что для всяких непустого множества $K$ и функции $(f_k)_{k\in K}\in\rr^K$
\beq{2.42}
\sum\limits_{k\in K}f_k\in\rr
\eeq
есть такое (единственное) число, что
\beq{2.43}
\sum\limits_{k\in K}f_k=\sum\limits_{j=1}^{|K|}f_{l(j)} \ \forall\,l\in(\pr{bi})[\overline{1,|K|};K].
\eeq
Тогда, в частности, для всякого непустого множества $A,$ функции $(f_\alpha)_{\alpha\in A}\in\rr^A$ и множества $K\in\pr{Fin}(A)$ определено
значение $$\sum\limits_{\alpha\in K}f_\alpha\in\rr;$$
мы учитываем здесь то, что $K$\,--- непустое конечное множество и $$(f_\alpha)_{\alpha\in K}\in\rr^K$$ (сужение исходной функции
 $(f_\alpha)_{\alpha\in A}$ на $K$ есть, очевидно,  функция $(f_\alpha)_{\alpha\in K}$).

Совсем кратко напомним свойства конечных сумм  (\ref{2.43}) (подробнее см. \cite[\S1.4]{30}). Прежде всего отметим, что для всякого множества $A,$ функции $(f_\alpha)_{\alpha\in A}\in\rr^A,$ множеств $P\in\pr{Fin}(A)$ и $Q\in\pr{Fin}(A)$
\beq{2.44}
(P\cap Q=\zer)\Rightarrow\biggl(\,\sum\limits_{\alpha\in P\cup Q}f_\alpha=\sum\limits_{\alpha\in P}f_\alpha+\sum\limits_{\alpha\in Q}f_\alpha\biggl).
\eeq
Здесь же отметим следующее простое обстоятельство: для всякого непустого конечного множества $K$ непременно имеет место $K\in\pr{Fin}(K)$ и, более того,
$$\pp'(K)=\pr{Fin}(K).$$

\section{Элементы топологии, 1}\setcounter{equation}{0} \setcounter{proposition}{0}

 \ \ \ \ \ Применяемые в основной части конструкции базируются на топологических свойствах пространства конечно-аддитивных (к.-а.) мер ограниченной вариации, определяемых, вообще говоря, на полуалгебре множеств. Существенно то, что используемые топологии оказываются при этом  неметризуемыми, хотя и удобными в некоторых других отношениях (например, в смысле условий компактности). Поэтому возникает необходимость в сводке основных топологических понятий, которая и будет приведена в настоящем разделе. при этом будем использовать известные положения, приведенные, например, в \cite{39}.

\begin{center} \textbf{Топологии фиксированного множества} \end{center}

Всюду в пределах настоящего пункта фиксируем множество $X.$ Следуя \cite[(1.7.1)]{30} введем в рассмотрение семейство
\beq{3.1}
\pi[X]\!\triangleq\!\{\mathcal{X}\!\in\!\pp'(\pp(X))\mid(\zer\in\mathcal{X})  \&  (X\in \mathcal{X})  \&  (A\cap B\in \mathcal{X} \ \forall\,A\in\mathcal{X} \ \forall\,B\in\mathcal{X})\}.
\eeq

 Элементы семейства (\ref{3.1}) называем, следуя \cite[c.14]{37}, $\pi$-системами множества $X$ с <<нулем>> (пустое множество) и <<единицей>> (множество $X$). Заметим, что $\pp(X)\in\pi[X]$  и $\{\zer;X\}\in\pi(X).$  Если $\mathcal{X}\in\pi[X]$ и $\mathcal{Y}\in\pp(\mathcal{X}),$  то, в частности, $\mathcal{Y}\subset\pp(X)$ и, в частности, $\mathcal{Y}$ есть семейство (множеств), а потому определено (см.  (\ref{2.2})) множество-объединение $$\bigcup\limits_{Y\in\mathcal{Y}}Y\in\pp(X),$$ которое, вообще говоря, может не принадлежать $\mathcal{X}.$ Если же для всякого семейства $\mathcal{Y}\in\pp(\mathcal{X})$ $$\bigcup\limits_{Y\in\mathcal{Y}}Y\in\mathcal{X},$$
 то $\pi$-систему $\mathcal{X}$ называем \emph{топологией} на $X.$  Поскольку при этом $\zer\in\mathcal{X}$ в силу (\ref{3.1}), то в последнем случае достаточно рассматривать вариант $\mathcal{Y}\in\pp'(\mathcal{X}).$ С учетом этого
 \beq{3.1'}
 (\pr{top})[X]\triangleq\biggl\{\tau\in\pi[X]\mid\bigcup\limits_{G\in\mathcal{G}}G\in \tau \ \forall\,\mathcal{G}\in\pp'(\tau)\biggl\}
 \eeq
 есть множество всех топологий на $X.$ Если же $\tau\in(\pr{top})[X],$ то пару $(X,\tau)$ называем \emph{топологическим пространством} (ТП); ясно, что
\beq{3.2}
(\pp(X)\in(\pr{top})[X]) \ \& \ (\{\zer;X\}\in(\pr{top})[X])
\eeq
(дискретная и антидискретная топологии на $X$). Для определения топологий (в том числе <<более интересных>> в сравнении с (\ref{3.2})) широко используется конструкция на основе (топологических) баз. Чтобы ввести данную конструкцию, рассмотрим сначала одно вспомогательное определение: если $\mathcal{U}$\,--- семейство (множеств), то полагаем
\beq{3.3}
\{\cup\}(\mathcal{U})\triangleq\biggl\{\bigcup\limits_{V\in\mathcal{V}}V: \ \mathcal{V}\in\pp(\mathcal{U})\biggl\}
\eeq
(подробнее см. в \cite[замечание 7.4.1]{38}).
Из определения семейства (\ref{3.3}) следует \cite[с.\,65]{38},  что \\
a') \ $\bigcup\limits_{V\in \mathcal{V}}V\in\{\cup\}(\mathcal{U})$ при всяком выборе семейства $\mathcal{V}\in\pp(\mathcal{U});$ \\
b') \ для каждого множества  $W\in\{\cup\}(\mathcal{U})$ существует $\mathcal{W}\in\pp(\mathcal{U})$ такое, что  $$W=\bigcup\limits_{U\in\mathcal{W}}U.$$

Введем теперь в рассмотрение семейства, которые могут использоваться в качестве баз того или иного ТП. Последнее в такой редакции конструируется по заданному (и, как правило, более простому семейству). Упомянутые семейства\,--- топологические базы\,--- составляют множество
\begin{multline}\label{3.4}
(\pr{BAS})[X]\triangleq\biggl\{\mathcal{B}\in\pp(\pp(X))\mid\biggl(X=\bigcup\limits_{B\in\mathcal{B}}B\biggl)\ \& \\ \&
 \left(\forall\,B_1\in\mathcal{B} \ \forall\,B_2\in\mathcal{B} \ \forall\,x\in B_1\cap B_2 \ \exists\,B_3\in\mathcal{B}: (x\in B_3) \ \& \ (B_3\subset B_1\cap B_2)\right)\biggl\};
\end{multline}
в наиболее естественном случае $X\neq\zer$ согласно (\ref{3.4}) имеем, что
\beq{3.5}
(\pr{BAS})[X]\subset\pp'(\pp(X)).
\eeq

Итак (см. (\ref{3.5})), при $X\neq\zer$ и $\mathcal{B}\in\pr{BAS})[X]$ имеем, что $\mathcal{B}$\,--- непустое семейство п/м $X,$ для которого
\begin{multline}\label{3.6}
\biggl(X=\bigcup\limits_{B\in\mathcal{B}}B\biggl)\ \& \
 \biggl(\forall\,B_1\in\mathcal{B} \ \forall\,B_2\in\mathcal{B} \ \forall\, x\in B_1\cap B_2  \ \exists\,B_3\in\mathcal{B}: \\ (x\in B_3) \ \& \ (B_3\subset B_1\cap B_2)\biggl)
\end{multline}
(в (\ref{3.6}) мы <<расшифровали>> (\ref{3.4})).

С помощью баз, определяемых в (\ref{3.4})--(\ref{3.6}), по универсальному правилу конструируются топологии: при $\mathcal{B}\!\in\!(\pr{BAS})[X]$
\begin{multline}\label{3.7}
\{\cup\}(\mathcal{B})=\biggl\{\bigcup\limits_{B\in\beta}B: \ \beta\in\pp(\mathcal{B})\biggl\}=\\=\left\{G\in\pp(X)\mid\forall\,x\in G \ \exists\, B\in\mathcal{\mathcal{B}}: (x\in B) \ \& \ (B\subset G)\right\}\in(\pr{top})[X]
\end{multline}
(читателю предлагается проверить (\ref{3.7}) самостоятельно) и при этом
\beq{3.8}
\mathcal{B}\subset\{\cup\}(\mathcal{B}).
\eeq
Итак (см. (\ref{3.8})), все множества\,--- элементы базы  $\mathcal{B}$\,--- открыты в топологии, порождаемой этой базой по правилу (\ref{3.7}). В дальнейшем мы рассмотрим конкретные примеры, в которых база задается достаточно просто, а схема на основе (\ref{3.7}) позволяет в значительной степени <<разглядеть>> порождаемую этой базой топологию. В свою очередь, с заданной топологией, а, точнее, с заданным ТП, можно связать целое семейство баз, порождающих данную топологию: если $\tau\in(\pr{top})[X],$ то
\beq{3.9}
(\tau-\pr{BAS})_0[X]\triangleq\{\mathcal{B}\in(\pr{BAS})[X]\mid \tau=\{\cup\}(\mathcal{B})\}\in\pp'((\pr{BAS})[X])
\eeq
(заметим, кстати, что в силу (\ref{3.1}), (\ref{3.1'}) и (\ref{3.7}) $\tau\in(\tau-\pr{BAS})_0[X];$  данный пример, однако, обычно  малоинтересен). Семейства\,--- элементы (\ref{3.9})\,--- базы фиксированной топологии $\tau;$ говорят также\,--- базы ТП $(X,\tau).$ Отметим полезное
\begin{proposition}
Если $X\neq\zer$ и $\tau\in(\pr{top})[X],$ то
\beq{3.10}
(\tau-\pr{BAS})_0[X]=\biggl\{\mathcal{B}\in\pp'(\tau)\mid\forall\,G\in\tau\setminus\{\zer\} \ \exists\,\beta\in\pp'(\mathcal{B}): \ G=\bigcup\limits_{B\in\beta}B\biggl\}.
\eeq
\end{proposition}
Д о к а з а т е л ь с т в о. Обозначим через $\mathbb{T}$ множество в правой части  (\ref{3.10}). Сравним $(\tau-\pr{BAS})_0[X]$ и $\mathbb{T}$ (имеется в виду сравнение двух упомянутых множеств). Пусть сначала $\mathfrak{B}\in(\tau-\pr{BAS})_0[X].$  Тогда в силу (\ref{3.5}) и (\ref{3.9})
\beq{3.11}
\mathfrak{B}\in\pp'(\pp(X))
\eeq
($\mathfrak{B}$\,---  непустое семейство п/м $X$)  и
\beq{3.12}
\tau=\{\cup\}(\mathfrak{B}).
\eeq
Если $\mathbb{G}\in\tau\setminus\{\zer\},$ то согласно (\ref{3.3}) и (\ref{3.12}) для некоторого $\mathcal{G}\in\pp(\mathfrak{B})$ имеем равенство
\beq{3.13}
\mathbb{G}=\bigcup\limits_{B\in\mathcal{G}}B;
\eeq
поскольку $\mathbb{G}\neq\zer,$ то и  $\mathcal{G}\neq\zer$ (в противном случае из (\ref{3.13}) следовало бы равенство $\mathbb{G}=\zer,$  что невозможно). Итак, $$\mathcal{G}\in\pp'(\mathfrak{B}).$$
Так как  $\mathbb{G}$ был произвольным, установлено, что
\beq{3.14}
\forall\, G\in\tau\setminus\{\zer\} \ \exists\,\beta\in\pp'(\mathfrak{B}): \ G=\bigcup\limits_{B\in \beta}B.
\eeq
С учетом (\ref{3.11}) и (\ref{3.14}) получаем, что $\mathfrak{B}\in\mathbb{T},$ чем завершается обоснование вложения
\beq{3.15}
(\tau-\pr{BAS})_0[X]\subset\mathbb{T}.
\eeq

Выберем произвольно $\mathcal{T}\in\mathbb{T}.$  Тогда
\beq{3.16}
\mathcal{T}\in\pp'(\tau)
\eeq
и, кроме того, справедливо свойство
\beq{3.17}
\forall\,G\in\tau\setminus\{\zer\} \ \exists\,\beta\in\pp'(\mathcal{T}): \ G=\bigcup\limits_{B\in \beta}B.
\eeq
Из (\ref{3.16}) имеем, в частности, что $\mathcal{T}$\,--- непустое семейство п/м $X.$  Покажем, что
\begin{equation}\label{3.18}
\mathcal{T}\in(\pr{BAS})[X].
\eeq
Для этого сначала отметим, что $X\in\tau\setminus\{\zer\},$ а потому в силу  (\ref{3.17})
\beq{3.19}
X=\bigcup\limits_{B\in \beta_0}B,
\eeq
где $\beta_0\in\pp'(\mathcal{T}).$  Из (\ref{3.19}) имеем, стало быть, цепочку вложений $$X\subset\bigcup\limits_{B\in \mathcal{T}}B\subset X.$$
(см. (\ref{3.16})), означающую справедливость равенства
\beq{3.20}
X=\bigcup\limits_{B\in \mathcal{T}}B.
\eeq

Пусть  $T_1\in\mathcal{T}, T_2\in\mathcal{T}$ и $x_*\in T_1\cap T_2.$ С учетом (\ref{3.16}) имеем $T_1\in\tau$ и $T_2\in\tau,$ а потому  (см. (\ref{3.1}), (\ref{3.1'}))
\beq{3.21}
T_1\cap T_2\in\tau\setminus\{\zer\}.
\eeq
С учетом (\ref{3.17}) и (\ref{3.21}) подберем $\beta_*\in\pp'(\mathcal{T}),$ для которого
\beq{3.22}
T_1\cap T_2=\bigcup\limits_{B\in\beta_*}B.
\eeq
В этом случае $x_*\in\bigcup\limits_{B\in\beta_*}B,$ а потому для некоторого $B_*\in\beta_*$ имеет место $x_*\in B_*.$ Поскольку, в частности, $B_*\in\mathcal{T},$  имеем, что $$\exists\,B\in\mathcal{T}: (x_*\in B) \ \& \ (B\subset T_1\cap T_2)$$ (учитываем (\ref{3.22})). Коль скоро $T_1      , T_2$ и $x_*$ выбирались произвольно, установлено, что
\beq{3.23}
\forall\,B_1\in\mathcal{T} \ \forall\,B_2\in\mathcal{T} \ \forall\,x\in B_1\cap B_2 \ \exists\,B_3\in\mathcal{T}:(x\in B_3)\ \& \ (B_3\subset B_1\cap B_2).
\eeq
Согласно  (\ref{3.16})  $\mathcal{T}\in\pp'(\pp(X)).$  Поэтому (см. (\ref{3.4}), (\ref{3.20}), (\ref{3.23})) $\mathcal{T}$\,--- топологическая база (см. (\ref{3.18})). При этом
\beq{3.24}
\tau=\{\cup\}(\mathcal{T}).
\eeq
В самом деле, пусть $\mathbf{G}\in \tau.$  Тогда в силу (\ref{3.17})
$$(\mathbf{G}\neq\zer)\Rightarrow\biggl(\exists\,\beta\in\pp'(\mathcal{T}):\mathbf{G}=\bigcup\limits_{B\in\beta}B\biggl).$$
Получаем, следовательно (см. (\ref{3.3})), импликацию
\beq{3.25}
(\mathbf{G}\neq\zer)\Rightarrow(\mathbf{G}\in\{\cup\}(\mathcal{T})).
\eeq
Если же $\mathbf{G}=\zer,$ то для $\zer\in\pp(\mathcal{T})$ имеем равенство
$\mathbf{G}=\bigcup\limits_{B\in\zer}B,$  а тогда $\mathbf{G}\in\{\cup\}(\mathcal{T});$ см. (\ref{3.3}). Установлена импликация
$$(\mathbf{G}=\zer)\Rightarrow(\mathbf{G}\in\{\cup\}(\mathcal{T})).$$
С учетом (\ref{3.25}) получаем, что $\mathbf{G}\in\{\cup\}(\mathcal{T})$ во всех возможных случаях. Следовательно, установлено вложение
\beq{3.26}
\tau\subset\{\cup\}(\mathcal{T}).
\eeq
Выберем произвольно $\Gamma\in\{\cup\}(\mathcal{T}).$ В силу (\ref{3.1}) и (\ref{3.1'}) имеем, что
\beq{3.27}
(\Gamma=\zer)\Rightarrow(\Gamma\in\tau).
\eeq
Пусть $\Gamma\neq\zer,$ т.е. $\Gamma\in\{\cup\}(\mathcal{T})\setminus\{\zer\}.$  Из (\ref{3.3}) получаем  для некоторого $\gamma\in\pp'(\mathcal{T})$
\beq{3.28}
\Gamma=\bigcup\limits_{B\in\gamma}B.
\eeq
Итак, $\gamma$\,--- непустое подсемейство $\mathcal{T}.$
Напомним, что  в силу (\ref{3.16}) $$\mathcal{T}\subset\tau.$$
Тогда,  в частности, $\gamma$\,--- непустое подсемейство $\tau,$ т.е. $\gamma\in\pp'(\tau),$ а тогда согласно (\ref{3.1'}) $$\bigcup\limits_{B\in\gamma}B\in\tau.$$
Используя  (\ref{3.28}), получаем, что $\Gamma\in\tau$ и в рассматриваемом случае непустого множества $\Gamma.$  Установлена  импликация
$$(\Gamma\neq\zer)\Rightarrow(\Gamma\in\tau).$$
С учетом  (\ref{3.27}) имеем, что во всех во всех возможных случаях справедливо включение $\Gamma\in\tau.$ Таким образом,  $\{\cup\}(\mathcal{T})\subset\tau,$ откуда (см.  (\ref{3.26})) вытекает равенство $\tau=\{\cup\}(\mathcal{T}).$
Из  (\ref{3.9}),   (\ref{3.18}) получаем  (см.  (\ref{3.24})) включение $$\mathcal{T}\in(\tau-\pr{BAS})_0[X].$$
Поскольку выбор $\mathcal{T}$ был произвольным, установлено, что $\mathbb{T}\subset(\tau-\pr{BAS})_0[X],$ откуда согласно  (\ref{3.15}) получаем требуемое равенство $(\tau-\pr{BAS})_0[X]=\mathbb{T}. \hfill\square$

Таким образом, в случае, когда $X$\,--- непустое множество, а $\tau$\,--- топология $X,$   базы ТП  $(X,\tau)$ суть непустые семейства открытых (в $(X,\tau)$) множеств, реализующие все непустые открытые множества в виде объединений своих непустых подсемейств.

\begin{center} \textbf{Окрестности} \end{center}

Если $\tau\in(\pr{top})[X]$ и $x\in X,$ то полагаем, что
\beq{3.29}
N_\tau^0(x)\triangleq\{G\in\tau\mid x\in G\}
\eeq
(получая непустое семейство открытых в ТП $(X,\tau)$ множеств; всегда $X\in N_\tau^0(x)$) и, кроме того,
\beq{3.30}
N_\tau(x)\triangleq\{H\in\pp(X)\mid\exists\,G\in N_\tau^0(x): G\subset H\};
\eeq
связь семейств (\ref{3.29}), (\ref{3.30})  определяется очевидным равенством
\beq{3.30'}
N_\tau^0(x)=N_\tau(x)\cap\tau.
\eeq

Отметим следующее полезное свойство \cite[гл.\,I]{45}:
\beq{3.1000}
\tau=\{G\in\pp(X)\mid G\in N_\tau(x) \ \forall\,x\in G\} \ \forall\,\tau\in(\pr{top})[X].
\eeq
%\end{document}
В целях полноты изложения приведем доказательство свойства (\ref{3.1000}), фиксируя топологию $\tau\in(\pr{top})[X].$ Обозначим через $\mathcal{T}$ семейство в правой части (\ref{3.1000}).

Если $\mathbb{G}\in\tau,$ то $\mathbb{G}\in N_\tau^0(x) \ \forall\,x\in\mathbb{G}$ (учитываем, что $\mathbb{G}\subset X$). В частности имеем с учетом (\ref{3.30'}), что $$\mathbb{G}\in N_\tau(x) \ \forall\,x\in\mathbb{G}.$$
Последнее означает, что $\mathbb{G}\in\mathcal{T},$ чем и завершается проверка вложения
\beq{3.1001}
\tau\subset\mathcal{T}.
\eeq
Выберем произвольно $T\in\mathcal{T}.$ Тогда, в частности, $T\in\pp(X).$ При этом
\beq{3.1002}
T\in N_\tau(x) \ \forall\,x\in T.
\eeq
Отметим сразу, что (см. (\ref{3.1}), (\ref{3.1'}))
\beq{3.1003}
(T=\zer)\Rightarrow (T\in\tau).
\eeq
%\end{document}
Пусть теперь $T\neq\zer.$  Согласно  (\ref{3.1002}) имеем, что
\beq{3.1004}
\mathcal{J}_x\triangleq\{G\in N_\tau^0(x)\mid G\subset T\}\in\pp'(\tau) \ \forall\,x\in T.
\eeq
%\end{document}
Как следствие, получаем непустое семейство открытых множеств
\beq{3.1005}
\mathcal{J}\triangleq\bigcup\limits_{x\in T}\mathcal{J}_x\in\pp'(\tau),\eeq
а потому в силу (\ref{3.1'}) имеем свойство
\beq{3.1006}
\mathbf{G}\triangleq\bigcup\limits_{G\in\mathcal{J}}G\in\tau.\eeq
Из (\ref{3.1004}) и (\ref{3.1005}) вытекает, что $$G\subset T \ \ \forall\,G\in\mathcal{J}.$$
Согласно  (\ref{3.1006}) имеем теперь с очевидностью, что
\beq{3.1007}
\mathbf{G}\subset T.
\eeq
Пусть $x_*\in T.$ Тогда в силу (\ref{3.1004}) получаем, что $$\mathcal{J}_{x_*}=\{G\in N_\tau^0(x_*)\mid G\subset T\}\neq\zer.$$
Теперь выберем окрестность $G_*\in\mathcal{J}_{x_*},$ получая для $G_*\in N_\tau^0(x_*)$ вложение $$G_*\subset T.$$

Вместе с тем, из (\ref{3.29}) следует, что $x_*\in G_*.$ При этом (см. (\ref{3.1005})) справедливо включение $G_*\in\mathcal{J},$ т.к. $\mathcal{J}_{x_*}\subset\mathcal{J}.$  Используя (\ref{3.1006}), получаем, что    $G_*\subset\mathbf{G},$ а тогда $x_*\in\mathbf{G}.$ Итак, $T\subset\mathbf{G},$ а потому (см. (\ref{3.1007})) $T=\mathbf{G}.$ В силу (\ref{3.1006}) имеем, что $T\in\tau$  и при $T\neq\zer,$  т.е. $(T\neq\zer)\Rightarrow(T\in\tau).$
Применяя (\ref{3.1003}),  имеем, что $T\in\tau$ во всех возможных случаях. Поскольку выбор $T$ был произвольным, установлено, что $\mathcal{T}\subset\tau,$ а тогда (см. (\ref{3.1001})) $\tau=\mathcal{T},$ что  и требовалось доказать. $\hfill\Box$

\begin{center} \textbf{Локальные базы (фундаментальные системы окрестностей)} \end{center}

Если $\tau\in(\pr{top})[X],$ $x\in X$ и $\mathcal{B}$\,--- подсемейство $N_\tau(x),$ то назовем $\mathcal{B}$ локальной базой ТП $(X,\tau)$  (\emph{фундаментальной системой окрестностей}) в точке $x$, если
\beq{3.31}
\forall\,A\in N_\tau(x) \ \exists\, B\in\mathcal{B}: \ B\subset A.
\eeq
 Разумеется, (\ref{3.31}) не определяет локальную базу однозначно; имеет смысл говорить о множестве всех локальных баз для той или иной фиксированной точки ТП: при $\tau\in(\pr{top})[X]$  и $x\in X$
 $$(x-\pr{bas})[\tau]\triangleq\{\mathcal{B}\in\pp(N_\tau(x))\mid\forall\,A\in N_\tau(x) \ \exists\,B\in\mathcal{B}: \ B\subset A\}\in\pp'(\pp(N_\tau(x)))$$
 есть множество всех локальных баз $(X,\tau)$ в точке  $x;$  разумеется, $N_\tau(x)\in (x-\pr{bas})[\tau]$ и  $N_\tau^0(x)\in (x-\pr{bas})[\tau]\cap\pp(\tau).$

\begin{center} \textbf{Замыкание множества в топологическом пространстве; замкнутые множества} \end{center}

Если $\tau\in(\pr{top})[X]$ и $A\in\pp(X),$  то множество
\beq{3.32}
\pr{cl}(A,\tau)\triangleq\{x\in X\mid A\cap H\neq\zer \ \forall\,H\in N_\tau(x)\}\in\pp(X)
\eeq
 называется \emph{замыканием} $A$ \emph{в ТП}  $(X,\tau).$ Из (\ref{3.30}) и (\ref{3.30'}) вытекает, что
 \beq{3.33}
 \pr{cl}(A,\tau)=\{x\in X\mid A\cap G\neq\zer \ \forall\,G\in N_\tau^0(x)\} \ \forall\,\tau\in(\pr{top})[X] \ \forall\,A\in\pp(X).
 \eeq
В (\ref{3.32}), (\ref{3.33}) реализуется <<окрестностное>>  представление замыкания.

Полезным оказывается другое по форме представление, а именно, представление, связанное с использованием замкнутых множеств:\\ если $\tau\in(\pr{top})[X],$ то
\beq{3.34}
\mathcal{F}[\tau]\triangleq\mathbf{C}_X[\tau]=\{X\setminus G: \ G\in\tau\}=\{F\in\pp(X)\mid X\setminus F\in\tau\}\in\pp'(\pp(X))
\eeq
есть семейство всех замкнутых в $(X,\tau)$  п/м множества $X.$ В (\ref{3.34}) учтено (см. \cite[(1.1.6)]{30}) свойство двойного дополнения:
$$X\setminus(X\setminus H)=H \ \forall\,H\in\pp(X);$$
 в частности, $G=X\setminus(X\setminus G)$ при $G\in\tau$.

\begin{zam} Условимся об одном определении общего характера.
Всюду в дальнейшем для каждого семейства $\mathcal{A}$ и множества $B$ полагаем, что
\beq{3.35}
[\mathcal{A}](B)\triangleq\{A\in\mathcal{A}\mid B\subset A\};
\eeq
в качестве $\mathcal{A}$ можно, разумеется, использовать любое подсемейство $\pp(X).$ \end{zam}

Заметим, что в (\ref{3.35}) мы можем рассматривать тот (важный для нас)  случай, когда $\mathcal{A}=\mathcal{F}[\tau],$ где $\tau\in(\pr{top})[X],$  и $B\in\pp(X),$ получая при этом
\beq{3.36}
[\mathcal{F}[\tau]](B)=\{F\in\mathcal{F}[\tau]\mid B\subset F\}\in\pp'(\mathcal{F}[\tau])
\eeq
(легко заметить, что $X\in[\mathcal{F}[\tau]](B)$). Введено семейство всех замкнутых в ТП $(X,\tau)$ множеств, содержащих $B.$
С упомянутым семейством естественным образом связывается замыкание множества  $B:$ замыкание есть наименьший по включению элемент семейства (\ref{3.36}).
Точнее: если $\tau\in(\pr{top})[X]$ и $B\in\pp(X),$ то
\beq{3.37}
(\pr{cl}(B,\tau)\in[\mathcal{F}[\tau]](B))\ \& (\pr{cl}(B,\tau)\subset F \ \forall\,F\in[\mathcal{F}[\tau]](B)).
\eeq
Дополняем (\ref{3.37}) следующим очевидным свойством: если $\tau\in(\pr{top})[X],$ то
\beq{3.38}
\mathcal{F}[\tau]=\{F\in\pp(X)\mid F=\pr{cl}(F,\tau)\}
\eeq
 (равенство (\ref{3.38}) предлагается проверить читателю самостоятельно, используя (\ref{3.37})); в частности,
 $$\pr{cl}(F,\tau)=F \ \forall\,F\in\mathcal{F}[\tau].$$

\begin{center} \textbf{Некоторые аксиомы отделимости} \end{center}

Важную роль в дальнейшем изложении играют аксиомы отделимости. Ограничимся сейчас рассмотрением лишь некоторых их них. Итак, если $\tau\in(\pr{top})[X],$ то в виде пары $(X,\tau)$ имеем ТП. Называем  $(X,\tau)$\\
 1) \ $T_1$-пространством, если $\forall\,x_1\in X \ \forall\,x_2\in X\setminus\{x_1\} \ \exists\,H\in N_\tau(x_1): \ x_2\notin H;$\\
 2) \ $T_2$-пространством, если $\forall\,x_1\!\in\! X \forall\,x_2\!\in\! X\setminus\{x_1\}  \exists\,H_1\!\in\! N_\tau(x_1)  \exists\,H_2\!\in\! N_\tau(x_2): H_1\cap H_2\!=\!~\zer;$\\
 3) \ $T_3$-пространством \cite[c.\,71]{39}, \cite{40} если $(X,\tau)$ есть $T_1$-пространство и при этом
 $$N_\tau(x)\cap\mathcal{F}[\tau]\in(x-\pr{bas})[\tau] \ \forall\,x\in X$$
 (последнее означает справедливость следующего свойства (см. (\ref{3.31})): при всяком выборе точки $x\in X$ и окрестности $H\in N_\tau(x)$ найдется замкнутая в ТП $(X,\tau)$ окрестность $F\in N_\tau(x)$ со свойством $F\subset H$).
Прочие аксиомы отделимости здесь не рассматриваем, отсылая читателя к \cite{39}.

 \newpage

\begin{center} \textbf{Компактность} \end{center}

Если $\tau\in(\pr{top})[X],$ то ТП $(X,\tau)$ называем \emph{компактным}, если $\forall\,\mathcal{G}\in\pp'(\tau)$
\beq{3.39}
\biggl(X=\bigcup\limits_{G\in\mathcal{G}}G\biggl)\Rightarrow\biggl(\exists\,\mathcal{K}\in\pr{Fin}(\mathcal{G}): X=\bigcup\limits_{G\in\mathcal{K}}G\biggl).
\eeq
Иными словами, компактность ТП $(X,\tau)$ есть следующее свойство: из каждого покрытия $X$  открытыми множествами можно <<извлечь>> конечное подпокрытие.

Отметим важное эквивалентное представление. Сначала (см. \cite[п.\,6]{38} ) введем понятие (непустого) центрированного семейства: если $\mathcal{U}$\,--- непустое семейство (множеств), то называем его \emph{центрированным}, если $$\bigcap\limits_{U\in\mathcal{K}}U\neq\zer \ \forall\,\mathcal{K}\in\pr{Fin}(\mathcal{U}).$$

Таким образом, для произвольного непустого семейства $\mathcal{T}$ в виде $$\mathbb{Z}[\mathcal{T}]\triangleq\biggl\{\mathcal{U}\in\pp'(\mathcal{T})\mid \bigcap\limits_{U\in\mathcal{K}}U\neq\zer \ \forall\,\mathcal{K}\in\pr{Fin}(\mathcal{U})\biggl\}$$
имеем семейство всех центрированных подсемейств $\mathcal{T}.$ Упомянутое важное  свойство (эквивалентное определение компактности) состоит в следующем: ТП $(X,\tau),$ где $\tau\in(\pr{top})[X],$ компактно тогда и только тогда, когда
\beq{3.40}
\bigcap\limits_{F\in\mathbb{F}}F\neq\zer \ \forall\,\mathbb{F}\in\mathbb{Z}[\mathfrak{F}[\tau]].
\eeq

Компактные $T_2$-пространства называют \emph{компактами}.

\begin{center} \textbf{Множества, компактные в заданном ТП} \end{center}

Наряду с компактностью исходного ТП во многих случаях важно рассматривать множества в произвольном ТП, которым также удается <<приписать>>  свойство, подобное в значительной степени свойству (\ref{3.39}): если  $\tau\in(\pr{top})[X]$ (т.е. если задано ТП  $(X,\tau)$ с <<единицей>> $X$), то полагаем, что
\begin{multline}\label{3.41}
(\tau-\pr{comp})[X]\triangleq\biggl\{K\in\pp(X)\mid \forall\,\mathcal{G}\in\pp'(\tau) \ (K\subset \bigcup\limits_{G\in\mathcal{G}}G)\Rightarrow\\ \Rightarrow(\exists\,\mathcal{K}\in\pr{Fin}(\mathcal{G}): K\subset\bigcup\limits_{G\in\mathcal{K}}G)\biggl\}.
\end{multline}

Позднее рассмотрим эквивалентное представление множества (\ref{3.41}) в терминах подпространств исходного ТП $(X,\tau).$

\begin{center} \textbf{Простейшие примеры} \end{center}
$1^*.$  \ Если $X$\,--- непустое конечное множество, то ТП $(X,\tau)$ компактно при всяком выборе топологии $\tau\in(\pr{top})[X].$  В самом деле, поскольку $X$ конечно, то можно указать $n\in\mathbb{N}$ и кортеж $$(x_i)_{i\in\overline{1,n}}:\overline{1,n}\rightarrow X,$$
для которых $X=\{x_i:i\in\overline{1,n}\}.$
Пусть теперь $\mathcal{G}\in\pp'(\tau)$ обладает свойством
\beq{3.42}
X=\bigcup\limits_{G\in\mathcal{G}}G.
\eeq
В этом случае $\forall\,j\in\overline{1,n} \ \exists\,G\in\mathcal{G}: \ x_j\in G.$ Поэтому для некоторого кортежа
$(\mathbb{G}_i)_{i\in\overline{1,n}}:\overline{1,n}\rightarrow \mathcal{G}$ справедливо следующее свойство: $x_j\in\mathbb{G}_j \ \forall\,j\in\overline{1,n}.$ Тогда $$\widetilde{\mathcal{K}}\triangleq\left\{\mathbb{G}_i:i\in\overline{1,n}\right\}\in\pr{Fin}(\mathcal{G})$$
и при этом справедлива цепочка равенств $$\bigcup\limits_{G\in\widetilde{\mathcal{K}}}G=\bigcup\limits_{i=1}^n\mathbb{G}_i=X;$$
мы учитываем, что $\mathbb{G}_j\subset X$ по выбору $\mathcal{G}.$ Поэтому (при условии (\ref{3.42}))
$$\exists\,\mathcal{K}\in\pr{Fin}(\mathcal{G}): \ X=\bigcup\limits_{G\in\mathcal{K}}G.$$
Коль скоро выбор семейства $\mathcal{G}$ открытых в $(X,\tau)$ п/м $X$ со свойством (\ref{3.42}) был произвольным, из (\ref{3.39}) вытекает требуемое свойство компактности.

$2^*.$ Если $\mathbb{K}\in\pr{Fin}(X),$ где $X$\,--- произвольное непустое множество, то $\mathbb{K}\in(\tau-\pr{comp})[X] \ \forall\,\tau\in(\pr{top})[X]$ (доказательство повторяет рассуждения в $1^*$).

\newpage

\section{Элементы топологии, 2}\setcounter{equation}{0} \setcounter{proposition}{0}\setcounter{zam}{0}

\ \ \ \ \ В настоящем разделе обсудим <<взаимодействия>> различных, вообще говоря, ТП. Данное рассмотрение начнем с обсуждения понятия подпространства. Затем совсем кратко коснемся вопросов, связанных с непрерывностью отображений.

\begin{center} \textbf{Подпространства} \end{center}

В  данном подразделе фиксируем сначала  множество $Y$ и топологию $\tau\in(\pr{top})[Y],$ получая ТП $(Y,\tau).$  Если теперь $Z\in\pp(Y),$ то определено семейство
\beq{4.0}
\tau|_Z=\{G\cap Z: \ G\in\tau\}\in\pp'(\pp(Z)),
\eeq
 являющееся, как легко проверить, топологией на множестве $Z:$
\beq{4.1}
\tau|_Z\in(\pr{top})[Z].
\eeq
Топологию (\ref{4.1}) называют \emph{индуцированной на множестве} $Z$ из ТП $(Y,\tau).$  Тогда определены семейства
$\mathcal{F}[\tau|_Z]=\{Z\setminus G: \ G\in\tau|_Z\}\in\pp'(\pp(Z)),$ \ $\mathcal{F}[\tau]|_Z=\{F\cap Z: \ F\in\mathcal{F}[\tau]\}.$
Более того, как легко видеть, данные семейства совпадают:
\beq{4.2}
\mathcal{F}[\tau|_Z]=\mathcal{F}[\tau]|_Z.
\eeq
Проверку весьма очевидного свойства (\ref{4.2}) предоставляем читателю. Еще одно полезное (и связанное с (\ref{4.2})) свойство касается операции замыкания (см. (\ref{3.37})): если $A\in\pp(Z),$ то, поскольку при этом также имеет место $A\in\pp(Y),$  определены замыкания
\beq{4.3}
\pr{cl}(A,\tau|_Z), \ \pr{cl}(A,\tau);
\eeq
в первом случае имеем замыкание в подпространстве $(Z,\tau|_Z)$ исходного ТП $(Y,\tau),$ а во втором\,--- замыкание в самом этом ТП. Связь замыканий
(\ref{4.3}) определяется просто:
\beq{4.4}
\pr{cl}(A,\tau|_Z)=\pr{cl}(A,\tau)\cap Z
\eeq
(напомним, что в (\ref{4.4}) $A$ есть произвольное п/м $Z$).

Отметим также естественную связь окрестностей в исходном ТП и в его подпространстве (проверку приводимых ниже свойств предоставляем читателю). Если $z\in Z,$ то, с одной стороны, определены семейства $N^0_{\tau|_Z}(z)$  и $N_{\tau|_Z}(z)$ (относительных) окрестностей $z,$ а с другой~--- семейства $N^0_\tau(z)$ и $N_\tau(z).$ Поскольку, в частности, $z\in Y,$ имеем
\beq{4.5}
N^0_{\tau|_Z}(z)=N^0_{\tau}(z)|_Z=\{G\cap Z: \ G\in N^0_{\tau}(z)\},
\eeq
\beq{4.6}
N_{\tau|_Z}(z)=N_{\tau}(z)|_Z=\{H\cap Z: \ H\in N_{\tau}(z)\}.
\eeq

Свойства   (\ref{4.5}), (\ref{4.6}) будут полезны при описании сходимости направленностей в исходном ТП и в его подпространстве.
В построениях, связанных с (\ref{4.0})~--~(\ref{4.4}), множество $Z\in\pp(Y)$ полагалось произвольным. Полезно отметить специальные случаи, когда $Z$ открыто, либо замкнуто в ТП $(Y,\tau):$  легко понять, что
$$(Z\in \tau)\Rightarrow\left(\tau|_Z=\{G\in \tau\mid G\subset Z\}\right),$$
$$(Z\in\mathcal{F}_\tau)\Rightarrow\left(\mathcal{F}[\tau|_Z]=\{F\in\mathcal{F}_\tau\mid F\subset Z\}\right).$$
Возвращаясь к (\ref{3.41}), отметим важное эквивалентное представление: если $(X,\tau)$ есть ТП и $K\in\pp(X),$ то $K\in(\tau-\pr{comp})[X]$  тогда и только тогда, когда компактно ТП
\beq{4.7}
\left(K,\tau|_K\right).
\eeq
Иными словами, $K$ компактно в ТП $(X,\tau)$ тогда и только тогда, когда подпространство (\ref{4.7})  пространства $(X,\tau)$ компактно как ТП.

В связи с аксиомами отделимости отметим только одно свойство, связанное с переходом к подпространству:  если $(X,\tau)$ есть $T_2$-пространство и
 $\mathbb{Y}\in\pp(X),$ то
 \beq{4.8}
 (\mathbb{Y},\tau|_\mathbb{Y})
 \eeq
  также является $T_2$-пространством (доказательство легко следует из (\ref{4.6}) и определения $T_2$-пространства).
 \begin{zam}
 Проверим последнее свойство. Пусть $y_1\in \mathbb{Y}$ и $y_2\in \mathbb{Y}\setminus\{y_1\}.$ Тогда, в частности, $y_1\in X$ и $y_2\in X\setminus\{y_1\},$ т.к. $\mathbb{Y}\subset X.$ Коль скоро $(X,\tau)$ есть $T_2$-пространство, то для некоторых $\mathbb{H}_1\in N_\tau(y_1)$ и $\mathbb{H}_2\in N_\tau(y_2)$
 \beq{4.9}
 \mathbb{H}_1\cap\mathbb{H}_2=\zer.
 \eeq

 При этом (см. (\ref{4.6})) $\mathbb{H}_1\cap \mathbb{Y}\in N_{\tau|_\mathbb{Y}}(y_1)$ и $\mathbb{H}_2\cap \mathbb{Y}\in N_{\tau|_\mathbb{Y}}(y_2);$ кроме того, из (\ref{4.9}) вытекает, что $$(\mathbb{H}_1\cap \mathbb{Y})\cap(\mathbb{H}_2\cap \mathbb{Y})=\zer.$$
 Поскольку выбор $y_1$ и $y_2$ был произвольным, установлено следующее свойство ТП (\ref{4.8}): $\forall\,y^{(1)}\in \mathbb{Y} \ \forall\,y^{(2)}\in \mathbb{Y}\setminus\{y^{(1)}\} \ \exists\,H^{(1)}\in N_{\tau|_\mathbb{Y}}(y^{(1)})$  $\exists\,H^{(2)}\in N_{\tau|_\mathbb{Y}}(y^{(2)}):$
  $$H^{(1)}\cap H^{(2)}=\zer.$$
  В силу 2) получили, что (\ref{4.8}) есть $T_2$-пространство.\end{zam}
  Напомним одно важное свойство $T_2$-пространств: если $(X,\tau)$ есть \\ $T_2$-пространство, то (см. (\ref{3.34}))
 \beq{4.10}
 (\tau-\pr{comp})[X]\subset\mathcal{F}[\tau].
 \eeq

 Итак, в $T_2$-пространствах (т.е. в хаусдорфовых ТП) компактные множества замкнуты.

 \begin{center} \textbf{Сходимость в топологическом пространстве} \end{center}

 В построениях основной части важную роль играет свойство сходимости (последовательностей, направленностей) в ТП. В настоящем пункте кратко обсудим упомянутое свойство. Начнем с более простого понятия сходимости последовательностей (полагаем, что $(X,\tau)$ есть фиксированное ТП).

 Итак, пусть $(x_i)_{i\in\mathbb{N}}:\mathbb{N}\rightarrow X$ и $x\in X;$ полагаем, что
 \beq{4.11}
 \left((x_i)_{i\in\mathbb{N}}\stackrel{\tau}{\rightarrow} x\right)\stackrel{\pr{def}}{\Leftrightarrow}\left(\forall\,H\in N_\tau(x) \ \exists\,n\in\mathbb{N}:x_j\in H \ \forall\,j\in\overrightarrow{n,\infty}\right);
 \eeq
 в случае, когда выполняется условие, определенное в  (\ref{4.11}), говорим, что \emph{последовательность} $(x_i)_{i\in\mathbb{N}}$ \emph{сходится к} $x$ \emph{в ТП} $(X,\tau).$

 Имея (\ref{4.11}), мы рассмотрим естественное обобщение, а именно, так называемую \emph{сходимость по Мору-Смиту}, которая будет в дальнейшем очень важна, поскольку нужные варианты топологий не описываются в терминах сходящихся последовательностей. Для использования упомянутой сходимости нам потребуются прежде всего понятия направления и направленности в том или ином множестве. В свою очередь, направление есть частный случай (бинарного) отношения, именуемого предпорядком.

 Если $\mathbb{D}$ есть непустое множество и дано отношение $\rho\in\pp(\mathbb{D}\times\mathbb{D}),$ то $\forall\,x\in\mathbb{D} \ \forall\,y\in\mathbb{D}$
 \beq{4.12}
 (x\rho y)\stackrel{\pr{def}}\Leftrightarrow((x,y)\in \rho).
 \eeq
 Поэтому, учитывая (\ref{4.12}), будем вместо условия $(x,y)\in \rho$ использовать запись $x\rho y.$
Отношение $\rho\in\pp(\mathbb{D}\times\mathbb{D})$ (на непустом множестве $\mathbb{D}$) называется \emph{предпорядком}, если
$$(x\rho x \ \forall\,x\in\mathbb{D})\ \& \ (\forall\,x\in\mathbb{D}\ \forall\,y\in\mathbb{D} \ \forall\,z\in\mathbb{D} \ \ ((x\rho y)\ \& \ (y\rho z))\Rightarrow(x\rho z)).$$

Через $(\pr{Ord})[\mathbb{D}]$ обозначаем множество всех предпорядков на непустом множестве $\mathbb{D}:$
\begin{multline}\label{4.12'}
(\pr{Ord})[\mathbb{D}]\triangleq \{\zeta\in\pp(\mathbb{D}\times\mathbb{D})\mid (x\zeta x \ \  \forall\,x\in\mathbb{D}) \ \& \ (\forall\,x\in\mathbb{D}\ \forall\,y\in\mathbb{D}  \\ \forall\,z\in\mathbb{D} \ \ ((x\zeta y)\ \& \  (y\zeta z))\Rightarrow(x\zeta z))\}.
\end{multline}

Если $\rho\in(\pr{Ord})[\mathbb{D}],$ то пару $(\mathbb{D},\rho)$ называем частично упорядоченным множеством (ЧУМ). Среди всевозможных ЧУМ выделяем направленные множества (НМ). Для введения НМ определяем сначала понятие направления: предпорядок $\rho\in(\pr{Ord})[\mathbb{D}]$ (см. (\ref{4.12'})) называем направлением на $\mathbb{D},$ если $\forall\,x\in\mathbb{D}\ \forall\,y\in\mathbb{D} \ \exists\,z\in\mathbb{D}: \ (x\rho z)\ \& \ (y\rho z).$
С учетом этого введем множество $(\pr{DIR})[\mathbb{D}]$ всех направлений на множестве $\mathbb{D}:$
\begin{multline}\label{4.13}
(\pr{DIR})[\mathbb{D}]\triangleq\{\zeta\in(\pr{Ord})[\mathbb{D}]\mid \forall\,x\in\mathbb{D}\ \forall\,y\in\mathbb{D} \ \exists\,z\in\mathbb{D}: \\ (x\zeta z)\ \& \ (y\zeta z)\}\in\pp((\pr{Ord})[\mathbb{D}]).
\end{multline}

Если $\rho\in (\pr{DIR})[\mathbb{D}],$ то пару $(\mathbb{D},\rho)$ называем НМ (ниже рассматриваются только непустые НМ). Заметим в связи с последним понятием, что вместо букв при обозначении того или иного направления часто используются специальные символы: $\preceq,\ \angle$ и т.п. Подразумевается, что $\preceq\in(\pr{DIR})[\mathbb{D}],$ $\angle\in(\pr{DIR})[\mathbb{D}]$ и тому подобные соглашения.

Если $(\mathbb{D},\preceq)$ есть (непустое) НМ, а $X$\,--- произвольное множество, то определено множество $X^\mathbb{D}$  всевозможных отображений из $\mathbb{D}$ в множество $X$ (данное понятие содержательно при $X\neq\zer$).  Если к тому же $f\in X^\mathbb{D},$ то триплет $(\mathbb{D},\preceq,f)=((\mathbb{D},\preceq),f)$ называется \emph{направленностью в множестве} $X.$ Если же $f\in M^\mathbb{D},$ где $M\in\pp(X)$ (тогда, в частности, $f\in X^\mathbb{D}$),  то $(\mathbb{D},\preceq,f)$ будем называть также направленностью в $M.$ С учетом этого имеем при всяком выборе $A\in\pp(X)$ свойство: направленности в $A$ автоматически являются направленностями в $X.$ Для наших целей понятие направленности существенно в связи со сходимостью в ТП, обобщающей (\ref{4.11}).

Итак, если $(\mathbb{D},\preceq,f)$ есть направленность в $X$  и $x\in X,$  то
\beq{4.14}
((\mathbb{D},\preceq,f)\stackrel{\tau}\rightarrow x)\stackrel{\pr{def}}\Leftrightarrow(\forall\,H\in N_\tau(x) \ \exists\,d_1\in\mathbb{D} \ \forall\,d_2\in\mathbb{D} \ \ (d_1\preceq d_2)\Rightarrow(f(d_2)\in H)).
\eeq
Заметим, что $(\mathbb{N},\leq),$ где $\leq$ есть обычная упорядоченность натурального ряда $\mathbb{N},$ есть НМ; если, к тому же
$(x_i)_{i\in\mathbb{N}}\in X^\mathbb{N},$ то триплет $(\mathbb{N},\leq,(x_i)_{i\in\mathbb{N}})$ есть направленность в $X,$ порожденная исходной последовательностью $(x_i)_{i\in\mathbb{N}}.$ Если при этом полагать в (\ref{4.14}) $X=\mathbb{N}, \ \preceq=\leq$ и $f=(x_i)_{i\in\mathbb{N}},$ то (\ref{4.14}) сводится к (\ref{4.11}). Итак, обычная секвенциальная сходимость, определяемая в (\ref{4.11}),  действительно является частным случаем сходимости (\ref{4.14}); последняя, как уже отмечалось,   называется сходимостью по Мору-Смиту (см. \cite[гл.\,2]{40}).  В связи с упомянутой сходимостью напомним представление оператора замыкания, связанное с теоремой Биркгофа: если $M\in\pp(X),$ то $\pr{cl}(M,\tau)$ есть множество всех точек $x\in X,$ для каждой из которых существует направленность $(D,\preceq,f)$ в множестве $M$ со свойством
$$(D,\preceq,f)\stackrel{\tau}\rightarrow x.$$

Очень важной для дальнейшего изложения представляется процедура изотонного <<прореживания>> направленности до сходящейся поднаправленности. Рассмотрим краткую схему данной процедуры, полагая для каждого НМ $(\mathbb{D},\preceq),$ что
\begin{multline}\label{4.15}
(\preceq-\pr{cof})[\mathbb{D}]\triangleq\{H\in\pp(\mathbb{D})\mid\forall\,d\in\mathbb{D} \ \exists\,\delta\in H:d\preceq \delta\}=\\ \{H\in\pp'(\mathbb{D})\mid\forall\,d\in\mathbb{D} \ \exists\,\delta\in H: d\preceq \delta\}
\end{multline}
(напомним, что мы рассматриваем только непустые НМ); в (\ref{4.15}) определено семейство всех конфинальных в $(\mathbb{D},\preceq)$ п/м $\mathbb{D}.$ Мы будем использовать  данные понятия в духе \cite[раздел 2.2]{41}.

Если теперь $A$ и $B$\,--- непустые множества, $\preceq\in(\pr{DIR})[A]$ и $\angle\in(\pr{DIR})[B]$ (иными словами, $(A,\preceq)$ и $(B,\angle)$ суть НМ), то полагаем
\begin{multline}\label{4.16}
(\pr{Isot})[A;\preceq;B;\angle\,]\triangleq \bigl\{\varphi\in B^A\mid(\varphi^1(A)\in(\angle-\pr{cof})[B]) \ \& \\ \& \  (\forall\,\alpha_1\in A \ \forall\,\alpha_2\in A \ (\alpha_1\preceq\alpha_2)\Rightarrow(\varphi(\alpha_1)\angle\varphi(\alpha_2)))\bigl\}.
\end{multline}
Из (\ref{4.15}) и (\ref{4.16}) вытекает, что (при условиях, обеспечивающих (\ref{4.16})) для $\varphi\in(\pr{Isot})[A;\preceq;B;\angle\,]$
\beq{4.17}
\forall\,b\in B \ \exists\,\alpha\in A \ \forall\,a\in A \ (\alpha\preceq a)\Rightarrow(b\angle \varphi(a)).
\eeq

Отметим, что если в дополнение к условиям, определяющим (\ref{4.16}), заданы множество $H$ и отображение $h\in H^B$  (т.е. $(B,\angle,h)$  есть направленность в $H$), то при всяком выборе $\varphi\in(\pr{Isot})[A;\preceq;B;\angle\,]$  в  виде
\beq{4.18}
(A,\preceq,h\circ \varphi)\eeq
имеем направленность в $H,$ именуемую поднаправленностью $(B,\angle,h).$  Используем здесь изотонный способ прореживания $(B,\angle,h)$ до поднаправленности, определяемой в (\ref{4.18}). В частности, в качестве $H$ может использоваться $X,$ оснащенное топологией $\tau\in(\pr{top})[X];$  в этом случае $\forall\,h\in X^B \ \forall\,x\in X$
\beq{4.19}
\left((B,\angle,h)\stackrel{\tau}\rightarrow x\right)\Rightarrow\left((A,\preceq,h\circ\varphi)\stackrel{\tau}\rightarrow x \ \forall\,\varphi\in(\pr{Isot})[A;\preceq;B;\angle\,]\right).
\eeq

Итак, операция изотонного прореживания направленности сохраняет сходимость исходной направленности, если таковая имела место (предлагаем читателю проверить это самостоятельно, используя (\ref{4.17})).

В связи с применением направленностей  и сходимости по Мору-Смиту отметим следующие два важных свойства:\\
1*) \  если $A\in\pp(X),$ то $\pr{cl}(A,\tau)$ есть множество всех $x\in X,$ для каждого из которых существует направленность $(\mathbf{D},\sqsubseteq,f)$  в множестве $A$ со свойством $$(\mathbf{D},\sqsubseteq,f)\stackrel{\tau}\rightarrow x$$
(в этом положении отражено существо известной и уже упоминаемой теоремы Биркгофа);\\
2*) \ $(\tau-\pr{comp})[X]$ совпадает с семейством всех множеств $K\in\pp(X)$ таких, что при всяком выборе направленности $(\mathbf{D},\preceq,f)$ в множестве $K$ существуют (непустое) НМ $(\mathbb{D},\sqsubseteq~),$ отображение $\psi\in(\pr{Isot})[\mathbb{D};\sqsubseteq;\mathbf{D};\preceq]$ и точка $x\in K$ со свойством
\beq{4.20}
(\mathbb{D},\sqsubseteq,f\circ \psi)\stackrel{\tau}\rightarrow x.
\eeq
Введем теперь следующее обозначение. Если $(U,\tau_1), \ U\neq\zer;$
 и $(V,\tau_2), \ V\neq\zer$ суть два ТП, то в виде
\beq{4.100}
C(U,\tau_1,V,\tau_2)\triangleq\left\{f\in V^U\mid f^{-1}(G)\in\tau_1 \ \forall\,G\in\tau_2\right\}
\eeq
имеем множество всех отображений из множества $U$ в множество $V,$ \emph{непрерывных} в смысле $(U,\tau_1)$  и $(V,\tau_2).$

\section{Отношение  вписанности на пространстве разбиений}\setcounter{equation}{0} \setcounter{proposition}{0}\setcounter{zam}{0}

\ \ \ \ \ В настоящем разделе фиксируем непустое множество $E.$ Будем рассматривать семейство $\pp'(\pp(E))$ всех непустых подсемейств $\pp(E).$ Пусть до конца настоящего раздела $\mathcal{L}\in\pi[E].$ Если $A\in\pp(E),$ то в соответствии с (\ref{2.22}) через $\mathbf{D}(A,\mathcal{L})$ обозначаем семейство всех $\mathcal{K}\in\pr{Fin}(\mathcal{L}),$ для каждого из которых
\beq{5.0}
\biggl(A=\bigcup\limits_{L\in\mathcal{K}}L\biggl)\& \ \left(\forall\,L_1\in\mathcal{K} \ \forall\,L_2\in\mathcal{K} \ (L_1\cap L_2\neq\zer)\Rightarrow(L_1=L_2)\right);\eeq
тем самым введено семейство всех конечных разбиений $A$ множествами из $\mathcal{L}.$ В частности, определено семейство $\mathbf{D}(E,\mathcal{L}).$ Следуя \cite[(4.3.1)]{36}, полагаем, что $\forall\,\mathcal{A}\in\mathbf{D}(E,\mathcal{L}) \ \forall\,\mathcal{B}\in\mathbf{D}(E,\mathcal{L})$
\beq{5.1}
(\mathcal{A}\prec\mathcal{B})\Leftrightarrow(\forall\,B\in\mathcal{B} \ \exists\,A\in\mathcal{A}: B\subset A).
\eeq

Таким образом, введено бинарное отношение в $\mathbf{D}(E,\mathcal{L}),$ характеризуемое вписанностью одного разбиения в другое. Отметим одно эквивалентное представление данного бинарного отношения.
\begin{proposition}
Если $\mathcal{A}\in\mathbf{D}(E,\mathcal{L})$ и $\mathcal{B}\in\mathbf{D}(E,\mathcal{L}),$ то
\beq{5.2}
(\mathcal{A}\prec\mathcal{B})\Leftrightarrow(\forall\,A\in\mathcal{A}\setminus\{\zer\} \ \exists\,\mathcal{G}\in\mathbf{D}(A,\mathcal{L}):\mathcal{G}\subset \mathcal{B}).
\eeq
\end{proposition}
Д о к а з а т е л ь с т в о. Пусть $\mathcal{A}\prec\mathcal{B}.$ Тогда согласно (\ref{5.1})
\beq{5.3}
\forall\,B\in\mathcal{B} \ \exists\,A\in\mathcal{A}: B\subset A.
\eeq
Выберем произвольно $\mathbb{A}\in\mathcal{A}\setminus\{\zer\}.$  Тогда, в частности, $\mathbb{A}$ есть непустое п/м $E.$  Введем в рассмотрение семейство
\beq{5.4}
\mathfrak{B}\triangleq\{B\in\mathcal{B}\mid B\subset \mathbb{A}\}\in\pp(\mathcal{B}).
\eeq
Заметим, что в силу непустоты $\mathbb{A}$ найдется $\mathbb{B}\in\mathcal{B}$ со свойством
\beq{5.5}
\mathbb{A}\cap\mathbb{B}\neq\zer
\eeq
(в самом деле, имеем цепочку равенств
$$\mathbb{A}=\mathbb{A}\cap E=\mathbb{A}\cap\biggl(\,\bigcup\limits_{B\in\mathcal{B}}B\biggl)=\bigcup\limits_{B\in\mathcal{B}}(\mathbb{A}\cap B);$$
 тогда в случае $\mathbb{A}\cap B=\zer$ при всех $B\in\mathcal{B}$ имели бы $\mathbb{A}=\zer,$ что противоречит выбору $\mathbb{A}$). С учетом (\ref{5.3}) подберем $\mathbf{A}\in\mathcal{A}$ со свойством
 \beq{5.6}
 \mathbb{B}\subset\mathbf{A}.
 \eeq
Поскольку (см. (\ref{5.6})) $\mathbb{A}\cap\mathbb{B}\subset\mathbb{A}\cap\mathbf{A},$  то в силу (\ref{5.5}) $\mathbb{A}\cap\mathbf{A}\neq\zer.$  Так как  $\mathcal{A}$\,--- разбиение, то (см. (\ref{5.0})) $\mathbb{A}=\mathbf{A},$ а это значит, что (см. (\ref{5.6})) $\mathbb{B}\subset\mathbb{A}$ и (см. (\ref{5.4})) $\mathbb{B}\in\mathfrak{B}.$ Таким образом, $\mathfrak{B}\in\pp'(\mathcal{B}).$   Коль скоро $\mathcal{B}\subset\mathcal{L},$ то, в частности,
\beq{5.7}
\mathfrak{B}\in\pp'(\mathcal{L}).
\eeq
Более того, поскольку $\mathcal{B}$\,--- конечно, то и $\mathfrak{B}$\,--- конечно, а тогда (см. (\ref{5.7}))
\beq{5.8}
\mathfrak{B}\in\mathrm{Fin}(\mathcal{L}).
\eeq
Сравним $\mathbb{A}$ и объединение всех множеств из $\mathfrak{B}.$ Ясно, что (см. (\ref{5.4}))
\beq{5.9}
\bigcup\limits_{B\in\mathfrak{B}}B\subset\mathbb{A}.
\eeq

Пусть $a\in\mathbb{A}.$ Тогда, в частности, $a\in E,$ а потому $a\in\mathbf{B}$ для некоторого $\mathbf{B}\in\mathcal{B}.$  Покажем, что
\beq{5.10}
\mathbf{B}\in\mathfrak{B}.
\eeq

Допустим противное: $\mathbf{B}\notin\mathfrak{B}.$ Тогда согласно (\ref{5.4})
\beq{5.11}
\mathbf{B}\setminus\mathbb{A}\neq\zer.
\eeq
Поскольку $\mathbf{B}\in\mathcal{B},$  непременно $\exists\,A\in\mathcal{A}: \ \mathbf{B}\subset A.$ Пусть $A_0\in\mathcal{A}$ таково, что $\mathbf{B}\subset A_0.$ Так как $\mathbf{B}\setminus\mathbb{A}\subset\mathbf{B},$  имеем  из (\ref{5.11}), что $\mathbf{B}\neq\zer,$ а, стало быть, и $A_0\neq\zer.$ Поэтому
\beq{5.12}
A_0\in\mathcal{A}\setminus\{\zer\}.
\eeq
При этом, поскольку $\mathbb{A}\in\mathcal{A}$ и $A_0\in\mathcal{A},$ то
\beq{5.13}
(\mathbb{A}\cap A_0=\zer)\vee(\mathbb{A}=A_0).
\eeq
Напомним, что $a\in\mathbb{A}\cap A_0,$ и, тем более, $a\in\mathbb{A}\cap A_0$ (т.к. $\mathbf{B}\subset A_0).$  Поэтому $\mathbb{A}\cap A_0\neq\zer,$ а тогда в силу (\ref{5.13}) получаем равенство $$\mathbb{A}=A_0,$$ из которого следует, что $\mathbf{B}\subset\mathbb{A}$ вопреки (\ref{5.11}). Полученное противоречие доказывает (\ref{5.10}). Итак, получаем, как следствие, что
\beq{5.14}
a\in\bigcup\limits_{B\in\mathfrak{B}}B.
\eeq

Поскольку выбор $a$ был произвольным, установлено, что $$\mathbb{A}\subset\bigcup\limits_{B\in\mathfrak{B}}B,$$
а тогда с учетом (\ref{5.9}) имеем очевидное равенство
\beq{5.15}
\mathbb{A}=\bigcup\limits_{B\in\mathfrak{B}}B.
\eeq
Напомним, что $\mathcal{B}\subset\mathcal{L}.$ Кроме того, по выбору $\mathcal{B}$ имеем в силу свойства $\mathfrak{B}\subset\mathcal{B},$ что $\forall\,B_1\in\mathfrak{B} \ \forall\,B_2\in\mathfrak{B}$
\beq{5.16}
(B_1\cap B_2\neq\zer)\Rightarrow(B_1=B_2).
\eeq
Заметим, что $\mathfrak{B}\in\pr{Fin}(\mathcal{L})$ и $\mathbb{A}\neq\zer,$ а тогда из (\ref{5.15}) вытекает, что $\mathfrak{B}\neq\zer.$ При этом, конечно,
\beq{5.17}
\mathfrak{B}\in\pp'(\mathcal{B})
\eeq
и, следовательно, $\mathfrak{B}\in\pr{Fin}(\mathcal{L})$ (учитываем, что $\mathfrak{B}\subset\mathcal{B}\subset\mathcal{L}$). Из (\ref{5.0}), (\ref{5.15})--(\ref{5.17}) имеем, что
\beq{5.18}
\mathfrak{B}\in\mathbf{D}(\mathbb{A},\mathcal{L}).
\eeq

Тем самым установлено, в частности, что $\exists\,\mathcal{G}\in\mathbf{D}(\mathbb{A},\mathcal{L}): \mathcal{G}\subset\mathcal{B}.$ Поскольку выбор $\mathbb{A}$ был произвольным, установлено, что
$$\forall\,A\in\mathcal{A}\setminus\{\zer\} \ \exists\,\mathcal{G}\in\mathbf{D}(A,\mathcal{L}):\mathcal{G}\subset\mathcal{B}.$$
Таким образом,  завершено обоснование импликации
\beq{5.19}
(\mathcal{A}\prec\mathcal{B})\Rightarrow(\forall\,A\in\mathcal{A}\setminus\{\zer\} \ \exists\,\mathcal{G}\in\mathbf{D}(A,\mathcal{L}):\mathcal{G}\subset\mathcal{B}).
\eeq
Пусть теперь, напротив, имеет место свойство:
\beq{5.20}
\forall\,A\in\mathcal{A}\setminus\{\zer\} \ \exists\,\mathcal{G}\in\mathbf{D}(A,\mathcal{L}): \ \mathcal{G}\subset\mathcal{B}.
\eeq
Надо показать, что $\mathcal{A}\prec\mathcal{B}.$  Выберем произвольно $B_*\in\mathcal{B}.$ Покажем, что
\beq{5.21}
\exists\,A\in\mathcal{A}: B_*\subset A.
\eeq
Если $B_*=\zer,$ то свойство (\ref{5.21}) имеет место, т.к. $\mathcal{A}\neq\zer.$  Итак,
\beq{5.22}
(B_*=\zer)\Rightarrow(\exists\,A\in\mathcal{A}: B_*\subset A).
\eeq
Пусть теперь $B_*\neq\zer.$  Допустим, что (\ref{5.21}) не имеет места:
\beq{5.23}
B_*\setminus A\neq\zer \ \forall\,A\in\mathcal{A}.
\eeq

С другой стороны, в силу непустоты $B_*$
\beq{5.24}
\exists\,A\in\mathcal{A}:A\cap B_*\neq\zer
\eeq
(в самом деле, имеем равенство $$E=\bigcup\limits_{A\in\mathcal{A}}A$$  и, кроме того, $B_*\subset E,$ а потому
$$B_*=E\cap B_*=\bigcup\limits_{A\in\mathcal{A}}(A\cap B_*),$$
откуда вытекает следующее свойство
$$\bigcup\limits_{A\in\mathcal{A}}(A\cap B_*)\neq\zer,$$
означающее справедливость (\ref{5.24}).) Используя (\ref{5.24}), выберем и зафиксируем $A_*\in\mathcal{A},$ для которого
\beq{5.25}
A_*\cap B_*\neq\zer.
\eeq
В силу  (\ref{5.23}) имеем, вместе с тем, что
\beq{5.26}
B_*\setminus A_*\neq\zer,
\eeq
Из (\ref{5.25}) следует, кроме того, что $A_*\neq\zer,$ а тогда
\beq{5.27}
A_*\in\mathcal{A}\setminus\{\zer\}.
\eeq
С учетом (\ref{5.20}) подберем разбиение

\beq{5.28}
\mathcal{G}_*\in\mathbf{D}(A_*,\mathcal{L}),
\eeq
 для которого имеет место вложение
 \beq{5.29}
 \mathcal{G}_*\subset\mathcal{B}.
 \eeq
 Из (\ref{5.0}) и (\ref{5.28}) получаем, в частности, что
 \beq{5.30}
 A_*=\bigcup\limits_{L\in\mathcal{G_*}}L.
 \eeq
 При этом, конечно, $\mathcal{G}_*\in\pr{Fin}(\mathcal{L}).$  Из (\ref{5.25}) и (\ref{5.30}) вытекает, что
 $$\bigcup\limits_{L\in\mathcal{G}_*}(L\cap B_*)\neq\zer,$$ и  тогда можно указать множество $C_*\in\mathcal{G}_*,$ для которого
 \beq{5.31}
 C_*\cap B_*\neq\zer.
 \eeq
 Из (\ref{5.29}) следует, в частности, что $C_*\in\mathcal{B}.$ Поэтому $$(C_*\in\mathcal{B})\ \& \ (B_*\in\mathcal{B}),$$
 а тогда с учетом (\ref{5.0}) истинна импликация
 $$(C_*\cap B_*\neq\zer)\Rightarrow(C_*=B_*).$$
 Используя (\ref{5.31}), получаем следующее равенство:
 \beq{5.32}
 C_*=B_*.
 \eeq

 С другой стороны, из (\ref{5.30}) вытекает, что $C_*\subset A_*.$ Стало быть (см. (\ref{5.32})),
 \beq{5.33}
 B_*\subset A_*,
 \eeq
 что противоречит (\ref{5.26}). Полученное при условии (\ref{5.23}) противоречие (см. (\ref{5.26}), (\ref{5.33})) показывает, что само (\ref{5.23}) невозможно, а потому справедливо (\ref{5.21}) и в случае $B_*\neq\zer.$
 Итак, $$(B_*\neq\zer)\Rightarrow(\exists\,A\in\mathcal{A}: \ B_*\subset A).$$
  С учетом (\ref{5.22}) имеем теперь (\ref{5.21}) во всех возможных случаях. Поскольку выбор $B_*$ был произвольным, установлено, что
 $$\forall\,B\in\mathcal{B} \ \exists\,A\in\mathcal{A}: B\subset A.$$
  В силу (\ref{5.1}) получаем свойство $\mathcal{A}\prec\mathcal{B},$ чем и завершается проверка истинности импликации
 \beq{5.34}
 (\forall\,A\in\mathcal{A}\setminus\{\zer\} \ \exists\,\mathcal{G}\in\mathbf{D}(A,\mathcal{L}): \ \mathcal{G}\subset\mathcal{B})\Rightarrow(\mathcal{A}\prec\mathcal{B}).
  \eeq

 Теперь из (\ref{5.19}) и (\ref{5.34}) следует (\ref{5.2}). $\hfill\square$
 \begin{proposition}
 Бинарное отношение $\prec$ есть направление на $\mathbf{D}(E,\mathcal{L}):$
 \beq{5.35}
 \prec\in(\pr{DIR})[\mathbf{D}(E,\mathcal{L})].
 \eeq
 \end{proposition}
Д о к а з а т е л ь с т в о. Отметим, что  отношение
$$\prec\in\pp(\mathbf{D}(E,\mathcal{L})\times\mathbf{D}(E,\mathcal{L}))$$
определяется правилом
\beq{5.35'}\prec\triangleq\{z\in\mathbf{D}(E,\mathcal{L})\times\mathbf{D}(E,\mathcal{L})\mid\forall\, B\in\pr{pr}_2(z) \ \exists\,A\in\pr{pr}_1(z): B\subset A\}.\eeq
Разумеется, согласно (\ref{5.1}) и (\ref{5.35}) $\forall\,\mathcal{A}\in\mathbf{D}(E,\mathcal{L}) \ \forall\,\mathcal{B}\in\mathbf{D}(E,\mathcal{L})$
$$(\mathcal{A}\prec\mathcal{B})\Leftrightarrow((\mathcal{A},\mathcal{B})\in\prec).$$

Покажем, что (см. (\ref{4.12})) $\prec\in(\pr{Ord})[\mathbf{D}(E,\mathcal{L})].$ Из (\ref{5.1}) вытекает, что
\beq{5.36}
\mathcal{M}\prec\mathcal{M} \ \forall\,\mathcal{M}\in\mathbf{D}(E,\mathcal{L}).\eeq
Пусть теперь заданы три разбиения
\beq{5.37}
(\mathcal{U}\in\mathbf{D}(E,\mathcal{L})) \ \& \ (\mathcal{V}\in\mathbf{D}(E,\mathcal{L})) \ \& \ (\mathcal{W}\in\mathbf{D}(E,\mathcal{L}))
\eeq
со следующими свойствами
\beq{5.38}
(\mathcal{U}\prec\mathcal{V})\ \& \ (\mathcal{V}\prec\mathcal{W}).
\eeq
Из (\ref{5.1}) и (\ref{5.38}) имеем, что
\beq{5.39}
(\forall\,V\in\mathcal{V} \ \exists\,U\in\mathcal{U}: V\subset U)\ \& \ (\forall\,W\in\mathcal{W} \ \exists\,V\in\mathcal{V}: W\subset V).
\eeq
Из (\ref{5.39}) непосредственно следует свойство
\beq{5.40}
\forall\,W\in\mathcal{W} \ \exists\,U\in\mathcal{U}: \ W\subset U.
\eeq
Согласно (\ref{5.1}) и (\ref{5.40}) $\mathcal{U}\prec\mathcal{W},$ чем завершается (см. (\ref{5.38})) проверка импликации
\beq{5.41}
((\mathcal{U}\prec\mathcal{V}) \ \& \ (\mathcal{V}\prec\mathcal{W}))\Rightarrow(\mathcal{U}\prec\mathcal{W}).
\eeq
Из (\ref{4.12}), (\ref{5.36}) и (\ref{5.41}) вытекает свойство
\beq{5.42}
\prec\in(\pr{Ord})[\mathbf{D}(E,\mathcal{L})], \eeq
а тогда $(\mathbf{D}(E,\mathcal{L}),\prec)$ есть ЧУМ. Пусть выбраны произвольно
\beq{5.43}
(\mathcal{S}\in\mathbf{D}(E,\mathcal{L}))\ \& \ (\mathcal{T}\in\mathbf{D}(E,\mathcal{L})).
\eeq

Тогда, в частности, $\mathcal{S}\in\pr{Fin}(\mathcal{L})$ и $\mathcal{T}\in\pr{Fin}(\mathcal{L}).$ Как следствие, $\mathcal{S}\times\mathcal{T}$ есть непустое конечное множество, а, точнее, (см. \cite[с.\,47]{30})
\beq{5.44}
\mathcal{S}\times\mathcal{T}\in\pr{Fin}(\mathcal{L}\times\mathcal{L}).
\eeq

С учетом (\ref{5.44}) введем в рассмотрение семейство
\beq{5.45}
\mathfrak{R}\triangleq\{\pr{pr}_1(z)\cap\pr{pr}_2(z): \ z\in\mathcal{S}\times\mathcal{T}\}\in\pp'(\mathcal{L})
\eeq
(учитываем при этом (\ref{3.1})). Оно, как легко видеть, конечно. В самом деле, с учетом (\ref{5.44}) имеем
 \beq{5.46}
 n\triangleq|\mathcal{S}\times\mathcal{T}|\in\mathbb{N},
 \eeq
 причем согласно (\ref{2.20}) и (\ref{5.46}) $(\pr{bi})[\overline{1,n};\mathcal{S}\times\mathcal{T}]\neq\zer.$  Выберем произвольно биекцию
\beq{5.47}
\varphi\in(\pr{bi})[\overline{1,n};\mathcal{S}\times\mathcal{T}].
\eeq
Тогда, в частности, $\varphi:\overline{1,n}\rightarrow\mathcal{S}\times\mathcal{T},$ причем (см. (\ref{2.15''}), (\ref{5.47}))
\beq{5.48}
\forall\,z\in\mathcal{S}\times\mathcal{T} \ \exists\,j\in\overline{1,n}:z=\varphi(j).
\eeq

Поэтому согласно (\ref{5.45}), (\ref{5.47}) и  (\ref{5.48})
\beq{5.49}
\mathfrak{R}=\{\pr{pr}_1(\varphi(j))\cap\pr{pr}_2(\varphi(j)): \ j\in\overline{1,n}\}\in\pp'(\mathcal{L}).
\eeq
Введем теперь в рассмотрение отображение
$$\psi\triangleq(\pr{pr}_1(\varphi(j))\cap\pr{pr}_2(\varphi(j)))_{j\in\overline{1,n}}\in\mathcal{L}^n.
$$
Это означает, в частности, что $\psi:\overline{1,n}\rightarrow\mathcal{L},$ причем (см. (\ref{5.49})) $\psi^1(\overline{1,n})~=\mathfrak{R},$ откуда следует включение $$\psi\in(\pr{su})[\overline{1,n};\mathfrak{R}].$$ Таким образом, (см. \cite[(2.50)]{31}) $\mathfrak{R}$ есть непустое конечное множество и, с учетом (\ref{5.45}),
\beq{5.50}
\mathfrak{R}\in\pr{Fin}(\mathcal{L}).
\eeq
Покажем, что на самом деле $\mathfrak{R}\in\mathbf{D}(E,\mathcal{L}).$ Для этого прежде всего заметим, что (см. (\ref{5.47}), (\ref{5.49}))
\beq{5.51}
\bigcup\limits_{L\in\mathfrak{R}}L=\bigcup\limits_{i=1}^n(\pr{pr}_1(\varphi(i))\cap\pr{pr}_2(\varphi(i))=\bigcup\limits_{z\in\mathcal{S}\times\mathcal{T}}
(\pr{pr}_1(z)\cap\pr{pr}_2(z)).
\eeq
Из (\ref{5.50}) следует вложение
\beq{5.52}
\bigcup\limits_{L\in\mathfrak{R}}L\subset E,
\eeq
где учтено, что $\mathcal{L}\subset\pp(E).$ Пусть $x_*\in E.$ Тогда согласно (\ref{5.0}) и (\ref{5.43}) имеем для некоторых $S_*\in\mathcal{S}$ и $T_*\in\mathcal{T}$ включение
\beq{5.53}
x_*\in S_*\cap T_*
\eeq
(учитываем, что из (\ref{5.0}) следует $\left(E=\bigcup\limits_{S\in\mathcal{S}}S\right) \ \& \ \left(E=\bigcup\limits_{T\in\mathcal{T}}T\right),$  вытекающее из (\ref{5.43})). При этом $z_*\triangleq(S_*,T_*)\in\mathcal{S}\times\mathcal{T}$ и
$$(\pr{pr}_1(z_*)=S_*) \ \& \ (\pr{pr}_2(z_*)=T_*),$$ а тогда имеем, что (см. (\ref{5.51}))
$$S_*\cap T_*=\pr{pr}_1(z_*)\cap\pr{pr}_2(z_*)\subset\bigcup\limits_{L\in\mathfrak{R}}L,$$
откуда в силу (\ref{5.53}) следует включение $$x_*\in\bigcup\limits_{L\in\mathfrak{R}}L.$$
Поскольку выбор $x_*$ был произвольным, установлено, что $$E\subset\bigcup\limits_{L\in\mathfrak{R}}L.$$ С учетом (\ref{5.52}) получаем  следующее равенство
\beq{5.54}
E=\bigcup\limits_{L\in\mathfrak{R}}L.
\eeq
Пусть теперь выбраны произвольно
\beq{5.55}
(\Lambda_1\in\mathfrak{R})\ \& \ (\Lambda_2\in\mathfrak{R}),\eeq
для которых имеет место свойство
\beq{5.56}
\Lambda_1\cap\Lambda_2\neq\zer.
\eeq

Тогда в силу (\ref{5.49}) и (\ref{5.55}) имеем для некоторых $$(\alpha\in\overline{1,n})\ \& \ (\beta\in\overline{1,n})$$
следующие два равенства
\beq{5.57}
\Lambda_1=\pr{pr}_1(\varphi(\alpha))\cap\pr{pr}_2(\varphi(\alpha)),
\eeq

\beq{5.58}
\Lambda_2=\pr{pr}_1(\varphi(\beta))\cap\pr{pr}_2(\varphi(\beta)).
\eeq
Из (\ref{5.56})--(\ref{5.58}) следует
\beq{5.59}
(\pr{pr}_1(\varphi(\alpha))\cap\pr{pr}_1(\varphi(\beta))\neq\zer)\ \& \ (\pr{pr}_2(\varphi(\alpha))\cap\pr{pr}_2(\varphi(\beta))\neq\zer),
\eeq
где согласно (\ref{5.47}) $\pr{pr}_1(\varphi(\alpha))\in\mathcal{S},$ $\pr{pr}_1(\varphi(\beta))\in\mathcal{S},$ $\pr{pr}_2(\varphi(\alpha))\in\mathcal{T}$ и $\pr{pr}_2(\varphi(\beta))\in\mathcal{T}.$ Из (\ref{5.0}) и (\ref{5.43}) имеем импликацию
\beq{5.60}
(\pr{pr}_1(\varphi(\alpha))\cap\pr{pr}_1(\varphi(\beta))\neq\zer)\Rightarrow(\pr{pr}_1(\varphi(\alpha))=\pr{pr}_1(\varphi(\beta)))
\eeq
и, аналогичным образом,
\beq{5.61}
(\pr{pr}_2(\varphi(\alpha))\cap\pr{pr}_2(\varphi(\beta))\neq\zer)\Rightarrow(\pr{pr}_2(\varphi(\alpha))=\pr{pr}_2(\varphi(\beta))).
\eeq
Из (\ref{5.59})--(\ref{5.61}) вытекает, что $(\pr{pr}_1(\varphi(\alpha))=\pr{pr}_1(\varphi(\beta))\ \& \ (\pr{pr}_2(\varphi(\alpha))=\pr{pr}_2(\varphi(\beta)),$ откуда по основному свойству упорядоченных пар (см. (\ref{2.6'})) получаем равенство  $$\varphi(\alpha)=\varphi(\beta).$$
Но в этом случае согласно (\ref{5.57}), (\ref{5.58}) $\Lambda_1=\Lambda_2,$  чем завершается (см. (\ref{5.56})) проверка истинности импликации
$$(\Lambda_1\cap\Lambda_2\neq\zer)\Rightarrow(\Lambda_1=\Lambda_2).$$
Поскольку выбор (\ref{5.55}) был произвольным, установлено свойство: $\forall\,L_1\in~\mathfrak{R}$ $ \forall\,L_2\in\mathfrak{R}$
\beq{5.62}
(L_1\cap L_2\neq\zer)\Rightarrow(L_1=L_2).
\eeq
Из (\ref{5.50}), (\ref{5.54}) и (\ref{5.62}) следует (см. (\ref{5.0})) требуемое свойство
\beq{5.63}
\mathfrak{R}\in\mathbf{D}(E,\mathcal{L}).
\eeq

Выберем произвольно $R_0\in\mathfrak{R}.$ Тогда согласно (\ref{5.49}) имеем для некоторого $j_0\in\overline{1,n}$ равенство
\beq{5.64}
R_0=\pr{pr}_1(\varphi(j_0))\cap\pr{pr}_2(\varphi(j_0)),
\eeq
где $\varphi(j_0)\in\mathcal{S}\times\mathcal{T}$  в силу (\ref{5.47}). При этом согласно (\ref{5.47})
\beq{5.65}
\left(S_0\triangleq\pr{pr}_1(\varphi(j_0))\in\mathcal{S}\right)\ \& \ \left(T_0\triangleq\pr{pr}_2(\varphi(j_0))\in\mathcal{T}\right).
\eeq
Из (\ref{5.64}) и (\ref{5.65}) вытекает равенство  $R_0=S_0\cap T_0.$ Это означает, в частности, что (см. (\ref{5.65})) $R_0\subset S_0,$ а потому $\exists\,S\in\mathcal{S}:R_0\subset S.$ Кроме того, имеем согласно  (\ref{5.65}), что $\exists\,T\in\mathcal{T}:R_0\subset T.$ Следовательно,
\beq{5.66}
(\exists\,S\in\mathcal{S}:R_0\subset S)\ \& \ (\exists\,T\in\mathcal{T}:R_0\subset T).
\eeq

Поскольку выбор $R_0$ был произвольным, получаем $\forall\,B\in\mathfrak{R}$
$$(\exists\,S\in\mathcal{S}:B\subset S)\ \& \ (\exists\,T\in\mathcal{T}:B\subset T).$$
С учетом (\ref{5.1}) имеем
\beq{5.67}
(\mathcal{S}\prec\mathfrak{R})\ \& \ (\mathcal{T}\prec\mathfrak{R}).
\eeq
Используя (\ref{5.63}) и (\ref{5.67}), получаем, в частности, следующее свойство:
$$\exists\,\mathcal{Z}\in\mathbf{D}(E,\mathcal{L}):(\mathcal{S}\prec\mathcal{Z})\ \& \ (\mathcal{T}\prec\mathcal{Z}).$$
Поскольку выбор $\mathcal{S}$ и $\mathcal{T}$  в (\ref{5.43}) был произвольным, установлено, что
\beq{5.68}
\forall\,\mathcal{X}\in\mathbf{D}(E,\mathcal{L}) \ \forall\,\mathcal{Y}\in\mathbf{D}(E,\mathcal{L}) \ \exists\,\mathcal{Z}\in\mathbf{D}(E,\mathcal{L}):(\mathcal{X}\prec\mathcal{Z})\ \& \ (\mathcal{Y}\prec\mathcal{Z}).
\eeq
  Как уже отмечалось, $\prec\in(\pr{Ord})[\mathbf{D}(E,\mathcal{L})],$ и, таким образом,  имеем из (\ref{4.13}) и (\ref{5.68})  свойство (\ref{5.35}). $\hfill\square$

Заметим в связи с определением $\mathbf{D}(E,\mathcal{L})$ (см. (\ref{5.0})) очевидное свойство $\{E\}\in\mathbf{D}(E,\mathcal{L}),$ а потому $\mathbf{D}(E,\mathcal{L})\neq\zer.$ В виде $$(\mathbf{D}(E,\mathcal{L}),\prec)$$ имеем теперь (непустое) НМ. В предложении 5.1 приведено полезное представление этого НМ, которое будет использоваться далее в связи со свойством конечной аддитивности функций множеств.

\section{Измеримые пространства с полуалгебрами и алгебрами множеств}\setcounter{equation}{0} \setcounter{proposition}{0}\setcounter{zam}{0}

\ \ \ \ \ Возвращаясь к (\ref{3.1}), отметим другие важные частные случаи данного определения. Ограничимся здесь полуалгебрами и алгебрами множеств, что  будет достаточным для наших построений, использующих к.-а. меры (случай стандартных измеримых пространств подробно освещается в литературе по теории меры; кроме того, см. \cite[\S 1.7]{30}).  Как и в предыдущем разделе,  фиксируем множество $E,$ получая из (\ref{3.1}),  что
\beq{5'.1}
\pi[E]=\{\mathcal{E}\in\pp'(\pp(E))\mid(\zer\in\mathcal{E})\ \& \ (E\in\mathcal{E}) \ \& \ (A\cap B\in\mathcal{E} \ \forall\,A\in\mathcal{E} \ \forall\,B\in\mathcal{E})\}.
\eeq

В (\ref{5'.1}) мы имеем семейство всех $\pi$-систем п/м $E$ с <<нулем>> и <<единицей>>. Среди данных $\pi$-систем выделяем алгебры и полуалгебры п/м $E.$ Итак, введем семейство всех алгебр п/м~$E:$
\beq{5'.2}
(\pr{alg})[E]\triangleq\{\mathcal{A}\in\pi[E]\mid E\setminus A\in\mathcal{A} \ \forall\, A\in\mathcal{A}\}.
\eeq

Если $\mathcal{A}\in(\pr{alg})[E],$ то называем $(E,\mathcal{A})$ \emph{измеримым пространством  (ИП) с алгеброй множеств}. Отметим  полезное следствие (\ref{5'.2}). Если $\mathcal{A}\in(\pr{alg})[E],$ то
\beq{5'.3}
A_1\cup A_2\in\mathcal{A} \ \forall\,A_1\in\mathcal{A} \ \forall\,A_2\in\mathcal{A}.
\eeq
\begin{zam}
Проверим (\ref{5'.3}), фиксируя $\mathcal{A}.$ Согласно (\ref{5'.1}) и (\ref{5'.2}) при $A_1\in\mathcal{A}$ и $A_2\in\mathcal{A}$
\beq{5'.4}
A_1\cup A_2=(E\setminus(E\setminus A_1))\cup(E\setminus(E\setminus A_2))=E\setminus((E\setminus A_1)\cap(E\setminus A_2)),
\eeq
где $E\setminus A_1\in\mathcal{A}$ и $E\setminus A_2\in\mathcal{A}.$ Поскольку, в частности, $\mathcal{A}\in\pi[E]$ (см. (\ref{5'.2})), то справедливо свойство
$$(E\setminus A_1)\cap(E\setminus A_2)\in\mathcal{A}$$ и, как следствие (снова используем (\ref{5'.2})),
 $$E\setminus((E\setminus A_1)\cap(E\setminus A_2))\in\mathcal{A},$$
 что в силу (\ref{5'.4}) означает тот факт, что $A_1\cup A_2\in\mathcal{A}.$ $\hfill\square$
\end{zam}
Разумеется, по индукции установливается, что
$$\bigcup\limits_{i=1}^n A^{(i)}\in\mathcal{A} \ \forall\,n\in\mathbb{N} \ \forall\,(A^{(i)})_{i\in\overline{1,n}}\in\mathcal{A}^n.$$
Учитывая (\ref{5'.1}), введем в рассмотрение семейство
\beq{5'.5}
\Pi[E]\triangleq\{\mathcal{E}\in\pi[E]\mid\forall\,\Sigma\in\mathcal{E} \ \exists \, n\in\mathbb{N}:\Delta_n(E\setminus\Sigma,\mathcal{E})\neq\zer\};
\eeq
элементы семейства (\ref{5'.5}), а это~---~ усовершенствованные $\pi$-системы, называются \emph{полуалгебрами п/м} $E.$ Если $\mathcal{E}\in\Pi[E],$ то пару $(E,\mathcal{E})$ называем \emph{ИП с полуалгеброй множеств}.
\begin{proposition}
Справедливо равенство
\beq{5'.6}
\Pi[E]=\{\mathcal{E}\in\pi[E]\mid\mathbf{D}(E\setminus\Sigma,\mathcal{E})\neq\zer \ \forall\,\Sigma\in\mathcal{E}\}.
\eeq
\end{proposition}
Д о к а з а т е л ь с т в о. Пусть $\Xi$ есть $\pr{def}$ семейство в правой части (\ref{5'.6}). Надо показать, что $\Pi[E]=\Xi.$ Пусть $\mathcal{I}\in\Pi[E].$ Тогда согласно (\ref{5'.5}) $\mathcal{I}\in\pi[E]$ и при этом
\beq{5'.7}
\forall\,\Sigma\in\mathcal{I} \ \exists\,n\in\mathbb{N}:\Delta_n(E\setminus\Sigma,\mathcal{I})\neq\zer.
\eeq

Выберем произвольно $\Lambda\in\mathcal{I},$ после чего с использованием (\ref{5'.7}) подберем $\mathbf{n}\in\mathbb{N}$ и
\beq{5'.8}
(\Lambda_i)_{i\in\overline{1,\mathbf{n}}}\in\Delta_{\mathbf{n}}(E\setminus\Lambda,\mathcal{I}).
\eeq
С учетом (\ref{2.23}) имеем теперь, что
\beq{5'.9}
\mathfrak{L}\triangleq\{\Lambda_i:i\in\overline{1,\mathbf{n}}\}\in\mathbf{D}(E\setminus\Lambda,\mathcal{I}).
\eeq
Заметим, что согласно (\ref{2.32}) $\mathfrak{L}\in\pr{Fin}(\mathcal{I})$ обладает свойствами
\beq{5'.10}
\biggl(E\setminus\Lambda=\bigcup\limits_{L\in\mathfrak{L}}L\biggl)\ \& \ (\forall\,A\in\mathfrak{L} \ \forall\,B\in\mathfrak{L} \  \ (A\cap B\neq\zer)\Rightarrow(A=B)).
\eeq
Нам, однако, достаточно (\ref{5'.9}): имеем, в частности, что $\mathbf{D}(E\setminus\Lambda,\mathcal{I})\neq\zer.$  Поскольку выбор $\Lambda$ был произвольным, получено свойство
$\mathbf{D}(E\setminus\Sigma,\mathcal{I})\neq\zer$ $\forall\,\Sigma\in\mathcal{I}.$ Это означает, что $\mathcal{I}\in\Xi.$ Таким образом, установлено вложение
\beq{5'.11}
\Pi[E]\subset\Xi.
\eeq
Выберем произвольно $\mathcal{J}\in\Xi.$ Тогда $\mathcal{J}\in\pi[E]$ и при этом
\beq{5'.12}
\mathbf{D}(E\setminus\Sigma,\mathcal{J})\neq\zer \ \forall\,\Sigma\in\mathcal{J}.
\eeq
Выберем произвольно $J\in\mathcal{J},$ после чего, используя (\ref{5'.12}), подберем
\beq{5'.13}
\mathfrak{M}\in\mathbf{D}(E\setminus\ J,\mathcal{J}).
\eeq
Тогда,  в частности, $\mathfrak{M}\in\pr{Fin}(\mathcal{J})$ и для $\mathbf{m}\triangleq|\mathfrak{M}|\in\mathbb{N}$  согласно (\ref{2.28})
  \beq{5'.14}
  (\pr{bi})[\overline{1,\mathbf{m}};\mathfrak{M}]\subset\Delta_{\mathbf{m}}(E\setminus J,\mathcal{J}).
  \eeq
  Вместе с тем согласно (\ref{2.20}) имеем, что
\beq{5'.15}
(\pr{bi})[\overline{1,\mathbf{m}};\mathfrak{M}]\neq\zer,
\eeq
а тогда (см. (\ref{5'.14}), (\ref{5'.15})) справедливо свойство
\beq{5'.16}
\Delta_{\mathbf{m}}(E\setminus J,\mathcal{J})\neq\zer.
\eeq
Это означает, в частности, что  $\exists\,n\in\mathbb{N}: \ \Delta_{n}(E\setminus J,\mathcal{J})\neq\zer.$ Поскольку выбор $J$ был произвольным, установлено, что $$\forall\,\Sigma\in\mathcal{J} \ \exists\,n\in\mathbb{N}: \ \Delta_{n}(E\setminus \Sigma,\mathcal{J})\neq\zer.$$
С учетом (\ref{5'.5}) имеем включение $\mathcal{J}\in\Pi[E].$  Тем самым установлено вложение $\Xi\subset\Pi[E],$
а потому (см. (\ref{5'.11})) $\Pi[E]=\Xi.$  Что и требовалось доказать. $\hfill\square$

Заметим, что при всяком выборе $\mathcal{L}\in\Pi[E]$
\beq{5'.17}
\mathbf{a}_E^0(\mathcal{L})\triangleq\{A\in\pp(E)\mid \exists\,n\in\mathbb{N}: \ \Delta_n(A,\mathcal{L})\neq\zer\}\in(\pr{alg})[E],
\eeq
причем $\mathcal{L}\subset\mathbf{a}_E^0(\mathcal{L})$ и, кроме того, $\forall\,\mathcal{A}\in(\pr{alg})[E]$
\beq{5'.18}
(\mathcal{L}\subset\mathcal{A})\Rightarrow\left(\mathbf{a}_E^0(\mathcal{L})\subset\mathcal{A}\right).
\eeq
В (\ref{5'.17}), (\ref{5'.18}) мы имеем простой способ построения алгебры п/м $E,$ порожденной произвольной полуалгеброй п/м $E.$ Отметим, наконец, что \cite[(1.7.5)]{30}
\begin{multline}\label{5'.19}
(\sigma-\pr{alg})[E]\triangleq\biggl\{\mathcal{L}\in(\pr{alg})[E]\mid\bigcup\limits_{i\in\mathbb{N}}L_i\in\mathcal{L} \ \forall\,(L_i)_{i\in\mathbb{N}}\in\mathcal{L}^{\mathbb{N}}\biggl\}=\\=\biggl\{\mathcal{L}\in(\pr{alg})[E]\mid\bigcap\limits_{i\in\mathbb{N}}L_i\in\mathcal{L} \ \forall\,(L_i)_{i\in\mathbb{N}}\in\mathcal{L}^{\mathbb{N}}\biggl\}
\end{multline}
 есть (всегда непустое) семейство всех $\sigma$-алгебр п/м $E.$

\newpage

\begin{center} \section*{Глава 2. Элементы конечно-аддитивной теории меры} \end{center}
\section{Конечно-аддитивные меры: общие свойства}\setcounter{equation}{0} \setcounter{proposition}{0}\setcounter{zam}{0}

\ \ \ \ \ В настоящем разделе сохраняем предположения пятого раздела: $E$\,--- непустое множество, $\mathcal{L}\in\pi[E].$ Тогда $\mathbb{R}^\mathcal{L}$ есть в силу  (\ref{2.8}) множество всех функций
\beq{6.1}
\mu:\mathcal{L}\rightarrow\mathbb{R}.
\eeq

В дальнейшем функции вида (\ref{6.1}) будем именовать \emph{функциями множеств}. Мы ограничиваемся при этом в/з функциями множеств. Среди всех таких функций будем выделять \emph{конечно-аддитивные (к.-а.) меры}.

Итак, функцию множеств  $\mu$ (\ref{6.1}) называем \emph{к.-а. мерой} (на $\mathcal{L}$), если
\beq{6.2}
\mu(L)=\sum\limits_{i=1}^n\mu(L_i)  \ \ \forall\,L\in\mathcal{L} \ \ \forall\,n\in\mathbb{N} \ \ \forall\,(L_i)_{i\in\overline{1,n}}\in\Delta_n(L,\mathcal{L}).
\eeq

В связи с (\ref{6.2}) получаем, что
\begin{multline}\label{6.3}
(\pr{add})[\mathcal{L}]\triangleq\biggl\{\mu\in\mathbb{R}^{\mathcal{L}}\mid \mu(L)=\sum\limits_{i=1}^n\mu(L_i) \ \ \forall\,L\in\mathcal{L} \ \\ \forall\,n\in\mathbb{N} \ \  \forall\,(L_i)_{i\in\overline{1,n}}\in\Delta_n(L,\mathcal{L})\biggl\}
\end{multline}
есть множество всех (в/з) к.-а. мер на $\mathcal{L}.$ Разумеется,  в (\ref{6.3}) речь идет о знакопеременных к.-а. мерах, иногда называемых зарядами. Здесь и ниже мы используем символику \cite[гл.\,2]{30}.  Отметим, что из (\ref{6.3}) легко следует свойство \cite[(2.2.9)]{30}
\beq{6.4}
\mu(\zer)=0 \ \ \forall\,\mu\in(\pr{add})[\mathcal{L}].
\eeq

Для некоторых преобразований (\ref{6.3}) будем использовать (\ref{2.43}).  Отметим, что при всяком выборе $\mu\in(\pr{add})[\mathcal{L}]$ и $\mathcal{K}\in(\pr{Fin})(\mathcal{L})$ определена функция $$(\mu(L))_{L\in\mathcal{K}}\in\mathbb{R}^{\mathcal{K}},$$  для которой (см. (\ref{2.42}), (\ref{2.43})) имеем следующее свойство
$$\sum\limits_{L\in\mathcal{K}}\mu(L)=\sum\limits_{j=1}^{|\mathcal{K}|}\mu(\Lambda_j)\in\mathbb{R} \ \ \ \forall\,(\Lambda_j)_{j\in\overline{1,|\mathcal{K}|}}\in(\pr{bi})[\overline{1,|\mathcal{K}|};\mathcal{K}];$$
к этому следует добавить то (см. (\ref{2.20})), что $$(\pr{bi})[\overline{1,|\mathcal{K}|};\mathcal{K}]\neq\zer.$$
В качестве $\mathcal{K}$  можно использовать разбиение из $\mathbf{D}(\mathbb{L},\mathcal{L}),$ где~ $\mathbb{L}\in~\mathcal{L}.$
\begin{proposition}
Справедливо вложение
\beq{6.5}
(\pr{add})[\mathcal{L}]\subset\biggl\{\mu\in\mathbb{R}^{\mathcal{L}}\mid \mu(\mathbb{L})=\sum\limits_{L\in\mathcal{K}}\mu(L) \ \ \forall\,\mathbb{L}\in\mathcal{L} \ \ \forall\,\mathcal{K}\in\mathbf{D}(\mathbb{L},\mathcal{L})\biggl\}.
\eeq
\end{proposition}
Д о к а з а т е л ь с т в о. Обозначим через $\mathbf{A}$ множество в правой части (\ref{6.5}). Требуется установить вложение  $(\pr{add})[\mathcal{L}]\subset\mathbf{A}.$ Данное вложение достаточно для всех наших последующих построений и мы ограничимся его проверкой, т.е. проверкой  (\ref{6.5}).
Выберем произвольно $\mu_*\in(\pr{add})[\mathcal{L}].$ Тогда $\mu_*:\mathcal{L}\rightarrow\mathbb{R}$ и согласно (\ref{6.3})
\beq{6.6}
 \mu_*(L)=\sum\limits_{i=1}^n\mu_*(L_i) \ \ \forall\,L\in\mathcal{L} \  \ \forall\,n\in\mathbb{N} \ \  \forall\,(L_i)_{i\in\overline{1,n}}\in\Delta_n(L,\mathcal{L}).
\eeq
Выберем произвольно множество $\mathbf{L}\in\mathcal{L}$ и разбиение $\mathfrak{L}\in\mathbf{D}(\mathbf{L},\mathcal{L}).$ Тогда определена функция
$(\mu_*(L))_{L\in\mathfrak{L}}\in\mathbb{R}^\mathfrak{L},$  т.е. $L\mapsto\mu_*(L):\mathfrak{L}\rightarrow\mathbb{R}.$ При этом определено значение $$\sum\limits_{L\in\mathfrak{L}}\mu_*(L)\in\mathbb{R},$$ для которого справедливо представление
\beq{6.7}
\sum\limits_{L\in\mathfrak{L}}\mu_*(L)=\sum\limits_{j=1}^{|\mathfrak{L}|}\mu_*(L_j) \ \ \forall\,(L_j)_{j\in\overline{1,|\mathfrak{L}|}}\in(\pr{bi})[\overline{1,|\mathfrak{L}|};\mathfrak{L}],
\eeq
где $(\pr{bi})[\overline{1,|\mathfrak{L}|};\mathfrak{L}]\neq\zer.$ Полагаем для краткости, что
\beq{6.8}
n\triangleq|\mathfrak{L}|,
\eeq
получая $n\in\mathbb{N}$ и при этом $(\pr{bi})[\overline{1,n};\mathfrak{L}]\neq\zer.$ С учетом этого выбираем произвольно
\beq{6.9}
(\Lambda_j)_{j\in\overline{1,n}}\in(\pr{bi})[\overline{1,n;}\mathfrak{L}].
\eeq
Тогда $\mu_*(\Lambda_j)\in\mathbb{R} \ \ \forall\,j\in\overline{1,n}.$  Согласно (\ref{6.7})--(\ref{6.9})
\beq{6.10}
\sum\limits_{L\in\mathfrak{L}}\mu_*(L)=\sum\limits_{j=1}^n\mu_*(\Lambda_j).
\eeq
Из (\ref{6.9}) с учетом свойства 2*) \  получаем, что $(\Lambda_j)_{j\in\overline{1,n}}\in\Delta_n(\mathbf{L},\mathcal{L}),$  а тогда из (\ref{6.6}) имеем равенство
\beq{6.11}
\mu_*(\mathbf{L})=\sum\limits_{j=1}^n\mu_*(\Lambda_j).
\eeq
 С учетом (\ref{6.10}) и (\ref{6.11}) получаем свойство $\mu_*(\mathbf{L})=\sum\limits_{L\in\mathfrak{L}}\mu_*(L).$

 Поскольку выбор $\mathbf{L}$ и $\mathfrak{L}$ был произвольным, установлено, что
 $$\mu_*(\mathbb{L})=\sum\limits_{L\in\mathcal{K}}\mu(L) \ \ \forall\,\mathbb{L}\in\mathcal{L} \ \ \forall\,\mathcal{K}\in\mathbf{D}(\mathbb{L},\mathcal{L}).$$
 Это означает, что $\mu_*\in\mathbf{A}.$ Итак, $(\pr{add})[\mathcal{L}]\subset\mathbf{A}.$ $\hfill\square$

Заметим, что утверждение предложения 7.1 может быть усилено, но для наших целей вполне достаточно очевидного следствия.
\begin{corollary}
Если $\mu\in(\pr{add})[\mathcal{L}],$ $\mathbb{L}\in\mathcal{L}$ и $\mathcal{K}\in\mathbf{D}(\mathbb{L},\mathcal{L}),$ то $$\mu(\mathbb{L})=\sum\limits_{L\in\mathcal{K}}\mu(L).$$
\end{corollary}
Доказательство непосредственно следует из предложения 7.1.

Введем в рассмотрение специальные множества в пространстве к.-а. мер:
\beq{6.12}
(\pr{add})_+[\mathcal{L}]\triangleq\{\mu\in(\pr{add})[\mathcal{L}]\mid 0\leq\mu(L) \ \forall\,L\in\mathcal{L}\},
\eeq
\beq{6.13}
\mathbb{P}(\mathcal{L})\triangleq\{\mu\in(\pr{add})_+[\mathcal{L}]\mid\mu(E)=1\},
\eeq
\beq{6.14}
\mathbb{T}(\mathcal{L})\triangleq\{\mu\in\mathbb{P}(\mathcal{L})\mid\forall\,L\in\mathcal{L} \ (\mu(L)=0)\vee(\mu(L)=1)\},
\eeq
\begin{multline}\label{6.15}
\mathbb{A}(\mathcal{L})\triangleq\biggl\{\mu\in(\pr{add})[\mathcal{L}]\mid\exists\,c\in[0,\infty[: \sum\limits_{i=1}^{n}|\mu(L_i)|\leqslant c \\ \forall\,n\in \mathbb{N} \ \forall\,(L_i)_{i\in\overline{1,n}}\in\Delta_n(E,\mathcal{L})\biggl\}.
\end{multline}
Из (\ref{6.3}), (\ref{6.12})--(\ref{6.15}) легко следует цепочка вложений (см. \cite[(2.2.11)]{30})
\beq{6.16}
\mathbb{T}(\mathcal{L})\subset\mathbb{P}(\mathcal{L})\subset(\pr{add})_+[\mathcal{L}]\subset\mathbb{A}(\mathcal{L})\subset(\pr{add})[\mathcal{L}]\subset\mathbb{R}^\mathcal{L}.
\eeq

\begin{zam}
В дальнейшем нам потребуются некоторые понятия, связанные с линейными пространствами. При этом будут рассматриваться только пространства в/з функций, а соответствующие линейные операции (а также умножение и порядок) в пространстве в/з функций, имеющих всякий раз общую область определения, определяются поточечно в духе \cite[\S\,1.6]{30}. Сейчас ограничимся совсем краткими напоминаниями, фиксируя в пределах настоящего замечания непустое множество $X$ (в качестве $X$ могут, в частности, использоваться множество $E$ и семейство~$\mathcal{L}$).

Если $\alpha\in\mathbb{R}$ и $f\in\mathbb{R}^X,$ то (в/з) функцию $\alpha f\in\mathbb{R}^X$ (иными словами, функцию $\alpha f:X\rightarrow \mathbb{R},$  называемую произведением $f$ на скаляр $\alpha$)   определяем условиями
\beq{6.17}
(\alpha f)(x)\triangleq\alpha f(x) \ \forall\,x\in X.
\eeq

Если же $f\in\mathbb{R}^X$ и $g\in\mathbb{R}^X,$ то $f+g\in\mathbb{R}^X$ определяется условиями
\beq{6.18}
(f+g)(x)\triangleq f(x)+g(x) \ \forall\,x\in X.
\eeq
Операции (\ref{6.17}) и (\ref{6.18}) можно комбинировать: при $\alpha\in\mathbb{R}, \beta\in\mathbb{R}, f\in\mathbb{R}^X$  и $g\in\mathbb{R}^X$
\beq{6.18'}
\alpha f+\beta g=(\alpha f(x)+\beta g(x))_{x\in X}\in\mathbb{R}^X;
\eeq
при этом, конечно, $(\alpha f+\beta g)(\overline{x})=\alpha f(\overline{x})+\beta g(\overline{x})$ при $\overline{x}\in X.$
Посредством (\ref{6.17}) определяется, в частности, функция, обратная (по сложению) по отношению к исходной: если $f\in\mathbb{R}^X,$ то $-f\in\mathbb{R}^X$ определяем правилом:
\beq{6.19}
-f\triangleq(-1)\cdot f;
\eeq
точку используем здесь в методических целях: функция $-f$ в (\ref{6.19}) есть $\alpha f$ при $\alpha=-1.$  При этом, конечно,
\beq{6.20}
-f=(-f(x))_{x\in X}\in\mathbb{R}^X.
\eeq

С помощью (\ref{6.18}), (\ref{6.19}) естественным образом определяется разность двух в/з функций: если $\varphi\in\mathbb{R}^X$ и $\psi\in\mathbb{R}^X,$ то $\varphi-\psi\in\mathbb{R}^X$ определяется правилом: $$\varphi-\psi\triangleq\varphi+(-\psi);$$ легко видеть (см. (\ref{6.18}), (\ref{6.20})), что $$(\varphi-\psi)(x)=\varphi(x)-\psi(x) \ \forall\,x\in X.$$
Операция (\ref{6.18}) естественным образом распространяется на произвольные упорядоченные конечные наборы в/з функций: если $n\in\mathbb{N}$ и
 $$(f_i)_{i\in\overline{1,n}}:\overline{1,n}\rightarrow\mathbb{R}^X,$$
 то, как обычно, полагаем, что
 \beq{6.21}
 \sum\limits_{i=1}^{n}f_i\triangleq\biggl(\sum\limits_{i=1}^{n}f_i(x)\biggl)_{x\in X}\in\mathbb{R}^X;
 \eeq
следовательно, имеем из (\ref{6.21}), что
$\left(\sum\limits_{i=1}^{n}f_i\right)(\overline{x})=\sum\limits_{i=1}^{n}f_i(\overline{x}) \ \forall\,\overline{x}\in X.$

На основе (\ref{6.21}) конструируется соответствующий обобщенный аналог (\ref{6.18'}), а именно: линейная комбинация соответствующего набора в/з функций. Итак, если $n\in\mathbb{N},$
$(\alpha_i)_{i\in\overline{1,n}}:\overline{1,n}\rightarrow\mathbb{R}, \ (f_i)_{i\in\overline{1,n}}:\overline{1,n}\rightarrow\mathbb{R}^X, \ \mbox{то}$
\beq{6.22}
\sum\limits_{i=1}^{n}\alpha_if_i\triangleq\sum\limits_{i=1}^{n}(\alpha_if_i)\in\mathbb{R}^X;
\eeq
в/з функция (\ref{6.22}), т.е. функция $\sum\limits_{i=1}^{n}\alpha_if_i: X\rightarrow\mathbb{R},$  имеет, очевидно, следующие значения:
$$\biggl(\sum\limits_{i=1}^{n}\alpha_if_i\biggl)(x)=\sum\limits_{i=1}^{n}\alpha_if_i(x) \ \forall\,x\in X.$$

Если в (\ref{6.22}) $\alpha_1\in[0,\infty[,\ldots,\alpha_n\in[0,\infty[$ и при этом $\sum\limits_{i=1}^{n}\alpha_i=1,$
то функцию (\ref{6.22}) называют \emph{выпуклой комбинацией} $f_1,\ldots,f_n.$

По аналогии с вышеупомянутыми определениями вводится произведение двух в/з функций: если $f\in\mathbb{R}^X$ и $g\in\mathbb{R}^X,$ то $f g\in\mathbb{R}^X$ определяется правилом
\beq{6.23}
(f g)(x)\triangleq f(x)g(x) \ \forall\,x\in X.
\eeq
Через $\mathcal{O}_X$ обозначаем в/з функцию на $X,$ тождественно равную нулю: $\mathcal{O}_X\!\in\!\mathbb{R}^X$ такова,~что $$\mathcal{O}_X(x)\triangleq 0 \ \forall\,x\in X.$$
Полагаем далее, что  \beq{6.23'}
\mathbb{B}(X)\triangleq\{f\in\mathbb{R}^X\mid\exists\,c\in[0,\infty[: \ |f(x)|\leqslant c \ \forall\,x\in X\}.
\eeq
Следуя \cite[(1.6.1)]{30}, полагаем при $H\in\pp'(\mathbb{R}^X),$ что
\begin{multline}\label{6.24}
(\pr{LIN})[H]\triangleq\{S\in\pp'(H)\mid (\alpha f\in S \ \forall\,\alpha\in\mathbb{R} \ \forall\,f\in S) \  \& \\ \& \ (f+g\in S \ \forall\,f\in S \ \forall\,g\in S)\};\end{multline}
элементы множества (возможно, пустого) (\ref{6.24})\,--- линейные подпространства $\mathbb{R}^X,$ содержащиеся в $H.$ В частности, определено $(\pr{LIN})[\mathbb{R}^X]$ (случай, когда в (\ref{6.24}) $H=\mathbb{R}^X$). Отметим, что
\beq{6.25'}
\mathbb{B}(X)\in (\pr{LIN})[\mathbb{R}^X].
\eeq

 Мы ввели линейное пространство ограниченных в/з функций на $X.$ Следуя \cite[(1.6.4)]{30},  при $H\in\pp'(\mathbb{R}^X)$ и $M\in\pp'(\mathbb{R}^X)$ рассматриваем множество
$$(\pr{LIN})[H|M]\triangleq\{S\in(\pr{LIN})[H]\mid M\subset S\},$$
тогда, в частности, $(\pr{LIN})[\mathbb{R}^X|M]\in\pp'((\pr{LIN})[\mathbb{R}^X]),$  поскольку выполняется условие $\mathbb{R}^X\in(\pr{LIN})[\mathbb{R}^X|M].$

С учетом этого определяем \emph{линейную оболочку} произвольного непустого п/м  $\mathbb{R}^X:$
\beq{6.25}
(\pr{sp})[M]\triangleq\bigcap\limits_{S\in(\pr{LIN})[\mathbb{R}^X| M]}S\in(\pr{LIN})[\mathbb{R}^X| M] \ \forall\,M\in\pp'(\mathbb{R}^X)
\eeq
(см. \cite[(1.6.5)--(1.6.7)]{30}). Хорошо известна структура линейной оболочки:
если $M\!\in\!\pp'(\mathbb{R}^X),$~то
\beq{6.26}
(\pr{sp})[M]=\biggl\{f\in\mathbb{R}^X\mid \exists\,n\in\mathbb{N} \ \exists\,(\alpha_i)_{i\in\overline{1,n}}\in\mathbb{R}^n \ \exists\,(f_i)_{i\in\overline{1,n}}\in M^n: f= \sum\limits_{i=1}^{n}\alpha_if_i\biggl\}.
\eeq
  Множество $K\!\in\!\pp'(\mathbb{R}^X)$ называется \emph{конусом}, если $\alpha f\!\in\! K  \forall\,\alpha\!\in\,]0,\infty[  \forall\,f\in~K.$ При этом $\pp'(H)\subset\pp'(\mathbb{R}^X) \ \forall\,H\in\pp'(\mathbb{R}^X).$ Полагаем при $H\in\pp'(\mathbb{R}^X),$ что
$$(\pr{cone})[H]\triangleq\{K\in\pp'(H)\mid\alpha f\in K \ \forall\,\alpha\in]0,\infty[ \ \forall\,f\in K\},$$
получая семейство всех конусов в $\mathbb{R}^X,$ содержащихся в $H.$ В частности, имеем
\beq{6.27}
(\pr{cone})[\mathbb{R}^X]=\{K\in\pp'(\mathbb{R}^X)\mid\alpha f\in K \ \forall\,\alpha\in]0,\infty[ \ \forall\,f\in K\}.
\eeq
\end{zam}

Напомним, что в замечании 7.1 в качестве $X$ могут использоваться множество $E$ и $\pi$-система  $\mathcal{L}.$
С учетом предложения 2.5.1 \cite{30}, имеем, что
$$(\pr{add})[\mathcal{L}]\in(\pr{LIN})[\mathbb{R}^\mathcal{L}]$$
и, кроме того, $\mathbb{A}(\mathcal{L})\in(\pr{LIN})[(\pr{add})[\mathcal{L}]],$ тогда, в частности,
\beq{6.27'}
\mathbb{A}(\mathcal{L})\in(\pr{LIN})[\mathbb{R}^\mathcal{L}]
\eeq
(см. \cite[следствие 2.5.1]{30}).  При этом (см. \cite[предложение 2.5.4]{30})
\beq{6.28}
(\pr{add})_+[\mathcal{L}]\in(\pr{cone})[\mathbb{A}(\mathcal{L})].
\eeq
С учетом (\ref{6.28}) получаем, что
\begin{multline}\label{6.29}
(\pr{add})^+[\mathcal{L};\mu]\triangleq\{\nu\in(\pr{add})_+[\mathcal{L}]\mid\forall\,L\in\mathcal{L} \ (\mu(L)=0)\Rightarrow(\nu(L)=0)\}\in \\ \in(\pr{cone})[(\pr{add})_+[\mathcal{L}]] \ \forall\,\mu\in(\pr{add})_+[\mathcal{L}].
\end{multline}

Заметим, что $\mathcal{O}_\mathcal{L}\in(\pr{add})^+[\mathcal{L};\mu] \ \forall\,\mu\in(\pr{add})_+[\mathcal{L}].$ Элементы множеств, определяемых в (\ref{6.29}) суть неотрицательные в/з к.-а. меры на $\mathcal{L}$ со свойством слабой абсолютной непрерывности относительно той или иной фиксированной к.-а. меры.  С учетом (\ref{6.15}) имеем, что при $\mu\in\mathbb{A}(\mathcal{L})$ корректно определяется (см. подробнее в \cite[\S 2.2]{30})
\beq{6.29'}
V_\mu\triangleq\sup\biggl(\bigl\{\,\sum\limits_{i=1}^n|\mu(L_i)|:n\in\mathbb{N}, \ (L_i)_{i\in\overline{1,n}}\in\Delta_n(E,\mathcal{L})\bigl\}\biggl)\in[0,\infty[.
\eeq

\begin{center} \textbf{Топологии на пространствах конечно-аддитивных мер} \end{center}

\ \ \ \ \ В настоящем разделе введем несколько характерных топологий на $\mathbb{A}(\mathcal{L})$ и на некоторых п/м $\mathbb{A}(\mathcal{L}).$ Отложим  рассмотрение наиболее важной $*$-слабой топологии $\mathbb{A}(\mathcal{L})$ и вернемся к нему после изложения конструкций интегрирования и представления линейных непрерывных функционалов на пространстве ярусных в/з функций. Сейчас рассматриваем $\mathbb{A}(\mathcal{L})$ как подпространство $\mathbb{R}^\mathcal{L}=\{\mathcal{L}\rightarrow\mathbb{R}\},$ что позволяет привлечь конструкции на основе тихоновских произведений (см. \cite[(4.2.8), (4.2.9)]{36}.) Но в настоящем пособии мы не будем воспроизводить все эти полезные сами по себе построения и ограничимся более краткими и по возможности непосредственными определениями (см. в этой связи также \cite[\S\S 2.6,\,4.6]{41}).

Сначала совсем кратко напомним два варианта тихоновской степени  (имеется в виду тихоновское произведение экземпляров одного и того же пространства). Рассматриваем при этом $\mathbb{R}^\mathcal{L}.$ Мы полагаем, что
$\tau_{\mathbb{R}}\triangleq\{G\in\pp(\mathbb{R})\mid\forall\,x\in G \ \exists\,\varepsilon\in]0,\infty[ \ : \ ]x-\varepsilon,x+\varepsilon[ \ \subset G\},$ $\tau_\partial\triangleq\pp(\mathbb{R});$ при этом, конечно, $\tau_{\mathbb{R}}\in(\pr{top})[\mathbb{R}]$ и $\tau_\partial\in(\pr{top})[\mathbb{R}].$

Полагаем также \cite[(2.6.27)]{41} при $\mu\in\rr^\mathcal{L}, \ \mathcal{K}\in\pr{Fin}(\ml)$ и $\varepsilon\bn,$ что
\beq{6.30}
\mathbb{N}_{\ml}(\mu,\mathcal{K},\varepsilon)\triangleq\left\{\nu\in\rr^\ml\mid|\mu(L)-\nu(L)|<\varepsilon \ \forall\,L\in\mathcal{K}\right\}.
\eeq
Тогда, как легко видеть, имеем \cite[(2.6.30)]{41}, что
\begin{multline}\label{6.31}
\otimes^\ml(\tau_{\rr})\triangleq\bigl\{G\in\pp(\rr^\ml)\mid\forall\,\mu\in G \ \exists\,\mathcal{K}\in\pr{Fin}(\ml) \ \exists\,\varepsilon\bn \ : \\ \mathbb{N}_\ml(\mu,\mathcal{K},\varepsilon)\subset  G\bigl\}\in (\pr{top})[\rr^\ml],
\end{multline}
причем $(\rr^\ml,\otimes^\ml(\tau_{\rr}))$ есть $T_2$-пространство. Заметим, что
\beq{6.32}
\forall\,\mu\in\rr^\ml \ \forall\,H\in N_{\otimes^\ml(\tau_{\rr})}(\mu) \ \exists\,\mathcal{K}\in\pr{Fin}(\ml) \ \exists\,\varepsilon\bn \ : \mathbb{N}_{\ml}(\mu,\mathcal{K},\varepsilon)\subset H;
\eeq
 см. \cite[c.\,59]{41}.  Разумеется, из  (\ref{6.30}) и (\ref{6.31}) вытекает, что $\forall\,\mu\in\rr^\mathcal{L}, \ \forall\,\mathcal{K}\in\pr{Fin}(\ml)$  $\forall\,\varepsilon\bn$
 \beq{6.33}
 \mathbb{N}_{\ml}(\mu,\mathcal{K},\varepsilon)\in\otimes^\ml(\tau_{\rr}).
 \eeq

В виде $(\rr^\ml,\otimes^\ml(\tau_{\rr}))$ имеем тихоновскую степень $(\rr,\tau_\rr)$ при использовании $\ml$ в качестве индексного множества;  (\ref{6.31})\,--- так называемая \emph{топология поточечной сходимости}.
Теперь введем другой вариант тихоновской степени, полагая сначала, что
\beq{6.34}
\mathbb{N}_{\ml}^{(\partial)}(\mu,\mathcal{K})\triangleq\left\{\nu\in\rr^\ml\mid\mu(L)=\nu(L) \ \forall\,L\in\mathcal{K}\right\} \ \forall\,\mu\in\rr^\mathcal{L} \ \forall\,\mathcal{K}\in\pr{Fin}(\ml).
\eeq
Тогда, как легко проверить, имеем, что
\beq{6.35}
\otimes^\ml(\tau_\partial)\triangleq\left\{G\in\pp(\rr^\ml)\mid\forall\,\mu\in G \ \exists\,\mathcal{K}\in\pr{Fin}(\ml): \ \mathbb{N}_{\ml}^{(\partial)}(\mu,\mathcal{K})\subset G \right\}\in(\pr{top})[\rr^\ml],
\eeq
причем $(\rr^\ml,\otimes^\ml(\tau_\partial))$  есть $T_2$-пространство.

Заметим, что из  (\ref{6.30}) и  (\ref{6.34}) следует, в частности, свойство
\beq{6.36}
\mathbb{N}_\ml^{(\partial)}(\mu,\mathcal{K})\subset\mathbb{N}_\ml(\mu,\mathcal{K},\varepsilon) \  \ \forall\,\mu\in\rr^\ml \  \ \forall\,\mathcal{K}\in\pr{Fin}(\ml) \ \ \forall\,\varepsilon\bn.
\eeq
Из  (\ref{6.31}),  (\ref{6.35}) и  (\ref{6.36}) получаем очевидное вложение
\beq{6.37}
\otimes^\ml(\tau_\rr)\subset\otimes^\ml(\tau_\partial).
\eeq
Из  (\ref{4.0}) и  (\ref{6.37}), в свою очередь, следует, что
\beq{6.38}
\otimes^\ml(\tau_\rr)|_Z\subset\otimes^\ml(\tau_\partial)|_Z \ \ \forall\,Z\in\pp(\rr^\ml).
\eeq
В  (\ref{6.38}) имеем свойство сравнимости подпространств тихоновских произведений. Для наших целей полезны, в частности, подпространства с <<единицей>>  $\mathbb{A}(\ml);$ им соответствуют топологии
\beq{6.39}
\tau_\otimes(\ml)\triangleq\otimes^\ml(\tau_\rr)|_{\mathbb{A}(\ml)}\in(\pr{top})[\mathbb{A}(\ml)],
\eeq
\beq{6.40}
\tau_0(\ml)\triangleq\otimes^\ml(\tau_\partial)|_{\mathbb{A}(\ml)}\in(\pr{top})[\mathbb{A}(\ml)].
\eeq
Из (\ref{6.38})--(\ref{6.40}) вытекает, что справедливо
\beq{6.41}
\tau_\otimes(\ml)\subset\tau_0(\ml).
\eeq

Итак,  топология (\ref{6.39}) слабее топологии (\ref{6.40}). Наряду с (\ref{6.41}) полезно отметить, что (см. (\ref{6.38})) при $Z\in\pp(\mathbb{A}(\ml))$
\beq{6.42}
\otimes^\ml(\tau_\rr)|_Z\subset\otimes^\ml(\tau_\partial)|_Z;
\eeq
в частности, (\ref{6.42}) справедливо при  $Z\in\pp((\pr{add})_+[\ml]).$ С учетом хорошо известного свойства транзитивности операции перехода к подпространству (см. \cite[2.1.2]{39})  при $H\in\pp(\mathbb{A}(\ml))$ имеем, что
$$\tau_\otimes(\ml)|_H=\otimes^\ml(\tau_\rr)|_H\in(\pr{top})[H],$$
$$\tau_0(\ml)|_H=\otimes^\ml(\tau_\partial)|_H\in(\pr{top})[H].$$
Из (\ref{6.42}) следуют свойства
\beq{6.43}
\tau_\otimes(\ml)|_H\subset\tau_0(\ml)|_H.
\eeq

В частности, для топологий \cite[с.\,80]{36}
\beq{6.44}
\tau_\otimes^+(\ml)\triangleq\tau_\otimes(\ml)|_{(\pr{add})_+[\ml]}=\otimes^\ml(\tau_\rr)|_{(\pr{add})_+[\ml]}\in(\pr{top})[(\pr{add})_+[\ml]],
\eeq
\beq{6.45}
\tau_0^+(\ml)\triangleq\tau_0(\ml)|_{(\pr{add})_+[\ml]}=\otimes^\ml(\tau_\partial)|_{(\pr{add})_+[\ml]}\in(\pr{top})[(\pr{add})_+[\ml]]
\eeq
имеем очевидное следствие (\ref{6.43})
\beq{6.46}
\tau_\otimes^+(\ml)\subset\tau_0^+(\ml).
\eeq

В дальнейшем свойства (\ref{6.43})--(\ref{6.46}) будут дополнены аналогами,  в которых используется так называемая $*$-слабая топология, для введения которой потребуются уже понятия из теории топологических векторных пространств. Эти понятия будут приведены ниже в краткой и предельно конкретизированной форме после изложения простейшей конструкции интеграла по к.-а. мере (подробнее см. в \cite[гл.\,3,4]{30}; общие вопросы интегрирования по к.-а. мере см. в \cite[гл.\,III]{43}).

\section{Ступенчатые, ярусные  и измеримые функции}\setcounter{equation}{0}\setcounter{proposition}{0}\setcounter{zam}{0}

\ \ \ \ \ Сохраняем предположения предыдущего раздела: фиксируем широко понимаемое ИП $(E,\ml),$ $E\neq\zer,$ $\ml\in\pi[E]$ (по мере надобности на выбор $\ml$ будут накладываться дополнительные условия).  Напомним, что $\rr^E$ есть множество всех в/з функций на множестве~$E:$ $$\rr^E=\{E\rightarrow\rr\}.$$

Тогда при $f\in\rr^E$ и $x\in E$ в виде $f(x)\in\rr$ имеем значение функции $f$ в точке $x.$  Будем использовать (\ref{6.23'}), (\ref{6.25'}), (\ref{6.26}) и (\ref{6.27}) при $X=E.$  В частности, с учетом (\ref{6.25'}) мы вводим обычную $\sup$-норму линейного пространства $\mathbb{B}(E),$  полагая, что (см. (\ref{6.23'}))
\beq{7.1}
\|f\|\triangleq\sup(\{|f(x)|:x\in E\}) \ \forall\,f\in\mathbb{B}(E).
\eeq

 Следовательно, $\|\cdot\|\triangleq(\|f\|)_{f\in\mathbb{B}(E)}$ является требуемой нормой: $\|\cdot\|$ есть неотрицательная в/з функция на $\mathbb{B}(E)$ (функционал на $\mathbb{B}(E)$), для которого\\
1) \ \ \ $\forall\,f\in\mathbb{B}(E) \ \ (\|f\|=0)\Leftrightarrow(f=\mathcal{O}_E);$ \\
2) \ \ \ $\|\alpha f\|=|\alpha|\cdot\|f\| \ \ \forall\,\alpha\in\rr \ \forall\,f\in\mathbb{B}(E);$\\
3) \ \ \ $\|f+g\|\leqslant\|f\|+\|g\| \ \forall\,f\in\mathbb{B}(E) \ \forall\,g\in\mathbb{B}(E).$

\begin{zam}
Строго говоря, норма $\|\cdot\|$ определена при фиксации множества $E,$ что для наших построений вполне достаточно.
\end{zam}
Отметим полезное следствие: если $f\in\mathbb{B}(E)$ и $g\in\mathbb{B}(E),$ то $$\bigl|\,\|f\|-\|g\|\,\bigl|\leqslant\|f-g\|.$$

Разумеется, пара $(\mathbb{B}(E),\|\cdot\|)$ есть нормированное пространство. Важно отметить полезную связь с понятием равномерной сходимости. Напомним последнее понятие, ограничиваясь случаем сходимости функций из $\mathbb{B}(E)$ (более общий случай см. в \cite[(2.6.9)]{30}):  если
\beq{7.2}
(f_i)_{i\in\mathbb{N}}:\mathbb{N}\rightarrow\mathbb{B}(E)
\eeq
и $f\in\mathbb{B}(E),$ то
$$\left((f_i)_{i\in\mathbb{N}}\rightrightarrows f\right)\Leftrightarrow\left(\forall\,\varepsilon\in]0,\infty[ \ \exists\,n\in\mathbb{N}: |f_j(x)-f(x)|<\varepsilon \ \ \forall\,j\in\overrightarrow{n,\infty} \ \ \forall\,x\in E\right);$$
при этом, с учетом (\ref{7.1}) имеем, что
\beq{7.3}
\left((f_i)_{i\in\mathbb{N}}\rightrightarrows f\right)\Leftrightarrow\left((\|f_i-f\|)_{i\in\mathbb{N}}\rightarrow 0\right).
\eeq
Отметим, используя (\ref{7.3}), важное свойство полноты $(\mathbb{B}(E),\|\cdot\|):$ если $(f_i)_{i\in\mathbb{N}}$ соответствует (\ref{7.2}), то
\begin{multline}\label{7.4}
\left(\forall\,\varepsilon\in]0,\infty[ \ \exists\,n\in\mathbb{N}: \ \|f_k-f_l\|<\varepsilon \ \forall\,k\in\overrightarrow{n,\infty} \ \forall\,l\in\overrightarrow{n,\infty}\right)\Rightarrow \\ \Rightarrow\left(\exists\,f\in\mathbb{B}(E):(f_i)_{i\in\mathbb{N}}\rightrightarrows f\right).
\end{multline}

Итак, $(\mathbb{B}(E),\|\cdot\|)$ есть банахово (полное, нормированное) пространство. К этому следует добавить следующее легкопроверяемое свойство
\cite[(2.6.26)]{30}:
$$fg\in\mathbb{B}(E) \ \forall\,f\in\mathbb{B}(E) \ \forall\,g\in\mathbb{B}(E).$$
Таким образом, при $f\in\mathbb{B}(E)$ и $g\in\mathbb{B}(E)$ определено значение $\|fg\|\in[0,\infty[,$ для которого \cite[(2.6.27)]{30}
$$\|fg\|\leqslant\|f\|\cdot\|g\|.$$
Условимся о следующем соглашении: если $A\in\pp(E),$ то полагаем, что
\beq{7.5}
\chi_A[E]:E\rightarrow\rr
\eeq
определяется условиями
\beq{7.6}
\left(\chi_A[E](x)\triangleq 1 \ \forall\,x\in A\right) \ \& \ \left(\chi_A[E](y)\triangleq 0 \ \forall\,y\in E\setminus A\right).
\eeq

Функцию  (\ref{7.5}), (\ref{7.6})  называем \emph{индикатором множества} $A.$ Учитывая то, что множество $E$ в последующих построениях настоящей главы будет фиксированным, принимаем соглашение
\beq{7.7}
\chi_A\triangleq\chi_A[E] \ \forall\,A\in\pp(E),
\eeq
имеющее целью упростить обозначения. С учетом (\ref{6.23'}), (\ref{7.5})--(\ref{7.7}) имеем следующее очевидное свойство
\beq{7.8}
\chi_A\in\mathbb{B}(E) \ \forall\,A\in\pp(E).
\eeq
При этом, как легко видеть, отображение $$\Sigma\mapsto\chi_\Sigma:\pp(E)\rightarrow\rr^E$$
взаимно однозначно (инъективно): $\forall\,A\in\pp(E) \ \forall\,B\in\pp(E)$
\beq{7.9}
(\chi_A=\chi_B)\Rightarrow(A=B).
\eeq
Здесь же отметим, что с учетом (\ref{6.23}), (\ref{7.6}) и (\ref{7.7}) легко проверяется свойство
\beq{7.10}
\chi_{A\cap B}=\chi_A\chi_B \ \forall\,A\in\pp(E) \ \forall\,B\in\pp(E).
\eeq
Напомним, что $\ml\subset\pp(E),$  а тогда при $L\in\ml$ функция $\chi_L:E\rightarrow \rr$ определяется условиями (см. (\ref{7.6}), (\ref{7.7}))
\beq{7.11}
\left(\chi_L(x)\triangleq 1 \ \forall\,x\in L\right) \ \& \ \left(\chi_L(y)\triangleq0 \ \forall\,y\in E\setminus L\right).
\eeq
При этом, конечно, $\chi_L\in\mathbb{B}(E) \ \forall\,L\in\ml.$  Отметим, что (см.  \cite[(2.7.2), предложение 2.7.2]{30})
\begin{multline}\label{7.12}
B_0(E,\ml)\triangleq \bigl\{f\in\rr^E\mid\exists\,n\in\mathbb{N} \ \exists\,(\alpha_i)_{i\in\overline{1,n}}\in\rr^n \ \exists\,(L_i)_{i\in\overline{1,n}}\in\Delta_n(E,\ml): \\ f=\sum\limits_{i=1}^n\alpha_i\chi_{L_i}\bigl\}\in (\pr{LIN})[\mathbb{B}(E)].
\end{multline}
Итак, $B_0(E,\ml)$ является (линейным) многообразием в $\mathbb{B}(E).$ Заметим, что согласно (\ref{7.12}) все постоянные на $E$ в/з функции содержатся в множестве  $B_0(E,\ml), $ элементы которого  называем \emph{ступенчатыми} (относительно $\ml$) \emph{функциями}. Из (\ref{7.12}) следует также, что
\beq{7.13}
fg\in B_0(E,\ml) \ \forall\,f\in B_0(E,\ml) \ \forall\,g\in B_0(E,\ml)
\eeq
(см. \cite[предложение 2.7.4]{30}).  Отметим здесь же, что
\beq{7.14}
B_0^+(E,\ml)\triangleq\{f\in B_0(E,\ml)\mid 0\leqslant f(x) \ \forall\,x\in E\}\in(\pr{cone})[B_0(E,\ml)]
\eeq
 есть конус неотрицательных ступенчатых в/з функций на  $E.$

 Теперь, учитывая то, что $\mathbb{B}(E)$ оснащено нормой (и, как следствие, метрикой), можно ввести обычную топологию равномерной сходимости на $\mathbb{B}(E),$ используя (\ref{7.3}). Такое построение последовательно проведено в \cite[c.\,109,110]{30} (см., в частности, \cite[(2.7.23), (2.7.24)]{30}). Однако сейчас мы, имея более конкретные цели и используя (\ref{7.3}), введем нужное пространство \emph{ярусных} (в смысле $(E,\ml)$) в/з функций, следуя \cite[(2.7.25)]{30}:  полагаем, что
 \begin{multline}\label{7.15}
 B(E,\ml)\triangleq\left\{f\in\mathbb{B}(E)\mid\exists\,(f_i)_{i\in\mathbb{N}}\in B_0(E,\ml)^{\mathbb{N}}:(f_i)_{i\in\mathbb{N}}\rightrightarrows f\right\}= \\ =\left\{f\in\mathbb{B}(E)\mid\exists\,(f_i)_{i\in\mathbb{N}}\in B_0(E,\ml)^{\mathbb{N}}:(\|f_i-f\|)_{i\in\mathbb{N}}\rightarrow 0\right\},
 \end{multline}
 получая, конечно, линейное пространство (см. \cite[предложение 2.7.5]{30}):

В виде $\mathbf{B}_*(\ml)\triangleq\{f\in B(E,\ml), \ \|f\|\leqslant 1\}$ имеем замкнутый единичный шар $B(E,\ml)$ с центром в начале координат (в <<нуле>>). Более того, $B(E,\ml)$  в оснащении нормой подпространства $(\mathbb{B}(E),\|\cdot\|)$ само является банаховым пространством (подробнее см. в \cite[c.\,111]{30}), для которого
$$B_0(E,\ml)\in(\pr{LIN})[B(E,\ml)]$$
и, в частности, $B_0(E,\ml)\subset B(E,\ml).$
Здесь же отметим, используя свойства конкретной нормы (\ref{7.1}), что
\beq{7.17}
fg\in B(E,\ml) \ \forall\,f\in B(E,\ml) \ \forall\,g\in B(E,\ml).
\eeq
Свойство (\ref{7.17}) наследуется от аналогичного свойства (\ref{7.13}), действующего в многообразии ступенчатых функций (см. \cite[предложение 2.7.7]{30}).
Сейчас мы отдельно коснемся свойств, имеющих место в случае, когда  $(E,\ml)$ есть ИП с полуалгеброй множеств.
Полагаем, если не оговорено противное, до конца настоящего раздела, что
\beq{7.18}
\ml\in\Pi[E],
\eeq
получая в виде $(E,\ml)$ ИП с полуалгеброй множеств. Тогда \cite[предложение 2.7.3]{30} имеем, что
\beq{7.19}
\chi_L\in  B_0(E,\ml) \ \forall\,L\in\ml.
\eeq

Доказательство, соответствующее \cite[предложение 2.7.3]{30}, легко следует из (\ref{5'.5});  рекомендуем читателю провести его самостоятельно. С учетом (\ref{6.26}) и (\ref{7.19}) получаем свойство (см. \cite[следствие 2.7.1]{30})
\beq{7.20}
B_0(E,\ml)=(\pr{sp})[\{\chi_L: L\in\ml\}];
\eeq
в связи с (\ref{7.20}) читателю рекомендуется ознакомиться с \cite[замечание 2.7.1]{30}. Итак, в рассматриваемом сейчас случае (\ref{7.18}) $B_0(E,\ml)$ есть линейная оболочка множества индикаторов всевозможных множеств из $\ml.$
Из (\ref{7.17}) и (\ref{7.19}) имеем (при условии (\ref{7.18})), что \cite[предложение 2.7.8]{30}
\beq{7.21}
f\chi_L\in B(E,\ml) \ \forall\,f\in B(E,\ml) \ \forall\,L\in\ml.
\eeq

Свойство (\ref{7.21}) будет неоднократно использоваться в дальнейшем. Поясним его: при $f\in B(E,\ml)$  и $L\in\ml$ функция $f\chi_L$ такова, что
\beq{7.22}
((f\chi_L)(x)=f(x) \ \forall\,x\in L) \ \& \ \left((f\chi_L)(\widetilde{x})=0 \ \forall\,\widetilde{x}\in E\setminus L\right).
\eeq
Из (\ref{7.22}) видно, что $f\chi_L$  является <<частью>> $f,$ выделяемой множеством $L,$ т.е. своеобразной <<срезкой>> $f.$
Отметим, что в стандартной теории меры важную роль играют измеримые функции (для наших целей достаточны функции ограниченные). Полагаем в этой связи
\beq{7.23'}
(\mathbb{B}-\pr{Meas})[E;\ml]\triangleq\left\{f\in\mathbb{E}\mid f^{-1}(\,]-\infty;c[)\in\ml \ \forall\,c\in\rr\right\}
\eeq
(здесь использовали простейший вариант определения свойства измеримости, а именно: измеримость множеств Лебега). Напомним, что \cite[теорема 2.8.1]{30}
\beq{7.24'}
(\ml\in(\sigma-\pr{alg})[E])\Rightarrow(B(E,\ml)=(\mathbb{B}-\pr{Meas})[E;\ml]).
\eeq
\begin{zam}
В связи с (\ref{7.24'}) отметим (см. (\ref{5'.19})), что уже при $\ml\in(\pr{alg})[E]$ возможен случай, когда  $B(E,\ml)\neq(\mathbb{B}-\pr{Meas})[E;\ml]$ (см. \cite[c.\,119-120]{30}), хотя $(\mathbb{B}-\pr{Meas})[E;\ml]\subset B(E,\ml),$ см. \cite[предложение 2.8.3]{30}.
\end{zam}

\begin{center} \textbf{Линейные ограниченные функционалы на пространстве ярусных функций} \end{center}

\ \ \ \ \ Относительно $(E,\ml)$ сохраняем предположение $\ml\in\pi[E],$ где $E\neq\zer$ (дополнительные условия на $\ml$ накладываются по мере надобности). Напомним, что (см. (\ref{7.14}))
\beq{7.24*}
B_0^+(E,\ml)\triangleq\{f\in B_0(E,\ml)\mid0\leqslant f(x) \ \forall\,x\in E\}.
\eeq
Полагаем в дальнейшем, что
\beq{7.24**}
B^+(E,\ml)\triangleq\{f\in B(E,\ml)\mid0\leqslant f(x) \ \forall\,x\in E\}.
\eeq
При этом  $B_0^+(E,\ml)\in(\pr{cone})[B_0^+(E,\ml)],$  $B^+(E,\ml)\in(\pr{cone})[B^+(E,\ml)] $ и  $$B_0^+(E,\ml)\subset B^+(E,\ml).$$
Отметим также, что
$$\rr^{B(E,\ml)}=\{B(E,\ml)\rightarrow\rr\}$$
есть множество всех функционалов на линейном пространстве $B(E,\ml).$ Среди последних выделяем линейные функционалы, образующие в своей совокупности пространство, алгебраически сопряженное к $B(E,\ml)$ (при определении линейных операций на $B(E,\ml)$ используем (\ref{6.17}), (\ref{6.18}), а также следствия этих соотношений). Пусть
\begin{multline}\label{7.23}
B'(E,\ml)\triangleq\bigl\{\varphi\in\rr^{B(E,\ml)}\mid(\varphi(\alpha f)=\alpha\varphi(f) \ \forall\,\alpha\in\rr \ \ \forall\,f\in B(E,\ml)) \ \& \\ \& \ (\varphi(f+g)=\varphi(f)+\varphi(g) \ \forall\,f\in B(E,\ml) \ \forall\,g\in B(E,\ml))\bigl\}.
\end{multline}

Таким образом, мы ввели пространство всевозможных линейных функционалов на $B(E,\ml).$ Ясно, что
\begin{multline*}
\varphi(\alpha f+\beta g)=\alpha\varphi(f)+\beta\varphi(g) \ \forall\,\varphi\in B'(E,\ml) \ \forall\,\alpha\in\rr \  \forall\,\beta\in\rr   \ \forall\,f\in B(E,\ml) \\ \forall\,g\in B(E,\ml).
\end{multline*}

Прочие очевидные следствия (\ref{7.23}) предоставляем читателю в качестве упражнения.

Учитывая теперь возможности, связанные с оснащением $B(E,\ml)$ нормой, индуцированной из $(\mathbb{B}(E),\|\cdot\|),$ полагаем, что в общем случае $\ml\in\pi[E]$
\beq{7.24}
B^*(E,\ml)\triangleq\{\varphi\in B'(E,\ml)\mid\exists\,c\in[0,\infty[:|\varphi(s)|\leqslant c\|s\| \ \forall\,s\in B(E,\ml)\},
\eeq
получая пространство линейных ограниченных (линейных непрерывных) функционалов на $B(E,\ml).$ Данное пространство (\ref{7.24}), однако, представляет для нас основной интерес при условии (\ref{7.18}), а потому в дальнейшем мы, как правило, ограничиваемся этим случаем. В общем  случае имеем в  (\ref{7.24}) пространство, топологически сопряженное к $B(E,\ml)$ в оснащении последнего $\sup$\,-нормой.
При этом  согласно (\ref{7.23}), (\ref{7.24}) имеем при $\ml\in\pi[E],$ что
$$B^*(E,\ml)\subset\rr^{B(E,\ml)};$$
в пространстве $\rr^{B(E,\ml)}$  линейные операции определены посредством (\ref{6.17}), (\ref{6.18}) при  $X=B(E,\ml)$ (получить очевидные следствия данного варианта (\ref{6.17}), (\ref{6.18}) предоставляем читателю в качестве упражнения), при этом \cite[предложение 3.5.1]{30}
\beq{7.25}
B^*(E,\ml)\in(\pr{LIN})[\rr^{B(E,\ml)}].
\eeq
Учитывая (\ref{7.25}), в упомянутом общем случае $\ml\in\pi[E]$ введем отображение
\beq{7.26}
\varphi\mapsto\|\varphi\|^*: B^*(E,\ml)\rightarrow[0,\infty[
\eeq
по следующему правилу, учитывающему (\ref{7.24}): если  $\psi\in B^*(E,\ml),$ то
\beq{7.27}
\|\psi\|^*\triangleq\sup(\{|\psi(s)|:s\in\mathbf{B}_*(\ml)\}).
\eeq
Напомним, что посредством (\ref{7.26}), (\ref{7.27}) определена \cite[предложение 3.5.4]{30}  традиционная норма $B^*(E,\ml)$ как пространства топологически сопряженного к $B(E,\ml).$
\begin{zam}
Отметим, что норма (\ref{7.26}), (\ref{7.27}) определена при условии фиксации  $\pi$-системы $\ml,$ что вполне достаточно для всех наших ближайших целей. В дальнейших общих теоретических конструкциях будем сохранять для рассматриваемого ИП обозначение $(E,\ml),$ что позволяет распространять (\ref{7.26}), (\ref{7.27}) и на последующие построения.
\end{zam}

Ставим своей ближайшей целью дать описание пространства $B^*(E,\ml)$ с упомянутой нормой.

\section{Интегрирование ярусных функций  по конечно-\\аддитивной мере ограниченной вариации}\setcounter{equation}{0}\setcounter{proposition}{0}\setcounter{zam}{0}\setcounter{corollary}{0}\setcounter{definition}{0}
   \ \ \ \ \  В настоящем разделе в краткой форме излагается конструкция интегрирования по  к.-а. мере ограниченной вариации, на основе которой далее будет построено представление пространства  (\ref{7.25}) для ИП с полуалгеброй множеств. Однако сначала мы рассмотрим несколько более общий случай ИП, фиксируя непустое множество $E$ и $\pi$-систему $\ml~\in~\pi[E];$  некоторые детали для случая ИП с полуалгеброй множеств обсудим в конце данного раздела.

\begin{center} \textbf{Элементарный интеграл на пространстве ступенчатых функций} \end{center}
  \ \ \ \ \ \ Рассмотрим вопросы интегрирования ступенчатых в/з функций по произвольной   к.-а. мере (ограниченность вариации здесь не требуется). Получающийся интеграл будем называть \emph{элементарным}. По сути дела такой интеграл~--- <<хорошо организованная>> конечная сумма. Последнее суждение поясняется предложением, доказательство которого приведено в \cite[c.\,122-123]{30}  Мы ориентируемся при этом на (\ref{7.12}),  имея в виду использование операции (\ref{6.22}).

\begin{proposition}
Если $m\in\mathbb{N},$ $(\alpha_i)_{i\in\overline{1,m}}\in\rr^m,$ $(L_i)_{i\in\overline{1,m}}\in\Delta_m(E,\ml),$  $n\in\mathbb{N},$ $(\beta_j)_{j\in\overline{1,n}}\in\rr^n$ и  $(\Lambda_j)_{j\in\overline{1,n}}\in\Delta_n(E,\ml)$  таковы, что $$\sum\limits_{i=1}^m\alpha_i\chi_{L_i}=\sum\limits_{j=1}^n\beta_j\chi_{\Lambda_j},$$
то справедливо следующее свойство: $$\sum\limits_{i=1}^m\alpha_i\mu(L_i)=\sum\limits_{j=1}^n\beta_j\mu(\Lambda_j) \ \forall\,\mu\in(\pr{add})[\ml].$$
\end{proposition}

\begin{corollary}
Если $f\in B_0(E,\ml)$ и $\mu\in(\pr{add})[L],$ то  $\exists\,!\,c\in\rr:\forall\,n\in\mathbb{N}$ $\forall\,(\alpha_i)_{i\in\overline{1,n}}\in\rr^n,$
$\forall\,(L_i)_{i\in\overline{1,n}}\in\Delta_n(E,\ml)$
$$\biggl(f=\sum\limits_{i=1}^n\alpha_i\chi_{L_i}\biggl)\Rightarrow\biggl(\sum\limits_{i=1}^n\alpha_i\mu(L_i)=c\biggl).$$
\end{corollary}
В связи с данным следствием см. \cite[c.\,123]{30}. Теперь корректно следующее
\begin{definition}
Если  $f\in B_0(E,\ml)$ и $\mu\in(\pr{add})[\ml],$ то полагаем, что
\beq{8.1}
\int\limits_E^{(\pr{el})} f\,d\mu\in\rr
\eeq
есть по определению  такое единственное число, что  $\forall\,n\in\mathbb{N}$ $\forall\,(\alpha_i)_{i\in\overline{1,n}}\in\rr^n$
$\forall\,(L_i)_{i\in\overline{1,n}}\in\Delta_n(E,\ml)$
$$\biggl(f=\sum\limits_{i=1}^n\alpha_i\chi_{L_i}\biggl)\Rightarrow\biggl(\,\int\limits_E^{(\pr{el})} f\,d\mu=\sum\limits_{i=1}^n\alpha_i\mu(L_i)\biggl);$$
мы называем число (\ref{8.1}) элементарным интегралом ступенчатой функции $f$ по к.-а. мере $\mu.$
\end{definition}

В связи с определением напомним, что согласно (\ref{7.12})
\beq{8.2}
\sum\limits_{i=1}^n\alpha_i\chi_{L_i}\!\!\in\! B_0(E,\ml)  \,\forall\,n\!\!\in\!\mathbb{N} \, \forall\,(\alpha_i)_{i\in\overline{1,n}}\in\rr^n \ \forall\,(L_i)_{i\in\overline{1,n}}\in\Delta_n(E,\ml).
\eeq
Поэтому (\ref{8.1}) применимо к ступенчатым функциям, определенным посредством (\ref{8.2}). Иными словами, при $n\in\mathbb{N}, (\alpha_i)_{i\in\overline{1,n}}\in\rr^n, (L_i)_{i\in\overline{1,n}}\in\Delta_n(E,\ml)$ и $\mu\in(\pr{add})[\ml]$  определен элементарный интеграл
$$\int\limits_E^{(\pr{el})}\sum\limits_{i=1}^n\alpha_i\chi_{L_i} \,d\mu\in\rr;$$
разумеется, в этом случае имеем равенство
\beq{8.3}
\int\limits_E^{(\pr{el})}\sum\limits_{i=1}^n\alpha_i\chi_{L_i} \,d\mu=\sum\limits_{i=1}^n\alpha_i\mu(L_i).
\eeq
Собственно говоря, (\ref{8.3}) есть другая форма записи определения 9.1. Поскольку $B_0(E,\ml)$\,--- линейное многообразие, то при $\mu\in(\pr{add})[\ml]$ определены элементарные интегралы
\begin{multline}\label{8.4}
\biggl(\,\int\limits_E^{(\pr{el})}\alpha f\,d\mu\in\rr \ \forall\,\alpha\in\rr \ \forall\,f\in B_0(E,\ml)\biggl)\ \& \\ \& \biggl(\,\int\limits_E^{(\pr{el})}( f+g)\,d\mu\in\rr \  \forall\,f\in B_0(E,\ml) \ \forall\,g\in B_0(E,\ml)\biggl).
\end{multline}
С учетом (\ref{7.13}) имеем также, что при $f\in B_0(E,\ml)$ определены (см. раздел~6) значения
\begin{multline}\label{8.5}
\biggl(\,\int\limits_E^{(\pr{el})}f\,d(\alpha\mu)\in\rr \ \forall\,\alpha\in\rr \  \forall\,\mu\in(\pr{add})[\ml]\biggl)\ \& \\ \& \
 \biggl(\int\limits_E^{(\pr{el})}f\,d(\mu+\nu)\in\rr \ \forall\,\mu\in(\pr{add})[\ml] \ \forall\,\nu\in(\pr{add})[\ml]\biggl).
\end{multline}
Несколько дополняя (\ref{8.4}), имеем (см. (\ref{8.1})), что $$\forall\,f\in B_0(E,\ml) \ \forall\,g\in  B_0(E,\ml) \ \forall\,\mu\in~(\pr{add})[\ml] \ \ \int\limits_E^{(\pr{el})}fg\,d\mu\in\rr.$$

С учетом \cite[предложение 3.2.2]{30} заметим, что элементарный интеграл как отображение
 \beq{8.6}
 (f,\mu)\mapsto\int\limits_E^{(\pr{el})}f\,d\mu: B_0(E,\ml)\times(\pr{add})[\ml]\rightarrow\rr
 \eeq
 является билинейным функционалом. Иными словами, при $\mu\in(\pr{add})[\ml]$
 \begin{multline}\label{8.7}
 \biggl(\,\int\limits_E^{(\pr{el})}\alpha f\,d\mu=\alpha\int\limits_E^{(\pr{el})}f\,d\mu \ \forall\,\alpha\in\rr \ \forall\,f\in B_0(E,\ml) \biggl)\ \& \\ \& \
 \biggl(\int\limits_E^{(\pr{el})}(f+g)\,d\mu= \int\limits_E^{(\pr{el})}f\,d\mu+\int\limits_E^{(\pr{el})}g\,d\mu  \ \forall\,f\in B_0(E,\ml) \  \forall\,g\in B_0(E,\ml)\biggl)
 \end{multline}
  и, кроме того, при $f\in B_0(E,\ml)$
\begin{multline}\label{8.8}
 \biggl(\int\limits_E^{(\pr{el})}f\,d(\alpha\mu)=\alpha\int\limits_E^{(\pr{el})}f\,d\mu \ \forall\,\alpha\in\rr \ \forall\,\mu\in(\pr{add})[\ml] \biggl)\ \& \\ \& \
 \biggl(\,\int\limits_E^{(\pr{el})}f\,d(\mu+\nu)= \int\limits_E^{(\pr{el})}f\,d\mu+\int\limits_E^{(\pr{el})}f\,d\nu  \ \forall\mu\in(\pr{add})[\ml] \
 \forall\,\nu\in(\pr{add})[\ml]\biggl).
 \end{multline}
В (\ref{8.7}), (\ref{8.8}) имеем реализацию упомянутого свойства билинейности функционала (\ref{8.6}).
 Разумеется, (\ref{8.7}) и (\ref{8.8}) распространяются на случаи, когда используются линейные комбинации ступенчатых функций и к.-а. мер; мы не будем на этом останавливаться, предоставляя эти случаи читателю в качестве  упражнения. Отметим только наиболее очевидные обстоятельства, полезные в дальнейшем. Так, при $\alpha\in\rr, \ \beta\in\rr,$ $f\in B_0(E,\ml),$ $g\in B_0(E,\ml)$ для (ступенчатой) в/з функции $\alpha f+\beta g\in B_0(E,\ml)$
\beq{8.9}
\int\limits_E^{(\pr{el})}(\alpha f+\beta g)\,d\mu=\alpha\int\limits_E^{(\pr{el})} f\,d\mu+\beta\int\limits_E^{(\pr{el})}g\,d\mu \ \forall\,\mu\in(\pr{add})[\ml].
\eeq
Далее  при $f\in B_0(E,\ml),$ для (ступенчатой) функции $$-f=(-1)f\in B_0(E,\ml)$$
имеем теперь очевидное свойство
\beq{8.10}
\int\limits_E^{(\pr{el})}(-f)\,d\mu=-\int\limits_E^{(\pr{el})}f\,d\mu \ \forall\,\mu\in(\pr{add})[\ml].
\eeq
Как следствие получаем, что при $f\in B_0(E,\ml)$ и $g\in B_0(E,\ml)$ для функции  $$f-g=f+(-g)\in B_0(E,\ml)$$
справедливо следующее свойство:
\beq{8.11}
\int\limits_E^{(\pr{el})}(f-g)\,d\mu=\int\limits_E^{(\pr{el})}f\,d\mu-\int\limits_E^{(\pr{el})}g\,d\mu \  \ \forall\,\mu\in(\pr{add})[\ml].
\eeq
 Предлагаем заинтересованному читателю выписать самостоятельно аналоги (\ref{8.9})--(\ref{8.11}) для случая, когда фиксируется функция из $B_0(E,\ml),$ а изменяется так или иначе к.-а. мера, по которой осуществляется интегрирование (используется линейная комбинация в $(\pr{add})[\ml],$ разность к.-а. мер и т.п.).

 Напомним ряд свойств (см. \cite[\S\,3.2]{30},(\ref{6.29'})) элементарного интеграла.
 \beq{8.12}
 \biggl|\int\limits_E^{(\pr{el})}f\,d\mu\biggl|\leqslant\|f\| V_{\mu} \ \forall\,f\in B_0(E,\ml) \ \forall\,\mu\in\mathbb{A}(\ml).
 \eeq
 С учетом (\ref{8.11}) и (\ref{8.12}) имеем, в частности, что
 \beq{8.13}
 \biggl|\int\limits_E^{(\pr{el})}f\,d\mu-\int\limits_E^{(\pr{el})}g\,d\mu\biggl|\leqslant\|f-g\| V_{\mu} \ \forall\,f\in B_0(E,\ml) \ \forall\,g\in B_0(E,\ml) \ \forall\,\mu\in\mathbb{A}(\ml).
 \eeq
Из (\ref{8.13}) и свойства полноты вещественной прямой с метрикой-модулем получаем (см. \cite[\S\,3.3]{30}), что
$\forall\,f\in B(E,\ml) \   \forall\,\mu\in\mathbb{A}(\ml) \ \exists\,!\,c\in\rr \ \forall\,(f_i)_{i\in\mathbb{N}}\in B_0(E,\ml)^\mathbb{N}$
\beq{8.14}
\left((f_i)_{i\in\mathbb{N}}\rightrightarrows f\right)\Rightarrow\left(\biggl(\int\limits_E^{(\pr{el})}f_i\,d\mu\biggl)_{i\in\mathbb{N}}\rightarrow c\right).
\eeq

Предлагаем читателю самостоятельно доказать свойство (\ref{8.14}), используя (\ref{8.13}), и сравнить данное рассуждение с \cite[предложение 3.3.1]{30}.   С учетом (\ref{8.14}) корректно (см. (\ref{7.25})) следующее
\begin{definition}
Если $f\in B(E,\ml)$ и $\mu\in\mathbb{A}(\ml),$ то полагаем, что
\beq{8.15}
\int\limits_E f\,d\mu\in\rr
\eeq
есть $\pr{def}$ такое единственное число, что $\forall\,(f_i)_{i\in\mathbb{N}}\in B_0(E,\ml)^\mathbb{N}$
$$\left((f_i)_{i\in\mathbb{N}}\rightrightarrows f\right)\Rightarrow\left(\biggl(\,\int\limits_E f_i\,d\mu\biggl)_{i\in\mathbb{N}}\rightarrow \int\limits_E f\,d\mu\right);$$
 число (\ref{8.15}) называем  ярусным интегралом функции $f$ по к.-а. мере $\mu.$
\end{definition}

В связи со свойствами ярусного интеграла (см. \cite[\S 3.4]{30}) следует учитывать (\ref{6.27'}). Имеем теперь отображение
\beq{8.16}
(f,\mu)\mapsto \int\limits_E f\,d\mu: B(E,\ml)\times\mathbb{A}(\ml)\rightarrow\rr.
\eeq
По аналогии с (\ref{8.4}), (\ref{8.5}) отметим, что \\
1) при $\mu\in\mathbb{A}(\ml)$ определены значения
\begin{multline*}
\biggl( \,\int\limits_E\alpha f\,d\mu\in\rr \ \forall\,\alpha\in\rr \  \forall\,f\in B(E,\ml)\biggl)\ \& \\ \& \
 \biggl(\,\int\limits_E (f+g)\,d\mu\in\rr \ \forall\,f\in B(E,\ml) \ \forall\,g\in B(E,\ml)\biggl);
\end{multline*}
 2) при $f\in B(E,\ml)$ определены значения
\begin{multline*}
\biggl(\,\int\limits_E f\,d(\alpha\mu)\in\rr \ \forall\,\alpha\in\rr \  \forall\,\mu\in\mathbb{A}(\ml)\biggl)\ \&  \\ \&
 \biggl(\,\int\limits_E f\,d(\mu+\nu)\in\rr \ \forall\,\mu\in\mathbb{A}(\ml) \ \forall\,\nu\in\mathbb{A}(\ml)\biggl).
\end{multline*}
Более того, предельным переходом из (\ref{8.7}), (\ref{8.8}) извлекается следующее свойство: отображение (\ref{8.16}) есть билинейный функционал, т.е. \\
1') если фиксирована  $\mu\in\mathbb{A}(\ml),$ то
\begin{multline}\label{8.17}
\biggl(\,\int\limits_E\alpha f\,d\mu=\alpha\int\limits_E  f\,d\mu \ \forall\,\alpha\in\rr \  \forall\,f\in B(E,\ml)\biggl)\ \& \\ \&
 \biggl(\,\int\limits_E (f+g)\,d\mu=\int\limits_E f\,d\mu+\int\limits_E g\,d\mu \ \forall\,f\in B(E,\ml) \ \forall\,g\in B(E,\ml)\biggl);
\end{multline}
 2') если фиксирована (ярусная) функция  $f\in B(E,\ml),$ то
\begin{multline}\label{8.18}
\biggl(\,\int\limits_E f\,d(\alpha\mu)=\alpha\int\limits_E  f\,d\mu \ \forall\,\alpha\in\rr \  \forall\,\mu\in\mathbb{A}(\ml)\biggl)\ \&  \\ \& \
 \biggl(\,\int\limits_E f\,d(\mu+\nu)=\int\limits_E  f\,d\mu+\int\limits_E  f\,d\nu \ \forall\,\mu\in\mathbb{A}(\ml) \ \forall\,\nu\in\mathbb{A}(\ml)\biggl).
\end{multline}
Отметим очевидные следствия (\ref{8.17}), (\ref{8.18}). Подобно (\ref{8.9}) имеем при $\alpha\!\in\!\rr,$ $ \beta\in\rr,$ $f\in B(E,\ml)$ и $g\in B(E,\ml)$  для функции $$\alpha f+\beta g\in B(E,\ml) $$ следующее важное свойство
$$\int\limits_E( \alpha f+\beta g)\,d\mu=\alpha\int\limits_E f\,d\mu+\beta\int\limits_E g\,d\mu \ \forall\,\mu\in\mathbb{A}(\ml).$$
Разумеется, из (\ref{6.19}) и (\ref{8.17}) получаем, что при $\mu\in\mathbb{A}(\ml)$ и  $f\in B(E,\ml)$
\beq{8.19}
\int\limits_E (-f)\,d\mu=-\int\limits_E f\,d\mu.
\eeq
Как следствие из (\ref{8.17}), (\ref{8.19}) вытекает, что
\beq{8.20}
\int\limits_E (f-g)\,d\mu=\int\limits_E f\,d\mu-\int\limits_E g\,d\mu \ \forall\,f\in B(E,\ml) \ \forall\,g\in B(E,\ml) \  \forall\,\mu\in\mathbb{A}(\ml).
\eeq
Аналогичным образом, из (\ref{6.19}) и (\ref{8.18})  при $f\in B(E,\ml)$   и $\mu\in\mathbb{A}(\ml)$ извлекается равенство
$$\int\limits_E f\,d(-\mu)=-\int\limits_E f\,d\mu.$$
Таким образом,  как частный случай (\ref{8.18}), имеем свойство  $$\int\limits_E f\,d(\mu-\nu)=\int\limits_E f\,d\mu-\int\limits_E f\,d\nu \ \forall\,f\in B(E,\ml) \ \forall\,\mu\in\mathbb{A}(\ml) \ \forall\,\nu\in\mathbb{A}(\ml).$$
Данные свойства используются ниже без дополнительных пояснений.

Отметим теперь, что из (\ref{7.15}) и (\ref{8.12}) предельным переходом получаем (см. \cite[предложение 3.4.2]{30}), что $ \forall\,f\in B(E,\ml)$ $\forall\,\mu\in\mathbb{A}(\ml)$
\beq{8.20'}
\left|\int\limits_E f\,d\mu\right|\leqslant\|f\| V_{\mu}.
\eeq
В свою очередь, из (\ref{8.20}) и (\ref{8.20'})  следует неравенство
\beq{8.21}
\biggl|\int\limits_E f\,d\mu-\int\limits_E g\,d\mu\biggl|\leqslant\|f-g\| V_{\mu} \ \forall\,f\in B(E,\ml) \ \forall\,g\in B(E,\ml) \  \forall\,\mu\in\mathbb{A}(\ml).
\eeq
Последнее свойство (см. (\ref{8.21})) характеризует при $\mu\in\mathbb{A}(\ml)$ функционал
$$f\mapsto\int\limits_Ef\,d\mu: B(E,\ml)\rightarrow\rr$$
как непрерывный в топологии $\sup$-нормы. В этой связи условимся о следующем обозначении:
\beq{8.22}
\int\limits_E\,d\mu\triangleq\biggl(\int\limits_Ef\,d\mu\biggl)_{f\in B(E,\ml)} \ \forall\,\mu\in\mathbb{A}(\ml).
\eeq
Из (\ref{7.23}), (\ref{8.17}) и  (\ref{8.22}) вытекает, что
$$\int\limits_E\,d\mu\in B'(E,\ml) \ \forall\,\mu\in\mathbb{A}(\ml).$$
С учетом (\ref{8.20'}) и последнего свойства имеем (см. (\ref{7.24})), что $$\int\limits_E\,d\mu\in B^*(E,\ml) \ \forall\,\mu\in\mathbb{A}(\ml).$$

Иными словами, установлено, что (в общем случае $\ml\in\pi[E]$) $\mathbb{A}(\ml)$ погружается в $B^*(E,\ml)$ посредством правила
\beq{8.23}
\mu\mapsto\int\limits_E\,d\mu: \mathbb{A}(\ml)\rightarrow B^*(E,\ml).
\eeq
В связи с (\ref{8.23}) напомним \cite[\S\,3.5]{30}, что
\beq{8.24}
V_{\mu}=\biggl\|\int\limits_E\,d\mu\biggl\|^* \  \forall\,\mu\in\mathbb{A}(\ml);
\eeq
отображение (\ref{8.23}) является линейным оператором
\begin{multline}\label{8.25}
\biggl(\,\int\limits_E \,d(\alpha\mu)=\alpha\int\limits_E \,d\mu \ \forall\,\alpha\in\rr \ \forall\,\mu\in\mathbb{A}(\ml)\biggl)\ \& \\ \&
\biggl(\,\int\limits_E \,d(\mu+\nu)=\int\limits_E \,d\mu+\int\limits_E \,d\nu\  \ \forall\,\mu\in\mathbb{A}(\ml) \ \forall\,\nu\in\mathbb{A}(\ml)\biggl).
\end{multline}
В  (\ref{8.24}) используем линейные операции в $\mathbb{A}(\ml)$ и в $ B^*(E,\ml),  B^*(E,\ml)\subset  B'(E,\ml);$ см. в этой связи (\ref{8.18}). Отметим, что в \cite[\S\,3.5]{30} приведены некоторые дополнения к (\ref{8.24}), (\ref{8.25}), которые сейчас обсуждаться не будут.

\section{Интегральное представление линейных непрерывных функционалов; оснащение $*$-слабой топологией}\setcounter{equation}{0}\setcounter{proposition}{0}\setcounter{zam}{0}\setcounter{corollary}{0}\setcounter{definition}{0}
   \ \ \ \ \ Всюду в дальнейшем будем полагать, что
   \beq{9.1}
   \ml\in\Pi[E],
   \eeq
   где $E$\,--- непустое множество. В виде $(E,\ml)$ имеем ИП с полуалгеброй множеств. Разумеется, все построения предыдущего раздела сохраняют свою силу и мы сделаем лишь некоторые добавления. Прежде всего напомним, что индикаторы множеств из $\ml$ являются ступенчатыми функциями (см. (\ref{7.19})). Поэтому при $L\in\ml$ и $\mu\in(\pr{add})[\mathcal{L}]$ определен элементарный интеграл $$\int\limits_E^{(\pr{el})}\chi_L\,d\mu\in\rr.$$ Из определения 9.1  следует, что (см. \cite[предложение 3.2.3]{30})
   \beq{9.2}
   \int\limits_E^{(\pr{el})}\chi_L\,d\mu=\mu(L)
   \eeq
   (при доказательстве (\ref{9.2}) полезно учитывать определение полуалгебры множеств). Напомним (\ref{7.21}), (\ref{7.22}) и следующее важное свойство \cite[предложение 3.6.1]{30}: функционал
   \beq{9.3}
   \mu\mapsto V_{\mu}:\mathbb{A}(\ml)\rightarrow[0,\infty[
   \eeq
является нормой на  $\mathbb{A}(\ml)$ (\ref{6.27'}). В дальнейшем норму (\ref{9.3}) называем \emph{сильной} или \emph{нормой-вариацией}:
\begin{multline}\label{9.4}
(V_{\alpha\mu}=|\alpha| V_{\mu} \ \forall\,\alpha\in\rr \ \forall\,\mu\in\mathbb{A}(\ml)) \ \& \\ \& (V_{\mu+\nu}\leqslant V_{\mu}+V_{\nu} \ \forall\,\mu\in\mathbb{A}(\ml) \ \forall\,\nu\in\mathbb{A}(\ml)) \ \& \\ \& \ (\forall\,\mu\in\mathbb{A}(\ml) \ (V_{\mu}=0)\Leftrightarrow(\mu=\mathcal{O}_{\ml})).
\end{multline}
Теперь при условии (\ref{9.1}) дополним свойства оператора (\ref{8.23}), а именно, в рассматриваемом случае справедлива следующая (см. \cite[теорема 3.6.1]{30})
\begin{theorem}
Оператор (\ref{8.23}) есть изометрический изоморфизм $\mathbb{A}(\ml)$ в сильной норме (\ref{9.3}) на $B^*(E,\ml)$ в оснащении традиционной нормой (\ref{7.26}): \\
1) \  $\left(\,\int\limits_E \,d\mu\right)_{\mu\in \mathbb{A}(\ml)}\in(\pr{bi})[\mathbb{A}(\mathcal{L});B^*(E,\ml)];$\\
2) \   (\ref{8.23}) есть линейный оператор (см. (\ref{8.25}));\\
3) оператор (\ref{8.23}) обладает изометричностью в смысле (\ref{8.24}).
\end{theorem}
Д о к а з а т е л ь с т в о \ \  теоремы см.  в  \cite[c.\,152-156]{30}.

Итак,  фактически $\mathbb{A}(\ml), \ B^*(E,\ml)$ отождествимы. Данное обстоятельство позволяет ввести в $\mathbb{A}(\ml)$ стандартную $*$-слабую топологию, что отвечает по существу использованию $B(E,\ml)$-топологии множества $B^*(E,\ml).$ Тем самым будут дополнены представления раздела 7.
Прежде всего введем  при $\mu\in\mathbb{A}(\ml), \ K\in\pr{Fin}(B(E,\mathcal{L}))$ и $\varepsilon\bn$ (непустое) множество
\beq{9.5}
N^*_{\ml}(\mu,K,\varepsilon)\triangleq\left\{\nu\in \mathbb{A}(\ml):\biggl|\int\limits_E f\,d\mu-\int\limits_E f\,d\nu\biggl|<\varepsilon \ \forall\,f\in K\right\}.
\eeq
Легко видеть, что в терминах (\ref{9.5}) определяется нужная нам $*$-слабая топология $\mathbb{A}(\ml):$
\begin{multline}\label{9.6}
\tau_*(\ml)\triangleq\{G\!\in\!\pp(\mathbb{A}(\ml))\mid\forall\,\mu\in\! G \ \exists\,K\!\in\!\pr{Fin}(B(E,\mathcal{L}))  \exists\,\varepsilon\!\in]0,\infty[: \\ N^*_{\ml}(\mu,K,\varepsilon)\!\subset\! G\}\!\in\!(\pr{top})[\mathbb{A}(\ml)];
\end{multline}
разумеется, в виде
\beq{9.7}
(\mathbb{A}(\ml),\tau_*(\ml))
\eeq
  имеем ТП с <<единицей>> $\mathbb{A}(\ml).$ Множества (\ref{9.5}) образуют (в совокупности) базу ТП (\ref{9.7}); см. \cite[(4.6.3)]{41}.
  В частности,
\beq{9.8}
  N^*_{\ml}(\mu,K,\varepsilon)\in\tau_*(\ml) \ \forall\,\mu\in \mathbb{A}(\ml) \ \forall\,K\in\pr{Fin}(B(E,\mathcal{L})) \ \forall\,\varepsilon\bn.
\eeq
  Как следствие, из (\ref{3.28}) и (\ref{9.8}) получаем, что
\beq{9.8'}
N^*_{\ml}(\mu,K,\varepsilon)\in N^0_{\tau_*(\ml)}(\mu)\ \forall\,\mu\in \mathbb{A}(\ml) \ \forall\,K\in\pr{Fin}(B(E,\mathcal{L})) \ \forall\,\varepsilon\bn.
\eeq

  Более того, при фиксации $\mu\in\mathbb{A}(\ml)$ имеем, что множества $ N^*_{\ml}(\mu,K,\varepsilon),$ $ K\in\pr{Fin}(B(E,\mathcal{L})),$ $ \varepsilon\bn$ образуют локальную базу ТП (\ref{9.7}) в точке $\mu;$ это означает, в частности, что
\beq{9.9}
\forall\,H\in N_{\tau_*(\ml)}(\mu) \ \exists\,K\in\pr{Fin}(B(E,\mathcal{L})) \ \exists\,\varepsilon\bn: N^*_{\ml}(\mu,K,\varepsilon)\subset H;
\eeq
   см. \cite[(4.6.2), (4.6.7)]{41}. В связи с (\ref{9.9}) отметим простое следствие (\ref{3.34}):
\begin{multline}\label{9.10}
\pr{cl}(S,\tau_*(\ml))=\{\mu\in\mathbb{A}(\ml)\mid S\cap N^*_{\ml}(\mu,K,\varepsilon)\neq\zer \\ \forall\,K\in\pr{Fin}(B(E,\mathcal{L})) \ \forall\,\varepsilon\bn\,\}\ \forall\,S\in\pp(\mathbb{A}(\ml))
\end{multline}
(читателю предлагается самостоятельно проверить (\ref{9.10})). Напомним эдесь же, что согласно (\ref{3.34}) и (\ref{9.6})
\begin{multline}\label{9.11}
\mathcal{F}[\tau_*(\ml)]=\{\mathbb{A}(\ml)\setminus G: G\in\tau_*(\ml)\}=\\=\{F\in\pp(\mathbb{A}(\ml))\mid\mathbb{A}(\ml)\setminus F\in\tau_*(\ml)\}\in \\ \in \pp'(\pp(\mathbb{A}(\ml))).
\end{multline}
В (\ref{9.11}) имеем семейство всех $*$-слабо замкнутых п/м $\mathbb{A}(\ml).$ При этом (см. (\ref{3.37}))
\beq{9.12}
\pr{cl}(S,\tau_*(\ml))\in\mathcal{F}[\tau_*(\ml)] \ \forall\,S\in\pp(\mathbb{A}(\ml)).
\eeq
Напомним, что $(\tau_*(\ml)-\pr{comp})[\mathbb{A}(\ml)]$ есть (см. (\ref{3.41})) семейство всех п/м $\mathbb{A}(\ml),$ компактных в ТП (\ref{9.7}). С другой стороны, данное ТП по самому способу построения обладает \cite[гл.\,V]{43} следующим важным свойством: п/м $\mathbb{A}(\ml)$ компактно в ТП (\ref{9.7}) тогда и только тогда, когда оно сильно ограничено и $*$-слабо замкнуто.

Для более точного представления упомянутого свойства введем
$$U_b(\ml)\triangleq\{\mu\in\mathbb{A}(\ml)\mid V_{\mu}\leqslant b\} \ \forall\,b\in[0,\infty[$$
(шар в сильной норме с центром в <<нуле>>). Кроме того, введем в рассмотрение
$$\mathbb{B}_*(\ml)\triangleq\{H\in\pp(\mathbb{A}(\ml))\mid\exists\,c\in[0,\infty[: H\subset U_c(\ml)\},$$
получая семейство всех сильно ограниченных п/м $\mathbb{A}(\ml)$. Тогда
\beq{9.13}
(\tau_*(\ml)-\pr{comp})[\mathbb{A}(\ml)]=\mathbb{B}_*(\ml)\cap\mathcal{F}[\tau_*(\ml)].
\eeq
В связи с (\ref{9.13}) отметим в дальнейшем ряд конкретных примеров п/м $\mathbb{A}(\ml),$ компактных в ТП (\ref{9.7}). Так, в частности, имеем (см. \cite[(4.6.11)]{41}), что при $c\in[0,\infty[$
\beq{9.14}
U_c(\ml)\in(\tau_*(\ml)-\pr{comp})[\mathbb{A}(\ml)].
\eeq
Читателю предлагается самостоятельно проверить свойство $U_c(\ml)\in\mathcal{F}[\tau_*(\ml)],$
используя определение полной вариации к.-а. меры и определение $*$-слабой топологии. Из (\ref{9.14}) вытекает, что ТП (\ref{9.7}) $\sigma$-компактно, т.е. является счетным объединением компактных~п/м:
\beq{9.15}
\mathbb{A}(\ml)=\bigcup\limits_{k\in\mathbb{N}}U_k(\ml).
\eeq

На самом деле  (\ref{9.7}) есть локально-выпуклый $\sigma$-компакт: \\
1) \  (\ref{9.7}) есть хаусдорфово ТП (см. в этой связи (\ref{9.8'}), (\ref{9.9})) со свойством  $\sigma$-компактности, следующим из (\ref{9.15});\\
2) \  ТП (\ref{9.7}) является топологическим векторным пространством, в котором к.-а. мера $\mathcal{O}_\ml\in\mathbb{A}(\ml)$ обладает фундаментальной системой (локальной базой) выпуклых окрестностей $N_\ml^*(\mathcal{O}_\ml,K,\varepsilon),$ $K\in\pr{Fin}(B(E,\ml)),$ $\varepsilon\bn.$
Из 2) следует, что (\ref{9.7}) есть локально выпуклое ТП. Итак, (\ref{9.7}) есть топологическое векторное пространство, причем <<хорошее>>, поскольку является локально выпуклым $\sigma$-компактом (см. (\ref{9.15})).

Отметим некоторые соотношения, связывающие (\ref{9.7}) с ТП раздела 7, используя положения  \cite[раздел 4.6]{41}.  Прежде всего заметим, что (см. \cite[(4.6.19), (4.6.21)]{41})
\beq{9.15'}
(\tau_\otimes(\ml)\subset\tau_0(\ml))\ \& \ (\tau_\otimes(\ml)\subset\tau_*(\ml));
\eeq
это означает справедливость утверждения:  $\tau_\otimes(\ml)$  слабее каждой из топологий $\tau_0(\ml), \tau_*(\ml).$ Вместе с тем \cite[предложение 4.6.1]{41}
\beq{9.16}
\tau_*(\ml)|_B=\tau_\otimes(\ml)|_B=\otimes^\ml(\tr)|_B \ \forall\,B\in\mathbb{B}_*(\ml).
\eeq
Из (\ref{9.13}) и (\ref{9.16})  вытекает, в частности, свойство
\beq{9.17}
\tau_*(\ml)|_K=\tau_\otimes(\ml)|_K=\otimes^\ml(\tr)|_K \ \forall\,K\in(\tau_*(\ml)-\pr{comp})[\mathbb{A}(\ml)].
\eeq
Отметим вместе с тем, что сами топологии $\tau_*(\ml)$ и $\tau_\otimes(\ml)$ могут различаться уже в достаточно простых случаях (см. пример в \cite[раздел 4.4]{36}).  Полезно отметить, однако (см. \cite[(4.2.12)]{36}), соотношения
\beq{9.18}
\tau_*^+(\ml)\triangleq\tau_*(\ml)|_{(\pr{add})_+[\mathcal{L}]}=\tau_\otimes^+(\ml)\subset\tau_0^+(\ml),
\eeq
чем существенно дополняется (\ref{9.17}). Заметим, что из (\ref{6.43}) и (\ref{9.16}) вытекает,  в частности, свойство
\beq{9.19}
\tau_*(\ml)|_B\subset\tau_0(\ml)|_B \ \forall\,B\in\mathbb{B}_*(\ml).
\eeq
Из (\ref{9.13}), (\ref{9.17}), (\ref{9.19}) получем, в частности,
\beq{9.20}
\tau_*(\ml)|_K\subset\tau_0(\ml)|_K \ \forall\,K\in(\tau_*(\ml)-\pr{comp})[\mathbb{A}(\ml)].
\eeq
Кроме того, напомним \cite[(4.6.13)]{41}, что
\beq{9.21}
(\pr{add})_+[\mathcal{L}]\in\mathcal{F}[\tau_*(\ml)].
\eeq
Итак (см. (\ref{9.21})), конус $(\pr{add})_+[\mathcal{L}]$ является $*$-слабо замкнутым п/м $\mathbb{A}(\ml).$ Здесь же отметим весьма очевидные свойства \cite[(4.6.14)]{41} и \cite[(4.6.15)]{41}
\beq{9.22}
\mathbb{P}(\ml)\in(\tau_*(\ml)-\pr{comp})[\mathbb{A}(\ml)], \ \mathbb{T}(\ml)\in(\tau_*(\ml)-\pr{comp})[\mathbb{A}(\ml)];
\eeq
в этой связи напомним цепочку вложений (\ref{6.16}).

Из (\ref{4.4}), (\ref{9.18}) и (\ref{9.21}) вытекает следующий  очевидный факт: при $H\in\pp((\pr{add})_+[\mathcal{L}])$ $$\pr{cl}(H,\tau_*(\ml))\subset(\pr{add})_+[\mathcal{L}],$$
а потому имеем равенство
\beq{9.23}
\pr{cl}(H,\tau_*^+(\ml))=\pr{cl}(H,\tau_*(\ml))\cap(\pr{add})_+[\mathcal{L}]=\pr{cl}(H,\tau_*(\ml)).\eeq
Аналогичным образом из (\ref{4.4}) и (\ref{6.44}) следует равенство
\beq{9.24}
\pr{cl}(H,\tau_\otimes^+(\ml))=\pr{cl}(H,\tau_\otimes(\ml))\cap(\pr{add})_+[\mathcal{L}]
\eeq
и, вместе с тем, $H\subset(\pr{add})_+[\mathcal{L}].$ Наряду с (\ref{9.21}) имеем следующее свойство замкнутости:
\beq{9.25}
(\pr{add})_+[\mathcal{L}]\in\mathcal{F}[\tau_\otimes(\ml)].
\eeq
\begin{zam}
Проверим справедливость (\ref{9.25}), для чего рассмотрим множество
\beq{9.26}
\mathbf{A}\triangleq\mathbb{A}(\ml)\setminus(\pr{add})_+[\mathcal{L}].
\eeq

 Покажем, что $\mathbf{A}\in\tau_\otimes(\ml),$ для чего выберем и зафиксируем
 \beq{9.27}
 \mu_*\in\mathbf{A}.
 \eeq
 Тогда, в частности, $\mu_*\in\rr^\ml,$  т.е. $\mu_*:\ml\rightarrow\rr$  и при этом $\mu_*\notin(\pr{add})_+[\mathcal{L}].$
Поскольку (см. (\ref{9.26})) $\mu_*\in\mathbb{A}(\ml),$ то непременно $\mu_*\in(\pr{add})[\mathcal{L}]$ (см. (\ref{6.15})), а тогда
\beq{9.28}
\mu_*\in(\pr{add})[\mathcal{L}]\setminus(\pr{add})_+[\mathcal{L}].
\eeq
В силу (\ref{6.12}) и (\ref{9.28}) получаем для некоторого множества
\beq{9.29}
\Lambda\in\ml
\eeq
неравенство $\mu_*(\Lambda)<0.$  При этом, конечно, имеем (см. (\ref{9.29})) включение   $\{\Lambda\}\in\pr{Fin}(\ml).$ Тогда, полагая $\varepsilon_*\triangleq-\mu(\Lambda),$ имеем $\varepsilon_*\bn$ и согласно (\ref{6.30})
\beq{9.29}
\mathbb{N}_\ml(\mu_*,\{\Lambda\},\varepsilon_*)=\left\{\nu\in\rr^\ml\mid|\mu_*(\Lambda)-\nu(\Lambda)|<\varepsilon_*\right\}\in\otimes^\ml(\tr).
\eeq
Если $\nu\in\mathbb{N}_\ml(\mu_*,\{\Lambda\},\varepsilon_*),$ то $\nu(\Lambda)<\mu_*(\Lambda)+\varepsilon_*=0.$ Поэтому
\beq{9.30}
\mathbb{N}_\ml(\mu_*,\{\Lambda\},\varepsilon_*)\cap\mathbb{A}(\ml)\subset\mathbb{A}(\ml)\setminus(\pr{add})_+[\mathcal{L}]=\mathbf{A}.
\eeq
Согласно (\ref{6.39}) и (\ref{9.29}) имеем свойство
\beq{9.31}
\mathbb{N}_\ml(\mu_*,\{\Lambda\},\varepsilon_*)\cap\mathbb{A}(\ml)\in\tau_\otimes(\ml);
\eeq
см. в этой связи также (\ref{4.0}). При этом (см. (\ref{9.29}))
$\mu_*\in\mathbb{N}_\ml(\mu_*,\{\Lambda\},\varepsilon_*)\cap\mathbb{A}(\ml),$
   поэтому в силу (\ref{4.5}), (\ref{9.31}) имеем, что $$\mathbb{N}_\ml(\mu_*,\{\Lambda\},\varepsilon_*)\cap\mathbb{A}(\ml)\in N_{\tau_\otimes(\ml)}^0(\mu_*).$$
  Это означает (см. (\ref{9.30})), что
 \beq{9.32}
 \exists\,G\in N_{\tau_\otimes(\ml)}^0(\mu_*): G\subset\mathbf{A}.
 \eeq

 Поскольку $\mu_*$ (\ref{9.27}) выбиралось произвольно, то $\mathbf{A}\in\tau_\otimes(\ml)$ (из (\ref{9.32}) следует, что $\mathbf{A}$  есть объединение открытых множеств, т.к. в силу произвольности выбора $\mu_*$ было установлено вложение $$\forall\,\mu\in\mathbf{A} \ \exists\,G\in N_{\tau_\otimes(\ml)}^0(\mu): G\subset\mathbf{A};$$ дальнейшее рассуждение очевидно). Но в этом случае $(\pr{add})_+[\mathcal{L}]=\mathbb{A}(\ml)\setminus\mathbf{A}$ есть замкнутое в ТП $(\mathbb{A}(\ml),\tau_\otimes(\ml))$ множество, чем и завершается проверка (\ref{9.25}). $\hfill\square$
\end{zam}
Как следствие мы получили  (см. (\ref{3.35})) свойство: при $H\in\pp((\pr{add})_+[\mathcal{L}])$  $$(\pr{add})_+[\mathcal{L}]\in[\mathcal{F}[\tau_\otimes(\ml)]](H),$$ поэтому согласно (\ref{3.37}) $$\pr{cl}(H,\tau_\otimes(\ml)])\subset(\pr{add})_+[\mathcal{L}].$$ В итоге (см. (\ref{9.24})) получено следующее свойство
\beq{9.33}
\pr{cl}(H,\tau_\otimes^+(\ml))=\pr{cl}(H,\tau_\otimes(\ml)) \ \forall\,H\in\pp((\pr{add})_+[\mathcal{L}]).
\eeq
Заметим здесь, что согласно (\ref{3.34}) и (\ref{9.15'})
\beq{9.34}
\mathcal{F}[\tau_\otimes(\ml)]\subset\mathcal{F}[\tau_0(\ml)],
\eeq
а потому (см. (\ref{9.25}), (\ref{9.34})) $(\pr{add})_+[\mathcal{L}]\in\mathcal{F}[\tau_0(\ml)].$
Если $\mathbb{H}\in\pp((\pr{add})_+[\mathcal{L}]),$ то согласно (\ref{3.36})
$(\pr{add})_+[\mathcal{L}]\in[\mathcal{F}[\tau_0(\ml)]](\mathbb{H}),$
поэтому (см. (\ref{3.38})) имеем вложение $\pr{cl}(\mathbb{H},\tau_0(\ml))\subset(\pr{add})_+[\mathcal{L}].$ Используя (\ref{4.4}) и (\ref{6.45}), получаем цепочку равенств  $$\pr{cl}(\mathbb{H},\tau_0^+(\ml))=\pr{cl}(\mathbb{H},\tau_0(\ml))\cap(\pr{add})_+[\mathcal{L}]=\pr{cl}(\mathbb{H},\tau_0(\ml)).$$

Итак, установлено следующее свойство:
\beq{9.35}
\pr{cl}(H,\tau_0^+(\ml))=\pr{cl}(H,\tau_0(\ml)) \ \forall\,H\in\pp((\pr{add})_+[\mathcal{L}]).
\eeq

\section{Слабо абсолютно непрерывные конечно-аддитивные меры и их приближение неопределенными интегралами, 1}\setcounter{equation}{0}\setcounter{proposition}{0}\setcounter{zam}{0}\setcounter{corollary}{0}\setcounter{definition}{0}
   \ \ \ \ \ В пределах настоящего раздела фиксируем неотрицательную к.-а. меру
\beq{10.1}
\eta\in(\pr{add})_+[\mathcal{L}],
\eeq
 получая в виде триплета
 \beq{10.2}
 (E,\ml,\eta)
 \eeq
аналог стандартного пространства с мерой \cite[с.\,64]{44}.  Будем называть триплет   (\ref{10.2}) \emph{к.-а. пространством с мерой}. Для случая, определяемого в (\ref{10.1}), (\ref{10.2}), конкретизируем~(\ref{6.29}):
\begin{multline}\label{10.3}
(\pr{add})^+[\mathcal{L};\eta]=\{\nu\in(\pr{add})_+[\mathcal{L}]\mid\forall\,L\in\ml \ (\eta(L)=0)\Rightarrow(\nu(L)=0)\}\in \\ \in(\pr{cone})[(\pr{add})_+[\mathcal{L}]].
\end{multline}
Элементы конуса (\ref{10.3}) называем неотрицательными к.-а. мерами на $\ml,$  слабо абсолютно непрерывными относительно $\eta$ (\ref{10.1}).  Из (\ref{10.3}), в частности, следует включение
$$(\pr{add})^+[\ml;\eta]\in\pp'( \mathbb{A}(\ml)).$$

Напомним, что $\mathcal{O}_{\ml}\in(\pr{add})^+[\mathcal{L};\eta].$ Более содержательным примером к.-а. меры из конуса (\ref{10.3}) является неопределенный интеграл неотрицательной ярусной функции по к.-а. мере $\eta.$ Определение и свойства неопределенного интеграла подробно рассматривались для случая (\ref{9.1}), (\ref{10.1}) в \cite[\S\,3.7]{30} (см. также (\ref{7.21}), (\ref{7.22})). Как обычно, полагаем, что
$$\int\limits_L f\,d\mu\triangleq\int\limits_E f\chi_L\,d\mu \ \ \forall\,f\in B(E,\ml) \ \forall\,\mu\in\mathbb{A}(\ml) \ \forall\,L\in\ml.$$
 Отметим в этой связи, что при  $L=E$ \ \ $f\chi_L=f \ \  \forall\,f\in B(E,\ml),$  поэтому в случае $L=E$   вновь определенный интеграл на множестве $L$ для произвольной ярусной функции совпадает с интегралом в смысле определения 9.2. В частности, имеем с учетом (\ref{10.1}) представление
 \beq{10.4}
 \int\limits_L f\,d\eta=\int\limits_E f\chi_L\,d\eta\in\rr \ \ \forall\,f\in B(E,\ml) \  \forall\,L\in\ml.
 \eeq
Далее следуем \cite[\S\,3.7]{30}, учитывая (\ref{7.17}). Итак, при $f\in B(E,\ml)$ имеем \cite[предложение 3.7.1]{30}
  $$\biggl(\,\int\limits_E uf\,d\eta\biggl)_{u\in B(E,\ml)}\in B^*(E,\ml),$$
  а потому согласно теореме 10.1
  $$\exists\,!\,\mu\in\mathbb{A}(\ml): \ \biggl(\int\limits_E uf\,d\eta\biggl)_{u\in B(E,\ml)}=\int\limits_E d\mu.$$
В этой связи полагаем, что $\forall\,f\in B(E,\ml)$ $\pr{def}$ $f*\eta\in\mathbb{A}(\ml):$
\beq{10.5}
\int\limits_E uf\,d\eta=\int\limits_E u\,d(f*\eta) \ \forall\,u\in B(E,\ml).
\eeq

Из  (\ref{10.4}) и  (\ref{10.5}) следует  (см. (\ref{9.2})) при $f\in B(E,\ml)$ и $L\in\ml$ цепочка равенств
\beq{10.5'}
\int\limits_L f\,d\eta=\int\limits_E f\chi_L\,d\eta=\int\limits_E \chi_Lf\,d\eta=\int\limits_E \chi_L\,d(f*\eta)=(f*\eta)(L).
\eeq

Иными словами, имеем следующее представление:
\beq{10.6}
f*\eta=\biggl(\int\limits_L f\,d\eta\biggl)_{L\in\ml}\in\mathbb{A}(\ml) \ \  \forall\,f\in B(E,\ml).
\eeq
Как отмечено в \cite[предложение 3.7.4]{30},
\beq{10.7}
\biggl|\int\limits_L f\,d\eta\biggl|\leqslant\|f\|\eta(L) \ \  \forall\,f\in B(E,\ml) \ \ \forall\,L\in\ml.
\eeq
Тогда из (\ref{10.6}) и (\ref{10.7}) вытекает, что $\forall\,L\in\ml$
\beq{10.8}
(\eta(L)=0)\Rightarrow((f*\eta)(L)=0 \ \ \forall\,f\in B(E,\ml)).
\eeq
Далее из (\ref{7.24*}) и (\ref{7.24**}) следуют свойства (см. \cite[(3.4.33)]{30})
\beq{10.9}
f*\eta\in(\pr{add})_+[\mathcal{L}] \ \ \forall\,f\in B^+(E,\ml).
\eeq
Но тогда, согласно (\ref{10.3}), (\ref{10.8}) и (\ref{10.9}) имеем
\beq{10.10}
f*\eta\in(\pr{add})^+[\mathcal{L};\eta] \ \ \forall\,f\in B^+(E,\ml).
\eeq
В частности, из (\ref{10.10}) легко следует, что
\beq{10.11}
f*\eta\in (\pr{add})^+[\mathcal{L};\eta] \ \ \forall\,f\in B_0^+(E,\ml).
\eeq
Итак, неопределенные $\eta$-интегралы неотрицательных ступенчатых и ярусных функций доставляют примеры слабо абсолютно непрерывных относительно $\eta$ к.-а. мер на $\ml.$ Следуя \cite[(4.3.18), (4.3.19)]{36},  введем при $b\in[0,\infty[$ в рассмотрение множества
\beq{10.12}
M_*^+[\,b\mid E;\ml;\eta]\triangleq\biggl\{f\in B_0^+(E,\ml)\mid \int\limits_E f\,d\eta=b\biggl\},
\eeq
\beq{10.13}
\mathcal{M}_*^+[\,b\mid E;\ml;\eta]\triangleq\biggl\{f\in B^+(E,\ml)\mid \int\limits_E f\,d\eta=b\biggl\},
\eeq
\beq{10.14}
\widetilde{M}_*^+[\,b\mid E;\ml;\eta]\triangleq\{f*\eta:f\in M_*^+[\,b\mid E;\ml;\eta]\},
\eeq
\beq{10.15}
\widetilde{\mathcal{M}}_*^+[\,b\mid E;\ml;\eta]\triangleq\{f*\eta:f\in \mathcal{M}_*^+[\,b\mid E;\ml;\eta]\}.
\eeq
В (\ref{10.12}) и (\ref{10.13}) имеем два множества в/з функций на $E,$ а в (\ref{10.11}) и (\ref{10.15})\,--- два множества неотрицательных к.-а. мер. В частности, имеем (при $b\in[0,\infty[$) свойства
$$\widetilde{M}_*^+[\,b\mid E;\ml;\eta]\in\pp(\mathbb{A}(\ml)), \ \ \widetilde{\mathcal{M}}_*^+[\,b\mid E;\ml;\eta]\in\pp(\mathbb{A}(\ml)).$$
Поскольку в силу (\ref{10.12}) и (\ref{10.13}) $M_*^+[\,b\mid E;\ml;\eta]\subset\mathcal{M}_*^+[\,b\mid E;\ml;\eta],$ то (см. (\ref{10.14}), (\ref{10.15})) справедливо вложение
\beq{10.16}
\widetilde{M}_*^+[\,b\mid E;\ml;\eta]\subset\widetilde{\mathcal{M}}_*^+[\,b\mid E;\ml;\eta].
\eeq

Заметим, что $b\in[0,\infty[$ и при $\tau\in(\pr{top})[\mathbb{A}(\ml)]$ корректно определяются следующие множества-замыкания
$$\pr{cl}\left(\widetilde{M}_*^+[\,b\mid E;\ml;\eta],\tau\right)\in\pp(\mathbb{A}(\ml)), \ \ \pr{cl}\left(\widetilde{\mathcal{M}}_*^+[\,b\mid E;\ml;\eta],\tau\right)\in\pp(\mathbb{A}(\ml)).$$
Вместе с тем, из (\ref{10.9}), (\ref{10.14}), (\ref{10.15}) вытекает, что при $b\in[0,\infty[$
$$\widetilde{M}_*^+[\,b\mid E;\ml;\eta]\in\pp((\pr{add})_+[\mathcal{L}]), \ \ \widetilde{\mathcal{M}}_*^+[\,b\mid E;\ml;\eta]\in\pp((\pr{add})_+[\mathcal{L}]).$$
С учетом (\ref{9.23}), (\ref{9.33}) получаем  $\forall\,b\in[0,\infty[$
\begin{multline*}
\biggl(\pr{cl}\left(\widetilde{M}_*^+[\,b\mid E;\ml;\eta],\tau_*^+(\ml)\right)=\pr{cl}\left(\widetilde{M}_*^+[\,b\mid E;\ml;\eta],\tau_*(\ml)\right)\in \\ \in\pp((\pr{add})_+[\mathcal{L}])\biggl)\ \& \end{multline*}
\begin{multline*} \&\, \biggl(\pr{cl}\left(\widetilde{\mathcal{M}}_*^+[\,b\mid E;\ml;\eta],\tau_*^+(\ml)\right)=\pr{cl}\left(\widetilde{\mathcal{M}}_*^+[\,b\mid E;\ml;\eta],\tau_*(\ml)\right)\in \\ \in\pp((\pr{add})_+[\mathcal{L}])\biggl)\ \& \end{multline*}
\begin{multline*}
\& \,\biggl(\pr{cl}\left(\widetilde{M}_*^+[\,b\mid E;\ml;\eta],\tau_\otimes^+(\ml)\right)=\pr{cl}\left(\widetilde{M}_*^+[\,b\mid E;\ml;\eta],\tau_\otimes(\ml)\right)\in \\ \in\pp((\pr{add})_+[\mathcal{L}])\biggl)\ \& \end{multline*}
\begin{multline}\label{10.15} \&\,
\biggl(\pr{cl}\left(\widetilde{\mathcal{M}}_*^+[\,b\mid E;\ml;\eta],\tau_\otimes^+(\ml)\right)=\pr{cl}\left(\widetilde{\mathcal{M}}_*^+[\,b\mid E;\ml;\eta],\tau_\otimes(\ml)\right)\in\\ \in\pp((\pr{add})_+[\mathcal{L}])\biggl).
\end{multline}
Полезно учесть (\ref{9.18}). Тогда из (\ref{10.15}) следует цепочка равенств
\begin{multline}\label{10.16'}
\pr{cl}\left(\widetilde{M}_*^+[\,b\mid E;\ml;\eta],\tau_*(\ml)\right)=\pr{cl}\left(\widetilde{M}_*^+[\,b\mid E;\ml;\eta],\tau_\otimes(\ml)\right), \\  \pr{cl}\left(\widetilde{\mathcal{M}}_*^+[\,b\mid E;\ml;\eta],\tau_*(\ml)\right)=\pr{cl}\left(\widetilde{\mathcal{M}}_*^+[\,b\mid E;\ml;\eta],\tau_\otimes(\ml)\right).
\end{multline}
Особо отметим последнее положение в (\ref{9.18}). Дело в том, что согласно (\ref{3.34}) и (\ref{9.18})
$$\mathcal{F}[\tau_\otimes^+(\ml)]\subset\mathcal{F}[\tau_0^+(\ml)],$$ а тогда (см. (\ref{3.36})) при $H\in\pp((\pr{add})_+[\mathcal{L}])$
$$[\mathcal{F}[\tau_\otimes^+(\ml)]](H)\subset[\mathcal{F}[\tau_0^+(\ml)]](H),$$
откуда согласно (\ref{3.37}) следует, что
\beq{10.17}
\pr{cl}(H,\tau_0^+(\ml))\subset\pr{cl}(H,\tau_\otimes^+(\ml)).
\eeq
Учитывая (\ref{9.23}), (\ref{9.35}) и (\ref{10.17}), получаем также
\beq{10.18}
\pr{cl}(H,\tau_0(\ml))\subset\pr{cl}(H,\tau_\otimes(\ml)) \ \forall\,H\in\pp((\pr{add})_+[\mathcal{L}]).
\eeq
В частности, из (\ref{10.18}) вытекает цепочка вложений
\begin{multline*}
\pr{cl}\left(\widetilde{M}_*^+[\,b\mid E;\ml;\eta],\tau_0(\ml)\right)\subset\pr{cl}\left(\widetilde{M}_*^+[\,b\mid E;\ml;\eta],\tau_\otimes(\ml)\right), \\
\pr{cl}\left(\widetilde{\mathcal{M}}_*^+[\,b\mid E;\ml;\eta],\tau_0(\ml)\right)\subset\pr{cl}\left(\widetilde{\mathcal{M}}_*^+[\,b\mid E;\ml;\eta],\tau_\otimes(\ml)\right).
\end{multline*}
С учетом (\ref{10.16'}) из последних соотношений получаем полезные свойства
\beq{10.19}
\pr{cl}\left(\widetilde{M}_*^+[\,b\mid E;\ml;\eta],\tau_0(\ml)\right)\subset\pr{cl}\left(\widetilde{M}_*^+[\,b\mid E;\ml;\eta],\tau_*(\ml)\right),
\eeq
\beq{10.20}
\pr{cl}\left(\widetilde{\mathcal{M}}_*^+[\,b\mid E;\ml;\eta],\tau_0(\ml)\right)\subset\pr{cl}\left(\widetilde{\mathcal{M}}_*^+[\,b\mid E;\ml;\eta],\tau_*(\ml)\right).
\eeq
Кроме того, из  (\ref{3.32}), (\ref{3.33}) и (\ref{10.16}) очевидно следует вложение
\beq{10.21}
\pr{cl}\left(\widetilde{M}_*^+[\,b\mid E;\ml;\eta],\tau_*(\ml)\right)\subset\pr{cl}\left(\widetilde{\mathcal{M}}_*^+[\,b\mid E;\ml;\eta],\tau_*(\ml)\right),
\eeq
\beq{10.22}
\pr{cl}\left(\widetilde{M}_*^+[\,b\mid E;\ml;\eta],\tau_0(\ml)\right)\subset\pr{cl}\left(\widetilde{\mathcal{M}}_*^+[\,b\mid E;\ml;\eta],\tau_0(\ml)\right).
\eeq
Введем в рассмотрение следующее множество в пространстве к.-а. мер ограниченной вариации \cite[с.\,88]{36}
\beq{10.23}
\Xi_+^*[\,b\mid E;\ml;\eta]\triangleq\{\mu\in(\pr{add})^+[\mathcal{L};\eta]\mid\mu(E)=b\}.
\eeq
\begin{proposition}
Имеет место следующая цепочка вложений:
\beq{10.24}
\widetilde{M}_*^+[\,b\mid E;\ml;\eta]\subset\widetilde{\mathcal{M}}_*^+[\,b\mid E;\ml;\eta]\subset\Xi_+^*[\,b\mid E;\ml;\eta].
\eeq
\end{proposition}
Д о к а з а т е л ь с т в о. \\  Пусть $\mu\in\widetilde{\mathcal{M}}_*^+[\,b\mid E;\ml;\eta].$  Тогда для некоторой функции $\widetilde{f}\in\mathcal{M}_*^+[\,b\mid E;\ml;\eta]$
\beq{10.25}
\widetilde{\mu}=\widetilde{f}*\eta.
\eeq
Согласно  (\ref{10.6}) и (\ref{10.25}) имеем равенство
\beq{10.26}
\widetilde{\mu}(E)=\int\limits_E\widetilde{f}\,d\eta.
\eeq
При этом согласно (\ref{10.9}) получаем по выбору $\widetilde{f},$ что $0\leqslant\widetilde{\mu}(E).$
Более того, из (\ref{10.10}) и (\ref{10.25}) вытекает, что
\beq{10.27}
\widetilde{\mu}\in(\pr{add})^+[\mathcal{L};\eta].
\eeq

В частности, $\widetilde{\mu}\in(\pr{add})_+[\mathcal{L}].$  Из (\ref{10.13}) имеем, что $\int\limits_E f\,d\eta=b,$   а тогда в силу (\ref{10.26}) $\widetilde{\mu}(E)=b.$ Следовательно (см. (\ref{10.23}), (\ref{10.27})) $\widetilde{\mu}\in\Xi_+^*[\,b\mid E;\ml;\eta].$ Коль скоро выбор $\widetilde{\mu}$ был произвольным, установлено вложение  $\widetilde{\mathcal{M}}_*^+[\,b\mid E;\ml;\eta]\subset\Xi_+^*[\,b\mid E;\ml;\eta].$
С учетом (\ref{10.16}) получаем требуемое утверждение. $\hfill\square$

Отметим также следующее весьма очевидное свойство:
\beq{10.28}
(\pr{add})^+[\mathcal{L};\eta]\in\mathcal{F}[\tau_\otimes^+(\ml)].
\eeq
\begin{zam}
Проверим справедливость свойства (\ref{10.28}), используя формулу  (\ref{10.3}). Кроме того, будем учитывать (\ref{9.25}).

 В связи с (\ref{10.3}) введем в рассмотрение множество
 \beq{10.29}
 \Omega\triangleq\{\mu\in\rr^{\ml}\mid\forall\,L\in\ml \ (\eta(L)=0)\Rightarrow(\mu(L)=0)\}.
 \eeq
Тогда с учетом (\ref{10.3}) получаем равенство
\beq{10.30'}
(\pr{add})^+[\mathcal{L};\eta]=(\pr{add})_+[\mathcal{L}]\cap\Omega.
\eeq

Покажем, что $\Omega\in\mathcal{F}[\otimes^{\ml}(\tau_{\rr})].$
В самом деле, пусть $$\zeta\in\rr^\ml\setminus\Omega.$$
Тогда в силу (\ref{10.29}) можно указать  $\mathbb{L}\in\ml,$ для которого
\beq{10.30}
(\eta(\mathbb{L})=0) \ \& \ (\zeta(\mathbb{L})\neq 0).
\eeq
С учетом того, что $|\zeta(\mathbb{L})|\bn$ и согласно (\ref{6.30}), (\ref{6.33}), имеем
$$\mathbb{N}_{\ml}(\zeta,\{\mathbb{L}\},|\zeta(\mathbb{L})|)=\{\nu\in\rr^\ml:|\zeta(\mathbb{L})-\nu(\mathbb{L})|<|\zeta(\mathbb{L})|\}\in\otimes^\ml(\tau_\rr).$$
Пусть теперь $\xi\in\mathbb{N}_{\ml}(\zeta,\{\mathbb{L}\},|\zeta(\mathbb{L})|).$ Тогда $\xi\in\rr^\ml$ и $$|\zeta(\mathbb{L})-\xi(\mathbb{L})|<|\zeta(\mathbb{L})|.$$

Последнее означает, в частности, что $\xi(\mathbb{L})\neq 0$ (в противном случае имели бы  $|\zeta(\mathbb{L})|<|\zeta(\mathbb{L})|,$ что невозможно). В силу (\ref{10.29}), (\ref{10.30})  $$\xi\in\rr^\ml\setminus\Omega.$$ Поскольку выбор $\xi$ был произвольным, установлено вложение
$$\mathbb{N}_{\ml}(\zeta,\{\mathbb{L}\},|\zeta(\mathbb{L})|)\subset\rr^\ml\setminus\Omega.$$

 Однако выбор $\zeta$ также был произвольным. В итоге
  $$\forall\,\mu\in\rr^\ml\setminus\Omega \ \exists\,\mathcal{K}\in(\pr{Fin})(\mathcal{L}) \ \exists\,\varepsilon\bn:\mathbb{N}_{\ml}(\mu,\mathcal{K},\varepsilon)\subset\rr^\ml\setminus\Omega. $$
Из (\ref{6.31}) имеем теперь, что $$\rr^\ml\setminus\Omega\in\otimes^\ml(\tau_\rr),$$ а потому
\beq{10.31}
\Omega=\rr^\ml\setminus(\rr^\ml\setminus\Omega)\in\mathcal{F}[\otimes^\ml(\tau_{\rr})],
\eeq
что и требовалось доказать.
Напомним, что согласно (\ref{6.44})
$$\tau_\otimes^+(\ml)=\otimes^\ml(\tau_{\rr})|_{(\pr{add})_+[\mathcal{L}]}=\left\{(\pr{add})_+[\mathcal{L}]\cap G:G\in\otimes^\ml(\tau_\rr)\right\}.$$
Теперь воспользуемся (\ref{4.2}): в самом деле, имеем
$$\mathcal{F}[\tau_\otimes^+(\ml)]=\mathcal{F}[\otimes^\ml(\tau_\rr)]|_{(\pr{add})_+[\mathcal{L}]}=\{(\pr{add})_+[\mathcal{L}]\cap F:F\in\mathcal{F}[\otimes^\ml(\tau_\rr)]\}.$$
Поэтому согласно (\ref{10.31}) $(\pr{add})_+[\mathcal{L}]\cap\Omega\in\mathcal{F}[\tau_\otimes^+(\ml)],$ а тогда из (\ref{10.30'}) имеем требуемое свойство (\ref{10.28}). $\hfill\square$
\end{zam}
Как следствие получаем (см. (\ref{9.18})) $*$-слабую замкнутость конуса (\ref{10.3}), а именно:
\beq{10.32}
(\pr{add})^+[\mathcal{L};\eta]\in\mathcal{F}[\tau_*^+(\ml)].
\eeq
Напомним, что согласно положениям раздела 5, (\ref{6.44}) и (\ref{9.25})
\beq{10.33}
\mathcal{F}[\tau_\otimes^+(\ml)]=\{F\in\mathcal{F}[\tau_\otimes(\ml)]\mid F\subset(\pr{add})_+[\mathcal{L}]\},
\eeq
а потому из (\ref{10.28}) вытекает, в частности свойство замкнутости в топологии поточечной сходимости
\beq{10.34}
(\pr{add})^+[\mathcal{L};\eta]\in\mathcal{F}[\tau_\otimes(\ml)].
\eeq
Тогда из (\ref{9.34}) и (\ref{10.34}) получаем полезное свойство
\beq{10.35}
(\pr{add})^+[\mathcal{L};\eta]\in\mathcal{F}[\tau_0(\ml)].
\eeq

Вместе с тем из  (\ref{3.34}), (\ref{9.15'}) и (\ref{9.25}) следует $*$-слабая замкнутость конуса (\ref{6.28}), т.е. $$(\pr{add})_+[\mathcal{L}]\in\mathcal{F}[\tau_*(\ml)];$$ следовательно  (см. раздел 5), из (\ref{9.18}) получаем равенство
\beq{10.36}
\mathcal{F}[\tau_*^+(\ml)]=\{F\in\mathcal{F}[\tau_*(\ml)]\mid F\subset(\pr{add})_+[\mathcal{L}]\}.
\eeq
Из (\ref{10.32}) вытекает (см. (\ref{10.36})), что
\beq{10.37}
(\pr{add})^+[\mathcal{L};\eta]\in\mathcal{F}[\tau_*(\ml)].
\eeq
Вернемся к (\ref{10.23}) и (\ref{10.34}).
\begin{proposition}
Множество $\Xi_+^*[\,b\mid E;\ml;\eta]$ является замкнутым  в топологии $\tau_\otimes(\ml):$
\beq{10.38}
\Xi_+^*[\,b\mid E;\ml;\eta]\in\mathcal{F}[\tau_\otimes(\ml)].
\eeq
\end{proposition}
Д о к а з а т е л ь с т в о.  Согласно  (\ref{2.14}), (\ref{4.2}) и (\ref{6.44})
\beq{10.39}
\mathcal{F}[\tau_\otimes^+(\ml)]=\mathcal{F}[\tau_\otimes(\ml)]|_{(\pr{add})_+[\mathcal{L}]}=\{(\pr{add})_+[\mathcal{L}]\cap G: G\in\tau_\otimes(\ml)\}.
\eeq
С учетом (\ref{10.23}) имеем следующее равенство
\beq{10.40}
\Xi_+^*[\,b\mid E;\ml;\eta]=(\pr{add})^+[\mathcal{L};\eta]\cap\Lambda,
\eeq
где $\Lambda\triangleq\{\mu\in\mathbb{A}(\ml)\mid\mu(E)=b\}.$ В связи с  (\ref{10.40}) полезно учитывать (\ref{10.34}). При этом
$$\Lambda\in\mathcal{F}[\tau_\otimes(\ml)].$$
В самом деле, рассмотрим множество-дополнение $$\mathbb{A}(\ml)\setminus\Lambda.$$
Пусть $\mu_0\in\mathbb{A}(\ml)\setminus\Lambda.$ Тогда $\mu_0(E)\neq b$ и
\beq{10.41}
\varepsilon_0\triangleq|\mu_0(E)-b|\bn.
\eeq
Согласно (\ref{3.29}), (\ref{6.30}), (\ref{6.33}), (\ref{10.41}) и учитывая, что $\mu_0\in\rr^\ml,$   получаем
\beq{10.42}
\mathbb{N}_\ml(\mu_0,\{E\},\varepsilon_0)=\{\nu\in\rr^\ml\mid|\mu_0(E)-\nu(E)|<\varepsilon_0\}\in N_{\otimes^\ml(\tau_\rr)}^0(\mu_0).
\eeq
Как следствие из (\ref{2.219}), (\ref{4.5}), (\ref{6.39}) и (\ref{10.42}) вытекает свойство
\beq{10.43}
\mathbb{N}_\ml(\mu_0,\{E\},\varepsilon_0)\cap\mathbb{A}(\ml)=\{\mu\in\mathbb{A}(\ml)\mid|\mu_0(E)-\mu(E)|<\varepsilon_0\}\in N_{\tau_\otimes(\ml)}^0(\mu_0).
\eeq
Пусть теперь $\mu^0\in\mathbb{N}_\ml(\mu_0,\{E\},\varepsilon_0)\cap\mathbb{A}(\ml).$  Тогда в силу (\ref{10.43}) $\mu^0\in\mathbb{A}(\ml)$
и при этом справедливо неравенство
\beq{10.44}
|\mu_0(E)-\mu^0(E)|<\varepsilon_0.
\eeq

  С учетом (\ref{10.41}) получаем, что $\mu^0(E)\neq b$ (в самом деле, при $\mu^0(E)= b$ имеем $|\mu_0(E)-\mu^0(E)|=\varepsilon_0$ и, в силу (\ref{10.44}),  $\varepsilon_0<\varepsilon_0,$ что невозможно), а потому $\mu^0\notin\Lambda.$  Таким образом,  $\mu^0\in\mathbb{A}(\ml)\setminus\Lambda.$ Однако  выбор $\mu^0$ был произвольным,  а потому получаем вложение
$$\mathbb{N}_\ml(\mu_0,\{E\},\varepsilon_0)\cap\mathbb{A}(\ml)\subset\mathbb{A}(\ml)\setminus\Lambda,$$ тогда
$\mathbb{A}(\ml)\setminus\Lambda\in N_{\tau_\otimes(\ml)}(\mu_0)$ (см. (\ref{3.30})). Поскольку выбор  $\mu_0$ также был произвольным, имеем свойство \beq{10.45}
\mathbb{A}(\ml)\setminus\Lambda\in N_{\tau_\otimes(\ml)}(\mu) \ \forall\,\mu\in \mathbb{A}(\ml)\setminus\Lambda.
\eeq
Из (\ref{3.1000}), (\ref{10.45}) вытекает, что $\mathbb{A}(\ml)\setminus\Lambda\in \tau_\otimes(\ml),
$ и, как следствие, имеем свойство замкнутости $\Lambda:$ $\Lambda\in\mathcal{F}[\tau_\otimes(\ml)].$ Таким образом (см. (\ref{10.35})) $(\pr{add})^+[\mathcal{L};\eta]$, и  $\Lambda$ суть замкнутые в топологии $\tau_\otimes(\ml)$ множества, а потому $$(\pr{add})^+[\mathcal{L};\eta]\cap\Lambda\in\mathcal{F}[\tau_\otimes(\ml)]$$  (см. (\ref{3.34})). С учетом (\ref{10.40}) получаем требуемое свойство (\ref{10.38}). $\hfill\square$

Из (\ref{3.34}), (\ref{9.15'}) и предложения 11.2 вытекает свойство замкнутости
\beq{10.46}
\Xi_+^*[\,b\mid E;\ml;\eta]\in\mathcal{F}[\tau_0(\ml)]
\eeq
и, кроме того,
\beq{10.47}
\Xi_+^*[\,b\mid E;\ml;\eta]\in\mathcal{F}[\tau_*(\ml)].
\eeq
Учтем теперь предложение 11.1. Тогда из (\ref{3.36}) и (\ref{10.46}) имеем
$$\Xi_+^*[\,b\mid E;\ml;\eta]\in[\mathcal{F}[\tau_0(\ml)]]\left(\widetilde{\mathcal{M}}_*^+[\,b\mid E;\ml;\eta]\right),$$
 а из (\ref{3.36}) и (\ref{10.47}) следует,  что
$$\Xi_+^*[\,b\mid E;\ml;\eta]\in[\mathcal{F}[\tau_*(\ml)]]\left(\widetilde{\mathcal{M}}_*^+[\,b\mid E;\ml;\eta]\right).$$
Отметим также (см. (\ref{3.36}), предложения 11.1 и 11.2) очевидное свойство
$$\Xi_+^*[\,b\mid E;\ml;\eta]\in[\mathcal{F}[\tau_\otimes(\ml)]]\left(\widetilde{\mathcal{M}}_*^+[\,b\mid E;\ml;\eta]\right).$$
С учетом (\ref{3.37}) получаем следующие вложения:
\beq{10.48}
\pr{cl}\left(\widetilde{\mathcal{M}}_*^+[\,b\mid E;\ml;\eta],\tau_0(\ml)\right)\subset\Xi_+^*[\,b\mid E;\ml;\eta],
\eeq
\beq{10.49}
\pr{cl}\left(\widetilde{\mathcal{M}}_*^+[\,b\mid E;\ml;\eta],\tau_*(\ml)\right)\subset\Xi_+^*[\,b\mid E;\ml;\eta],
\eeq
\beq{10.50}
\pr{cl}\left(\widetilde{\mathcal{M}}_*^+[\,b\mid E;\ml;\eta],\tau_\otimes(\ml)\right)\subset\Xi_+^*[\,b\mid E;\ml;\eta].
\eeq
Полезно также иметь в виду свойство (\ref{10.16}). Из (\ref{10.21}), (\ref{10.49}) имеем цепочку вложений
\beq{10.50'}
\pr{cl}\left(\widetilde{M}_*^+[\,b\mid E;\ml;\eta],\tau_*(\ml)\right)\subset\pr{cl}\left(\widetilde{\mathcal{M}}_*^+[\,b\mid E;\ml;\eta],\tau_*(\ml)\right)\subset\Xi_+^*[\,b\mid E;\ml;\eta].
\eeq
Аналогичным образом, из (\ref{10.22}) и (\ref{10.48}) получаем
\beq{10.51}
\pr{cl}\left(\widetilde{M}_*^+[\,b\mid E;\ml;\eta],\tau_0(\ml)\right)\subset\pr{cl}\left(\widetilde{\mathcal{M}}_*^+[\,b\mid E;\ml;\eta],\tau_0(\ml)\right)\subset\Xi_+^*[\,b\mid E;\ml;\eta].
\eeq
Для дальнейших целей полезно отметить цепочку вложений (см. (\ref{10.19}), (\ref{10.50'}))
\beq{10.52}
\pr{cl}\left(\widetilde{M}_*^+[\,b\mid E;\ml;\eta],\tau_0(\ml)\right)\subset\pr{cl}\left(\widetilde{\mathcal{M}}_*^+[\,b\mid E;\ml;\eta],\tau_*(\ml)\right)\subset\Xi_+^*[\,b\mid E;\ml;\eta].
\eeq

Заметим, что соотношения (\ref{10.48})--(\ref{10.52}) являются вспомогательными; в следующем разделе  будет показано, что они обращаются в равенства (см. \cite[гл.\,4]{36}).

\section{Слабо абсолютно непрерывные конечно-аддитивные меры и их приближение неопределенными интегралами, 2}\setcounter{equation}{0}\setcounter{proposition}{0}\setcounter{zam}{0}\setcounter{corollary}{0}\setcounter{definition}{0}
   \ \ \ \ \  Сохраняем предположения (\ref{10.1}), (\ref{10.2}) предыдущего раздела, обозначения которого используем далее без дополнительных пояснений. Напомним предложение 5.2:
\beq{11.1}
(\mathbf{D}(E,\ml),\prec)
\eeq
есть непустое направленное множество (НМ). Всюду в дальнейшем полагаем, что
\beq{11.2}
\mathfrak{N}\triangleq\{L\in\ml\mid\eta(L)=0\};
\eeq
   при этом ясно, что $\zer\in\mathfrak{N},$  а потому $\mathfrak{N}\in\pp'(\ml).$  Введем теперь специальные функционалы  на $\ml,$ учитывая, что при $\mu\in(\pr{add})^+[\mathcal{L};\eta]$ и $L\in\ml\setminus\mathfrak{N}$ определена величина $$\frac{\mu(L)}{\eta(L)}\in[0,\infty[.$$

Итак, следуя \cite[раздел 4.3]{38}, полагаем, если $\mu\in(\pr{add})^+[\mathcal{L};\eta],$ то
\beq{11.3}
\theta_+[\mu]:\ml\rightarrow[0,\infty[
\eeq
определяется следующими правилами
\beq{11.4}
\left(\theta_+[\mu](L)\triangleq 0 \ \forall\,L\in\mathfrak{N}\right)\ \& \ \left(\theta_+[\mu](\Lambda)\triangleq\frac{\mu(\Lambda)}{\eta(\Lambda)} \ \forall\,\Lambda\in\ml\setminus\mathfrak{N}\right).
\eeq
Из (\ref{10.3}) и (\ref{11.4}) следует полезное свойство
\beq{11.5}
\mu(L)=\eta(L)\,\theta_+[\mu](L) \ \  \forall\,\mu\in(\pr{add})^+[\mathcal{L};\eta] \ \ \forall\,L\in\ml.
\eeq
\begin{zam}
Проверим (\ref{11.5}), фиксируя $\mu\in(\pr{add})^+[\mathcal{L};\eta]$ и $L\in\ml.$ Если $L\in\ml,$ то $\eta(L)=0,$ поэтому согласно (\ref{10.3}) $\mu(L)=0.$ В итоге значения $\mu(L)$ и $\eta(L)\,\theta_+[\mu](L)$  совпадают с нулем и, стало быть, между собой.  Если же $L\in\ml\setminus\mathfrak{N},$ то $\eta(L)\neq 0,$ а значит  (см. (\ref{11.4}))
$$\eta(L)\,\theta_+[\mu](L) =\eta(L)\,\frac{\mu(L)}{\eta(L)}=\mu(L).$$
Требуемое равенство $\mu(L)=\eta(L)\,\theta_+[\mu](L) $ выполнено во всех возможных случаях. $\hfill\square$
\end{zam}

При   $\mu\in(\pr{add})^+[\mathcal{L};\eta]$ и $L\in\ml$ определено (см. (\ref{11.3})) неотрицательное число $\theta_+[\mu](L).$  Если  теперь $\mathcal{K}\in\mathbf{D}(E,\ml),$ то и при $L\in\mathcal{K}$ имеем $\theta_+[\mu](L)\in[0,\infty[.$
Отметим здесь же, что согласно (\ref{5.0})
$$\forall\,\mathcal{K}\in\mathbf{D}(E,\ml) \ \forall\,x\in E \ \exists\,! L\in\mathcal{K}: x\in L.$$

С учетом этого корректно следующее определение (см. \cite[гл.\,4]{25}): если $\mu\in(\pr{add})^+[\mathcal{L};\eta]$ и $\mathcal{K}\in\mathbf{D}(E,\ml),$ то функция
\beq{11.6}
\Theta_{\mu}^+[\mathcal{K}]:E\rightarrow\rr
\eeq
 определяется посредством  условий
 \beq{11.7}
 \Theta_{\mu}^+[\mathcal{K}](x)\triangleq\theta_+[\mu](L) \ \forall\,L\in\mathcal{K} \ \forall\,x\in L.
 \eeq
 Отметим (см. \cite[раздел\,4.3]{25}), что (\ref{11.6}), (\ref{11.7}) определяют всякий раз неотрицательную ступенчатую функцию
\beq{11.7'}
\Theta_{\mu}^+[\mathcal{K}]\in B_0^+(E,\ml) \ \ \forall\,\mu\in(\pr{add})^+[\mathcal{L};\eta] \ \ \forall\,\mathcal{K}\in\mathbf{D}(E,\ml).
\eeq
\begin{zam}
Проверим (\ref{11.7'}), фиксируя $\mu\in(\pr{add})^+[\mathcal{L};\eta]$ и $\mathcal{K}\in\mathbf{D}(E,\ml).$ Тогда, в частности, $\mathcal{K}$ есть непустое конечное множество, мощность которого определяет
$$n\triangleq|\mathcal{K}|\in\mathbb{N,}$$
причем согласно (\ref{2.20})
\beq{11.8}
(\pr{bi})[\overline{1,n};\mathcal{K}]\neq\zer.
\eeq
С учетом (\ref{11.8}) выберем и зафиксируем биекцию
$(L_i)_{i\in\overline{1,n}}\in(\pr{bi})[\overline{1,n};\mathcal{K}].$
Тогда из (\ref{2.28}) получаем по определению $n,$ что
\beq{11.8'}
(L_i)_{i\in\overline{1,n}}\in\Delta_n(E,\ml).
\eeq
Поскольку при $j\in\overline{1,n},$ в частности, $L_j\in\mathcal{K}$ и $L_j\in\ml,$  определено $\theta_+[\mu](L_j)\in[0,\infty[.$
Тем самым определен кортеж
$(\theta_+[\mu](L_i))_{i\in\overline{1,n}}:\overline{1,n}\rightarrow[0,\infty[.$
Поэтому, в частности, $n\in\mathbb{N},$ $(\theta_+[\mu](L_i))_{i\in\overline{1,n}}\in\rr^n$ и $(L_i)_{i\in\overline{1,n}}\in\Delta_n(E,\ml),$ а также согласно (\ref{7.12}) определена ступенчатая функция
\beq{11.8''}
\varphi\triangleq\sum\limits_{i=1}^n\theta_+[\mu](L_i)\chi_{L_i}\in B_0(E,\ml).
\eeq
Более того, с учетом (\ref{7.14}) получаем
\beq{11.9}
\varphi\in B_0^+(E,\ml).
\eeq

Покажем, что $\Theta_\mu^+[\mathcal{K}]=\varphi.$ В самом деле, пусть $x_*\in E.$ Тогда согласно (\ref{2.22}) $x_*\in\mathbb{L}$ для некоторого $\mathbb{L}\in\mathcal{K}.$ Если $\widetilde{\mathbb{L}}\in\mathcal{K}$ таково, что $x_*\in\widetilde{\mathbb{L}},$ то $\mathbb{L}\cap\widetilde{\mathbb{L}}\neq\zer$ и (см. (\ref{2.22})) $\mathbb{L}=\widetilde{\mathbb{L}}.$
Поэтому имеем очевидное свойство: $x_*\notin L \ \forall\,L\in\mathcal{K}\setminus\{\mathbb{L}\}.$
Итак, $\mathbb{L}\in\mathcal{K}$ и при этом $x_*\in\mathbb{L},$ а тогда согласно (\ref{11.7'})
\beq{11.10}
\Theta_{\mu}^+[\mathcal{K}](x_*)=\theta_+[\mu](\mathbb{L}).
\eeq
С другой стороны, по выбору $(L_i)_{i\in\overline{1,n}}$ имеем, что $\mathbb{L}=L_r$ для некоторого $r\in\overline{1,n}.$ Итак, $x_*\in L_r;$ как следствие, $x_*\notin L_j \ \forall\,j\in\overline{1,n}\setminus\{r\}$ (см. (\ref{2.21})). В итоге  имеем равенство $\varphi(x_*)=\theta_+[\mu](L_r)=\theta_+[\mu](\mathbb{L})$ (см. определение $\varphi$).  С учетом (\ref{11.10}) получаем, что
\beq{11.11}
\Theta_{\mu}^+[\mathcal{K}](x_*)=\varphi(x_*).
\eeq
Поскольку $x_*$ выбиралось произвольно, то  из (\ref{11.11}) следует, что $$\Theta_{\mu}^+[\mathcal{K}](x)=\varphi(x) \ \forall\,x\in E.$$
Это  означает, что $\Theta_{\mu}^+[\mathcal{K}]=\varphi,$ а тогда (см. (\ref{11.9})) $\Theta_{\mu}^+[\mathcal{K}]\in B_0^+(E,\ml).$ $\hfill\square$
\end{zam}

Свойство (\ref{11.7'}) дополним полезным следствием, формулируемым в виде отдельного предложения.

\begin{proposition}
Если $\mu\in\Xi_+^*[\,b\mid E;\ml;\eta]$ и $\mathcal{K}\in\mathbf{D}(E,\ml),$ то
\beq{11.12}
\Theta_{\mu}^+[\mathcal{K}]\in M_*^+[\,b\mid E;\ml;\eta].
\eeq
\end{proposition}
Д о к а з а т е л ь с т в о. Из (\ref{10.23}) и (\ref{11.7'})  следует, что для выбранных $\mu$ и~$\mathcal{K}$
$$\Theta_{\mu}^+[\mathcal{K}]\in B_0^+(E,\ml).$$
Поэтому (см. (\ref{10.12})) осталось доказать равенство
\beq{11.13}
\int\limits_E \Theta_{\mu}^+[\mathcal{K}]\,d\eta=b.
\eeq

Действуя по аналогии с рассуждением в замечании 12.2 для $n\triangleq|\mathcal{K}|\in\mathbb{N}$ подбираем $(L_i)_{i\in\overline{1,n}}\in(\pr{bi})[\overline{1,n};\mathcal{K}]$ и получаем упорядоченное конечное разбиение $(L_i)_{i\in\overline{1,n}}$ (\ref{11.8'}), после чего определяем кортеж $(\theta_+[\mu](L_i))_{i\in\overline{1,n}}$ в $[0,\infty[,$  в терминах которого конструируем $\varphi\in B_0^+(E,\ml)$ (\ref{11.8''}), (\ref{11.9}). Как и в замечании 12.2 имеем равенство $\Theta_{\mu}^+[\mathcal{K}]=\varphi,$ откуда в силу определения 9.1 (см. также (\ref{8.3}))
$$\int\limits_E\Theta_{\mu}^+[\mathcal{K}]\,d\eta=\int\limits_E^{(el)}\Theta_{\mu}^+[\mathcal{K}]\,d\eta=\sum\limits_{i=1}^n\theta_+[\mu](L_i)\eta(L_i).$$

Теперь используем (\ref{11.5}), получая как следствие, что
\beq{11.14}
\int\limits_E\Theta_{\mu}^+[\mathcal{K}]\,d\eta=\sum\limits_{i=1}^n\mu(L_i).
\eeq
Однако в силу конечной аддитивности (см. (\ref{10.23})) имеем с учетом (\ref{11.8'}) равенство
$$\sum\limits_{i=1}^n\mu(L_i)=\mu(E)=b.$$
В силу (\ref{11.14}) получаем требуемое равенство (\ref{11.13}), чем и завершаем доказательство свойства (\ref{11.12}). $\hfill\square$

В дальнейших рассуждениях будем использовать представление конечной аддитивности в терминах неупорядоченных конечных разбиений в разделе 7 (см. следствие 7.1). В этой связи заметим, что при $\mu\in\Xi_+^*[\,b\mid E;\ml;\eta]$ определен оператор
\beq{11.15}
\Theta_{\mu}^+[\cdot]\triangleq\left(\Theta_{\mu}^+[\mathcal{K}]\right)_{\mathcal{K}\in\mathbf{D}(E,\ml)}\in M_*^+[\,b\mid E;\ml;\eta]^{\mathbf{D}(E,\ml)},
\eeq
 действующий из $\mathbf{D}(E,\ml)$ в  $M_*^+[\,b\mid E;\ml;\eta].$ Данный оператор легко преобразуется в отображение из $\mathbf{D}(E,\ml)$ в $\widetilde{M}_*^+[\,b\mid E;\ml;\eta],$  а именно: при $\mu\in\Xi_+^*[\,b\mid E;\ml;\eta]$  в виде
 \beq{11.16}
 \Theta_{\mu}^+[\cdot]*\eta\triangleq\left(\Theta_{\mu}^+[\mathcal{K}]*\eta\right)_{\mathcal{K}\in\mathbf{D}(E,\ml)}\in\widetilde{M}_*^+[\,b\mid E;\ml;\eta]^{\mathbf{D}(E,\ml)}
 \eeq
  имеем требуемую модификацию оператора (\ref{11.15}).   С учетом предложения 5.2 получаем важное свойство: если $\mu\in\Xi_+^*[\,b\mid E;\ml;\eta],$  то
   \beq{11.17}
   (\mathbf{D}(E,\ml),\prec,\Theta_{\mu}^+[\cdot]*\eta)
   \eeq
   есть  направленность в множестве $\widetilde{M}_*^+[\,b\mid E;\ml;\eta].$ Конструкция (\ref{11.16}), (\ref{11.17}) в сочетании со следствием 7.1 будет существенно использоваться для доказательства того, что последнее множество всюду плотно в $\Xi_+^*[\,b\mid E;\ml;\eta].$
\begin{proposition}
Если $\mu\in\Xi_+^*[\,b\mid E;\ml;\eta],$ то $$\left(\mathbf{D}(E,\ml),\prec,\Theta_{\mu}^+[\cdot]*\eta\right)\stackrel{\tau_0(\ml)}\longrightarrow \mu.$$
\end{proposition}
Д о к а з а т е л ь с т в о . Фиксируем $\mu\in\Xi_+^*[\,b\mid E;\ml;\eta],$ получая, в частности, $\mu\in\mathbb{A}(\ml).$  С учетом (\ref{3.29}) и (\ref{6.40}) имеем, что
 $$N^0_{\tau_0(\ml)}(\mu)=\{G\in\tau_0(\ml)\mid\mu\in G\}$$
 есть непустое семейство открытых (в смысле $\tau_0(\ml)$) окрестностей $\mu.$ Кроме того,
\beq{11.18}
N_{\tau_0(\ml)}(\mu)=\left\{H\in\pp(\mathbb{A}(\ml))\mid\exists\,G\in N^0_{\tau_0(\ml)}(\mu): G\subset H\right\}
\eeq
есть семейство всевозможных (не обязательно открытых) окрестностей $\mu.$ Воспользуемся представлением (\ref{4.14}). Пусть $\mathbb{H}\in N_{\tau_0(\ml)}(\mu).$  С учетом (\ref{11.18}) подберем открытую окрестность
\beq{11.19}
\mathbb{G}\in N^0_{\tau_0(\ml)}(\mu)
\eeq
со следующим свойством вложенности:
\beq{11.20}
\mathbb{G}\subset\mathbb{H}.
\eeq
Тогда в силу (\ref{11.19})  $\mathbb{G}\in\tau_0(\ml)$  и при этом $\mu\in\mathbb{G}.$ Вместе с тем, согласно (\ref{4.0}) и (\ref{6.40})
\beq{11.21}
\mathbb{G}=\mathbf{G}\cap\mathbb{A}(\ml),
\eeq
где $\mathbf{G}\in\otimes^\ml(\tau_\partial).$ С учетом  (\ref{6.35}) и (\ref{11.21}) для некоторого $\mathcal{K}\in(\pr{Fin})[\mathcal{L}]$ и множества
\beq{11.22}
\Gamma\triangleq\mathbb{N}_{\ml}^{(\partial)}(\mu,\mathcal{K})\cap\mathbb{A}(\ml)=\{\nu\in\mathbb{A}(\ml)\mid\mu(L)=\nu(L) \ \forall\,L\in\mathcal{K}\}
\eeq
 реализуется следующее вложение:
 \beq{11.23}
 \Gamma\subset\mathbb{G}.
 \eeq
Напомним, что в данном рассуждении мы интересуемся сходимостью направленности (\ref{11.17}) к фиксированной к.-а. мере $\mu.$

Выберем произвольно множество  $\Lambda\in\mathcal{K}\setminus\{\zer\}.$ В силу  (\ref{9.1}) имеем, что для некоторого  $m\in\mathbb{N}$
$$\Delta_m(E\setminus\Lambda,\ml)\neq\zer.$$
С учетом этого выберем разбиение $$(\Lambda_i)_{i\in\overline{1,m}}\in\Delta_m(E\setminus\Lambda,\ml).$$
Тогда $(\Lambda_i)_{i\in\overline{1,m}}:\overline{1,m}\rightarrow\ml$ и при этом
\beq{11.24}
\left(E\setminus\Lambda=\bigcup\limits_{i=1}^m\Lambda_i\right)\ \& \ \left(\Lambda_{i_1}\cap\Lambda_{i_2}=\zer \ \ \forall\,i_1\in\overline{1,m} \ \ \forall\,i_2\in\overline{1,m}\setminus\{i_1\}\right).
\eeq

Поскольку (по выбору $\mathcal{K}$ и $\Lambda$), в частности, $\Lambda\in\ml,$  введем в рассмотрение кортеж $$\left(\Lambda^{(i)}\right)_{i\in\overline{1,m+1}}:\overline{1,m+1}\rightarrow\ml$$ по следующему правилу:
\beq{11.25}
\left(\Lambda^{(j)}\triangleq\Lambda_j \ \forall\,j\in\overline{1,m}\right)\ \& \ \left(\Lambda^{(m+1)}\triangleq\Lambda\right).
\eeq
Легко понять (см. (\ref{2.35})), что
\beq{11.26}
\left(\Lambda^{(i)}\right)_{i\in\overline{1,m+1}}\in\Delta_{m+1}(E,\ml),
\eeq
а потому (см. (\ref{2.23}))
\beq{11.27}
\mathfrak{L}\triangleq\left\{\Lambda^{(i)}:i\in\overline{1,m+1}\right\}\in\mathbf{D}(E,\ml),
\eeq
причем $\Lambda\in\mathfrak{L}.$ Выберем произвольно разбиение
\beq{11.28}
\mathcal{G}\in\mathbf{D}(E,\ml),
\eeq
для которого
\beq{11.29}
\mathfrak{L}\prec\mathcal{G}.
\eeq
Отметим, что разбиение $\mathcal{G}$ (\ref{11.28}) со свойством (\ref{11.29}) непременно существует. Таковым, в частности, является <<исходное>> разбиение $\mathfrak{L}$ (\ref{11.27}). По выбору $\Lambda$ имеем, что
$$\Lambda\in\mathfrak{L}\setminus\{\zer\},$$
откуда согласно предложению 5.1 следует  (см. (\ref{11.27})--(\ref{11.29})) для некоторого разбиения
\beq{11.30}
\mathcal{D}\in\mathbf{D}(\Lambda,\ml)
\eeq
 вложение
\beq{11.31}
\mathcal{D}\subset\mathcal{G}.
\eeq

Поскольку $\Lambda\in\ml,$ то в силу (\ref{11.30}) получаем (см. следствие 7.1) свойство
\beq{11.32}
\nu(\Lambda)=\sum\limits_{L\in\mathcal{D}}\nu(L) \ \ \forall\,\nu\in(\pr{add})[\mathcal{L}].
\eeq
Из (\ref{11.32}), в частности, следует, что
\beq{11.33}
\mu(\Lambda)=\sum\limits_{L\in\mathcal{D}}\mu(L).
\eeq
С другой стороны, согласно (\ref{11.16}) и (\ref{11.32}) имеем также равенство
\beq{11.34}
\left(\Theta_{\mu}^+[\mathcal{G}]*\eta\right)(\Lambda)=\sum\limits_{L\in\mathcal{D}}\left(\Theta_{\mu}^+[\mathcal{G}]*\eta\right)(L).
\eeq
При этом $\Theta_{\mu}^+[\mathcal{G}]\in B_0^+(E,\ml)$ согласно (\ref{11.7}). В силу (\ref{11.7}) и (\ref{11.30}) реализуется система равенств
\beq{11.35}
\Theta_{\mu}^+[\mathcal{G}](x)=\theta_+[\mu](L) \ \forall\,L\in\mathcal{G} \ \forall\,x\in L.
\eeq
В силу  (\ref{10.6}) получаем равенства
$$\left(\Theta_{\mu}^+[\mathcal{G}]*\eta\right)(\Lambda)=\int\limits_\Lambda\Theta_{\mu}^+[\mathcal{G}]\,d\eta$$
и, кроме того,
$$\left(\Theta_{\mu}^+[\mathcal{G}]*\eta\right)(L)=\int\limits_L\Theta_{\mu}^+[\mathcal{G}]\,d\eta\ \forall\,L\in\mathcal{G}.$$
С учетом (\ref{10.4})
\beq{11.36}
\left(\Theta_{\mu}^+[\mathcal{G}]*\eta\right)(\Lambda)=\int\limits_E\Theta_{\mu}^+[\mathcal{G}]\chi_\Lambda\,d\eta;
\eeq
с другой стороны, имеем систему равенств
\beq{11.37}
\left(\Theta_{\mu}^+[\mathcal{G}]*\eta\right)(L)=\int\limits_E\Theta_{\mu}^+[\mathcal{G}]\chi_L\,d\eta \ \ \forall\,L\in\mathcal{G}.
\eeq

Напомним также (\ref{11.5}). С учетом (\ref{11.30}) получаем, что
\beq{11.38}
\mu(L)=\eta(L)\theta_+[\mu](L) \ \forall\,L\in\mathcal{D}.
\eeq
 Вместе с тем, что согласно (\ref{7.6}) и (\ref{11.35}) $$\Theta_{\mu}^+[\mathcal{G}]\chi_L=\left(\Theta_{\mu}^+[\mathcal{G}](x)\chi_L(x)\right)_{x\in E}=\theta_+[\mu](L)\chi_L \ \forall\,L\in\mathcal{G}.$$
Поэтому из (\ref{8.18}), (\ref{9.2}) и (\ref{11.37}) следует  при $L\in\mathcal{G}$
$$\left(\Theta_{\mu}^+[\mathcal{G}]*\eta\right)(L)=\int\limits_E\theta_+[\mu](L)\chi_L\,d\eta=\theta_+[\mu](L)\int\limits_E\chi_L\,d\eta=
\theta_+[\mu](L)\eta(L).$$
С учетом (\ref{11.31}) и (\ref{11.38}) получаем систему равенств $$\left(\Theta_{\mu}^+[\mathcal{G}]*\eta\right)(L)=\mu(L) \ \forall\,L\in\mathcal{D}.$$ Из (\ref{11.33}) и (\ref{11.34}) имеем теперь следующую цепочку равенств:
\beq{11.39}
\left(\Theta_{\mu}^+[\mathcal{G}]*\eta\right)(\Lambda)=\sum\limits_{L\in\mathcal{D}}\mu(L)=\mu(\Lambda).
\eeq
Итак (см. (\ref{11.29}), (\ref{11.39})), установлена импликация $$(\mathfrak{L}\prec\mathcal{G})\Rightarrow\left(\left(\Theta_{\mu}^+[\mathcal{G}]*\eta\right)(\Lambda)=\mu(\Lambda)\right).$$

Поскольку разбиение  $\mathcal{G}$ (со свойством  (\ref{11.29})) выбиралось произвольно, то  $\forall\,\mathcal{S}\in\mathbf{D}(E,\ml)$
$$(\mathfrak{L}\prec\mathcal{S})\Rightarrow\left(\left(\Theta_{\mu}^+[\mathcal{S}]*\eta\right)(\Lambda)=\mu(\Lambda)\right).$$

Учитывая (\ref{11.27}), получаем, как следствие, что
\beq{11.40}
\exists\,\mathcal{R}\in\mathbf{D}(E,\ml) \ \forall\,\mathcal{S}\in\mathbf{D}(E,\ml) \ \ (\mathcal{R}\prec\mathcal{S})\Rightarrow\left(\left(\Theta_{\mu}^+[\mathcal{S}]*\eta\right)(\Lambda)=\mu(\Lambda)\right).
\eeq
Однако выбор $\Lambda$ был произвольным. Поэтому (см. (\ref{11.40})) установлено свойство
\begin{multline}\label{11.41}
\forall\,L\!\in\!\mathcal{K}\setminus\{\zer\} \ \exists\,\mathcal{R}\!\in\!\mathbf{D}(E,\ml) \ \forall\,\mathcal{S}\!\in\!\mathbf{D}(E,\ml) \ \
(\mathcal{R}\prec\mathcal{S})\Rightarrow \\ \Rightarrow\left(\left(\Theta_{\mu}^+[\mathcal{S}]*\eta\right)(L)=\mu(L)\right).
\end{multline}
Из (\ref{6.2}) легко следует (см. \cite[(2.2.19)]{30})
$\nu(\zer)=0 \ \forall\,\nu\in(\pr{add})[\mathcal{L}].$
Тогда, в частности, $\mu(\zer)=0.$ Вместе с тем согласно (\ref{11.16}) $$\left(\Theta_{\mu}^+[\mathcal{S}]*\eta\right)(\zer)=0 \ \forall\,\mathcal{S}\in\mathbf{D}(E,\ml).$$

Коль скоро $\mathbf{D}(E,\ml)\neq\zer,$  последнее означает, в частности, справедливость свойства  $$\exists\,\mathcal{R}\in\mathbf{D}(E,\ml) \ \forall\,\mathcal{S}\in\mathbf{D}(E,\ml) \ \ (\mathcal{R}\prec\mathcal{S})\Rightarrow\left(\left(\Theta_{\mu}^+[\mathcal{S}]*\eta\right)(\zer)=\mu(\zer)\right).$$
С учетом (\ref{11.41}) получаем, что
\beq{11.42}
\forall\,L\in\mathcal{K} \ \exists\,\mathcal{R}\in\mathbf{D}(E,\ml) \ \forall\,\mathcal{S}\in\mathbf{D}(E,\ml) \ \
(\mathcal{R}\prec\mathcal{S})\Rightarrow\left(\left(\Theta_{\mu}^+[\mathcal{S}]*\eta\right)(L)=\mu(L)\right).
\eeq
Вместе с тем, как легко видеть,  из (\ref{4.13}) и предложения 5.2 рассуждением по индукции устанавливается положение
$$\forall\,n\in\mathbb{N} \ \forall\,(\mathcal{R}_i)_{i\in\overline{1,n}}\in\mathbf{D}(E,\ml)^n \ \exists\,\widetilde{\mathcal{R}}\in\mathbf{D}(E,\ml):\mathcal{R}_j\prec\mathcal{R} \ \forall\,j\in\overline{1,n}.$$
Тогда из (\ref{11.42}) с учетом конечности семейства $\mathcal{K}$  имеем следующее свойство:
\beq{11.43}
\exists\,\mathcal{R}\in\mathbf{D}(E,\ml) \ \forall\mathcal{S}\in\mathbf{D}(E,\ml) \ \
(\mathcal{R}\prec\mathcal{S})\Rightarrow\left((\Theta_{\mu}^+[\mathcal{S}]*\eta)(L)=\mu(L) \ \forall\,L\in\mathcal{K}\right).
\eeq
Из (\ref{11.22}), (\ref{11.23}) и (\ref{11.43}) получаем, что
\beq{11.44}
\exists\,\mathcal{R}\in\mathbf{D}(E,\ml) \ \forall\mathcal{S}\in\mathbf{D}(E,\ml) \ \
(\mathcal{R}\prec\mathcal{S})\Rightarrow\left(\Theta_{\mu}^+[\mathcal{S}]*\eta\in\mathbb{G}\right).
\eeq
В силу  (\ref{11.20}), (\ref{11.44})
$\exists\,\mathcal{R}\in\mathbf{D}(E,\ml) \ \forall\mathcal{S}\in\mathbf{D}(E,\ml) \ \
(\mathcal{R}\prec\mathcal{S})\Rightarrow\left(\Theta_{\mu}^+[\mathcal{S}]*\eta\in\mathbb{H}\right).$
Поскольку выбор $\mathbb{H}$ был произвольным, установлено положение
\beq{11.45}
\forall\,H\in N_{\tau_0(\ml)}(\mu) \ \exists\,\mathcal{R}\in\mathbf{D}(E,\ml) \ \forall\mathcal{S}\in\mathbf{D}(E,\ml) \ \
(\mathcal{R}\prec\mathcal{S})\Rightarrow\left(\Theta_{\mu}^+[\mathcal{S}]*\eta\in H\right).
\eeq
Учитывая (\ref{4.14}) и (\ref{11.45}), получаем требуемое свойство: $$\left(\mathbf{D}(E,\ml),\prec,\Theta_{\mu}^+[\cdot]*\eta\right)\stackrel{\tau_0(\ml)}\longrightarrow \mu.$$ $\hfill\square$
\begin{proposition}
Справедливо равенство $$\Xi_+^*[\,b\mid E;\ml;\eta]=\pr{cl}\left(\widetilde{M}_*^+[\,b\mid E;\ml;\eta],\tau_0(\ml)\right)=\pr{cl}\left(\widetilde{\mathcal{M}}_*^+[\,b\mid E;\ml;\eta],\tau_0(\ml)\right).$$
\end{proposition}
Д о к а з а т е л ь с т в о. Напомним (\ref{10.51}). С другой стороны,  из свойства 1*) раздела~4 вытекает, что $\pr{cl}\left(\widetilde{M}_*^+[\,b\mid E;\ml;\eta],\tau_0(\ml)\right)$ есть множество всех $\mu\in\mathbb{A}(\ml),$ для каждой из которых существует направленность $(\mathbf{D},\sqsubseteq,f)$ в $\widetilde{M}_*^+[\,b\mid E;\ml;\eta]$ со свойством  $\left(\mathbf{D},\sqsubseteq,f\right)\stackrel{\tau_0(\ml)}\longrightarrow \mu.$

С учетом (\ref{11.16}) и предложения 12.2 имеем для всякой к.-а. меры $\mu\in\Xi_+^*[\,b\mid E;\ml;\eta]$ свойство  сходимости направленности (\ref{11.17}) к $\mu,$ т.е.
$$\left(\mathbf{D}(E,\ml),\prec,\Theta_{\mu}^+[\cdot]*\eta\right)\stackrel{\tau_0(\ml)}\longrightarrow \mu.$$
Поэтому $\Xi_+^*[\,b\mid E;\ml;\eta]\subset\pr{cl}\left(\widetilde{M}_*^+[\,b\mid E;\ml;\eta],\tau_0(\ml)\right).$ Используя (\ref{10.51}), получаем требуемую систему равенств. $\hfill\square$
\begin{proposition}
Если $\mu\in \Xi_+^*[\,b\mid E;\ml;\eta],$ то $$\left(\mathbf{D}(E,\ml),\prec,\Theta_{\mu}^+[\cdot]*\eta\right)\stackrel{\tau_*(\ml)}\longrightarrow \mu.$$
\end{proposition}
Д о к а з а т е л ь с т в о.  Отметим, что
\beq{11.46}
\mu\in(\pr{add})_+[\mathcal{L}].
\eeq
Кроме того, $\tau_*^+(\ml)\subset\tau_0^+(\ml)$ в силу (\ref{9.18}). В виде (\ref{11.17}) имеем (см. предложение 12.2) направленность в $\widetilde{M}_*^+[\,b\mid E;\ml;\eta],$  сходящуюся к $\mu$ в топологии $\tau_0(\ml).$
С учетом (\ref{4.5}), (\ref{4.6}), (\ref{4.16}) получаем поэтому, что $$\forall\,H\in N_{\tau_0^+(\ml)}(\mu) \ \exists\,\mathcal{R}\in\mathbf{D}(E,\ml) \ \forall\mathcal{S}\in\mathbf{D}(E,\ml) \ \ (\mathcal{R}\prec\mathcal{S})\Rightarrow\left(\Theta_{\mu}^+[\mathcal{S}]*\eta\in H\right).$$
Из (\ref{9.18}) как следствие
\beq{11.47}
\forall\,H\in N_{\tau_*^+(\ml)}(\mu) \ \exists\,\mathcal{R}\in\mathbf{D}(E,\ml) \ \forall\mathcal{S}\in\mathbf{D}(E,\ml) \ \ (\mathcal{R}\prec\mathcal{S})\Rightarrow\left(\Theta_{\mu}^+[\mathcal{S}]*\eta\in H\right).
\eeq
(необходимо учесть (\ref{3.29}), (\ref{3.30}): из (\ref{3.29}) и (\ref{9.18}) вытекает вложение $$N_{\tau_*^+(\ml)}^0(\mu)\subset N_{\tau_0^+(\ml)}^0(\mu)$$ и, как следствие, $N_{\tau_*^+(\ml)}(\mu)\subset N_{\tau_0^+(\ml)}(\mu)$). С учетом (\ref{4.6}), (\ref{9.18}) получаем
$$N_{\tau_*^+(\ml)}(\mu)=\left\{H\cap(\pr{add})_+[\mathcal{L}]: \ H\in N_{\tau_*(\ml)}(\mu)\right\},$$
а потому  (см. (\ref{11.47})) имеем, что
\beq{11.48}
\forall\,H\in N_{\tau_*(\ml)}(\mu) \ \exists\,\mathcal{R}\in\mathbf{D}(E,\ml) \ \forall\mathcal{S}\in\mathbf{D}(E,\ml) \ \ (\mathcal{R}\prec\mathcal{S})\Rightarrow\left(\Theta_{\mu}^+[\mathcal{S}]*\eta\in H\right).
\eeq

Из (\ref{4.14}) и (\ref{11.48}) получаем требуемое утверждение. $\hfill\square$
\begin{proposition}
Справедливо равенство
\beq{11.49}
\Xi_+^*[\,b\mid E;\ml;\eta]=\pr{cl}\left(\widetilde{M}_*^+[\,b\mid E;\ml;\eta],\tau_*(\ml)\right)=\pr{cl}\left(\widetilde{\mathcal{M}}_*^+[\,b\mid E;\ml;\eta],\tau_*(\ml)\right).
\eeq
\end{proposition}
Д о к а з а т е л ь с т в о.  Напомним (\ref{10.50}). С другой стороны, согласно предложению 12.4 (см. (\ref{11.17})) при всяком выборе $\mu\in\Xi_+^*[\,b\mid E;\ml;\eta]$  можно указать направленность $(\mathbf{D},\preceq,g),$ в $\widetilde{M}_*^+[\,b\mid E;\ml;\eta],$ для которой
$$\left(\mathbf{D},\preceq,g\right)\stackrel{\tau_*(\ml)}\longrightarrow \mu;$$
это означает, (см. $1^*$) в разделе 4), что
$$\mu\in\pr{cl}\left(\widetilde{M}_*^+[\,b\mid E;\ml;\eta],\tau_*(\ml)\right).$$
Итак, $\Xi_+^*[\,b\mid E;\ml;\eta]\subset\pr{cl}\left(\widetilde{M}_*^+[\,b\mid E;\ml;\eta],\tau_*(\ml)\right),$  что в сочетании с (\ref{10.50}) и (\ref{10.50'}) означает справедливость (\ref{11.49}). $\hfill\square$

Из предложений 12.3 и 12.5 вытекает цепочка равеств
\begin{multline}\label{11.50}
\Xi_+^*[\,b\mid E;\ml;\eta]=\pr{cl}\left(\widetilde{M}_*^+[\,b\mid E;\ml;\eta],\tau_*(\ml)\right)=\pr{cl}\left(\widetilde{\mathcal{M}}_*^+[\,b\mid E;\ml;\eta],\tau_*(\ml)\right)=\\=\pr{cl}\left(\widetilde{M}_*^+[\,b\mid E;\ml;\eta],\tau_0(\ml)\right)=\pr{cl}\left(\widetilde{\mathcal{M}}_*^+[\,b\mid E;\ml;\eta],\tau_0(\ml)\right).
\end{multline}
Из (\ref{10.23}), (\ref{11.47}) и (\ref{11.50}) следует, в частности, что $\Xi_+^*[\,b\mid E;\ml;\eta]$  $*$-слабо замкнуто и сильно ограничено, а потому (см. (\ref{9.13}))
\beq{11.51}
\Xi_+^*[\,b\mid E;\ml;\eta]\in(\tau_*(\ml)-\pr{comp})[\mathbb{A}(\ml)].
\eeq
Тогда топология
\begin{multline}\label{11.52}
\tau_\eta^*[\ml\mid b]\triangleq\tau_*(\ml)|_{\Xi_+^*[\,b\mid E;\ml;\eta]}=\{\Xi_+^*[\,b\mid E;\ml;\eta]\cap G: G\in\tau_*(\ml)\}\in \\ \in(\pr{top})[\Xi_+^*[\,b\mid E;\ml;\eta]]
\end{multline}
превращает непустое множество (\ref{10.23}) в компактное ТП. Поскольку (\ref{9.7}) есть, как легко  видеть, $T_2$-пространство, то (см. (\ref{4.5}), (\ref{4.6}), (\ref{10.52}))
\beq{11.53}
\left(\Xi_+^*[\,b\mid E;\ml;\eta],\tau_\eta^*[\ml\mid b]\right)
\eeq
также является $T_2$-пространством. Иными словами, (\ref{11.53})\,--- непустой компакт. С учетом (\ref{11.50}) и (\ref{11.52}) имеем (см. (\ref{4.4})) равенства
\begin{multline}\label{11.54}
\pr{cl}\left(\widetilde{M}_*^+[\,b\mid E;\ml;\eta],\tau_\eta^*[\ml\mid b]\right)=\pr{cl}\left(\widetilde{M}_*^+[\,b\mid E;\ml;\eta],\tau_*(\ml)\right)\cap\Xi_+^*[\,b\mid E;\ml;\eta]=\\=\Xi_+^*[\,b\mid E;\ml;\eta],
\end{multline}
\begin{multline}\label{11.55}
\pr{cl}\left(\widetilde{\mathcal{M}}_*^+[\,b\mid E;\ml;\eta],\tau_\eta^*[\ml\mid b]\right)=\pr{cl}\left(\widetilde{\mathcal{M}}_*^+[\,b\mid E;\ml;\eta],\tau_*(\ml)\right)\cap\Xi_+^*[\,b\mid E;\ml;\eta]=\\=\Xi_+^*[\,b\mid E;\ml;\eta].
\end{multline}

Таким образом (см. (\ref{10.14}), (\ref{10.15}), предложение 11.1),  в виде
\beq{11.56}
f\mapsto f*\eta: M_*^+[\,b\mid E;\ml;\eta]\rightarrow\Xi_+^*[\,b\mid E;\ml;\eta],
\eeq
\beq{11.57}
f\mapsto f*\eta: \mathcal{M}_*^+[\,b\mid E;\ml;\eta]\rightarrow\Xi_+^*[\,b\mid E;\ml;\eta]
\eeq
имеем отображения, реализующие погружение соответствующих множеств в пространстве ярусных функций в (один и тот же) компакт (\ref{11.53}). Далее ограничимся работой с  (\ref{11.56}), предлагая заинтересованному читателю самостоятельно рассмотреть аналогичное использование (\ref{11.57}). Для краткости условимся отображение (\ref{11.56}) обозначать через $\mathcal{I.}$ Итак, $$\mathcal{I}:M_*^+[\,b\mid E;\ml;\eta]\rightarrow\Xi_+^*[\,b\mid E;\ml;\eta]$$  и при этом справедливо $$\mathcal{I}(f)=f*\eta \ \forall\,f\in M_*^+[\,b\mid E;\ml;\eta].$$
Тогда из (\ref{11.54}) следует с очевидностью, что
\beq{11.58}
\pr{cl}\left(\mathcal{I}^{\,1}(M_*^+[\,b\mid E;\ml;\eta]),\tau_\eta^*[\ml\mid b]\right)=\Xi_+^*[\,b\mid E;\ml;\eta],
\eeq
поскольку (см. (\ref{10.14})) справедливо равенство
\beq{11.59}
\widetilde{M}_*^+[\,b\mid E;\ml;\eta]=\mathcal{I}^{\,1}(M_*^+[\,b\mid E;\ml;\eta]).
\eeq
Вместе с тем из (\ref{11.50}) и (\ref{11.59}) вытекает цепочка равенств
\beq{11.60}
\Xi_+^*[\,b| E;\ml;\eta]=\pr{cl}\left(\mathcal{I}^{\,1}(M_*^+[\,b| E;\ml;\eta]),\tau_*(\ml)\right)=\pr{cl}\left(\mathcal{I}^{\,1}(M_*^+[\,b| E;\ml;\eta]),\tau_0(\ml)\right).
\eeq
Свойства (\ref{11.59}), (\ref{11.60}) будут существенно использоваться при построении и последующем исследовании расширения одной естественной задачи о достижимости, рассматриваемой в следующем разделе.

В заключении раздела коснемся одной возможности, связанной с предложением 12.2 и касающейся вопросов интегрирования ступенчатых функций. В этой связи напомним, что (\ref{11.17}) определяет направленность в $\widetilde{M}_*^+[\,b\mid E;\ml;\eta].$ Разумеется, в частности, (\ref{11.17}) является направленностью в $\mathbb{A}(\ml).$ В этой связи уместно рассмотреть вопросы, связанные со сходимостью в $(\mathbb{A}(\ml),\tau_0(\ml))$  уже для любых направленностей упомянутого типа.
\begin{proposition}

Если $(\mathbb{D},\preceq,g)$ есть направленность в $\mathbb{A}(\ml)$ и $\mu\in\mathbb{A}(\ml),$ то
\begin{multline}\label{11.61}
\left(\left(\mathbb{D},\preceq,g\right)\stackrel{\tau_0(\ml)}\longrightarrow\mu\right)\Rightarrow\\ \Rightarrow\left(\!\forall\,f\!\in\! B_0(E,\ml) \ \exists\,\delta_1\in\mathbb{D} \  \forall\,\delta_2\!\in\!\mathbb{D} \ (\delta_1\preceq\delta_2)\Rightarrow\biggl(\int\limits_E f\,d g(\delta_2)=\int\limits_E f\,d\mu\biggl)\right).
\end{multline}
\end{proposition}
Д о к а з а т е л ь с т в о.  Фиксируем направленность $(\mathbb{D},\preceq,g)$ в $\mathbb{A}(\ml)$  и $\mu\in\mathbb{A}(\ml).$ Пусть истинна посылка доказываемой импликации (\ref{11.61}), т.е.
\beq{11.61'}
\left(\mathbf{D},\preceq,g\right)\stackrel{\tau_0(\ml)}\longrightarrow \mu.
\eeq
  Тогда в силу (\ref{4.14}) имеем следующее свойство:
  \beq{11.62}
  \forall\,H\in N_{\tau_0(\ml)}(\mu) \ \exists\,d_1\in\mathbb{D} \ \forall\,d_2\in\mathbb{D} \ \ (d_1\preceq d_2)\Rightarrow(g(d_2)\in H).
  \eeq
  Напомним, что согласно (\ref{6.35}), (\ref{6.40})  $\otimes^{\ml}(\tau_\partial)\in(\pr{top})[\rr^\ml],$ $\mathbb{A}(\ml)\in\pp'(\rr^\ml)$ и при этом $ \tau_0(\ml)=\otimes^{\ml}(\tau_\partial)|_{\mathbb{A}(\ml)}.$ Кроме того, отметим, что $\mu\in\mathbb{A}(\ml);$ поэтому из (\ref{2.219}) и (\ref{4.6}) получаем
  \beq{11.63}
  N_{\tau_0(\ml)}(\mu)=N_{\otimes^{\ml}(\tau_\partial)}(\mu)|_{\mathbb{A}(\ml)}=\{\mathbb{A}(\ml)\cap H: \ H\in N_{\otimes^{\ml}(\tau_\partial)}(\mu)\}.
  \eeq
  Возвращаясь к (\ref{6.34}), отметим, что при $\mathcal{K}\in\pr{Fin}(\ml)$ и $\nu\in\mathbb{N}_\ml^{(\partial)}(\mu,\mathcal{K})$ $$\mathbb{N}_\ml^{(\partial)}(\mu,\mathcal{K})=\mathbb{N}_\ml^{(\partial)}(\nu,\mathcal{K}).$$
    Это означает, в частности, что (см. (\ref{6.35})) при $\mathcal{K}\in\pr{Fin}(\ml)$ $$\mathbb{N}_\ml^{(\partial)}(\mu,\mathcal{K})\in\otimes^\ml(\tau_\partial)$$  и при этом $\mu\in\mathbb{N}_\ml^{(\partial)}(\mu,\mathcal{K}).$ Как следствие (\ref{3.29}) $\mathbb{N}_\ml^{(\partial)}(\mu,\mathcal{K})\in N_{\otimes^{\ml}(\tau_\partial)}^0(\mu) \ \forall\,\mathcal{K}\in\pr{Fin}(\ml).$
    С~учетом (\ref{11.63}) получаем с очевидностью, что
\beq{11.64}
\mathbb{N}_\ml^{(\partial)}(\mu,\mathcal{K})\cap\mathbb{A}(\ml)\!=\!\{\nu\in\mathbb{A}(\ml)\mid\mu(L)=\nu(L) \, \forall\,L\!\in\!\mathcal{K}\}\!\in\! N_{\tau_0(\ml)}(\mu) \ \forall\,\mathcal{K}\in\pr{Fin}(\ml).
\eeq
Выберем произвольно $\varphi\in B_0(E,\ml).$ Тогда в силу  (\ref{7.12}) для некоторых $n\in\mathbb{N}, (\alpha_i)_{i\in\overline{1,n}}\in\rr^n$ и $(L_i)_{i\in\overline{1,n}}\in\Delta_n(E,\ml)$  справедливо равенство
\beq{11.65}
\varphi=\sum\limits_{i=1}^n\alpha_i\chi_{L_i}.
\eeq
Согласно определению 9.1
\beq{11.66}
\int\limits_E^{(\pr{el})} \varphi\,d\nu=\sum\limits_{i=1}^n\alpha_i\nu(L_i) \ \forall\,\nu\in\mathbb{A}(\ml).
\eeq
При этом, в частности, $\varphi\in B(E,\ml)$ и по определению 9.2 имеем  $\int\limits_E^{(\pr{el})} \varphi\,d\nu=\int\limits_E \varphi\,d\nu \ \forall\,\nu\in\mathbb{A}(\ml).$ В итоге получаем следующие равенства
\beq{11.67}
\int\limits_E \varphi\,d\nu=\sum\limits_{i=1}^n\alpha_i\nu(L_i) \ \forall\,\nu\in\mathbb{A}(\ml).
\eeq
В частности, из (\ref{11.67}) извлекаются свойства
\beq{11.68}
\biggl(\,\int\limits_E \varphi\,d\mu=\sum\limits_{i=1}^n\alpha_i\mu(L_i)\biggl)\ \& \biggl(\,\int\limits_E \varphi\,d g(\delta)=\sum\limits_{i=1}^n\alpha_ig(\delta)(L_i) \ \forall\,\delta\in\mathbb{D}\biggl).
\eeq
Теперь используем свойство 1*) раздела 3: $$\mathfrak{L}\triangleq\{L_i:i\in\overline{1,n}\}\in\mathbb{D}(E,\ml)$$
и,  в частности, (см. (\ref{2.22})) имеет место свойство
\beq{11.69}
\mathfrak{L}\in\pr{Fin}(\ml).
\eeq
Поэтому (см. (\ref{11.64}), (\ref{11.69})) получаем важное положение
\begin{multline}\label{11.70}
\mathbb{N}_{\ml}^{(\partial)}(\mu,\mathfrak{L})\cap\mathbb{A}(\ml)=  \{\nu\in\mathbb{A}(\ml)\mid\mu(\ml)=\nu(\ml) \ \forall\,L\in\mathfrak{L}\}=\\=\{\nu\in\mathbb{A}(\ml)\mid\mu(L_i)=\nu(L_i) \ \forall\,i\in\overline{1,n}\,\}\in N_{\tau_0(\ml)}(\mu).
\end{multline}
Согласно (\ref{11.62}) имеем, как следствие, для некоторого $\mathbf{d}\in\mathbb{D},$ что $\forall\,\delta\in\mathbb{D}$
$$(\mathbf{d}\preceq\delta)\Rightarrow\left(g(\delta)\in\mathbb{N}_{\ml}^{(\partial)}(\mu,\mathfrak{L})\cap\mathbb{A}(\ml)\right).$$

С учетом (\ref{11.70}) при $\delta\in\mathbb{D}$ со свойством $\mathbf{d}\preceq \delta$ реализуется следующая система равенств
$$\mu(L_i)=g(\delta)(L_i) \ \forall\,i\in\overline{1,n};$$
поэтому из (\ref{11.68}) извлекается равенство $$\int\limits_E \varphi\,d\mu=\int\limits_E \varphi\,d g(\delta).$$
Таким образом, установлена, в частности, следующая импликация $$(\mathbf{d}\preceq\delta)\Rightarrow\biggl(\,\int\limits_E \varphi\,d\mu=\int\limits_E \varphi\,d g(\delta)\biggl).$$ Поскольку выбор $\varphi$ был произвольным, справедливо, что $$\forall\,f\in B_0(E,\ml) \ \exists\,\delta_1\in\mathbb{D} \ \forall\,\delta_2\in\mathbb{D} \  \ (\delta_1\preceq\delta_2)\Rightarrow\biggl(\,\int\limits_E f\,d g(\delta_2)=\int\limits_E f\,d\mu\biggl).$$

Тем самым (см. (\ref{11.61'})) доказана импликация (\ref{11.61}). $\hfill\square$
\begin{corollary}
Если $\mu\in\Xi_+^*[\,b\mid E;\ml;\eta],$ то $\forall\,h\in B_0(E,\ml)$ $ \exists\,\mathcal{K}_1\!\in\!\mathbb{D}(E,\ml)$ $\forall\,\mathcal{K}_2\in\mathbb{D}(E,\ml)$
\beq{11.71}
(\mathcal{K}_1\prec\mathcal{K}_2)\Rightarrow\biggl(\,\int\limits_E h\,\Theta_{\mu}^+[\mathcal{K}_2]\,d\eta=\int\limits_E h\,d\mu\biggl).
\eeq
\end{corollary}
Д о к а з а т е л ь с т в о. Напомним, что согласно (\ref{11.16}), предложениям 12.2 и 12.6
\begin{multline}\label{11.72}
\forall\,h\in B_0(E,\ml) \ \exists\,\mathcal{K}_1\in\mathbb{D}(E,\ml) \ \forall\,\mathcal{K}_2\in\mathbb{D}(E,\ml) \\
(\mathcal{K}_1\prec\mathcal{K}_2)\Rightarrow\biggl(\int\limits_E h\,d(\Theta_{\mu}^+[\mathcal{K}_2]*\eta)=\int\limits_E h\,d\mu\biggl).
\end{multline}
Вместе с тем с учетом (\ref{10.5}), (\ref{10.12}) и (\ref{11.7'}) имеем, что $$\int\limits_E h\,d(\Theta_{\mu}^+[\mathcal{K}]*\eta)=\int\limits_E h\,\Theta_{\mu}^+[\mathcal{K}]\,d \eta \ \forall\,h\in B(E,\ml) \ \forall\,\mathcal{K}\in\mathbb{D}(E,\ml).$$
Поэтому (см. (\ref{11.72})) получаем с очевидностью (\ref{11.71}).$\hfill\square$
\begin{corollary}
Пусть $\mu\in\Xi_+^*[\,b\mid E;\ml;\eta]$  и $\mathbb{K}\in\pr{Fin}(B_0(E,\ml)).$ Тогда $\exists\,\mathcal{K}_1\in\mathbb{D}(E,\ml) \ \forall\,\mathcal{K}_2\in\mathbb{D}(E,\ml)$
\beq{11.73}
(\mathcal{K}_1\prec\mathcal{K}_2)\Rightarrow\biggl(\,\int\limits_E h\,\Theta_{\mu}^+[\mathcal{K}_2]\,d\eta=\int\limits_E h\,d\mu \ \forall\,h\in\mathbb{K}\biggl).
\eeq
\end{corollary}
Д о к а з а т е л ь с т в о. Из (\ref{4.13}) и предложения 6.2 имеем по индукции, что
  $\forall\,m\in~\mathbb{N} \ \forall\,(\mathcal{K}^{(i)})_{i\in\overline{1,m}}\in\mathbf{D}(E,\ml)^m \ \ \exists\,\widetilde{\mathcal{K}}\in\mathbf{D}(E,\ml):$ \beq{11.73'} \mathcal{K}^{(j)}\prec\widetilde{\mathcal{K}} \ \forall\,j\in\overline{1,m}.
 \eeq
  В силу (\ref{2.20})  для $\mathbf{k}\triangleq|\mathbb{K}|\in\mathbb{N}$ получаем свойство $(\pr{bi})[\overline{1,\mathbf{k}};\mathbb{K}]\neq\zer.$
 С учетом этого выберем (см. (\ref{2.15'}), (\ref{2.15''})) биекцию $(\mathbf{h}_i)_{i\in\overline{1,\mathbf{k}}}\in(\pr{bi})[\overline{1,\mathbf{k}};\mathbb{K}].$ Тогда, в частности, $(\mathbf{h}_i)_{i\in\overline{1,\mathbf{k}}}\in(\pr{su})[\overline{1,\mathbf{k}};\mathbb{K}],$ а потому
 \beq{11.74}
 \mathbb{K}=\{\mathbf{h}_i:i\in\overline{1,\mathbf{k}}\}.
 \eeq
  Заметим, что по следствию 12.1
 $\forall j\in\overline{1,\mathbf{k}} \ \exists\,\mathcal{K}_1\in\mathbb{D}(E,\ml) \ \forall\,\mathcal{K}_2\in\mathbb{D}(E,\ml)$
 \beq{11.75}
(\mathcal{K}_1\prec\mathcal{K}_2)\Rightarrow\biggl(\int\limits_E \mathbf{h}_j \Theta_{\mu}^+[\mathcal{K}_2]\,d\eta=\int\limits_E \mathbf{h}_j
\,d\mu\biggl).
\eeq
 Используя  конечность множества $\overline{1,\mathbf{k}},$ получаем, что
$\exists\,(\widehat{\mathcal{K}}_i)_{i\in\overline{1,\mathbf{k}}}\in\mathbf{D}(E,\ml)^{\mathbf{k}} \ \forall\,j\in\overline{1,\mathbf{k}} \ \forall\,\mathcal{K}\in\mathbf{D}(E,\ml)$
\beq{11.76}
\left(\widehat{\mathcal{K}}^{(j)}\prec\mathcal{K}\right)\Rightarrow\biggl(\int\limits_E \mathbf{h}_j \Theta_{\mu}^+[\mathcal{K}]\,d\eta=\int\limits_E \mathbf{h}_j\,d\mu\biggl).
\eeq
С учетом (\ref{11.73'}) подберем теперь
\beq{11.77}
\mathfrak{K}\in\mathbf{D}(E,\ml),
\eeq
  для которого имеет место свойство $\widehat{\mathcal{K}}^{(j)}\prec\mathfrak{K} \ \forall\,j\in\overline{1,\mathbf{k}}.$
Из (\ref{11.76}), (\ref{11.77}) получаем, что $\forall\,\mathcal{K}\in\mathbf{D}(E,\ml)$
$$\left(\mathfrak{K}\prec\mathcal{K}\right)\Rightarrow\biggl(\,\int\limits_E \mathbf{h}_j \Theta_{\mu}^+[\mathcal{K}]\,d\eta=\int\limits_E \mathbf{h}_j\,d\mu \ \forall\,j\in\overline{1,\mathbf{k}}\biggl).$$
Как следствие (см. (\ref{11.74})) имеем, что $\forall\,\mathcal{K}\in\mathbf{D}(E,\ml)$ $$\left(\mathfrak{K}\prec\mathcal{K}\right)\Rightarrow\biggl(\,\int\limits_E h\, \Theta_{\mu}^+[\mathcal{K}]\,d\eta=\int\limits_E h\,d\mu \ \forall\,h\in\mathbb{K}\biggl);$$ это означает (см. (\ref{11.77})) справедливость (\ref{11.73}). $\hfill\square$

\newpage

\begin{center} \section*{Глава 3. Абстрактная задача о достижимости с ограничениями асимптотического характера } \end{center}
\section{Вариант задачи  о достижимости при ограничениях моментного характера}\setcounter{equation}{0}\setcounter{proposition}{0}\setcounter{zam}{0}\setcounter{corollary}{0}\setcounter{definition}{0}
   \ \ \ \ \ Следуем далее обозначениям и определениям разделов 10--12, фиксируя натуральные числа
   \beq{12.1}
   (n\in\mathbb{N})\ \& \ (N\in\mathbb{N}),
   \eeq
а также следующие два кортежа ярусных функций:
\beq{12.2}
(\pi_i)_{i\in\overline{1,n}}:\overline{1,n}\rightarrow B(E,\ml),
\eeq
\beq{12.3}
(s_j)_{j\in\overline{1,N}}:\overline{1,N}\rightarrow B(E,\ml).
\eeq
Итак, $(\pi_i)_{i\in\overline{1,n}}\in B(E,\ml)^n$ и  $(s_j)_{j\in\overline{1,N}}\in B(E,\ml)^N$ (можно также говорить о ярусных вектор-функциях). Пусть, кроме того, задано множество
\beq{12.4}
Y\in\pp'(\rr^N).
\eeq
Множество $Y$ задает следующее ограничение на выбор $f\in M_*^+[\,b\mid E;\ml;\eta]:$
\beq{12.5}
\biggl(\int\limits_E s_j f\,d\eta\biggl)_{j\in\overline{1,N}}\in Y.
\eeq
 Какими могут быть векторы (числовые кортежи) $\left(\int\limits_E \pi_i f\,d\eta\right)_{i\in\overline{1,n}}\in\rr^n$ при условии, что $f$ удовлетворяет (\ref{12.5})? Данный вопрос можно сформулировать и несколько иначе, полагая, что $S$ есть отображение (на самом деле $N$-вектор-функционал)
  \beq{12.6}
  f\mapsto\biggl(\int\limits_E s_j f\,d\eta\biggl)_{j\in\overline{1,N}}: M_*^+[\,b\mid E;\ml;\eta]\rightarrow\rr^N,
  \eeq
  а $\Pi$ есть  соответственно отображение ($n$-вектор-функционал)
 \beq{12.6}
  f\mapsto\biggl(\int\limits_E \pi_i f\,d\eta\biggl)_{i\in\overline{1,n}}: M_*^+[\,b\mid E;\ml;\eta]\rightarrow\rr^n.
  \eeq
   Требуется определить  множество
     \beq{12.8}
     \Pi^1(S^{-1}(Y))\in\pp(\rr^n),
     \eeq
играющее объективно роль, аналогичную области достижимости в задачах управления. Интересен вопрос о том, что будет происходить с  (\ref{12.8}), если $Y$-ограничение (\ref{12.5}) будет ослабляться, т.е. $Y$ будет заменяться каким-то другим множеством $\widetilde{Y}\in\pp'(\rr^N)$ со свойством $Y\subset\widetilde{Y}.$ Разумеется, процедура ослабления $Y$-ограничения может быть весьма различной.
Так, например, можно представить себе равномерное ослабление $Y$-ограничения по всем направлениям, что формализуется посредством использования  $\varepsilon$-~окрестностей $Y$ (при $\varepsilon>0$), определяемых в той или иной норме пространства $\rr^n.$ С другой стороны, можно запретить ослабление данного ограничения по части координат. Сравнение двух упомянутых вариантов может представлять не только теоретический, но и определенный практический интерес.

В то же  время само множество   (\ref{12.8}) может оказаться и <<неинтересным>> с практической точки зрения. Это отвечает ситуации, когда отсутствует устойчивость задачи при ослаблении $Y$-ограничения (\ref{12.5}).

Условимся, что при $k\in\mathbb{N}$ через $\|\cdot\|_k$ обозначается следующая норма пространства $\rr^k:$ $\|\cdot\|_k\triangleq(\|x\|_k)_{x\in\rr^k},$  где при  всяком выборе $y=(y_i)_{i\in\overline{1,k}}\in\rr^k$
$$\|y\|_k\triangleq\max\limits_{1\leqslant i\leqslant k}|y(i)|\in[0,\infty[.$$
Разумеется, норма $\|\cdot\|_k$  порождает метрику $$\left(y^{(1)},y^{(2)}\right)\mapsto\left\|y^{(1)}-y^{(2)}\right\|_k:\rr^k\times\rr^k\rightarrow[0,\infty[,$$
которая, в свою очередь, порождает обычную топологию $\tau_{\rr}^{(k)}\in(\pr{top})[\rr^k]$ (покоординатной сходимости), базу которой составляют открытые шары, соответствующие упомянутой метрике, т.е. множества
$$\left\{\widetilde{y}\in\rr^k\mid\|\widetilde{y}-y\|_k<\varepsilon\right\}, \ y\in\rr^k, \ \varepsilon\bn.$$

В частности,  $\|\cdot\|_N$\,--- норма $\rr^N$ и при $\varepsilon\bn$ полагаем
\beq{12.9}
\mathbb{O}_N(Y,\varepsilon)\triangleq\left\{z\in\rr^N\mid\exists\,y\in Y: \ \|y-z\|_N<\varepsilon\right\}.
\eeq

С учетом (\ref{12.9}) определяем естественный способ <<равномерного>>  ослабления $Y$-ограничения, при котором требование (\ref{12.5}) на выбор  $f\in M_*^+[\,b\mid E;\ml;\eta]$ заменяется следующим:
\beq{12.10}
\biggl(\int\limits_E s_j f\,d\eta\biggl)_{j\in\overline{1,N}}\in\mathbb{O}_N(Y,\varepsilon),
\eeq
где $\varepsilon\bn.$ Число $\varepsilon,$ $\varepsilon>0$ в (\ref{12.10}) зачастую зафиксировать не удается, а потому рассматривается вся серия условий упомянутого типа. В этой связи полагаем, что
\beq{12.11}
\mathbb{Y}_\varepsilon\triangleq\left\{f\in M_*^+[\,b\mid E;\ml;\eta]\mid\biggl(\,\int\limits_E s_j f\,d\eta\biggl)_{j\in\overline{1,N}}\in\mathbb{O}_N(Y,\varepsilon)\right\} \ \forall\,\varepsilon\bn.
\eeq
Семейство $\mathfrak{Y}\triangleq\{\mathbb{Y}_\varepsilon:\varepsilon\bn\}$  (всех множеств (\ref{12.11})) рассматриваем в качестве объекта, определяющего \emph{ограничения асимптотического характера} (ОАХ): обращаясь к $\mathfrak{Y},$ мы ориентируемся на соблюдение условия (\ref{12.5}) <<с нарастающей точностью>>. При этом взамен одного множества (\ref{12.8}) естественным образом возникает непустое семейство множеств
\beq{12.12}
\Pi^1(S^{\,-1}(\mathbb{O}_N(Y,\varepsilon))), \ \varepsilon\bn.
\eeq

Учитывая то, что параметр $\varepsilon,$ $\varepsilon>0,$ по смыслу следует устремлять к нулю, вполне логичным представляется характеризация зависимости $$\varepsilon\mapsto\Pi^1(S^{\,-1}(\mathbb{O}_N(Y,\varepsilon))): \ \ ]0,\infty[\rightarrow\pp(\rr^n)$$
обобщенным пределом вида
\beq{12.13}
\bigcap\limits_{\varepsilon\bn}\pr{cl}\left(\Pi^1(S^{\,-1}(\mathbb{O}_N(Y,\varepsilon))),\tau_{\rr}^{(n)}\right)\in\pp(\rr^n).
\eeq
Мы будем называть (\ref{12.13}) \emph{множеством притяжения} (МП) в случае ОАХ, определяемых семейством $\mathfrak{Y}.$ Точный смысл этого будет пояснен ниже, а сейчас рассмотрим другой вариант ослабления $Y$-ограничения.

Итак, пусть $J\in\pp(\overline{1,N}).$ Индексы $j\in J$ предполагаются особыми: в <<направлениях>>, ими определяемых, требуется теперь точное совмещение компонент вектора в левой части (\ref{12.5}) с соответствующими компонентами вектора из $Y.$
Полагаем, что при $\varepsilon\bn$
\begin{multline}\label{12.14}
\widehat{\mathbb{O}}_N(Y,\varepsilon|J)\triangleq\bigl\{(z_i)_{i\in\overline{1,N}}\in\rr^N\mid\exists\,(y_i)_{i\in\overline{1,N}}\in Y:  (y_j=z_j \ \forall\,j\in J) \ \& \\ \& \  (|y_j-z_j|<\varepsilon \ \forall\,j\in\overline{1,N}\setminus J)\bigl\}.
\end{multline}
Как и прежде, параметр $\varepsilon;$ $\varepsilon>0,$ используемый в  (\ref{12.14}), не фиксируется, а потому определяется бесконечная серия условий
$$\biggl(\,\int\limits_E s_j f\,d\eta\biggl)_{j\in\overline{1,N}}\in\widehat{\mathbb{O}}_N(Y,\varepsilon| J), \ \  \ \varepsilon\bn.$$
В этой связи полагаем далее, что
\beq{12.15}
\widehat{\mathbb{Y}}_\varepsilon[J]\triangleq\left\{f\in M_*^+[\,b\mid E;\ml;\eta]\mid \biggl(\,\int\limits_E s_j f\,d\eta\biggl)_{j\in\overline{1,N}}\in\widehat{\mathbb{O}}_N(Y,\varepsilon| J)\right\} \ \ \forall\,\varepsilon\bn.
\eeq
Семейство $\widehat{\mathfrak{Y}}[J]\triangleq\{\widehat{\mathbb{Y}}_\varepsilon[J]:\varepsilon\bn\}$ определяет другой вариант ОАХ. При этом множество (\ref{12.8}) заменяется семейством $\Pi^1(S^{\,-1}(\mathbb{O}_N(Y,\varepsilon| J))),  \ \ \varepsilon\in ]0,\infty[\rightarrow\pp(\rr^n).$

Учитывая естественное желание устремить параметр $\varepsilon$ к $0,$ для характеризации (многозначной) зависимости
$$\varepsilon\mapsto\Pi^1(S^{\,-1}(\widehat{\mathbb{O}}_N(Y,\varepsilon|J))):  \ ]0,\infty[\rightarrow\pp(\rr^n)$$
логично использовать обобщенный предел
\beq{12.16}
\bigcap\limits_{\varepsilon\bn}\pr{cl}\left(\Pi^1(S^{\,-1}(\widehat{\mathbb{O}}_N(Y,\varepsilon|J))),\tau_{\rr}^{(n)}\right)\in\pp(\rr^n),
\eeq
который будем именовать МП, соответствующим ОАХ, определяемым семейством $\widehat{\mathfrak{Y}}[J].$
Сопоставляя (\ref{12.9}) и (\ref{12.14}), заметим, что
\beq{12.17}
\widehat{\mathbb{O}}_N(Y,\varepsilon|J)\subset \mathbb{O}_N(Y,\varepsilon) \ \ \forall\,\varepsilon\bn.
\eeq
Как следствие получаем, что
$$\Pi^1(S^{\,-1}(\widehat{\mathbb{O}}_N(Y,\varepsilon| J)))\subset\Pi^1(S^{\,-1}(\mathbb{O}_N(Y,\varepsilon))) \ \ \forall\,\varepsilon\bn.$$
Последнее свойство наследуется соответствующими замыканиями, а тогда из (\ref{12.13}) и (\ref{12.16}) получаем весьма очевидную оценку для МП:
\beq{12.18}
\bigcap\limits_{\varepsilon\bn}\pr{cl}\left(\Pi^1(S^{\,-1}(\widehat{\mathbb{O}}_N(Y,\varepsilon|J))),\tau_{\rr}^{(n)}\right)\subset
\bigcap\limits_{\varepsilon\bn}\pr{cl}\left(\Pi^1(S^{\,-1}(\mathbb{O}_N(Y,\varepsilon))),\tau_{\rr}^{(n)}\right).
\eeq
Вопрос о совпадении МП (\ref{12.13}), (\ref{12.16}) представляет не только теоретический, но, в ряде случаев, и определенный практический интерес: при его положительном решении можно говорить об асимптотической нечувствительности задачи при ослаблении $Y$-ограничения в <<направлениях>>, определяемых  индексами из $J.$ В связи с (\ref{12.13}), (\ref{12.16}) в следующем разделе напомним одну общую конструкцию, касающуюся применения расширений абстрактных задач о достижимости. Рассматриваемая в настоящем разделе задача является частным случаем, к рассмотрению которого мы вернемся позднее.

В заключении раздела отметим одно очевидное свойство: из (\ref{12.9}), (\ref{12.14}) имеем по определению $\|\cdot\|_N,$ что
\beq{12.18'}
(J=\zer)\Rightarrow\left(\mathbb{O}_N(Y,\varepsilon)=\widehat{\mathbb{O}}_N(Y,\varepsilon|J) \ \forall\,\varepsilon\bn\right).
\eeq
Как следствие имеем при $J=\zer,$ что (см. (\ref{12.11}), (\ref{12.15})) $$\mathbb{Y}_\varepsilon=\widehat{\mathbb{Y}}_\varepsilon[J] \ \ \forall\,\varepsilon\bn.$$
В итоге получаем импликацию
\beq{12.19}
(J=\zer)\Rightarrow\left(\mathfrak{Y}=\widehat{\mathfrak{Y}}[J]\right).
\eeq
Кроме того, при $J=\zer$ множества (\ref{12.13}) и (\ref{12.16}) совпадают (см. (\ref{12.18'})). Свойство (\ref{12.19}) будет использоваться при построении конкретной конструкции расширения.

Совсем кратко обсудим другой крайний случай, полагая до конца настоящего раздела $J=\overline{1,N}.$ Тогда  в силу  (\ref{12.14}) имеем
\beq{12.20}
\widehat{\mathbb{O}}_N(Y,\varepsilon|J)=\left\{z\in\rr^n\mid\exists\,(y_i)_{i\in\overline{1,N}}\in Y: y_j=z_j \ \forall\,j\in\overline{1,N}\right\}=Y \ \forall\,\varepsilon\bn.
\eeq
Согласно (\ref{12.20}) наша вторая схема введения ОАХ приводит к невозмущенной задаче: $\widehat{\mathfrak{Y}}[J]=\{Y\}$ и множество (\ref{12.16}) совпадает с $$\pr{cl}\left(\Pi^1(S^{\,-1}(Y)),\tau_{\rr}^{(n)}\right).$$

Поэтому сравнение множеств (\ref{12.13}), (\ref{12.16}) может быть полезным с точки зрения исследования условий устойчивости множества достижимости при ослаблении полной системы ограничений.

\section{Абстрактная задача о достижимости и ее <<конечно-аддитивное>> расширение}\setcounter{equation}{0}\setcounter{proposition}{0}\setcounter{zam}{0}\setcounter{corollary}{0}\setcounter{definition}{0}
   \ \ \ \ \ В этом разделе фиксируем непустое множество $X$  произвольной природы, ТП  $(\mathbb{X},\tau),$ $\mathbb{X}\neq\zer,$ и отображение
 \beq{13.1}
 \mathbf{f}:X\rightarrow\mathbb{X}.
 \eeq
   Если $\mathcal{X}\in\pp'(\pp(X)),$ то полагаем, что  $(\pr{\mathbf{as}})[X;\mathbb{X};\tau;\mathbf{f};\mathcal{X}]$  есть $\pr{def}$ множество всех $\mathbf{x}\in\mathbb{X},$ для каждого из которых можно указать такую направленность $(D,\preceq,g)$ в множестве $X,$ что
   \beq{13.1'}
   (\forall\,\Sigma\in\mathcal{X} \ \exists\,d\in D \ \forall\,\delta\in D \ (d\preceq\delta)\Rightarrow(g(\delta)\in\Sigma)) \ \& \ \left((D,\preceq,\mathbf{f}\circ g)\stackrel{\tau}\rightarrow \mathbf{x}\right).
   \eeq
   Множество $(\pr{\mathbf{as}})[X;\mathbb{X};\tau;\mathbf{f};\mathcal{X}]$ рассматриваем как МП на значениях $\mathbf{f}$ (\ref{13.1}) при ОАХ, определенных посредством $\mathcal{X};$ см. в этой связи \cite[раздел 4]{49}.
Полезно отметить некоторые важные частные случаи. Так
$$\beta[X]\triangleq\left\{\mathcal{B}\in\pp'(\pp(X))\mid\forall\,B_1\in\mathcal{B} \ \forall\,B_2\in\mathcal{B} \ \exists\,B_3\in \mathcal{B}: B_3\subset B_1\cap B_2\right\}$$
определяет семейство, элементами которого    являются всевозможные направленные подсемейства $\pp(X).$ Если $\mathcal{X}\in\beta[X],$ то \cite[(3.7)]{49}
\beq{13.2}
(\pr{\mathbf{as}})[X;\mathbb{X};\tau;\mathbf{f};\mathcal{X}]=\bigcap\limits_{\Sigma\in\mathcal{X}}\pr{cl}\left(f^1(\Sigma),\tau\right).
\eeq

Данный случай будет достаточным для всех  последующих построений (отметим очевидную аналогию (\ref{12.13}), (\ref{12.16}) и (\ref{13.2})).

Еще один естественный вопрос касается использования в (\ref{13.1'}) направленностей. Точнее, нельзя ли ограничиться применением последовательностей, что идейно соответствует подходу Дж. Варги \cite[гл.\,III]{44}? Ответ в общем случае отрицателен (см., в частности, пример в \cite{50}). Однако известны весьма общие условия (см. \cite{48}), при которых исчерпывающая реализация МП достигается в классе последовательностей. Все же введем в качестве самостоятельного объекта секвенциальное МП в общем случае; напомним в этой связи, что $X^\mathbb{N}$ есть множество всех последовательностей со значениями в $X.$

Итак, при $\mathcal{X}\in\pp'(\pp(X))$ полагаем, что $(\pr{\mathbf{sas}})[X;\mathbb{X};\tau;\mathbf{f};\mathcal{X}]$  есть $\pr{def}$ множество всех $\mathbf{x}\in\mathbb{X},$ для каждого из которых существует такая последовательность $(x_i)_{i\in\mathbb{N}}\in X^\mathbb{N},$ что
\beq{13.3}
\left(\forall\,\Sigma\in\mathcal{X} \ \exists\,k\in\mathbb{N}: \ x_j\in\Sigma \ \forall\,j\in\overrightarrow{k,\infty}\right)\ \& \ \left((\mathbf{f}(x_i))_{i\in\mathbb{N}}\stackrel{\tau}\rightarrow x\right).
\eeq
 Ясно, что (\ref{13.3}) является по существу частным случаем (\ref{13.1'}); при этом используется (\ref{4.11}). Поскольку последовательность является вариантом направленности (см. раздел 4), то
 \beq{13.4}
 (\pr{\mathbf{sas}})[X;\mathbb{X};\tau;\mathbf{f};\mathcal{X}]\subset(\pr{\mathbf{as}})[X;\mathbb{X};\tau;\mathbf{f};\mathcal{X}] \ \forall\,\mathcal{X}\in\pp'(\pp(X)).
 \eeq
Отметим условия \cite{48}, при которых МП и секвенциальное МП совпадают. Речь идет о случае, когда семейство $\mathcal{X}$ обладает счетной базой, а $(\mathbb{X},\tau)$ есть ТП с первой аксиомой счетности. Отметим следующее \\
С в о й с т в о. Если $\mathcal{X}\in\pp'(\pp(X))$ таково, что $\exists\,(\Sigma_i)_{i\in\mathbb{N}}\in\mathcal{X}^\mathbb{N} \ \forall\,\Sigma\in\mathcal{X}$
\beq{13.5}
\{k\in\mathbb{N}\mid\Sigma_k\subset\Sigma\}\neq\zer,
\eeq
 а также $\forall\,\mathbf{x}\in\mathbb{X} \ \exists\,(H_i)_{i\in\mathbb{N}}\in N_\tau(x)^\mathbb{N} \ \forall\,\mathbb{H}\in N_\tau(x)$
 \beq{13.6}
 \{k\in\mathbb{N}\mid H_k\subset\mathbb{H}\}\neq\zer,
 \eeq
то справедливо следующее равенство
\beq{13.7}
(\pr{\mathbf{as}})[X;\mathbb{X};\tau;\mathbf{f};\mathcal{X}]=(\pr{\mathbf{sas}})[X;\mathbb{X};\tau;\mathbf{f};\mathcal{X}].
\eeq
Заметим, кстати, что для всех наших целей достаточно ограничиться случаем $\mathcal{X}\in\beta[X].$ Ясно, что при $k\in \mathbb{N}$ ТП $\left(\rr^k,\tau_{\rr}^{(k)}\right)$  обладает свойством (\ref{13.6}), а семейства $\mathfrak{Y}$ и $\widehat{\mathfrak{Y}}[J],$ рассматриваемые в разделе 13, ---  обладают каждое свойством, подобным (\ref{13.5}).
Мы предлагаем читателю убедиться в этом самостоятельно.

Сейчас возвращаемся к общему случаю МП (см. (\ref{13.1'})). Заметим, что конструкцию на основе (\ref{13.2}) можно применить и в данном (общем) случае: если $\mathcal{X}\in\pp'(\pp(X)),$ то введем семейство
$$\mathfrak{X}\triangleq\left\{\bigcap\limits_{\Sigma\in\mathcal{K}}\Sigma:\mathcal{K}\in\pr{Fin}(\mathcal{X})\right\},$$
для которого очевидным образом имеют место свойства $\mathfrak{X}\in\beta[X]$ и $$(\pr{\mathbf{as}})[X;\mathbb{X};\tau;\mathbf{f};\mathcal{X}]=(\pr{\mathbf{as}})[X;\mathbb{X};\tau;\mathbf{f};\mathfrak{X}],$$
а потому согласно (\ref{13.2})
\beq{13.8}
(\pr{\mathbf{as}})[X;\mathbb{X};\tau;\mathbf{f};\mathcal{X}]=\bigcap\limits_{\Sigma\in\mathfrak{X}}\pr{cl}(f^1(\Sigma),\tau).
\eeq

Прием, связанный с (\ref{13.8}) (подробнее см. в \cite{48}), ниже не используется, но полезен с точки зрения вопросов общей теории расширений.

Заметим, что конкретное построение МП на основе определений (см. (\ref{13.1'}), (\ref{13.2})) возможно лишь в простейших случаях. Общий метод исследования задач такого рода связан с применением аппарата расширений (компактификаций). Возвращаясь к общему случаю (см. (\ref{13.1'})), наметим упомянутый подход, полагая, что так или иначе удалось указать компактное ТП $(\mathbf{K},\mathbf{t}),$ $\mathbf{K}\neq\zer,$ а также отображения $$p: X\rightarrow\mathbf{K}, \ \ q:\mathbf{K}\rightarrow\mathbb{X,}$$ для которых $q\in C(\mathbf{K},\mathbf{t},\mathbb{X},\tau)$ и, кроме того, $\mathbf{f}=q\circ p.$ Условимся каждый набор $(\mathbf{K},\mathbf{t},p,q)$ с упомянутыми свойствами называть \emph{компактификатором}; подчеркнем, что по аналогии с (\ref{13.1'}) (см. также (\ref{13.2})) определено МП
$$(\pr{\mathbf{as}})[X;\mathbf{K};\mathbf{t};p;\mathcal{X}]\in\pp(\mathbf{K}) \ \forall\,\mathcal{X}\in\pp'(\pp(X)).$$

Для наших последующих целей важно следующее известное  \cite{46}-- \cite{50} свойство: если   $(\mathbf{K},\mathbf{t},p,q)$\,--- компактификатор, то
\beq{13.9}
(\pr{\mathbf{as}})[X;\mathbb{X};\tau;\mathbf{f};\mathcal{X}]=q^1((\pr{\mathbf{as}})[X;\mathbf{K};\mathbf{t};p;\mathcal{X}]) \ \forall\,\mathcal{X}\in\pp'(\pp(X)).
\eeq
Итак, искомое МП в (\ref{13.9}) представляется в виде непрерывного образа вспомогательного МП, определяемого в компактном ТП. Данное представление (явно или неявно используемое в теории управления, см. \cite[гл.\,III,IV]{44}, \cite{22}) позволяет в ряде случаев существенно продвинуться в вопросах построения и исследования МП в исходной задаче. Подчеркнем в этой связи следующее важное обстоятельство: согласно (\ref{13.9}) один и тот же компактификатор может использоваться при построении МП, отвечающих различным семействам $\mathcal{X}\in\pp'(\pp(X)).$ Мы, в частности, можем, выбирая  $\mathcal{X}_1\in\pp'(\pp(X))$ и $\mathcal{X}_2\in\pp'(\pp(X)),$  сравнивать МП $(\pr{\mathbf{as}})[X;\mathbb{X};\tau;\mathbf{f};\mathcal{X}_1]$ и
$(\pr{\mathbf{as}})[X;\mathbb{X};\tau;\mathbf{f};\mathcal{X}_2].$ По многим причинам представляют интерес случаи, когда два упомянутых МП совпадают. Такое равенство очевидным образом выполняется, если для некоторого компактификатора $(\mathbf{K},\mathbf{t},p,q)$
\beq{13.10}
(\pr{\mathbf{as}})[X;\mathbf{K};\mathbf{t};p;\mathcal{X}_1]=(\pr{\mathbf{as}})[X;\mathbf{K};\mathbf{t};p;\mathcal{X}_2].
\eeq

Прием, связанный с реализацией  (\ref{13.10}), обсудим далее применительно к задаче раздела 13, имея в виду идею отождествления при некоторых условиях на множества (\ref{12.13}) и (\ref{12.16}).

\section{Множества  притяжения и условия асимптотической нечувствительности при ослаблении части ограничений}\setcounter{equation}{0}\setcounter{proposition}{0}\setcounter{zam}{0}\setcounter{corollary}{0}\setcounter{definition}{0}
\setcounter{condition}{0}\setcounter{theorem}{0}
   \ \ \ \ \ В настоящем разделе будем рассматривать следующую конкретизацию общих определений предыдущего раздела:
   \beq{14.1}
   X=M_*^+[\,b\mid E;\ml;\eta], \ \mathbb{X}=\rr^n, \ \tau=\tau_{\rr}^{(n)}, \ \mathbf{f}=\Pi.
   \eeq
Семейство $\mathcal{X}$ раздела 14 будет конкретизироваться двумя способами:
\beq{14.2}
\mathcal{X}=\mathfrak{Y}, \ \mathcal{X}=\widehat{\mathfrak{Y}}[J].
\eeq

Обстоятельство, связанное с (\ref{14.2}), важно не только с теоретической, но и с практической точки зрения, поскольку на этой основе удается получить условия асимптотической нечувствительности (а, по сути дела, --- грубости) при ослаблении чпсти ограничений. Однако прежде всего отметим, что
\beq{14.3}
\left(\mathfrak{Y}\in\beta[M_*^+[\,b\mid E;\ml;\eta]]\right)\ \& \ \left(\widehat{\mathfrak{Y}}[J]\in\beta[M_*^+[\,b\mid E;\ml;\eta]]\right)
\eeq
(действительно, каждое из семейств в (\ref{14.3}) получается как множество значений изотонных многозначных отображений, определенное каждое на $]0,\infty[$). Поэтому для множества (\ref{12.13}) имеем в силу (\ref{13.2}) равенство
\beq{14.4}
(\pr{\mathbf{as}})\left[M_*^+[\,b\mid E;\ml;\eta];\rr^n;\tau_{\rr}^{(n)};\Pi;\mathfrak{Y}\right]=
\bigcap\limits_{\varepsilon\bn}\pr{cl}\left(\Pi^1\left(S^{\,-1}(\mathbb{O}_N(Y,\varepsilon))\right),\tau_{\rr}^{(n)}\right).
\eeq
Аналогичным образом из (\ref{13.2}) и (\ref{14.3}) вытекает, что для множества (\ref{12.16})
\begin{multline}\label{14.5}
(\pr{\mathbf{as}})\left[M_*^+[\,b\mid E;\ml;\eta];\rr^n;\tau_{\rr}^{(n)};\Pi;\widehat{\mathfrak{Y}}[J]\right]=\\=
\bigcap\limits_{\varepsilon\bn}\pr{cl}\left(\Pi^1\left(S^{\,-1}\left(\widehat{\mathbb{O}}_N(Y,\varepsilon\mid J)\right)\right),\tau_{\rr}^{(n)}\right).
\end{multline}
В  построениях учитывается следующее очевидное свойство: при $\varepsilon\bn$ по определению~$S$
\beq{14.6}
\left(\mathbb{Y}_\varepsilon=S^{-1}(\mathbb{O}_N(Y,\varepsilon))\right)\ \& \ \left(\widehat{\mathbb{Y}}_\varepsilon[J]= S^{-1}\left(\widehat{\mathbb{O}}_N(Y,\varepsilon\mid J)\right)\right).
\eeq

В этой связи отметим, также, что
\begin{multline}\label{14.7}
(\pr{\mathbf{as}})\left[M_*^+[\,b\mid E;\ml;\eta];\rr^n;\tau_{\rr}^{(n)};\Pi;\mathfrak{Y}\right]=
\bigcap\limits_{\Sigma\in\mathfrak{Y}}\pr{cl}\left(\Pi^1\left(\Sigma,\tau_{\rr}^{(n)}\right)\right)=\\=
\bigcap\limits_{\varepsilon\bn}\pr{cl}\left(\Pi^1\left(\mathbb{Y}_\varepsilon\right),\tau_{\rr}^{(n)}\right).
\end{multline}

Из (\ref{14.6}) и (\ref{14.7}) легко следует (\ref{14.4}). Аналогичным образом имеем, наконец, что (см. (\ref{14.5}), (\ref{14.6}))
\begin{multline}\label{14.8}
(\pr{\mathbf{as}})\left[ M_*^+[\,b\mid E;\ml;\eta];\rr^n;\tau_{\rr}^{(n)};\Pi;\widehat{\mathfrak{Y}}[J]\right]=\!
\bigcap\limits_{\varepsilon\bn}\pr{cl}\left(\Pi^1\left(\widehat{\mathbb{Y}}_\varepsilon[J]\right),\tau_{\rr}^{(n)}\right)=\\=
\bigcap\limits_{\Sigma\in\widehat{\mathfrak{Y}}[J]}\pr{cl}\left(\Pi^1\left(\Sigma\right),\tau_{\rr}^{(n)}\right).
\end{multline}
Напомним, что согласно (\ref{12.18}), (\ref{14.4}) и (\ref{14.5})
\beq{14.9}
(\pr{\mathbf{as}})\left[ M_*^+[\,b\mid E;\ml;\eta];\rr^n;\tau_{\rr}^{(n)};\Pi;\widehat{\mathfrak{Y}}[J]\right]\subset(\pr{\mathbf{as}})\left[ M_*^+[\,b\mid E;\ml;\eta];\rr^n;\tau_{\rr}^{(n)};\Pi;\mathfrak{Y}\right].
\eeq

В дальнейшем нас будут интересовать условия, при которых (\ref{14.9}) превращается в равенство. Однако предварительно отметим еще несколько полезных свойств.

Прежде всего заметим, что для рассматриваемого конкретного случая выполняются условия, подобные (\ref{13.5}), (\ref{13.6}), а именно:\\
$1') \ \ \ \exists\,(\Sigma_i)_{i\in\mathbb{N}}\in\mathfrak{Y}^\mathbb{N} \ \forall\,\Sigma\in\mathfrak{Y}$ \ \ \ $\{k\in\mathbb{N}\mid\Sigma_k\subset\Sigma\}\neq\zer;$ \\
$2') \  \exists\,(\Sigma_i)_{i\in\mathbb{N}}\in\widehat{\mathfrak{Y}}[J]^\mathbb{N} \ \forall\,\Sigma\in\widehat{\mathfrak{Y}}[J]$ \ \ \ $\{k\in\mathbb{N}\mid\Sigma_k\subset\Sigma\}\neq\zer;$ \\
$3') \ \ \ \forall\,z\in\rr^n \ \exists\,(H_i)_{i\in\mathbb{N}}\in N_{\tau_{\rr}^{(n)}}(z)^\mathbb{N} \ \forall\,\mathbf{H}\in N_{\tau_{\rr}^{(n)}}(z)$ \ \ \ $\{k\in\mathbb{N}\mid H_k\subset\mathbb{H}\}\neq\zer.$
Читателю предлагается проверить положения $1')$--$3')$ самостоятельно, используя определения раздела 11 и метризуемость топологии $\tau_{\rr}^{(n)}.$
Теперь с учетом  (\ref{13.7}) получаем, что
\beq{14.10}
(\pr{\mathbf{as}})\left[ M_*^+[\,b\mid E;\ml;\eta];\rr^n;\tau_{\rr}^{(n)};\Pi;\mathfrak{Y}\right]=(\pr{\mathbf{sas}})\left[ M_*^+[\,b\mid E;\ml;\eta];\rr^n;\tau_{\rr}^{(n)};\Pi;\mathfrak{Y}\right],
\eeq
\beq{14.11}
(\pr{\mathbf{as}})\!\left[ M_*^+[\,b|E;\ml;\eta];\rr^n;\tau_{\rr}^{(n)};\!\Pi;\widehat{\mathfrak{Y}}[J]\right]\!=\!(\pr{\mathbf{sas}})\!\left[ M_*^+[\,b| E;\ml;\eta];\rr^n;\tau_{\rr}^{(n)};\!\Pi;\widehat{\mathfrak{Y}}[J]\right].
\eeq

Итак, в виде (\ref{12.13}) и (\ref{12.16}) имеем (см. (\ref{14.4}), (\ref{14.10}) и (\ref{14.5}), (\ref{14.11})) секвенциальные МП (см. (\ref{13.3})). Читателю предоставляется аккуратное оформление высказываний, связанных с (\ref{14.10}), (\ref{14.11}) в терминах конкретных определений раздела 11 (при этом следует учитывать (\ref{13.3})).

Ориентируясь на использование (\ref{13.9}), рассмотрим один из возможных вариантов компактификатора. В этой связи напомним некоторые положения разделов 9, 10. В частности, далее через $I$ будет обозначаться отображение
$$f\mapsto f*\eta: \ M_*^+[\,b\mid E;\ml;\eta]\rightarrow \widetilde{M}_*^+[\,b\mid E;\ml;\eta].$$
Ясно, что $I=\mathcal{I}.$ С учетом предложения 9.1 получаем также, что
\beq{14.12}
I=(f*\eta)_{f\in M_*^+[\,b\mid E;\ml;\eta]}\in\Xi_+^*[\,b\mid E;\ml;\eta]^{M_*^+[\,b\mid E;\ml;\eta]}
\eeq
(данный вариант записи более удобен для наших целей). Кроме того, введем в рассмотрение следующий оператор $\widetilde{\Pi}:$ по определению $\widetilde{\Pi}$ есть отображение
\beq{14.13}
\mu\mapsto\biggl(\,\int\limits_E \pi_i\,d\mu\biggl)_{i\in\overline{1,n}}:\Xi_+^*[\,b\mid E;\ml;\eta]\rightarrow\rr^n.
\eeq
Заметим, что $\widetilde{\Pi}(\mu)=\left(\,\int\limits_E \pi_i\,d\mu\right)_{i\in\overline{1,n}}\in\rr^n;$  в частности, имеем при $f\in M_*^+[\,b\mid E;\ml;\eta]$ с учетом  (\ref{10.5}) и (\ref{14.12}) цепочку равенств
$$\widetilde{\Pi}(f*\eta)=\left(\widetilde{\Pi}\circ I\right)(f)=\biggl(\,\int\limits_E \pi_i\,d(f*\eta)\biggl)_{i\in\overline{1,n}}=\biggl(\,\int\limits_E \pi_i f\,d\eta\biggl)_{i\in\overline{1,n}}=\Pi(f).$$
 Это означает справедливость равенства
 \beq{14.14}
 \Pi=\widetilde{\Pi}\circ I.
 \eeq
Возвращаясь к (\ref{14.13}), заметим, что
\beq{14.15}
\widetilde{\Pi}\in C\left(\Xi_+^*[\,b\mid E;\ml;\eta],\tau_\eta^*[\ml|b],\rr^n,\tau_{\rr}^{(n)}\right).
\eeq

В самом деле, $\widetilde{\Pi}$  есть сужение на  $\Xi_+^*[\,b\mid E;\ml;\eta]$ отображения
  \beq{14.16}\mu\mapsto\biggl(\,\int\limits_E \pi_i\,d\mu\biggl)_{i\in\overline{1,n}}: \mathbb{A}(\ml)\rightarrow\rr^n,\eeq
  которое обозначим сейчас через $\widehat{\Pi}.$ Итак,  $\widehat{\Pi}$ есть оператор (\ref{14.16}), а $\widetilde{\Pi}=\left(\widehat{\Pi}\mid\Xi_+^*[\,b\mid E;\ml;\eta]\right).$
При этом
\beq{14.17}
\widehat{\Pi}\in C\left(\mathbb{A}(\ml),\tau_*(\ml),\rr^n,\tau_{\rr}^{(n)}\right),
\eeq
что следует из самого определения $*$-слабой топологии. Однако в целях полноты изложения мы все же проверим (\ref{14.7}). Фиксируем
\beq{14.18}
\mathbb{G}\in \tau_{\rr}^{(n)}.
\eeq
Тогда $\widehat{\Pi}^{-1}(\mathbb{G})=\left\{\widetilde{\mu}\in\mathbb{A}(\ml)\mid\widehat{\Pi}(\widetilde{\mu})\in\mathbb{G}\right\}\in\pp(\mathbb{A}(\ml)).$ Выберем произвольно
\beq{14.19}
\zeta\in\widehat{\Pi}^{-1}(\mathbb{G}).\eeq Получаем, что  $\zeta\in\mathbb{A}(\ml)$ и, кроме того, $$\widehat{\Pi}(\zeta)=\biggl(\,\int\limits_E \pi_i\,d\zeta\biggl)_{i\in\overline{1,n}}\in\mathbb{G}.$$
С учетом (\ref{14.18}) имеем теперь, что для некоторого $\xi\bn$
\beq{14.20}
Z_\xi\triangleq\left\{z\in\rr^n\mid\biggl|\int\limits_E \pi_j\,d\zeta-z(j)\biggl|<\xi \ \ \ \forall\,j\in\overline{1,n}\right\}\subset\mathbb{G}.
\eeq
Заметим, что $\mathbb{K}\triangleq\{\pi_i:i\in\overline{1,n}\,\}\in\pr{Fin}(B(E,\ml)),$ а потому (см. (\ref{9.5}), (\ref{9.8'}))
\begin{multline}\label{14.21}
N_{\ml}^*(\zeta,\mathbb{K},\xi)=\left\{\nu\in\mathbb{A}(\ml)\mid \biggl|\int\limits_E h\,d\zeta-\int\limits_E h\,d\nu\biggl|<\xi \ \ \forall\,h\in\mathbb{K}\right\}=\\=\left\{\nu\in\mathbb{A}(\ml)\mid \biggl|\int\limits_E \pi_j\,d\zeta-\int\limits_E \pi_j\,d\nu\biggl|<\xi  \ \ \forall\,j\in\overline{1,n}\right\}\in N_{\tau_*(\ml)}^0(\zeta).
\end{multline}

Пусть выбрана произвольно к.-а. мера $\widetilde{\zeta}\in N_{\ml}^*(\zeta,\mathbb{K},\xi).$ Тогда $\widetilde{\zeta}\in\mathbb{A}(\ml)$ и при этом справедлива система неравенств
\beq{14.22}
\biggl|\int\limits_E \pi_j\,d\,\zeta-\int\limits_E \pi_j\,d\,\widetilde{\zeta}\biggl|<\xi  \ \ \forall\,j\in\overline{1,n}.
\eeq
Как следствие имеем с учетом (\ref{14.20}) свойство
$$\widehat{\Pi}\left(\widetilde{\zeta}\right)=\biggl(\int\limits_E \pi_i\,d\,\widetilde{\zeta}\biggl)_{i\in\overline{1,n}}\in Z_\xi,$$
а, стало быть, $\widehat{\Pi}\left(\widetilde{\zeta}\right)\in\mathbb{G}.$
В итоге $$\widetilde{\zeta}\in\widehat{\Pi}^{-1}(\mathbb{G}).$$
Поскольку выбор $\widetilde{\zeta}$ был произвольным, установлено вложение
\beq{14.23}
N_{\ml}^*(\zeta,\mathbb{K},\xi)\subset\widehat{\Pi}^{-1}(\mathbb{G}).
\eeq
Коль скоро выбор  $\zeta$ (\ref{14.19}) был произвольным,  установлено, что множество  $\widehat{\Pi}^{-1}(\mathbb{G})\!\!\in\!\pp(\mathbb{A}(\ml))$  обладает свойством
$\forall\,\widetilde{\mu}\in\widehat{\Pi}^{-1}(\mathbb{G}) \ \exists\,K\in\pr{Fin}(B(E,\ml)) \ \exists\,\varepsilon\bn:$ $$N_{\ml}^*\left(\widetilde{\mu},K,\varepsilon\right)\subset\widehat{\Pi}^{-1}(\mathbb{G}).$$
В силу (\ref{9.6}) получаем важное свойство
\beq{14.24}
\widehat{\Pi}^{-1}(\mathbb{G})\in\tau_*(\ml).
\eeq

Итак (см. (\ref{14.18}), (\ref{14.24})), имеет место
\beq{14.25}
\widehat{\Pi}^{-1}(G)\in\tau_*(\ml) \ \forall\,G\in\tr^{(n)}.
\eeq
Последнее означает справедливость требуемого свойства (\ref{14.17}): $\widehat{\Pi}$  есть непрерывная вектор-функция на $ \mathbb{A}(\ml).$ Напомним теперь, что
\beq{14.26}
\widetilde{\Pi}: \ \Xi_+^*[\,b\mid E;\ml;\eta]\rightarrow\rr^n
\eeq
  и при этом $\widetilde{\Pi}(\mu)=\widehat{\Pi}(\mu) \ \forall\,\mu\in\Xi_+^*[\,b\mid E;\ml;\eta].$
Тогда имеем при $G\in\tr^{(n)}$
\begin{multline*}\widetilde{\Pi}^{-1}(G)=\left\{\mu\in\Xi_+^*[\,b\mid E;\ml;\eta]\mid\widehat{\Pi}(\mu)\in G\right\}=\\=\widehat{\Pi}^{-1}(G)\cap\Xi_+^*[\,b\mid E;\ml;\eta]\in\tau_\eta^*[\ml|b],\end{multline*}
где учитывается (\ref{14.25}) и определение топологии $\tau_\eta^*[\ml|b].$ Следовательно,
$$\widetilde{\Pi}^{-1}(G)\in\tau_\eta^*[\ml|b] \ \  \forall\,G\in\tr^{(n)}.$$
С учетом (\ref{14.26}) получаем теперь справедливость свойства (\ref{14.15}).

Таким образом, мы располагаем набором
\beq{14.27}
\left(\Xi_+^*[\,b\mid E;\ml;\eta],\tau_\eta^*[\ml|b],I,\widetilde{\Pi}\right),
\eeq
для которого (см. (\ref{14.14}), (\ref{14.15})) справедливы следующие положения:
\beq{14.28}
\left(\widetilde{\Pi}\in C\left(\Xi_+^*[\,b\mid E;\ml;\eta],\tau_\eta^*[\ml|b],\rr^n,\tau_{\rr}^{(n)}\right)\right)\ \& \left(\Pi=\widetilde{\Pi}\circ I\right).
\eeq

Из (\ref{14.28}) вытекает, что (\ref{14.27}) есть компактификатор (здесь имеется в виду реализация (\ref{14.1}) для абстрактной задачи о достижимости раздела 14). Поэтому согласно (\ref{13.9}) имеем
\begin{multline}\label{14.29}
(\pr{\mathbf{as}})\left[M_*^+[\,b\mid E;\ml;\eta];\rr^n;\tau_{\rr}^{(n)};\Pi;\mathcal{X}\right]=\\=\widetilde{\Pi}^1\left((\pr{\mathbf{as}})\left[M_*^+[\,b\mid E;\ml;\eta];\Xi_+^*[\,b\mid E;\ml;\eta];\tau_\eta^*[\ml|b];I;\mathcal{X}\right]\right)  \\ \forall\,\mathcal{X}\in\pp'\left(\pp\left(M_*^+[\,b\mid E;\ml;\eta]\right)\right).
\end{multline}
Представление (\ref{14.29}) вполне применимо в случаях
$$\mathcal{X}=\mathfrak{Y}, \ \ \mathcal{X}=\widehat{\mathfrak{Y}}[J].$$

Таким образом,  имеем из (\ref{14.11}), (\ref{14.29}) следующие две цепочки равенств
\begin{multline}\label{14.30}
(\pr{\mathbf{as}})\!\left[M_*^+[\,b\mid E;\ml;\eta];\rr^n;\!\tau_{\rr}^{(n)};\!\Pi;\mathfrak{Y}\right]\!=\!(\pr{\mathbf{sas}})\!\left[M_*^+[\,b\mid E;\ml;\eta];\rr^n;\!\tau_{\rr}^{(n)};\!\Pi;\mathfrak{Y}\right]=\\=\widetilde{\Pi}^1\left((\pr{\mathbf{as}})\left[M_*^+[\,b\mid E;\ml;\eta];\Xi_+^*[\,b\mid E;\ml;\eta];\tau_\eta^*[\ml|b];I;\mathfrak{Y}\right]\right),
\end{multline}
\begin{multline}\label{14.31}
(\pr{\mathbf{as}})\left[M_*^+[\,b\mid E;\ml;\eta];\rr^n;\tau_{\rr}^{(n)};\Pi;\widehat{\mathfrak{Y}}[J]\right]=\\=(\pr{\mathbf{sas}})\left[M_*^+[\,b\mid E;\ml;\eta];\rr^n;\tau_{\rr}^{(n)};\Pi;\widehat{\mathfrak{Y}}[J]\right]=\\=\widetilde{\Pi}^1\left((\pr{\mathbf{as}})\left[M_*^+[\,b\mid E;\ml;\eta];\Xi_+^*[\,b\mid E;\ml;\eta];\tau_\eta^*[\ml|b];I;\widehat{\mathfrak{Y}}[J]\right]\right).
\end{multline}
Следовательно, (см. (\ref{14.30}), (\ref{14.31})), для построения двух МП в основной задаче имеет смысл определить вспомогательные МП
\beq{14.32}
\mathfrak{M}\triangleq(\pr{\mathbf{as}})\left[M_*^+[\,b\mid E;\ml;\eta];\Xi_+^*[\,b\mid E;\ml;\eta];\tau_\eta^*[\ml|b];I;\mathfrak{Y}\right]\in\pp\left(\Xi_+^*[\,b\mid E;\ml;\eta]\right),
\eeq
\begin{multline}\label{14.33}
\widehat{\mathfrak{M}}_J\triangleq(\pr{\mathbf{as}})\left[M_*^+[\,b\mid E;\ml;\eta];\Xi_+^*[\,b\mid E;\ml;\eta];\tau_\eta^*[\ml|b];I;\widehat{\mathfrak{Y}}[J]\right]\in \\ \in\pp\left(\Xi_+^*[\,b\mid E;\ml;\eta]\right).
\end{multline}
Учитывая  (\ref{13.2}) и (\ref{14.3}),  получаем для МП (\ref{14.32}), (\ref{14.33}) следующие представления
\beq{14.33'}
\mathfrak{M}=\bigcap\limits_{\varepsilon\bn}\pr{cl}\left(I^1\left(\mathbb{Y}_\varepsilon\right),\tau_{\eta}^{*}[\ml|b]\right),\ \
\widehat{\mathfrak{M}}_J=\bigcap\limits_{\varepsilon\bn}\pr{cl}\left(I^1\left(\widehat{\mathbb{Y}}_\varepsilon[J]\right),\tau_{\eta}^{*}[\ml|b]\right).
\eeq
С учетом замкнутости $\Xi_+^*[\,b\mid E;\ml;\eta]$  в топологии $\tau_*(\ml),$ что уже отмечалось ранее, получаем при $\varepsilon\bn,$ что
\begin{multline*}
\pr{cl}\left(I^1\left(\mathbb{Y}_\varepsilon\right),\tau_{\eta}^{*}[\ml|b]\right)=
\pr{cl}\left(I^1\left(\mathbb{Y}_\varepsilon\right),\tau_*(\ml)\right)\cap \Xi_+^*[\,b\mid E;\ml;\eta]=\pr{cl}\left(I^1\left(\mathbb{Y}_\varepsilon\right),\tau_*(\ml)\right),\\
\pr{cl}\left(I^1\left(\widehat{\mathbb{Y}}_\varepsilon[J]\right),\tau_{\eta}^{*}[\ml|b]\right)=
\pr{cl}\left(I^1\left(\widehat{\mathbb{Y}}_\varepsilon[J]\right),\tau_*(\ml)\right)\cap \Xi_+^*[\,b\mid E;\ml;\eta]=\\=\pr{cl}\left(I^1\left(\widehat{\mathbb{Y}}_\varepsilon[J]\right),\tau_*(\ml)\right).
\end{multline*}
Заметим, кстати, что согласно (\ref{12.11}) $$\mathbb{Y}_\varepsilon=S^{-1}(\mathbb{O}_N(Y,\varepsilon)) \ \forall\,\varepsilon\bn.$$
Аналогичным образом из (\ref{12.15}) извлекается свойство
$$\widehat{\mathbb{Y}}_\varepsilon[J]=S^{-1}(\widehat{\mathbb{O}}_N(Y,\varepsilon|J)) \ \forall\,\varepsilon\bn.$$

В двух последних представлениях учтено определение $S$ (см. (\ref{12.6})). Таким образом,
\beq{14.34}
\mathfrak{M}=\bigcap\limits_{\varepsilon\bn}\pr{cl}\left(I^1\left(S^{-1}(\mathbb{O}_N(Y,\varepsilon))\right),\tau_{*}(\ml)\right)=
\bigcap\limits_{\varepsilon\bn}\pr{cl}\left(I^1\left(\mathbb{Y}_\varepsilon\right),\tau_{*}(\ml)\right),
\eeq
\beq{14.35}
\widehat{\mathfrak{M}}_J\!=\!\bigcap\limits_{\varepsilon\bn}\!\pr{cl}\left(I^1\left(S^{-1}\left(\widehat{\mathbb{O}}_N(Y,\varepsilon|J)\right)\right),
\tau_*(\ml)\right)\!=\!\bigcap\limits_{\varepsilon\bn}\!\pr{cl}\left(I^1\left(\widehat{\mathbb{Y}}_\varepsilon[J]\right),\tau_*(\ml)\right).
\eeq
С учетом (\ref{12.17}), (\ref{14.34}) и (\ref{14.35}) получаем, что всегда
\beq{14.36}
\widehat{\mathfrak{M}}_J\subset\mathfrak{M}.
\eeq
Отметим еще одно полезное обстоятельство, полагая, что $\widetilde{S}$  есть $\pr{def}$ отображение
\beq{14.37}
\mu\mapsto\biggl(\,\int\limits_E s_j\,d\,\mu\biggl)_{j\in\overline{1,n}}:\Xi_+^*[\,b\mid E;\ml;\eta]\rightarrow\rr^N.
\eeq
Тогда
\beq{14.38}
\widetilde{S}\in C\left(\Xi_+^*[\,b\mid E;\ml;\eta],\tau_\eta^{*}[\ml|b],\rr^N,\tau_{\rr}^{(n)}\right).
\eeq

Проверка (\ref{14.38}) аналогична обоснованию свойства (\ref{14.15}).
Всюду в дальнейшем предполагается выполненным следующее
\begin{condition} Функции $s_j, \ j\in J,$ являются  ступенчатыми:
$s_j\in B_0(E,\ml)$ $\forall\,j\in J.$
\end{condition}

Итак, мы рассматриваем далее случай, когда в/з функции $s_j, \ j\in J,$ являются ступенчатыми относительно ИП $(E,\ml).$ Кроме того, постулируем в дальнейшем следующее
 \begin{condition}
Множество $Y$ замкнуто в  $\left(\rr^N,\tau_{\rr}^{(n)}\right).$
\end{condition}
\begin{proposition}
Справедлива следующая цепочка равенств
\beq{14.39}
\widetilde{S}^{\,-1}(Y)=\mathfrak{M}=\widehat{\mathfrak{M}}_J.
\eeq
\end{proposition}
Д о к а з а т е л ь с т в о. Учитываем (\ref{14.36}). Пусть $\mu_0\in\widetilde{S}^{\,-1}(Y).$  Тогда
\beq{14.40}
\mu_0\in\Xi_+^*[\,b\mid E;\ml;\eta]
\eeq
 и при этом
 \beq{14.41}
 \widetilde{S}(\mu_0)=\biggl(\int\limits_E s_j\,d\,\mu_0\biggl)_{j\in\overline{1,N}}\in Y.
 \eeq

Покажем, что $\mu_0\in\widehat{\mathfrak{M}}_J.$ Для этого отметим, что согласно (\ref{14.33}) $\widehat{\mathfrak{M}}_J$ есть (см. (\ref{13.1'})) множество всех таких $\mu\in\Xi_+^*[\,b\mid E;\ml;\eta],$ что для некоторой направленности
$(D,\preceq,g)$ в $M_*^+[\,b\mid E;\ml;\eta]$
$$\left(\forall\,\Sigma\in\widehat{\mathfrak{Y}}[J] \ \exists\,d\in D \ \forall\,\delta\in D \ (d\preceq\delta)\!\Rightarrow\!(g(\delta)\in\Sigma)\right) \&  \left((D,\preceq,I\circ g)\stackrel{\tau_\eta^*[\ml|b]}\longrightarrow \!\mu\right).$$

С учетом определения $\widehat{\mathfrak{Y}}[J]$ получаем, что $\widehat{\mathfrak{M}}_J$ является множеством всех к.-а. мер $\mu\in\Xi_+^*[\,b\mid E;\ml;\eta],$ для каждой из которых существует направленность $(D,\preceq,g)$ в $M_*^+[\,b\mid E;\ml;\eta]$ со свойствами
\begin{multline}\label{14.42}
\left(\forall\,\varepsilon\bn \ \exists\,d\in D \ \forall\,\delta\in D \ (d\preceq\delta)\Rightarrow\left(g(\delta)\in\widehat{\mathbb{Y}}_\varepsilon[J]\right)\right)\ \& \\ \& \
\left((D,\preceq,I\circ g)\stackrel{\tau_\eta^*[\ml|b]}\longrightarrow \mu\right).
\end{multline}
При этом согласно (\ref{14.12}) для произвольной направленности $(D,\preceq,g)$ в $M_*^+[\,b\mid E;\ml;\eta]$ оператор $I\circ g$ направленности-композиции $(D,\preceq,I\circ g)$ имеет вид
$$I\circ g=(g(\delta)*\eta)_{\delta\in D}\in\Xi_+^*[\,b\mid E;\ml;\eta]^D.$$

Рассмотрим теперь введенную ранее направленность $\left(\mathbf{D}(E,\ml),\prec, \Theta_{\mu_0}^+[\cdot]\right)$ в множестве $M_*^+[\,b\mid E;\ml;\eta].$
Тогда из (\ref{11.6}) имеем, в частности, что $$I\circ\Theta_{\mu_0}^+[\cdot]=\Theta_{\mu_0}^+[\cdot]*\eta=\left(\Theta_{\mu_0}^+[\mathcal{K}]*\eta\right)_{\mathcal{K}\in\mathbf{D}(E,\ml)}\in
\widetilde{M}_*^+[\,b\mid E;\ml;\eta]^{\mathbf{D}(E,\ml)}.$$
С учетом предложения 11.1 получаем
$$I\circ\Theta_{\mu_0}^+[\cdot]=\Theta_{\mu_0}^+[\cdot]*\eta\in\Xi_+^*[\,b\mid E;\ml;\eta]^{\mathbf{D}(E,\ml)}.$$

Таким образом,  (\ref{11.17}) есть направленность в $\Xi_+^*[\,b\mid E;\ml;\eta];$ при этом согласно предложению 12.2 получаем свойство сходимости
\beq{14.43}
\left(\mathbf{D}(E,\ml),\prec,I\circ\Theta_{\mu_0}^+[\cdot]\right)\stackrel{\tau_0(\ml)}\longrightarrow\mu_0.
\eeq
Из условия 15.1 следует очевидное положение:  $$\widetilde{\mathbb{K}}\triangleq\{s_j:j\in J\}, \ \ \widetilde{\mathbb{K}}\subset B_0(E,\ml),$$
есть конечное множество. При этом
\beq{14.44}
(J\neq\zer)\Rightarrow\left(\widetilde{\mathbb{K}}\in\pr{Fin}(B_0(E,\ml))\right).
\eeq
В силу следствия 12.2, учитывая  определение $\widetilde{\mathbb{K}}$ и (\ref{14.44}), имеем импликацию
\begin{multline*}(J\neq\zer)\Rightarrow\biggl(\exists\,\mathcal{K}_1\in\mathbf{D}(E,\ml) \ \forall\,\mathcal{K}_2\in\mathbf{D}(E,\ml) \ \ (\mathcal{K}_1\prec\mathcal{K}_2)\Rightarrow\\ \Rightarrow\bigl(\int\limits_E s_j\Theta_{\mu_0}^+[\mathcal{K}_2]\,d\eta=\int\limits_E s_j\,d\mu_0 \ \forall\,j\in J\bigl)\biggl).\end{multline*}
 Если же $J=\zer,$ то $\widetilde{\mathbb{K}}=\zer$ и, поскольку $\mathbf{D}(E,\ml)\neq\zer,$ то непременно $\exists\,\mathcal{K}_1\in\mathbf{D}(E,\ml) \ \forall\,\mathcal{K}_2\in\mathbf{D}(E,\ml)$
\beq{14.45}
(\mathcal{K}_1\prec\mathcal{K}_2)\Rightarrow\biggl(\,\int\limits_E s_j\Theta_{\mu_0}^+[\mathcal{K}_2]\,d\eta=\int\limits_E s_j\,d\mu_0 \ \forall\,j\in J\biggl);
\eeq
в рассматриваемом случае в качестве $\mathcal{K}_1$ в (\ref{14.45}) достаточно выбрать любое разбиение из $\mathbf{D}(E,\ml).$
Таким образом, во всех возможных случаях
\begin{multline}\label{14.46}
\exists\,\mathcal{K}_1\in\mathbf{D}(E,\ml) \ \forall\,\mathcal{K}_2\in\mathbf{D}(E,\ml) \\ (\mathcal{K}_1\prec\mathcal{K}_2)\Rightarrow\biggl(\,\int\limits_E s_j\Theta_{\mu_0}^+[\mathcal{K}_2]\,d\eta=\int\limits_E s_j\,d\mu_0 \ \forall\,j\in J\biggl).
\end{multline}

С другой стороны, имеем свойство, определяемое в предложении 12.4:
$$\left(\mathbf{D}(E,\ml),\prec,\Theta_{\mu_0}^+[\cdot]*\eta\right)\stackrel{\tau_*(\ml)}\longrightarrow \mu_0.$$
Из (\ref{4.14}) получаем, следовательно, что $\forall\,H\in N_{\tau_*(\ml)}(\mu_0) \ \exists\,\mathcal{K}'\in\mathbf{D}(E,\ml) \ \forall\,\mathcal{K}''\in\mathbf{D}(E,\ml)$
\beq{14.47}
\left(\mathcal{K}'\prec\mathcal{K}''\right)\Rightarrow\left(\Theta_{\mu_0}^+[\mathcal{K}'']*\eta\in H\right).
\eeq

Из (\ref{11.52}) вытекает (см. (\ref{4.6})), что $\forall\,\widetilde{H}\in N_{\tau_\eta^*[\ml|b]}(\mu_0) \ \exists\,\mathcal{K}'\in\mathbf{D}(E,\ml) \ \forall\,\mathcal{K}''\in\mathbf{D}(E,\ml)$
$$\left(\mathcal{K}'\prec\mathcal{K}''\right)\Rightarrow\left(\Theta_{\mu_0}^+[\mathcal{K}'']*\eta\in \widetilde{H}\right)$$
(учли здесь то, что $\Theta_{\mu_0}^+[\mathcal{K}]*\eta\in \Xi_+^*[\,b\mid E;\ml;\eta]$ при $\mathcal{K}\in\mathbf{D}(E,\ml),$ как уже отмечалось ранее). Это означает согласно (\ref{4.14}), что
$$\left(\mathbf{D}(E,\ml),\prec,I\circ\Theta_{\mu_0}^+[\cdot]\right)\stackrel{\tau_\eta^*[\ml|b]}\longrightarrow\mu_0. $$
С учетом непрерывности $\widetilde{S}$ (см. (\ref{4.100})) получаем свойство сходимости
\beq{14.48}
\left(\mathbf{D}(E,\ml),\prec,\widetilde{S}\circ I\circ\Theta_{\mu_0}^+[\cdot]\right)\stackrel{\tau_\rr^{(N)}}\longrightarrow\widetilde{S}(\mu_0)
\eeq
(мы использовали \cite[(2.5.4)]{41}, но, в данном случае (\ref{14.48}) можно проверить проще, используя определение $*$-слабой топологии: предлагаем читателю установить (\ref{14.48}) самостоятельно). Однако из (\ref{10.5}), (\ref{12.6}), (\ref{14.37}) вытекает, что
\begin{multline*}\left(\widetilde{S}\circ I\right)(f)=\widetilde{S}(I(f))=\widetilde{S}(f*\eta)=\\=\biggl(\int\limits_E s_j\,d(f*\eta)\biggl)_{j\in\overline{1,N}}=
\biggl(\int\limits_E s_j f\,d\eta\biggl)_{j\in\overline{1,N}}=S(f) \ \ \ \forall\,f\in M_*^+[\,b\mid E;\ml;\eta].\end{multline*}
Иными словами, $\widetilde{S}\circ I=S,$ а потому в силу (\ref{14.48})
 \beq{14.48'}
\left(\mathbf{D}(E,\ml),\prec,S\circ\Theta_{\mu_0}^+[\cdot]\right)\stackrel{\tau_\rr^{(N)}}\longrightarrow\widetilde{S}(\mu_0)
\eeq
 (учли хорошо известное свойство ассоциативности композиций; читателю предлагается извлечь (\ref{14.48'}) из (\ref{14.48}) самостоятельно). Рассмотрим теперь комбинацию (\ref{14.41}), (\ref{14.46}) и (\ref{14.48}). Итак, выберем произвольно $\widehat{\varepsilon}\bn$  и рассмотрим множество
 \beq{14.49}
 \widehat{\mathbb{Y}}_{\widehat{\varepsilon}}[J]=\left\{f\in M_*^+[\,b\mid E;\ml;\eta]\mid S(f)\in\widehat{\mathbb{O}}_N(Y,\widehat{\varepsilon}\mid J)\right\}= S^{-1}\left(\widehat{\mathbb{O}}_N(Y,\widehat{\varepsilon}\mid J)\right)
 \eeq
 (см. (\ref{12.15})), где согласно (\ref{12.14}) имеет место равенство
 \begin{multline}\label{14.50}
 \widehat{\mathbb{O}}_N(Y,\widehat{\varepsilon}\,| J)=\bigl\{(z_i)_{i\in\overline{1,n}}\in\rr^N|\exists\,(y_i)_{i\in\overline{1,N}}\in Y: (y_j=z_j \ \forall\,j\in J)\ \& \\ \& \ (|y_j-z_j|<\widehat{\varepsilon} \ \forall\,j\in\overline{1,N}\setminus J)\bigl\}.
 \end{multline}
  С учетом (\ref{14.46}) подберем $\mathcal{K}_*\in\mathbf{D}(E,\ml)$ из условия $\forall\,\mathcal{K}\in\mathbf{D}(E,\ml)$
 $$(\mathcal{K}_*\prec\mathcal{K})\Rightarrow\biggl(\,\int\limits_E s_j\Theta_{\mu_0}^+[\mathcal{K}]\,d\eta=\int\limits_E s_j\,d\mu_0 \ \forall\,j\in J\biggl).$$

 Кроме того, используя (\ref{14.48'}) и представление топологии $\tau_{\rr}^{(N)}$ в терминах нормы $\|\cdot\|_N,$ подберем $\mathcal{K}^*\in\mathbf{D}(E,\ml)$ так, что $\forall\,\mathcal{K}\in\mathbf{D}(E,\ml)$
 $$\left(\mathcal{K}^*\prec\mathcal{K}\right)\Rightarrow\left(\,\left\|\left(S\circ\Theta_{\mu_0}^+[\cdot]\right)(\mathcal{K})-\widetilde{S}
 (\mu_0)\right\|_N<\widehat{\varepsilon}\,\right).$$
 Тогда $\mathcal{K}_*\in\mathbf{D}(E,\ml)$ и $\mathcal{K}^*\in\mathbf{D}(E,\ml).$ Поскольку $\prec$ есть направление на $\mathbf{D}(E,\ml),$ то для некоторого $\widetilde{\mathfrak{K}}\in \mathbf{D}(E,\ml)$ имеют место свойства
 \beq{14.51}
 \left(\mathcal{K}_*\prec\widetilde{\mathfrak{K}}\right)\ \& \ \left(\mathcal{K}^*\prec\widetilde{\mathfrak{K}}\right).
 \eeq
 Выберем произвольно разбиение $\mathfrak{K}\in \mathbf{D}(E,\ml),$ для которого
 \beq{14.52}
 \widetilde{\mathfrak{K}}\prec\mathfrak{K}.
 \eeq

 Тогда в силу (\ref{14.51}) и (\ref{14.52})
 \beq{14.53}
 \left(\mathcal{K}_*\prec\mathfrak{K}\right)\ \& \ \left(\mathcal{K}^*\prec\mathfrak{K}\right).
 \eeq
  При этом  $\Theta_{\mu_0}^+[\mathfrak{K}]\in M_*^+[\,b\mid E;\ml;\eta]$ и, в частности, $\Theta_{\mu_0}^+[\mathfrak{K}]\!\!\in\! B_0^+[E;\ml]$ Из (\ref{14.53}) имеем, в частности, по выбору $\mathcal{K}_*$ систему равенств
\beq{14.54}
\int\limits_E s_j\Theta_{\mu_0}^+[\mathfrak{K}]\,d\eta=\int\limits_E s_j\,d\mu_0 \ \forall\,j\in J.
\eeq
С другой стороны, по выбору $\mathcal{K}^*$ получаем  из (\ref{14.53}) неравенство
\beq{14.55}
\left\|\left(S\circ\Theta_{\mu_0}^+[\cdot]\right)(\mathfrak{K})-\widetilde{S}(\mu_0)\right\|_N<\widehat{\varepsilon}.
\eeq
С учетом определения  $S$ (см. (\ref{12.6})) реализуется цепочка равенств
$$\left(S\circ\Theta_{\mu_0}^+[\cdot]\right)(\mathfrak{K})=S\left(\Theta_{\mu_0}^+[\mathfrak{K}]\right)=\biggl(\,\int\limits_E s_j\Theta_{\mu_0}^+[\mathfrak{K}]\,d\eta\biggl)_{j\in\overline{1,N}}.$$

В то же время по определению $\widetilde{S}$ получаем, что (см. (\ref{14.37}))
$$\widetilde{S}(\mu_0)=\biggl(\int\limits_E s_j\,d\mu_0\biggl)_{j\in\overline{1,N}}.$$
Тогда из (\ref{14.55}) вытекает справедливость системы неравенств $$\biggl|\int\limits_E s_j\Theta_{\mu_0}^+[\mathfrak{K}]\,d\eta-\int\limits_E s_j\,d\mu_0\biggl|<\widehat{\varepsilon} \ \ \  \forall\,j\in\overline{1,N}.$$
Таким образом, с учетом (\ref{14.41}), (\ref{14.54}) и (\ref{14.55}) получаем для  вектора $$\biggl(\int\limits_E s_j\Theta_{\mu_0}^+[\mathfrak{K}]\,d\eta\biggl)_{j\in\overline{1,N}}\in\rr^N$$ очевидное свойство
\begin{multline*}
 \exists\,(y_i)_{i\in\overline{1,N}}\in Y:\\
\biggl(y_j=\int\limits_E s_j\Theta_{\mu_0}^+[\mathfrak{K}]\,d\eta \ \ \forall\,j\in J\biggl)\ \& \left(\,\biggl|\int\limits_E s_j\Theta_{\mu_0}^+[\mathfrak{K}]\,d\eta-y_j\biggl|<\widehat{\varepsilon} \ \ \forall\,j\in\overline{1,N}\setminus J\right).
\end{multline*}

С учетом (\ref{14.50}) получаем следующее свойство
$$S\left(\Theta_{\mu_0}^+[\mathfrak{K}]\right)=\biggl(\,\int\limits_E s_j\Theta_{\mu_0}^+[\mathfrak{K}]\,d\eta\biggl)_{j\in\overline{1,N}}\in\widehat{\mathbb{O}}_N(Y,\widehat{\varepsilon}\mid J),$$
откуда в силу (\ref{14.49}) вытекает при условии (\ref{14.52}) включение
$$\Theta_{\mu_0}^+[\mathfrak{K}]\in\widehat{\mathbb{Y}}_{\widehat{\varepsilon}}[J].$$

Итак (см. (\ref{14.52})), истинна импликация
$$\left(\widetilde{\mathfrak{K}}\prec\mathfrak{K}\right)\Rightarrow\left(\Theta_{\mu_0}^+[\mathfrak{K}]\in
\widehat{\mathbb{Y}}_{\widehat{\varepsilon}}[J]\right).$$
Поскольку выбор $\mathfrak{K}$ был произвольным, установлено, что $\exists\,\widehat{\mathcal{K}}\!\!\in\!\mathbf{D}(E,\ml) \ \forall\,\mathcal{K}\!\!\in\!\mathbf{D}(E,\ml)$
$$\left(\widehat{\mathcal{K}}\prec\mathcal{K}\right)\Rightarrow\left(\Theta_{\mu_0}^+[\mathcal{K}]
\in\widehat{\mathbb{Y}}_{\widehat{\varepsilon}}[J]\right).$$
Напомним, что и выбор $\widehat{\varepsilon}$ \  был произвольным. Таким образом, направленность  $\left(\mathbf{D}(E,\ml),\prec,\Theta_{\mu_0}^+[\cdot]\right)$ в $M_*^+[\,b\mid E;\ml;\eta]$  такова, что
\begin{multline*}
\left(\left(\mathbf{D}(E,\ml),\prec,I\circ\Theta_{\mu_0}^+[\cdot]\right)\stackrel{\tau_\eta^*[\ml|b]}\longrightarrow\mu_0\right)\ \& \\ \left(\forall\,\varepsilon\bn \ \exists\,\widehat{\mathcal{K}}\in\mathbf{D}(E,\ml)\ \forall\,\mathcal{K}\in\mathbf{D}(E,\ml) \ \left(\widehat{\mathcal{K}}\prec\mathcal{K}\right)\Rightarrow\left(\Theta_{\mu_0}^+[\mathcal{K}]
\in\widehat{\mathbb{Y}}_{\varepsilon}[J]\right)\right).
\end{multline*}
С учетом (\ref{14.40}) получаем теперь (см. (\ref{14.42})), что $\mu_0\in\widehat{\mathfrak{M}}_J.$ Таким образом, установлена цепочка вложений
\beq{14.56}
\widetilde{S}^{-1}(Y)\subset\widehat{\mathfrak{M}}_J\subset\mathfrak{M}.
\eeq
Выберем теперь произвольно $\mu_*\in\mathfrak{M}.$  Тогда, в частности, имеем (см. (\ref{14.32}))
\beq{14.500}
\mu_*\in\Xi_+^*[\,b\mid E;\ml;\eta].
\eeq

При этом согласно (\ref{14.33'}) справедливо следующее свойство:
\beq{14.501}
\mu_*\in\pr{cl}\left(I^1(\mathbb{Y}_\varepsilon),\tau_\eta^*[\ml|b]\right)\ \ \forall\,\varepsilon\bn.
\eeq
Покажем, что $\mu_*\in\widetilde{S}^{-1}(Y).$  В самом деле, пусть ${\ae}\bn.$ Тогда в силу (\ref{14.501})
\beq{14.502}
\mu_*\in\pr{cl}\left(I^1\left(\mathbb{Y}_{\frac{{\ae}}{2}}\right),\tau_\eta^*[\ml|b]\right),
\eeq
где (см. (\ref{12.11})) справедливо равенство
 \begin{multline*}
 \mathbb{Y}_{\frac{{\ae}}{2}}=\left\{f\in M_*^+[\,b\mid E;\ml;\eta]\mid \biggl(\,\int\limits_E s_j f\,d\eta\biggl)_{j\in\overline{1,N}}\in\mathbb{O}_N\left(Y,\frac{{\ae}}{2}\right)\right\}=\\=S^{\,-1}\left(\mathbb{O}_N\left(Y,\frac{{\ae}}{2}\right)\right).\end{multline*}

Кроме того, в силу (\ref{14.12}) имеем равенство $$I^1\left(Y_{\frac{{\ae}}{2}}\right)=\left\{f*\eta:f\in Y_{\frac{{\ae}}{2}}\right\}.$$ Свойство (\ref{14.502}) означает, что
\beq{14.503}
G\cap I^1\left(Y_{\frac{{\ae}}{2}}\right)\neq\zer \ \ \forall\,G\in N_{\tau_\eta^*[\ml|b]}^0(\mu_*).
\eeq
  Вместе с тем $\mu_*\in\mathbb{A}(\ml),$ а тогда для $\mathbf{K}\triangleq\{s_j:j\in\overline{1,N}\}\in\pr{Fin}(B(E,\ml))$ получаем, в частности представление
 \begin{multline}\label{14.503'}
N_\ml^*\left(\mu_*,\mathbf{K},\frac{{\ae}}{2}\right)=\left\{\nu\in\mathbb{A}(\ml)\mid\biggl|\int\limits_E s_j\,d\mu_*-\int\limits_E s_j\,d\nu\biggl|<\frac{{\ae}}{2} \ \forall\,j\in\overline{1,N}\right\}\in \\ \in N_{\tau_*(\ml)}^{0}(\mu_*).
\end{multline}
Поэтому имеем очевидное свойство
\begin{multline*}
G^*\triangleq N_\ml^*\left(\mu_*,\mathbf{K},\frac{{\ae}}{2}\right)\cap\Xi_+^*[\,b\mid E;\ml;\eta]=\\=\left\{\nu\in\Xi_+^*[\,b\mid E;\ml;\eta]\mid\biggl|\int\limits_E s_j\,d\mu_*-\int\limits_E s_j\,d\nu\biggl|<\frac{{\ae}}{2} \ \forall\,j\in\overline{1,N}\right\}\in N_{\tau_\eta^*[\ml|b]}^{0}(\mu_*);
\end{multline*}
следует учесть определение топологии  $\tau_\eta^*[\ml|b].$ Тогда в силу (\ref{14.503})
$G^*\cap I^1\left(Y_{\frac{{\ae}}{2}}\right)\neq\zer. $
С учетом последнего свойство выберем и зафиксируем
\beq{14.504}
\mu^*\in G^*\cap I^1\left(Y_{\frac{{\ae}}{2}}\right).
\eeq
Тогда, в частности, имеем включения $\mu^*\in N_\ml^*\left(\mu_*,\mathbf{K},\frac{{\ae}}{2}\right)$ и $\mu^*\!\!\in\!\Xi_+^*[\,b\mid E;\ml;\eta].$
Из (\ref{14.503'}) вытекает с очевидностью система неравенств
\beq{14.505}
\biggl|\int\limits_E s_j\,d\mu_*-\int\limits_E s_j\,d\mu^*\biggl|<\frac{{\ae}}{2} \ \forall\,j\in\overline{1,N}.
\eeq

С другой стороны, $\mu^*\in  I^1\left(\mathbb{Y}_{\frac{{\ae}}{2}}\right).$ Это означает, что $\mu^*=\varphi*\eta,$ где $\varphi\in\mathbb{Y}_{\frac{{\ae}}{2}}.$ Тогда $$\left(\int\limits_E s_j \varphi\,d\eta\right)_{j\in\overline{1,N}}\in\mathbb{O}_N\left(Y,\frac{{\ae}}{2}\right),$$
где $\varphi\in M_*^+[\,b\mid E;\ml;\eta]$ и, в частности, $\varphi\in B_0^+(E,\ml).$ При этом по свойствам неопределенного интеграла $\int\limits_E s_j \varphi\,d\eta=\int\limits_E s_j\,d\mu^* \ \ \forall\,j\in\overline{1,N}.$
В итоге имеем следующее включение $\left(\int\limits_E s_j \,d\mu^*\right)_{j\in\overline{1,N}}\in\mathbb{O}_N\left(Y,\frac{{\ae}}{2}\right).$
С учетом (\ref{12.9}) получаем, что для некоторого вектора
\beq{14.506}
(y_j^*)_{j\in\overline{1,N}}\in Y
\eeq
  справедлива система неравенств $$\left|y_j^*-\int\limits_E s_j\,d\mu^*\right|<\frac{{\ae}}{2} \ \forall\,j\in\overline{1,N}.$$
Как следствие из (\ref{14.505}) имеем теперь
$$\biggl|y_j^*-\int\limits_E s_j\,d\mu_*\biggl|<\ae \ \forall\,j\in\overline{1,N}.$$
Мы получили, следовательно, неравенство
$$\left\|(y_j^*)_{j\in\overline{1,N}}-\widetilde{S}(\mu_*)\right\|_N<{\ae}.$$

Поскольку ${\ae}$ выбиралось произвольно, то  (см. (\ref{14.506}))
 $$\forall\,\varepsilon\bn \ \exists\,y\in Y: \ \left\|y-\widetilde{S}(\mu_*)\right\|_N<\varepsilon.$$
 С учетом условия 15.2 получаем включение  $\widetilde{S}(\mu_*)\in Y,$ а тогда $\mu_*\in\widetilde{S}^{\,-1}(Y),$ что и требовалось доказать. Итак,
  $$\mathfrak{M}\subset\widetilde{S}^{\,-1}(Y).$$
  Последнее вложение в сочетании с (\ref{14.56}) завершает доказательство (\ref{14.39}). $\hfill\square$

Из (\ref{14.30})--(\ref{14.33}) и предложения 15.1 вытекает основная
\begin{theorem}
Справедлива следующая цепочка равенств
\begin{multline*}
(\pr{\mathbf{as}})\left[M_*^+[\,b\mid E;\ml;\eta];\rr^n;\tau_{\rr}^{(n)};\Pi;\mathfrak{Y}\right]=\\=(\pr{\mathbf{as}})\left[M_*^+[\,b\mid E;\ml;\eta];\rr^n;\tau_{\rr}^{(n)};\Pi;\widehat{\mathfrak{Y}}[J]\right]=\\=(\pr{\mathbf{sas}})\left[M_*^+[\,b\mid E;\ml;\eta];\rr^n;\tau_{\rr}^{(n)};\Pi;\mathfrak{Y}\right]=\\=(\pr{\mathbf{sas}})\left[M_*^+[\,b\mid E;\ml;\eta];\rr^n;\tau_{\rr}^{(n)};\Pi;\widehat{\mathfrak{Y}}[J]\right]=\widetilde{\Pi}^1\left(\widetilde{S}^{\,-1}(Y)\right).
\end{multline*}
\end{theorem}

Мы получили следующий важный факт: $\widetilde{\Pi}^1\left(\widetilde{S}^{\,-1}(Y)\right)$ есть <<универсальное>> МП. Упомянутая  <<универсальность>> означает, в частности, справедливость свойства асимптотической нечувствительности  при ослаблении части ограничений. Данная часть определяется индексами $j\in J$ (см. (\ref{12.9}), (\ref{12.11}) и (\ref{12.14}), (\ref{12.15})).
\begin{zam}
Рассмотрим случай, когда среди функций $s_j,$ $j\in\overline{1,N},$ нет ступенчатых. Итак, пусть в пределах настоящего замечания $$s_j\notin B_0(E,\ml) \ \forall\,j\in\overline{1,N}.$$

Тогда  полагаем $J=\zer,$ что позволяет говорить о выполнении условия 15.1 (конечно, данное условие не является в нашем случае содержательным, но с формальной точки зрения оно соблюдается при упомянутой конкретизации множества $J$). С учетом (\ref{12.19}) получаем равенство $\mathfrak{Y}=\widehat{\mathfrak{Y}}[J].$ В этих условиях утверждение теоремы 15.1 сводится к цепочке равенств:
  \begin{multline*}
 (\pr{\mathbf{as}})\left[M_*^+[\,b\mid E;\ml;\eta];\rr^n;\tau_{\rr}^{(n)};\Pi;\mathfrak{Y}\right]=\\=(\pr{\mathbf{sas}})\left[M_*^+[\,b\mid E;\ml;\eta];\rr^n;\tau_{\rr}^{(n)};\Pi;\mathfrak{Y}\right]=\widetilde{\Pi}^1\left(\widetilde{S}^{\,-1}(Y)\right).\end{multline*}
\end{zam}
\begin{zam}
Рассмотрим теперь другой крайний случай, полагая в пределах настоящего замечания, что
\beq{14.507}
s_j\in B_0(E,\ml) \ \forall\,j\in\overline{1,N}.
\eeq

   Полагаем, что  $J=\overline{1,N},$ а условие 15.2 выполнено и, стало быть, верно утверждение теоремы 15.1.  Как уже отмечалось в разделе 13, в этом случае $\widehat{\mathfrak{Y}}[J]=\{Y\},$ а потому согласно (\ref{12.16}) и (\ref{14.5})
$$(\pr{\mathbf{as}})\left[M_*^+[\,b\mid E;\ml;\eta];\rr^n;\tau_{\rr}^{(n)};\Pi;\widehat{\mathfrak{Y}}[J]\right]=
\pr{cl}\left(\widetilde{\Pi}^1\left(\widetilde{S}^{\,-1}(Y)\right),\tau_{\rr}^{(n)}\right).$$
 Таким образом,  при условии (\ref{14.507}) из теоремы 15.1 следует  устойчивость множества достижимости в $\rr^n$ при ослаблении  $Y$-ограничения:
 \begin{multline}\label{14.507'}
 (\pr{\mathbf{as}})\left[M_*^+[\,b\mid E;\ml;\eta];\rr^n;\tau_{\rr}^{(n)};\Pi;\mathfrak{Y}\right]=\\=(\pr{\mathbf{sas}})\left[M_*^+[\,b\mid E;\ml;\eta];\rr^n;\tau_{\rr}^{(n)};\Pi;\mathfrak{Y}\right]=\pr{cl}\left(\widetilde{\Pi}^1\left(\widetilde{S}^{\,-1}(Y)\right),\tau_{\rr}^{(n)}\right).
 \end{multline}

С учетом определений раздела 13 (см., в частности, (\ref{12.9}), (\ref{12.11})) цепочка равенств (\ref{14.507'}) означает, что наши возможности по реализации значений вектор-функционала $\Pi$ практически не изменяются при замене $Y$-ограничения (\ref{12.5}) подобным условием, где в правой части (\ref{12.5}) множество $Y$ заменено достаточно малой окрестностью вида (\ref{12.9}). Точнее, имеется в виду замена $Y$-ограничения условием (\ref{12.10}) при достаточно малом значении $\varepsilon,$ $\varepsilon>0.$ Правда, данное суждение об устойчивости так же, как и ранее сформулированное в теореме 15.1 утверждение об асимптотической нечувствительности при ослаблении части ограничений, дано в терминах МП, т.е. в терминах предельных множеств. Однако в \cite[гл.\,4]{46} и \cite[гл.\,3]{41}  указано, как именно данные, асимптотические по смыслу положения, перевести в суждения о близости конкретных множеств достижимости при достаточно малых значениях $\varepsilon, \ \varepsilon>0.$ Мы предлагаем читателю проделать это самостоятельно, ориентируясь на предложение 3.6.3 и следствие 3.6.3 в \cite{41}.

Отметим еще одно полезное следствие теоремы 15.1, реализующееся при условии (\ref{14.507}):
$$\widetilde{\Pi}^1\left(\widetilde{S}^{\,-1}(Y)\right)=\pr{cl}\left(\Pi^1\left(S^{\,-1}(Y)\right),\tau_{\rr}^{(n)}\right).$$
В этом последнем равенстве имеем следующее свойство: действие обобщенных элементов, точно соблюдающих  $Y$-ограничение (имеется в виду условие
  $$\biggl(\int\limits_E s_j \,d\mu\biggl)_{j\in\overline{1,N}}\in Y$$
  на выбор $\mu\in\Xi_+^*[\,b\mid E;\ml;\eta]$), получается в случае (\ref{14.507}) замыканием действия обычных управлений из  $ M_*^+[\,b\mid E;\ml;\eta],$ которые также соблюдают $Y$-ограничение точно (см. в этой связи (\ref{12.5})), в то время как при отказе от (\ref{14.507}) могут возникать эффекты, не реализуемые в классе точных обычных решений.
Отметим, наконец, один естественный вариант (\ref{12.5}), отвечающий ситуации (\ref{14.507}):
$$\biggl(\,\int\limits_{\mathbb{L}_j}f\,d\eta\biggl)_{j\in\overline{1,N}}\in Y$$
(здесь $f\in M_*^+[\,b\mid E;\ml;\eta]$), где $\mathbb{L}_1\in\ml,\ldots,\mathbb{L}_N\in\ml.$  В самом деле, кортеж (\ref{12.3}) в данном случае  задается условиями $$s_j=\chi_{\mathbb{L}_j} \ \forall\,j\in\overline{1,N};$$
условие (\ref{14.507}) выполняется очевидным образом.\end{zam}

\newpage

\begin{center} \section*{Глава 4. Линейные задачи управления с ограничениями импульсного характера} \end{center}
\section{Линейные управляемые системы с разрывностью в коэффициентах при управляющих воздействиях (случай ограничений импульсного характера)}\setcounter{equation}{0}\setcounter{proposition}{0}\setcounter{zam}{0}\setcounter{corollary}{0}\setcounter{definition}{0}
   \ \ \ \ \ В настоящем и последующих разделах будем рассматривать применение весьма общей конструкции расширения, приведенной  в предыдущей главе,  в  задаче о достижимости для одного класса управляемых систем. Речь идет об управлении линейной системой
\beq{16.1}
\dot{x}(t)=A(t)x(t)+f(t)c(t),
\eeq
      функционирующей в $n$-мерном фазовом пространстве  на конечном промежутке времени $[t_0,\theta_0],$ где $t_0\in\rr,$ $\theta_0\in\rr,$ $t_0<\theta_0.$ В  (\ref{16.1}) $A(\cdot)=(A(t),t_0\leqslant t\leqslant\theta_0)$ есть матричнозначная функция, а $c=c(\cdot)=(c(t),t_0\leqslant t\leqslant\theta_0)$\,~--- $n$-вектор-функция, не являющаяся, вообще говоря, непрерывной. Имея в виду последующее согласование содержательной задачи с абстрактной схемой предыдущей главы, будем  полагать далее, что $A(t)$ есть $n\times n$-матрица; полагаем для простоты все компоненты отображения $A=A(\cdot)$ непрерывными в/з функциями на $[t_0,\theta_0].$

   Полагаем далее, что у системы (\ref{16.1}) задано начальное состояние: $x(t_0)=x_0\in\rr^n,$ а выбор $f$ находится в нашем распоряжении. Будем полагать, что $f$\,---  неотрицательная в/з функция, имеющая смысл интенсивности подачи топлива в двигательную установку, ориентация которой задана посредством вектор-функции $c(\cdot).$ Иными словами, ориентация двигательной установки реализуется по заданной программе. Мы будем рассматривать постановку, в которой весь имеющийся запас топлива должен быть израсходован до момента $\theta_0$ (отметим один вполне жизненный аналог такого требования: самолет с пассажирами, идущий на вынужденную посадку, должен израсходовать весь  запас топлива в целях безопасности посадки); заметим, однако, что другие ограничения ресурсного характера также исследовались (см., например, \cite{25,46,47}), но мы этим заниматься не будем, предоставляя соответствующие (идейно аналогичные) построения заинтересованному читателю.

Теперь мы уточним условия на вектор-функцию с компонентами $$c_1=c_1(\cdot),\ldots,c_n=c_n(\cdot).$$

Условимся в дальнейшем множество $E$ (см. предыдущую главу) отождествлять с отрезком $[t_0,\theta_0].$  Итак,  в дальнейших построениях
\beq{16.2}
E=[t_0,\theta_0]
\eeq
является промежутком управления. В качестве $\ml$ рассматриваем алгебру п/м $E$ (\ref{16.2}), порожденную полуалгеброй
\beq{16.2'}
\mathfrak{I}\triangleq\left\{L\in\pp(E)\mid\exists\,c\in E \, \exists\,d\in E: \ \left(\,]c,d\subset L\right)\ \& \left(L\subset[c,d]\right)\right\}\in \Pi[E].
\eeq

Далее $\ml$ определяется следующим образом (см. \cite[(1.7.9)]{30}):
\beq{16.2''}
\ml=\left\{\Lambda\in\pp(E)\mid\exists\,m\in\mathbb{N}:\Delta_m(\Lambda,\mathfrak{I})\neq\zer\right\};
\eeq
разумеется, в частности, $\ml\in\Pi[E],$ а потому наше конкретное ИП $(E,\ml)$ удовлетворяет всем условиям, введенным  в построениях предыдущей главы.

Будем предполагать в дальнейшем, что $\eta$ есть след меры Лебега на алгебру $\ml$ ($\eta$ есть мера \cite[(8.5.1)]{38}  при условии, что $a=t_0$  и $b=\theta_0$), т.е. мера, получаемая продолжением на $\ml$ функции длины, определенной на полуалгебре $\mathfrak{I}.$ Заметим, что мера $\eta$ в данном случае счетно-аддитивна. Таким образом,  мы используем ниже конкретный вариант пространства $(E,\ml,\eta),$  связанный с (\ref{16.2}). При этом, конечно,
$$\int\limits_E h\,d\mu=\int\limits_{[t_0,\theta_0]} h\,d\mu \ \ \forall\,h\in B(E,\ml) \ \forall\,\mu\in\mathbb{A}(\ml).$$
Особенно простую природу имеют здесь функции из $B_0(E,\ml).$ По сути дела это кусочно-постоянные функции, которые, правда, могут обладать скачками в отдельных точках (точнее, в конечном числе таковых). Так, например, функция $f:E\rightarrow \rr,$
для которой
 \begin{multline*}
\left(f(t)\triangleq 1 \ \forall\,t\in\left[t_0,\frac{t_0+\theta_0}{2}\right[\,\right)\ \& \ \left(f\left(\frac{t_0+\theta_0}{2}\right)\triangleq 3\right) \ \& \\ \& \left(f(t)\triangleq 2 \ \forall\,t\in\left]\frac{t_0+\theta_0}{2},\theta_0\right]\,\right)\end{multline*}
является элементом $B_0^+(E,\ml).$ При интегрировании  функции $f$ по \\ $\mu\!\!\in\! (\pr{add})^+[\ml;\eta]$ значение в точке $\frac{t_0+\theta_0}{2}$ оказывается, конечно, несущественным, а, точнее, не влияющим на значение интеграла. Мы будем называть здесь элементы  $$B_0(E,\ml)$$ кусочно-постоянными функциями, игнорируя упомянутую и несущественную в интегральном смысле особенность.

Возвращаясь к (\ref{16.1}), полагаем, что все компоненты вектор-функции  $c,$ т.е. функции $c_1,\ldots,c_n$  являются ярусными  относительно нашего конкретного варианта $(E,\ml).$  Иными словами, требуем, чтобы каждая из функций $c_j, \ j\in\overline{1,n},$ являлась равномерным пределом последовательности кусочно-постоянных в/з функций. В этих условиях при всяком выборе $f\in B_0(E,\ml)$ (т.е. при всяком выборе кусочно-постоянной управляющей программы)
\beq{16.3}
\mathbf{x}(\theta_0)\triangleq\Phi(\theta_0,t_0)x_0+\int_Ef(t)\Phi(\theta_0,t)c(t)\eta(dt)
\eeq
есть терминальное состояние системы (\ref{16.1}) при использовании в/з функции $f$ в качестве управляющего воздействия (интегрирование вектор-функции осуществляется покомпонентно). В (\ref{16.3}) имеем простое следствие формулы Коши (см. \cite{23}). Ясно, что $\mathbf{x}(\theta_0)$ зависит от $f,$ т.е.
\beq{16.3'}
\mathbf{x}(\theta_0)=\mathbf{x}_f(\theta_0).
\eeq
Теперь введем  ограничения на выбор $f.$ Здесь используем прежде всего требование
\beq{16.4}
f\in M_*^+[\,b\mid E;\ml;\eta].
\eeq
 В рассматриваемом случае условие (\ref{16.4}) означает, что $f$ есть неотрицательная кусочно-постоянная в/з функция, у которой
\beq{16.5}
\int_Ef\,d\eta=\int\limits_{[t_0,\theta_0]} f\,d\eta=b.
\eeq

С учетом  (\ref{16.5}) мы получаем следующее толкование условия (\ref{6.4}): в качестве допустимого (реализуемого) программного управления мы можем использовать любую кусочно-постоянную  неотрицательную функцию $f=(f(t),$ $t_0\leqslant t\leqslant\theta_0)$ с импульсом силы, равным $b.$ Константа  $b$ определяет по смыслу имеющийся запас топлива, которое в процессе управления должно быть полностью израсходовано. Таким образом, множество $M_*^+[\,b\mid E;\ml;\eta]$  имеет в нашем случае вполне понятный с практической точки зрения смысл: мы рассматриваем реализуемые программные управления (т.е. кусочно-постоянные функции), регулирующие тягу двигателя и такие, что при этом весь запас топлива должен быть израсходован полностью,  что как раз и формализуется посредством (\ref{16.5}).

Теперь рассмотрим вопрос, связанный с введением $Y$-ограничения. В принципе данное ограничение может иметь достаточно произвольный характер, что и было отражено в предыдущей главе. В \cite[\S 1.3]{25} подробно рассматривался случай, когда упомянутое $Y$-ограничение связано с обеспечением нужного режима работы двигателя (чередование импульсов и пауз). Сейчас мы остановимся на другой возможности: упомянутое ограничение может возникать из соображений, связанных с краевыми и промежуточными условиями. Обсудим простейший вариант (возможные обобщения предоставляются читателю). Итак, предположим, что $N=n$ и задан некоторый момент $t_*\in[t_0,\theta_0].$ Требуется обеспечить условие
\beq{16.6}
x(t_*)\in\widetilde{Y},
\eeq
где $\widetilde{Y}$\,--- непустое замкнутое множество в $\rr^N.$ Наша система должна <<побывать>> в множестве $\widetilde{Y}$ в момент $t_*.$ Отметим, кстати, что само множество $\widetilde{Y}$ в ряде случаев является цилиндрическим и, стало быть, <<бесконечно протяженным>> по части координат. Данный случай оставляем на рассмотрение читателю. С использованием формулы Коши приходим к следующей эквивалентной форме (\ref{16.6}):
\beq{16.7}
\int\limits_{[t_0,t_*]}f(t)\Phi(t_*,t)c(t)\,\eta(dt)=\int\limits_{E}f(t)\Phi(t_*,t)c(t)\chi_{[t_0,t_*]}(t)\,\eta(dt)\in Y,
\eeq
где $Y$ (в данной конкретизации) имеет вид
$$Y=\{y-\Phi(t_*,t_0)x_0:\,y\in\widetilde{Y}\};$$
множество $Y$ соответствует предположениям предыдущей главы. В свою очередь, (\ref{16.6}) является, как легко видеть, вариантом уже рассмотренного $Y$-ограничения. Предлагаем читателю самостоятельно подобрать ярусные функции $s_1,\ldots,s_N,$ для которых (\ref{16.6}) сводится к виду, рассматриваемому в предыдущей главе.

Заметим, что само условие (\ref{16.7}) может быть усложнено. Так, например, оно может соответствовать условиям на реализацию фазового состояния  не в один момент времени; таких моментов может быть несколько. Мы, однако, сохраним требование к реализации $Y$-ограничения в том виде, как это было сформулировано в предыдущей главе, имея в виду конкретный вариант $(E,\ml,\eta),$ связанный с (\ref{16.2}).

Обсудим далее проблему построения и исследования области достижимости, а также (и даже в большей степени) ее регуляризации в виде МП. При этом по смыслу данной задачи следует рассматривать отображение
\beq{16.8}
f\mapsto x_f(\theta_0):\, M_*^+[\,b\mid E;\ml;\eta]\rightarrow\rr^n.
\eeq
Однако в силу (\ref{16.3}) имеем, что для всех наших целей достаточно рассматривать отображение (а, точнее, вектор-функционал)
\beq{16.9}
f\mapsto\int\limits_{E}f(t)\Phi(\theta_0,t)c(t)\,\eta(dt):\,M_*^+[\,b\mid E;\ml;\eta]\rightarrow\rr^n,
\eeq
поскольку значения отображений (\ref{16.8}), (\ref{16.9}) различаются только сдвигом на фиксированный вектор $\Phi(\theta_0,t_0)x_0.$ Что же касается
(\ref{16.9}), то данное отображение является вариантом $\Pi$ (см. предыдущую главу). В самом деле, вектор-функция $\Psi_f$ вида
 $$t\mapsto f(t)\Phi(\theta_0,t)c(t): E\rightarrow\rr^n,$$  где $f\in M_*^+[\,b\mid E;\ml;\eta],$  имеет следующие компоненты
 $\Psi_{f,1},\ldots,\Psi_{f,n}$
 (при этом $\Psi_f(t)=\left(\Psi_{f,i}(t)\right)_{i\in\overline{1,n}} \ \forall\,t\in E$).  Напомним, что $n$-мерные векторы есть кортежи вещественных чисел <<длины>> $n;$  при $j\in\overline{1,n}$ и $t\in E$
 $$\Psi_{f,j}(t)=f(t)\sum\limits_{k=1}^n\Phi_{j,k}(\theta_0,t)c_k(t).$$

 Разумеется, каждая функция
 \beq{16.10}
 t\mapsto\sum\limits_{k=1}^n\Phi_{j,k}(\theta_0,t)c_k(t):\,E\rightarrow\rr,
 \eeq
 где $j\in\overline{1,n}$ является ярусной в смысле исследуемого сейчас варианта ИП $(E,\ml;)$ см. (\ref{16.2}). Поэтому мы можем при каждом  $j\in\overline{1,n}$ отождествить функцию $\pi_j$ используемую в заключительной части предыдущей главы с (\ref{16.10}), тогда (\ref{16.9}) превращается в  вектор-функционал $\Pi$ предыдущей главы. Важно, что в этих терминах мы получаем возможность изучать область достижимости, определяемую в терминах (\ref{16.8}), что представляет не только теоретический, но и практический интерес (здесь имеем в виду также и МП, являющееся по сути дела своеобразной регуляризацией области достижимости). После сведения (\ref{16.8}) к виду (\ref{16.9}) и отождествления функций (\ref{16.10}) с $\pi_1,\ldots,\pi_n$ получаем возможность исследовать исходную задачу методами предыдущей главы; при этом отображение $\Pi$ будет конкретизироваться в виде (\ref{16.9}).

Полезно, однако, отметить одну интерпретацию в терминах самой задачи управления. Для этого вернемся к  (\ref{16.3}), где указано представление траектории в последний момент времени; при этом предполагалось, что данная траектория порождена обычным управлением $f.$ Сейчас, однако, мы можем определить траектории системы (\ref{16.1}), порожденные обобщенными управлениями-мерами. В качестве последних будем использовать элементы множества
  $\Xi_+^*[\,b\mid E;\ml;\eta]$ в той конкретизации $(E,\ml,\eta),$ которая использовалась в настоящем разделе (заметим, что в (\ref{16.3}) следует учитывать (\ref{16.3'}), т.е. на самом деле вектор в правой части (\ref{16.3}) определяет $x_f(\theta_0)$).

Итак, мы при $\mu\in\Xi_+^*[\,b\mid E;\ml;\eta]$ введем вектор-функцию $\widetilde{x}_\mu=\widetilde{x}_\mu(\cdot)$ на $E,$ для которой
\beq{16.11}
\widetilde{x}_\mu(t)\triangleq\Phi(t,t_0)x_0+\int\limits_{[t_0,t]}\Phi(t,\xi)c(\xi)\,\mu(d\xi) \ \forall\,t\in E,
\eeq
где интеграл следует вычислять покомпонентно; в частности (и, по сути дела, только это для нас существенно)
\beq{16.12}
\widetilde{x}_\mu(\theta_0)\triangleq\Phi(\theta_0,t_0)x_0+\int\limits_{E}\Phi(\theta_0,t)c(t)\,\mu(dt).
\eeq

Прием такого рода был использован в \cite[\S\,6.5]{25}. Вектор-функцию (\ref{16.11}) можно толковать как обобщенную траекторию системы (\ref{16.1}). Ясно, что при $\mu=f*\eta,$ где $f$ удовлетворяет (\ref{16.4}),
$$\widetilde{x}_\mu(\cdot)=\mathbf{x}_f(\cdot);$$
в частности, имеем в данном случае равенство $\widetilde{x}_\mu(\theta_0)=\mathbf{x}_f(\theta_0).$

\section{Управление материальной точкой при  моментных ограничениях }\setcounter{equation}{0}\setcounter{proposition}{0}\setcounter{zam}{0}\setcounter{corollary}{0}\setcounter{definition}{0}
   \ \ \ \ \ Начиная с этого раздела, мы обращаемся к простейшей механической системе~--- материальной точке. Таким образом, мы конкретизируем систему  (\ref{16.1}) следующим образом:
   \beq{17.1}
   \dot{x}_1(t)=x_2(t), \ \dot{x}_2(t)=\mathbf{c}(t)f(t).
   \eeq
Полагаем  для простоты $t_0=0, \ \theta_0=1, \ x_1(0)=0, \ x_2(0)=0.$ Итак, постулируем, что в начальный момент $t_0=0$ материальная точка покоится <<в нуле>>. Сохраняя предположения относительно ИП $(E,\ml),$ введенные в предыдущем разделе, имеем здесь
\beq{17.1'}
E=[0,1],
\eeq
а  $\ml$ определяется посредством (\ref{16.2''}), где $\mathfrak{I}$ конкретизируется в виде (\ref{16.2'}) при условиях, что $t_0=0$ и $\theta_0=1.$

Мы полагаем, что $\mathbf{c}\in B(E,\ml)$ (получаем вариант   $c=c(\cdot)$ предыдущего раздела, где $c_1=c_1(\cdot)$ есть функция, тождественно равная нулю, а $c_2=c_2(\cdot)$ отождествляется с функцией $\mathbf{c},$ т.е. $c_2=\mathbf{c}$). В наших последующих построениях $n=N=2$ (функции $s_1$ и $s_2$ будут введены позже).

Сейчас мы заметим, что при выборе всякой управляющей функции $f\in M_*^+[\,b\mid E;\ml;\eta]$ траектория $\mathbf{x}_f=\mathbf{x}_f(\cdot)$ системы (\ref{17.1}) имеет следующие компоненты
\beq{17.200}
\mathbf{x}_{f,1}(t)=\int\limits_{[0,t]}(t-\tau)\mathbf{c}(\tau)f(\tau)\,\eta(d\tau),
\eeq
\beq{17.201}
\mathbf{x}_{f,2}(t)=\int\limits_{[0,t]}\mathbf{c}(\tau)f(\tau)\,\eta(d\tau).
\eeq
В целях согласования с конструкциями предыдущей главы мы используем здесь интегралы по мере $\eta,$ определенные в разделе 9. В частности,
$$\mathbf{x}_{f,1}(1)=\int\limits_{E}(1-t)\mathbf{c}(t)f(t)\,\eta(dt), \ \mathbf{x}_{f,2}(1)=\int\limits_{E}\mathbf{c}(t)f(t)\,\eta(dt).
$$
В соответствии с (\ref{16.8}), (\ref{16.9}) отображение
\beq{17.1''}
f\mapsto\mathbf{x}_f(1): M_*^+[\,b\mid E;\ml;\eta]\rightarrow\rr^2
\eeq
задает вектор-функционал $\Pi$ в конкретизации настоящей главы. Итак,
$$\Pi(f)=(\mathbf{x}_{f,i}(1))_{i\in\overline{1,2}}\in\rr^2$$
есть вектор, у которого первая компонента
$$\Pi(f)(1)=\int\limits_{E}(1-t)\mathbf{c}(t)f(t)\,\eta(dt),$$
а вторая компонента имеет вид
$$\Pi(f)(2)=\int\limits_{E}\mathbf{c}(t)f(t)\,\eta(dt);$$
здесь $f\in M_*^+[\,b\mid E;\ml;\eta].$ С учетом этих представлений получаем, что $\pi_1\!\!\in\! B(E,\ml)$ определяется условием
\beq{17.2}
\pi_1(t)=(1-t)\mathbf{c}(t) \ \forall\,t\in E,
\eeq
а $\pi_2\in B(E,\ml)$ есть функция $\mathbf{c},$ т.е.
\beq{17.3}
\pi_2=\mathbf{c}.
\eeq
Тем самым у нас конкретизирован кортеж (\ref{12.2}). Данная конкретизация приводит к весьма важному варианту $\Pi$ (\ref{17.1'}), связанному с построением и исследованием области достижимости системы (\ref{17.1}).

Теперь рассмотрим одну конкретизацию $Y$-ограничения, связанную с управлением системой (\ref{17.1}) и определяемую т.н. промежуточными условиями. Итак, фиксируем два произвольных момента времени:
$$t_1\in E, \ \ t_2\in E,$$
где $E$ соответствует (\ref{17.1'}). Таким образом, $$t_1\in \rr, \ t_2\in\rr, \ 0\leqslant t_1\leqslant1, \ 0\leqslant t_2\leqslant 1. $$

Напомним, что в настоящем разделе $N=2,$ а $Y$ является, следовательно (см. (\ref{12.4})), непустым замкнутым множеством на плоскости. Рассматриваем следующие ограничения на выбор $f\in M_*^+[\,b\mid E;\ml;\eta]:$
\beq{17.4}
\left(\mathbf{x}_{f,i}(t_i)\right)_{i\in\overline{1,2}}\in Y.
\eeq
\begin{zam}
Условие (\ref{17.4}) является более понятным в том частном случае, когда $Y$ определяется на содержательном уровне выражением
\beq{17.5}
Y=Y_1\times Y_2,
\eeq
где $Y_1$ и $Y_2$  суть непустые замкнутые п/м вещественной прямой $\rr.$ Тогда (\ref{17.4}) сводится к требованию о совместном выполнении условий
\beq{17.6}
\mathbf{x}_{f,1}(t_1)\in Y_1, \ \mathbf{x}_{f,2}(t_2)\in Y_2.
\eeq

Итак, наша материальная точка должна иметь в момент времени $t_1$ координату из множества $Y_1$  и в момент времени $t_2$  скорость из множества $Y_2.$ \end{zam}

Возвращаясь к общему случаю условия (\ref{17.4}) (имеется в виду обобщение (\ref{17.5}), (\ref{17.6})), с учетом формулы Коши и определения $S$ раздела 13 получаем следующую конкретизацию $s_1,s_2$ в (\ref{12.3}):  функция $s_1\in B(E,\ml)$  имеет вид
\beq{17.7}
\left(s_1(t)\triangleq(t_1-t)\,\mathbf{c}(t) \ \forall\,t\in[0,t_1]\right) \ \& \ \left(s_1(t)\triangleq 0 \ \forall\,t\in]t_1,1]\right),
\eeq
а функция $s_2\in B(E,\ml)$ такова, что
\beq{17.8}
\left(s_2(t)\triangleq \mathbf{c}(t) \ \forall\,t\in[0,t_2]\right) \ \& \ \left(s_2(t)\triangleq 0 \ \forall\,t\in]t_2,1]\right).
\eeq

В связи с (\ref{17.7}), (\ref{17.8}) отметим, что в левой части (\ref{12.5}) используются интегралы на всем множестве $E$ (\ref{17.1'}). Иными словами,
\beq{17.8'}
s_1=((t_1-t)\,\mathbf{c}(t))_{t\in E}\cdot\chi_{[0,t_1]},
\eeq
\beq{17.8''}
s_2=\mathbf{c}\chi_{[0,t_2]}.
\eeq

Упомянутые конкретизации кортежа (\ref{12.3}) легко извлекаются из (\ref{17.200}), (\ref{17.201}) и мы рекомендуем читателю провести все выкладки самостоятельно.

В связи с вопросами ослабления $Y$-ограничения напомним, что в рассматриваемом сейчас случае
\begin{multline}\label{17.9}
\mathbb{O}_N(Y,\varepsilon)=\mathbb{O}_2(Y,\varepsilon)=\bigl\{(z_i)_{i\in\overline{1,2}}\in\rr^2\mid\exists\,(y_i)_{i\in\overline{1,2}}\in Y:(\,|y_1-z_1|<\varepsilon)\ \& \\ \& \ (\,|y_2-z_2|<\varepsilon)\bigl\} \ \forall\,\varepsilon\bn.
\end{multline}

Посредством (\ref{17.9}) конкретизируется семейство $\mathfrak{Y}$ раздела 13. В связи с выбором $J$ допускаем следующие два возможных случая:\\
\noindent (1*) \ \  если $\mathbf{c}\in B_0(E,\ml),$ то полагаем, что $J=\{2\}$ (синглетон, содержащий индекс 2);\\
\noindent (2*) \ \  если $\mathbf{c}\notin B_0(E,\ml),$ то полагаем $J=\zer.$

Случай (1*) представляется более интересным в связи с вопросом об условиях асимптотической нечувствительности при ослаблении  части ограничений (последнее наиболее понятно в случае ограничений (\ref{17.6})). В случае  (2*) получаем ситуацию, подобную обсуждаемой в замечании 15.1.

Рассмотрим конкретизацию обобщенной задачи о достижимости, рассматривая к.-а. меры из $\Xi_+^*[\,b\mid E;\ml;\eta]$ в качестве обобщенных управлений. Прежде всего отметим, что в рассматриваемом сейчас случае при $\mu\in\Xi_+^*[\,b\mid E;\ml;\eta]$ обобщенная траектория
$$\widetilde{x}_\mu: E\rightarrow\rr^2$$
имеет следующие компоненты:
$$\widetilde{x}_{\mu,1}(t)=\int\limits_{[0,t]}(t-\xi)\,\mathbf{c}(\xi)\mu(d\,\xi),
\widetilde{x}_{\mu,2}(t)=\int\limits_{[0,t]}\mathbf{c}(\xi)\mu(d\,\xi),$$
где $t\in E.$  В частности, определены значения
$$\widetilde{x}_{\mu,1}(t_1), \ \widetilde{x}_{\mu,2}(t_2), \ \widetilde{x}_{\mu}(1)=\left(\widetilde{x}_{\mu,i}(1)\right)_{i\in\overline{1,2}}.$$
Условие  (\ref{17.4}) имеет следующий обобщенный аналог
\beq{17.10}
\left(\widetilde{x}_{\mu,i}(t_i)\right)_{i\in\overline{1,2}}\in Y.
\eeq
В свою очередь терминальное состояние, определяемое к.-а. мерой  $\mu,$ имеет вид $$\widetilde{x}_{\mu}(1)=\left(\widetilde{x}_{\mu,i}(1)\right)_{i\in\overline{1,2}}\in \rr^2,$$
где
\beq{17.11}
\widetilde{x}_{\mu,1}(1)=\int\limits_E(1-t)\,\mathbf{c}(t)\mu(dt),
\eeq
\beq{17.12}
\widetilde{x}_{\mu,2}(1)=\int\limits_E\mathbf{c}(t)\mu(dt).
\eeq
Теперь используем  конкретизацию $\pi_1,\pi_2$  в виде  (\ref{17.2}), (\ref{17.3}). Итак, имея в виду (\ref{17.11}) и (\ref{17.12}), получаем, что  $\widetilde{\Pi}$ (\ref{14.3}) в рассматриваемом  сейчас случае есть отображение
\beq{17.13}
\mu\mapsto\widetilde{x}_{\mu}(1):\Xi_+^*[\,b\mid E;\ml;\eta]\rightarrow\rr^2.
\eeq
Из (\ref{17.11})--(\ref{17.13}) следует, конечно, что
$$\widetilde{\Pi}(\mu)(1)=\int\limits_E(1-t)\,\mathbf{c}(t)\mu(dt),$$
$$\widetilde{\Pi}(\mu)(2)=\int\limits_E\mathbf{c}(t)\mu(dt),$$
где $\mu\in\Xi_+^*[\,b\mid E;\ml;\eta].$ Аналогичными рассуждениями конкретизируется отображение $\widetilde{S}$  (\ref{14.37}). Здесь мы используем (\ref{17.8'}), (\ref{17.8''}). Итак, в нашем случае
\beq{17.14}
\widetilde{S}:\Xi_+^*[\,b\mid E;\ml;\eta]\rightarrow\rr^2
\eeq
имеет следующие компоненты:
$$\widetilde{S}(\mu)(1)=\int\limits_{E}s_1\,d\mu=\int\limits_{[0,t_1]}(t_1-t)\,\mathbf{c}(t)\mu(dt),$$
$$\widetilde{S}(\mu)(2)=\int\limits_{E}s_2\,d\mu=\int\limits_{[0,t_2]}\mathbf{c}(t)\mu(dt),$$
где $\mu\!\!\in\!\Xi_+^*[b|\,E;\ml;\eta].$ С учетом определения $\widetilde{x}_{\mu,1}\!\!=\!\widetilde{x}_{\mu,1}(\cdot)$ и  $\widetilde{x}_{\mu,2}\!\!=\!\widetilde{x}_{\mu,2}(\cdot)$ получаем, что
$$\widetilde{S}(\mu)(1)=\widetilde{x}_{\mu,1}(t_1), \ \widetilde{S}(\mu)(2)=\widetilde{x}_{\mu,2}(t_2).$$
Таким образом, имеем очевидное равенство
\beq{17.15}
\widetilde{S}^{\,-1}(Y)=\left\{\mu\in\Xi_+^*[\,b\mid E;\ml;\eta]\mid\left(\widetilde{x}_{\mu,j}(t_j)\right)_{j\in\overline{1,2}}\in Y\right\};
\eeq
иными словами,  (\ref{17.15}) есть множество всех обобщенных управлений, допустимых в смысле соблюдения промежуточных условий. Наконец, с учетом этого множество $\widetilde{\Pi}^1\left(\widetilde{S}^{\,-1}(Y)\right)$ имеет следующий вид:
\beq{17.16}
\widetilde{\Pi}^1\left(\widetilde{S}^{\,-1}(Y)\right)=\left\{\widetilde{x}_\mu(1):\mu\in\widetilde{S}^{\,-1}(Y)\right\},
\eeq
т.е. (\ref{17.16}) есть область достижимости в классе обобщенных управлений, для  которых допустимость (помимо принадлежности множеству $ \Xi_+^*[\,b\mid E;\ml;\eta]$) определяется в терминах точного соблюдения промежуточных условий (мы не исключаем, впрочем, тот случай, когда $t_1=t_2=1,$ при котором реализуются краевые условия). Согласно теореме 15.1 упомянутое множество (\ref{17.6}) определяет весьма универсальное МП на фазовой плоскости системы (\ref{17.1}). Упомянутая универсальность проявляется, впрочем,  в рассматриваемом случае при $\mathbf{c}\in B_0(E,\ml).$ Этот случай обсудим в следующем разделе, а сейчас сосредточимся на обсуждении варианта МП, связанного с семейством $\mathfrak{Y},$ используемым для введения ОАХ. Итак, речь пойдет о МП вида (\ref{12.3}), где  ослабление  $Y$-ограничения осуществляется посредством (\ref{17.9}): вместо (\ref{17.4}) требуем при заданном $\varepsilon\bn,$  чтобы $f\in M_*^+[\,b\mid E;\ml;\eta]$ было таким, что
\beq{17.17}
\exists\,(y_i)_{i\in\overline{1,2}}\in Y:\left(\,|\mathbf{x}_{f,1}(t_1)-y_1|<\varepsilon\right)\ \& \ \left(\,|\mathbf{x}_{f,2}(t_2)-y_2|<\varepsilon\right).
\eeq

Множество всех управлений $f\in M_*^+[\,b\mid E;\ml;\eta]$ со свойством (\ref{17.17}) в нашем конкретном случае совпадает с $\mathbb{Y}_\varepsilon$  и характеризует ослабленный в сравнении с (\ref{17.4}) режим управления материальной точкой. Данному условию (\ref{17.17}) отвечает определяемая обычным способом область достижимости
$$\Pi^1(\mathbb{Y}_\varepsilon)=\Pi^1(S^{\,-1}(\mathbb{O}_2(Y,\varepsilon)))=\{\mathbf{x}_f(1): f\in\mathbb{Y}_\varepsilon\}.$$

Операция (\ref{12.13}) определяет предел таких конкретных областей достижимости при $\varepsilon\downarrow 0.$ Данный предел согласно теореме 15.1
совпадает с множеством (\ref{17.16}). Полезно отметить, что даже при условии $\mathbf{c}(t)\equiv 1$ упомянутое МП, определяемое в виде (\ref{17.16}), может не совпадать с множеством $$\pr{cl}\left(\Pi^1(S^{\,-1}(Y)),\tau_\rr^{(2)}\right),$$ которое является замыканием области достижимости в классе обычных управлений, хотя всегда $$\pr{cl}\left(\Pi^1(S^{\,-1}(Y)),\tau_\rr^{(2)}\right)\subset\widetilde{\Pi}^1\left(\widetilde{S}^{\,-1}(Y)\right).$$

В связи с этим наиболее простым случаем системы (\ref{17.1}) отметим работу \cite{50}; см. также примеры в учебном пособии  \cite{26}.

\section{Одно достаточное  условие асимптотической нечувствительности при ослаблении части ограничений}\setcounter{equation}{0}\setcounter{proposition}{0}\setcounter{zam}{0}\setcounter{corollary}{0}\setcounter{definition}{0}
   \ \ \ \ \ В настоящем разделе сосредоточимся на частном, но важном случае построения МП для системы (\ref{17.1}), когда
   \beq{18.1}
   \mathbf{c}\in B_0(E,\ml).
   \eeq

Условие (\ref{18.1}) в случае (\ref{17.1'}) (и при условии, что $\ml$ получается в виде (\ref{16.2''}) при $t_0=0$  и $\theta_0=1$) может отвечать весьма естественным вариантам управляемой динамики. Сейчас отметим следующие два совсем простых примера, в которых реализуется
 только один разрыв функции $\mathbf{c}=\mathbf{c}(\cdot),$ возникающий в момент $t^0\in]0,1[.$

 1)   Предположим, что функция $\mathbf{c}$ получается склеиванием двух постоянных <<кусков>>, отвечающих константам $c_1\bn$ и $c_2\bn:$
 \beq{18.2}
 \left(\mathbf{c}(t)=c_1 \ \ \forall\,t\in\left[0,t^0\right[\,\right)\ \& \ \left(\mathbf{c}(t)=c_2 \ \ \forall\,t\in\left[t^0,1\right]\right).
 \eeq
Будем при этом полагать, что $c_1<c_2.$ Учитывая, что $c_1,c_2$\,--- величины, обратно пропорциональные значениям массы, вариант (\ref{18.2}) можно связать с моделью, характеризующей сброс какой-то части массы (некоторого груза) в момент $t^0.$

2) Рассмотрим вариант, когда $|\mathbf{c}(t)|\equiv 1,$ но в момент $t^0$ изменяется знак функции $\mathbf{c}.$ Так, например, можно рассматривать случай, когда
\beq{18.3}
 \left(\mathbf{c}(t)=1 \ \ \forall\,t\in\left[0,t^0\right[\,\right)\ \& \ \left(\mathbf{c}(t)=-1 \ \ \forall\,t\in\left[t^0,1\right]\right).
 \eeq

В (\ref{18.3}) фактически говорится о реверсе двигателя. Точнее, мы, формируя то или иное управление $f\in M_*^+[\,b\mid E;\ml;\eta],$ определяем режим подачи топлива в двигатель, ориентация которого в момент $t^0$ меняется на противоположную.

Заметим в связи с 1), 2), что мы в своей задаче считаем функцию $\mathbf{c}=\mathbf{c}(\cdot)$  заданной и занимаемся только выбором управления $f.$ На самом же деле возможна ситуация, когда $\mathbf{c}$ формируется какой-то другой системой управления при соблюдении тех или иных ограничений; здесь  такую возможность не рассматриваем.

Итак, обратимся к условию (\ref{18.1}) в общем его виде: $\mathbf{c}=\mathbf{c}(\cdot)$   является далее ступенчатой функцией. Отметим, что функция
$$\tau\mapsto(t-\tau)\,\mathbf{c}(\tau):E\rightarrow\rr$$
ступенчатой уже, как правило, не является (данная функция будет ступенчатой  в не представляющем интереса случае $\mathbf{c}(\tau)\equiv 0$). С учетом (\ref{18.1}) и последнего замечания логично полагать
\beq{18.2}
J=\{2\};
\eeq
имеется в виду обеспечение условия 15.1 (см. также представление (\ref{17.8''})).

Условие 15.1 выполнено (условие 15.2 также полагаем выполненным: $Y$ есть непустое замкнутое множество на плоскости). Тогда $\widetilde{\Pi}^1\left(\widetilde{S}^{\,-1}(Y)\right)$  есть (см. теорему 15.1) МП, универсальное в следующем смысле: данное множество есть МП и в случае, когда ОАХ определяются семейством $\mathfrak{Y,}$ и в случае, когда такие ограничения определяются семейством $\widehat{\mathfrak{Y}}[J],$  где $J$\,--- синглетон   (\ref{18.2}). С учетом (\ref{14.7}), (\ref{14.8}) и теоремы 15.1 получаем совпадение множеств (\ref{12.13}) и (\ref{12.15}). При этом в рассматриваемом плоском случае при $\varepsilon\bn$
 \begin{multline}\label{18.4}
\widehat{\mathbb{O}}_N(Y,\varepsilon|J)=\widehat{\mathbb{O}}_2(Y,\varepsilon|\{2\})=
\bigl\{(z_i)_{i\in\overline{1,2}}\in\rr^2\mid\exists\,(y_i)_{i\in\overline{1,2}}\in Y:\\ (|y_1-z_1|<\varepsilon)\ \& \ (y_2=z_2)\bigl\};
\end{multline}
  полезно сравнить (\ref{17.9}) и (\ref{18.4}). В результате такого сравнения из теоремы 15.1  получаем следующий вывод: с точки зрения <<истинной>> асимптотики множества терминальных состояний, грубо говоря, все равно, ослаблялось ли условие совмещения скоростной координаты МП со второй компонентой подходящего вектора из $Y$ или нет. Данный вывод, формализуемый в терминах совпадения двух МП, допускает также и <<окрестностную>> интерпретацию, подобную \cite[(4.3.19)]{46}.

  Итак, наша задача при истолковании ее результата на языке МП или в терминах <<вилки>>  окрестностей, подобной \cite[(4.3.19)]{46} груба <<в направлении>> скоростной координаты при ослаблении соответствующих условий. Таким образом, для управляемой материальной точки получили естественный пример непосредственной реализации единого МП теоремы 15.1, имеющей достаточно понятный механический смысл.

\section{Пример задачи с <<асимптотическими>> ограничениями}\setcounter{equation}{0}\setcounter{proposition}{0}\setcounter{zam}{0}\setcounter{corollary}{0}\setcounter{definition}{0}
   \ \ \ \ \ Рассмотрим на примере задачи управления материальной точкой другой тип ОАХ (общая постановка задачи и соответствующее решение приведены в \cite{51,52}). Речь пойдет об использовании реализуемых импульсов управления исчезающе малой протяженности. Такой режим управления оказывается полезным в некоторых задачах космической навигации при использовании модели в виде ньютоновой точки в поле тяготения. Здесь мы ограничимся простейшим иллюстративным примером. Речь пойдет об управлении системой  (\ref{17.1}) на единичном промежутке времени.
Сохраняем предположение о том, что конкретное программное управление является элементом $M_*^+[\,b\mid E;\ml;\eta],$  однако дополнительные ограничения изначально имеют асимптотический характер. Полагаем далее для простоты, что $b=1.$ Для краткости введем обозначения, согласующиеся с \cite{51}: пусть
\beq{19.1}
\mathbf{F}\triangleq M_*^+[\,b\mid E;\ml;\eta]=M_*^+[\,1\mid E;\ml;\eta]=\biggl\{f\in B_0^+(E,\ml)\mid\int\limits_Ef\,d\eta=1\biggl\}.
\eeq
С каждым управлением $f\in\mathbf{F}$ связываем два момента времени из $E=[0,1],$ характеризующие его длительность: при условии, что
\beq{19.2}
\pr{supp}(f)\triangleq\{t\in E\mid f(t)\neq 0\}
\eeq
 полагаем, что
\beq{19.3}
\mathbf{t}_0(f)\triangleq\inf(\pr{supp(f)})\in E, \ \ \mathbf{t}^0(f)\triangleq\sup(\pr{supp(f)})\in E.
\eeq
Отметим, что определение (\ref{19.3}) корректно, поскольку (\ref{19.2}) является всякий раз непустым множеством, а, точнее, непустым подмножеством  $E$ (\ref{17.1'}). Разумеется, не исключается тот случай, что $f(t)\neq 0 \ \forall\,t\in E;$ тогда $\mathbf{t}_0(f)=0,\ \mathbf{t}^0(f)=1.$  Но он не представляет интереса. Введем при $\varepsilon\bn$ непустое множество
\beq{19.4}
\mathbf{F}_\varepsilon\triangleq\{f\in\mathbf{F}\mid \mathbf{t}^0(f)-\mathbf{t}_0(f)<\varepsilon\}.
\eeq
Непустота (\ref{19.4}) весьма очевидна: при $\varepsilon\bn$ введем $\overline{\varepsilon}\triangleq\inf\left(\left\{\frac{\varepsilon}{2};1\right\}\right)\bn$ и рассмотрим функцию $$f_\varepsilon=\frac{1}{\overline{\varepsilon}}\,\chi_{[0,\overline{\varepsilon}]}\in\mathbf{F};$$
тогда согласно  (\ref{19.3}) имеем $$\mathbf{t}_0(f_\varepsilon)=0, \ \mathbf{t}^0(f_\varepsilon)=\overline{\varepsilon}\leqslant\frac{\varepsilon}{2}<\varepsilon.$$
Последнее означает, что $f_\varepsilon\in\mathbf{F}_\varepsilon.$ Итак, (\ref{19.4})\,--- непустое подмножество $\mathbf{F}.$ Тогда, как легко заметить,
\beq{19.5}
\mathfrak{F}\triangleq\{\mathbf{F}_\varepsilon:\varepsilon\bn\,\}\in\beta_0[\mathbf{F}],
\eeq
где $\beta_0[\mathbf{F}]\triangleq\{\mathcal{B}\in\beta[\mathbf{F}]\mid\zer\notin\mathcal{B}\}$ есть семейство, элементами которого также являются семейства, именуемые базами фильтров (непустые направленные семейства, состоящие из непустых множеств). Заметим, что $\mathfrak{F}$ можно использовать для введения ОАХ подобно тому, как было сделано в разделе 14. По сути мы следуем соглашениям (\ref{14.1}) при условиях (\ref{17.1'}) и $b=1.$ Однако в качестве семейства подмножеств $\mathbf{F}$ мы теперь используем $\mathfrak{F}$ (\ref{19.5}) (напомним в этой связи (\ref{19.1})). Нас интересует построение МП
\beq{19.6}
(\pr{\mathbf{as}})\left[\mathbf{F};\rr^n;\tau_\rr^{(n)};\Pi;\mathfrak{F}\right]\in\pp(\rr^n),
\eeq
где $n=2,$ как и в предыдущем разделе. Заметим в связи с (\ref{19.6}), что построение соответствующего аналога множества (\ref{12.8}) в нашем случае по сути дела лишено смысла.
В самом деле, если попытаться построить аналог невозмущенной задачи, т.е. задачи, решаемой в классе точных решений --- управлений, понимаемых в духе Дж. Варги, то в качестве множества всех точных (соблюдающих все ограничения точно) управлений следовало бы использовать пересечение всех множеств $\mathbf{F}_\varepsilon,$ $\varepsilon>0.$ Управляющие импульсы $f,$ являющиеся элементами такого пересечения, обязаны (см. (\ref{19.1})--(\ref{19.5})) иметь нулевую протяженность, т.е. удовлетворять условию $$\mathbf{t}_0(f)=\mathbf{t}^0(f,)$$ и, вместе с тем,  $\eta$-интеграл, равный единице (см. (\ref{19.1})). Таких управлений в множестве $\mathbf{F}$ не существует (предлагаем читателю убедиться в этом самостоятельно). Итак, мы должны признать тогда, что обычная область достижимости в момент $t=1,$ есть пустое множество (действительно, в пересечении всех множеств из $\mathfrak{F}$ нет ничего).

С практической же точки зрения интерес представляют множества $\Pi^1(\mathbf{F}_\varepsilon)$ при малых~$\varepsilon,$ $\varepsilon>0,$ а также предел этих множеств при $\varepsilon\downarrow 0.$ В связи с представлением МП (\ref{19.6}) полезно заметить, что в нашем случае выполняются условия, конкретизирующие (\ref{13.5}), (\ref{13.6}).    Отметим (см. (\ref{13.5})), что последовательность
$$k\mapsto\mathbf{F}_{\frac{1}{k}}:\mathbb{N}\rightarrow\mathfrak{F}$$
такова, что $\forall\,\Sigma\in\mathfrak{F} \ \exists\,j\in\mathbb{N}:$ $$\mathbf{F}_{\frac{1}{j}}\subset\Sigma.$$

Условие (\ref{13.6}) (в нужном конкретном варианте) выполнено в силу метризуемости $\left(\rr^n,\tau_\rr^{(n)}\right)=\left(\rr^2,\tau_\rr^{(2)}\right)$ (см. (\ref{19.6})). Поэтому согласно (\ref{13.7})
\beq{19.8}
(\pr{\mathbf{as}})\left[\mathbf{F};\rr^n;\tau_\rr^{(n)};\Pi;\mathfrak{F}\right]=(\pr{\mathbf{sas}})\left[\mathbf{F};\rr^n;\tau_\rr^{(n)};\Pi;\mathfrak{F}\right].
\eeq
Для МП (\ref{19.8}) имеем свойство совпадения с пересечением всех множеств
$$\pr{cl}\left(\Pi^1(\mathbf{F}_\varepsilon),\tau_\rr^{(n)}\right), \ \varepsilon\bn;$$
конкретное представление (\ref{19.8}) установлено в \cite[Теорема 5.1]{51}
(при этом существенно использовалась конструкция на основе (\ref{13.9})). Ограничимся здесь обсуждением данного представления для рассматриваемого плоского случая (при $n=2$).

Если $x\in\rr^2$ и $y\in\rr^2,$ то полагаем $[x;y]_n=[x;y]_2\triangleq\{\alpha x+(1-\alpha)y:\alpha\in[0,1]\}$, получая выпуклую оболочку $\{x;y\}.$  При этом, как отмечено в \cite{53}, функции $\pi_1$ и $\pi_2$ (см. (\ref{17.2}), (\ref{17.3})) обладают следующими важными свойствами: при $j\in\overline{1,2}$ и $t\in\,]0,1]$ определен предел слева функции $\pi_j$ в точке $t,$ обозначаемый через $$\lim\limits_{\theta\uparrow t}\pi_j(\theta);$$
при $j\in\overline{1,2}$ и $\overline{t}\in[0,1[$ определен предел справа функции $\pi_j$  в точке $\overline{t},$ который обозначаем
$$\lim\limits_{\theta\downarrow \overline{t}}\pi_j(\theta).$$
С учетом этого введем вектор-функции
\beq{19.9}
\widehat{\pi}_\uparrow: \, ]0,1]\rightarrow\rr^2, \ \ \widehat{\pi}_\downarrow: \, [0,1[\rightarrow\rr^2,
\eeq
 полагая в дальнейшем, что
 \beq{19.10}
 \left(\widehat{\pi}_\uparrow(t)\triangleq\left(\lim\limits_{\theta\uparrow t}\pi_j(\theta)\right)_{i\in\overline{1,2}} \ \forall\,t\in\,]0,1]\right)\ \& \ \left(\widehat{\pi}_\downarrow(\overline{t})\triangleq\lim\limits_{\theta\downarrow\overline{t}}\pi_j(\theta) \ \forall\,\overline{t}\in\,[0,1[\right)
 \eeq
 В терминах вектор-функций (\ref{19.9}), (\ref{19.10}) конкретизируем сегменты
 $$\left[\widehat{\pi}_\uparrow(\tau);\widehat{\pi}_\downarrow(\tau)\right]_2, \ \tau\in\,]0,1[.$$

 Тогда (см. \cite{51}) множество (\ref{19.8}) совпадает с
 \beq{19.11}
 \left(\bigcup\limits_{t\in]0,1[}\left[\widehat{\pi}_\uparrow(t);\widehat{\pi}_\downarrow(t)\right]_2\right)\cup\left\{\widehat{\pi}_\downarrow(0);
 \widehat{\pi}_\uparrow(1)\right\}.
 \eeq
 Итак, (\ref{19.11}) есть искомое МП, отвечающее применению $\mathfrak{F}$ в качестве ОАХ. В связи с определением односторонних пределов в/з функций $\pi_1$ и $\pi_2$ заметим, что указанные пределы определяются в силу (\ref{17.2}), (\ref{17.3}) аналогичными пределами функции $\mathbf{c}=\mathbf{c}(\cdot),$ которая, как отмечено в \cite{53}, также обладает такими пределами: при $t\in]0,1]$ определен предел слева
 $$\lim\limits_{\theta\uparrow t}\mathbf{c}(\theta)\in\rr, $$
 а при $\overline{t}\in[0,1[$ определен аналогичный предел справа
 $$\lim\limits_{\theta\downarrow\overline{t}}\mathbf{c}(\theta)\in\rr. $$
 При этом с учетом (\ref{17.2}), (\ref{17.3}) получаем, что
 \beq{19.12}
 \lim\limits_{\theta\uparrow t}\pi_1(\theta)=(1-t)\lim\limits_{\theta\uparrow t}\mathbf{c}(\theta) \ \ \forall\,t\in]0,1],
 \eeq
 \beq{19.13}
 \lim\limits_{\theta\uparrow t}\pi_2(\theta)=\lim\limits_{\theta\uparrow t}\mathbf{c}(\theta) \ \ \forall\,t\in]0,1],
 \eeq
 \beq{19.14}
 \lim\limits_{\theta\downarrow\overline{t}}\pi_1(\theta)=(1-\overline{t})\lim\limits_{\theta\downarrow\overline{t}}\mathbf{c}(\theta) \ \ \forall\,\overline{t}\in[0,1[,
 \eeq
 \beq{19.15}
 \lim\limits_{\theta\downarrow \overline{t}}\pi_2(\theta)=\lim\limits_{\theta\downarrow\overline{t}}\mathbf{c}(\theta) \ \ \forall\,t\in[0,1[.
 \eeq
 Итак, располагая конкретной ярусной функцией и используя последние четыре соотношения, мы в виде (\ref{19.11}) получаем МП (\ref{19.8}).

В заключение раздела рассмотрим простейший пример, полагая выполненным (\ref{18.3}). Рассмотрим вариант системы (\ref{17.1}), когда на промежутке управления осуществляется однократный реверс двигательной установки. Пусть  для определенности  $t^0=\frac12.$ Полагаем, что $\mathbf{c}(t)\triangleq 1$ при $t\in[0,t^0[$ и $\mathbf{c}(t)\triangleq -1$ при $t\in[t^0,1].$
Определим прежде всего $\widehat{\pi}_\downarrow(0), \widehat{\pi}_\uparrow(1).$  Используя непрерывность $\mathbf{c}=\mathbf{c}(\cdot)$ в точках $0$ и $1,$ получаем из (\ref{19.12})~--~(\ref{19.15}),  что
$$\lim\limits_{\theta\uparrow 1}\pi_1(\theta)=0, \ \ \lim\limits_{\theta\uparrow 1}\pi_2(\theta)=-1; \ \ \lim\limits_{\theta\downarrow 0}\pi_1(\theta)=1; \ \ \lim\limits_{\theta\downarrow 0}\pi_2(\theta)=1.$$

Тем самым полностью определены два упомянутых вектора $\widehat{\pi}_\uparrow(1);\widehat{\pi}_\downarrow(0),$ которые реализуются в виде УП следующим образом:
\beq{19.16}
(0,-1), \ (1,1);
\eeq
напомним, что согласно  соглашениям раздела 3 каждый вектор из $\rr^2,$ строго говоря, является отображением из $\overline{1,2}$ в $\rr,$ т.е. кортежем <<длины>> 2; однако, понятно, что вариант (\ref{19.16}) является всего лишь несущественной редакцией векторов $\widehat{\pi}_\uparrow(1);\widehat{\pi}_\downarrow(0).$

Далее рассмотрим представление $\widehat{\pi}_\uparrow(t);\widehat{\pi}_\downarrow(t)$ при $t\in]0,\frac12[.$ Согласно (\ref{18.3}) в каждой такой точке $t$ наша функция $\mathbf{c}=\mathbf{c}(\cdot)$ непрерывна и  $\mathbf{c}(t)=1.$ Поэтому
$$\lim\limits_{\theta\uparrow t}\pi_1(\theta)=1-t, \ \ \lim\limits_{\theta\uparrow t}\pi_2(\theta)=1; \ \ \lim\limits_{\theta\downarrow t}\pi_1(\theta)=1-t; \ \ \lim\limits_{\theta\downarrow t}\pi_2(\theta)=1.$$
Тогда для каждого из векторов $\widehat{\pi}_\uparrow(t);\widehat{\pi}_\downarrow(t),$ где $t\in]0,\frac12[$  имеем следующие представления в виде УП:
\beq{19.17}
(1-t,1), \ \ (1-t,1).
\eeq
Заметим, что в терминах (\ref{19.11}) имеем, конечно, равенства
$$\left[\widehat{\pi}_\uparrow(t);\widehat{\pi}_\downarrow(t)\right]_2 =\left\{\widehat{\pi}_\uparrow(t)\right\}= \left\{\widehat{\pi}_\downarrow(t)\right\} \ \forall\,t\in\left]0,\frac12\right[,$$
где элементы, формирующие синглетоны наглядно представляются УП (\ref{19.17}).

Особый интерес представляет множество
\beq{19.18}
\left[\widehat{\pi}_\uparrow\left(\frac12\right);\widehat{\pi}_\downarrow\left(\frac12\right)\right]_2\in\pp'(\rr^2).
\eeq
Для его представления определим конкретный вид каждого из векторов
\beq{19.19}
\widehat{\pi}_\uparrow\left(\frac12\right);\widehat{\pi}_\downarrow\left(\frac12\right).
\eeq
Для этого мы прежде всего заметим, что в силу (\ref{18.3})
$$\lim\limits_{\theta\uparrow \frac12}\mathbf{c}(\theta)=1, \ \ \lim\limits_{\theta\downarrow\frac12}\mathbf{c}(\theta)=-1.$$
Теперь учтем (\ref{19.12})--(\ref{19.15}); получаем, что
$$\lim\limits_{\theta\uparrow \frac12}\pi_1(\theta)=\frac12, \ \ \lim\limits_{\theta\uparrow\frac12}\pi_2(\theta)=1;\ \
\lim\limits_{\theta\downarrow \frac12}\pi_1(\theta)=-\frac12, \ \ \lim\limits_{\theta\downarrow\frac12}\pi_2(\theta)=-1.$$
Таким образом, на языке УП (\ref{19.19}) принимает следующий вид:
\beq{19.20}
\left(\frac 12,1\right), \ \ \left(-\frac12,-1\right).
\eeq
Если теперь выбрать произвольное число $\alpha\in[0,1],$ то вектор $\mathbf{x}_\alpha\triangleq\alpha\widehat{\pi}_\uparrow\left(\frac12\right)+(1-\alpha)\widehat{\pi}_\downarrow\left(\frac12\right)$
имеет, очевидно, следующие компоненты $$\mathbf{x}_\alpha(1)=\frac{\alpha-(1-\alpha)}{2}=\frac{2\alpha-1}{2}=\alpha-\frac12,$$
$$\mathbf{x}_\alpha(2)=\alpha-(1-\alpha)=2\alpha-1.$$
Сегмент $\left[\widehat{\pi}_\uparrow\left(\frac12\right);\widehat{\pi}_\downarrow\left(\frac12\right)\right]_2=\{\mathbf{x}_\alpha:\alpha\in[0,1]\}$ обладает следующим представлением в терминах УП:
\beq{19.21}
\left\{\left(\alpha-\frac12,2\alpha-1\right):\alpha\in[0,1]\right\}
\eeq
Полезно отметить, кстати, что
\beq{19.22}
\mathbf{x}_\alpha(1)|_{\alpha=0}=-\frac12, \ \ \mathbf{x}_\alpha(2)|_{\alpha=0}=-1;
\eeq
\beq{19.23}
\mathbf{x}_\alpha(1)|_{\alpha=1}=\frac12, \ \ \mathbf{x}_\alpha(2)|_{\alpha=1}=1.
\eeq
Для нас сейчас существенно (\ref{19.23}), поскольку с учетом (\ref{19.17}) имеем в виде $\mathbf{x}_\alpha(1)|_{\alpha=1}$ предел значений $\pi_\uparrow(t)=\pi_\downarrow(t)$ при $t\uparrow\frac12$ (предел слева). Выражение  (\ref{19.22}) учтем позже.

Для $t\in]\frac12,1[$ векторы $\widehat{\pi}_\uparrow\left(t\right);\widehat{\pi}_\downarrow\left(t\right)$ допускают следующие представления:
$$\widehat{\pi}_\uparrow(t)\left(1\right)=\lim\limits_{\theta\uparrow t}\pi_1(\theta)=-(1-t)=t-1, \  \ \widehat{\pi}_\uparrow(t)\left(2\right)=\lim\limits_{\theta\uparrow t}\pi_2(\theta)=-1,$$
$$\widehat{\pi}_\downarrow(t)\left(1\right)=\lim\limits_{\theta\downarrow t}\pi_1(\theta)=t-1, \  \ \widehat{\pi}_\downarrow(t)\left(2\right)=\lim\limits_{\theta\downarrow t}\pi_2(\theta)=-1.$$
В результате каждый из таких векторов  $\widehat{\pi}_\uparrow\left(t\right), \ \widehat{\pi}_\downarrow\left(t\right)$ имеет на языке УП следующее представление:
\beq{19.24}
(t-1,-1).
\eeq
Отметим  справедливость следующего свойства (см. (\ref{19.24})): $\mathbf{x}_\alpha(1)|_{\alpha=0}$ есть предел значений $$\pi_\uparrow\left(t\right)=\pi_\downarrow\left(t\right)\ \text{при} \  t\downarrow\frac{1}{2}.$$

Ясно, что для $t\in\left]\frac{1}{2},1\right[$ имеет место равенство
\beq{19.25}
\left[\,\pi_\uparrow\left(t\right);\pi_\downarrow\left(t\right)\right]_2=\left\{\pi_\uparrow\left(t\right)\right\}
=\left\{\pi_\downarrow\left(t\right)\right\}.
\eeq
Отметим, наконец, что $\widehat{\pi}_\downarrow\left(0\right)$ есть предел значений $$\widehat{\pi}_\downarrow\left(t\right)=\widehat{\pi}_\uparrow\left(t\right)\ \text{при} \  t\downarrow 0,$$ а $\widehat{\pi}_\uparrow\left(1\right)$\,--- предел значений  $$\widehat{\pi}_\downarrow\left(t\right)=\widehat{\pi}_\uparrow\left(t\right) \text{при} \ t\uparrow 1;$$  см. (\ref{19.16}), (\ref{19.17}), (\ref{19.24}).

Для построения требуемого МП осталось воспользоваться представлением (\ref{19.11}), в котором теперь конкретизированы все его компоненты. Приводим ниже соответствующий рисунок и предоставляем читателю убедиться в том, что полученный <<зигзаг>> и есть искомое~МП.
\begin{center}
\includegraphics[width=85mm]{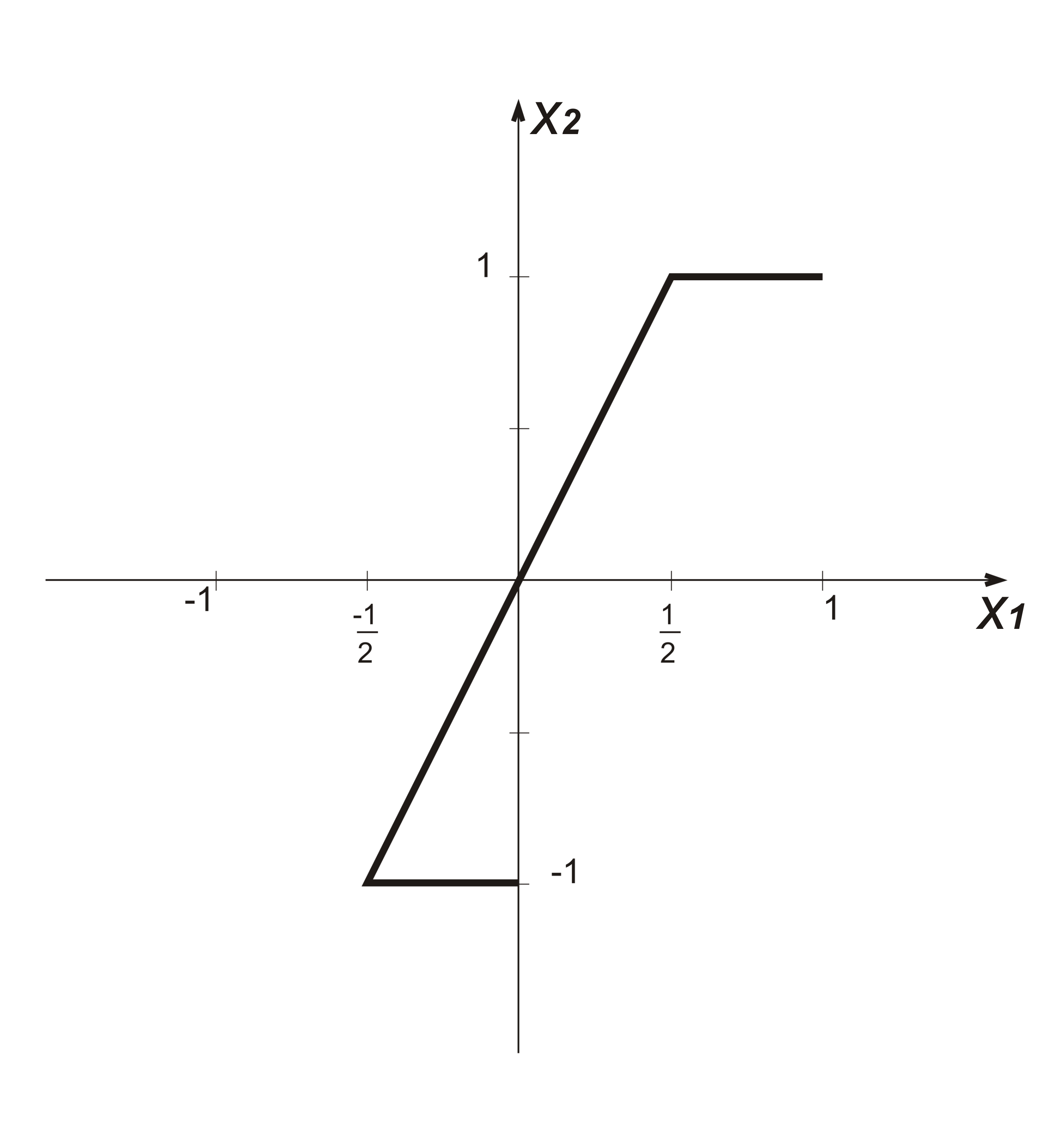}
\end{center}

\newpage

\section{Заключение}\setcounter{equation}{0}\setcounter{proposition}{0}\setcounter{zam}{0}\setcounter{corollary}{0}\setcounter{definition}{0}
   \ \ \ \ \ В пособии рассмотрена одна естественная задача о достижимости при наличии ограничений асимптотического характера. В наиболее понятном варианте, описанном в последней главе, ее можно связывать с построением асимптотического аналога области достижимости линейной импульсно управляемой системы с разрывностью в коэффициентах при управляющих воздействиях. Известный эффект произведения разрывной функции на обобщенную в данном конкретном случае получает свое исчерпывающее описание в терминах конечно-аддитивных мер со свойством слабой абсолютной непрерывности относительно соответствующего сужения меры Лебега. В настоящем пособии рассмотрен только один тип ограничений импульсного характера, отвечающий (в задачах управления)  директивному требованию полного расходования имеющегося энергоресурса. Однако в имеющейся к настоящему моменту обширной серии публикаций (отметим здесь только монографии \cite{25},\cite{46},\cite{47},\cite{41}) указаны и другие варианты упомянутых ограничений, причем конструкция, излагаемая в данном пособии, при минимальной корректировке позволит заинтересованному читателю разобраться и в тех постановках, которые здесь не рассматривались.

   Полезно отметить и естественные аналогии построений последней главы со ставшими уже классическими конструкциями расширений нелинейных задач управления с геометрическими ограничениями (см. \cite{44}, \cite{22}). В обоих случаях основные усилия по построению решения соответствующей задачи с ограничениями, соблюдаемыми с высокой, но всё же конечной степенью точности, <<перекладываются>> на обобщенную задачу, а, точнее, на задачу о достижимости в классе обобщенных управлений (последние в условиях геометрических ограничений порождают так называемые скользящие режимы).

   Особенности, связанные с применением конечно-аддитивных мер в качестве обобщенных элементов потребовали, однако, использование целого ряда нетрадиционных (для расширений задач теории управления) конструкций. Так, в частности, основные топологии, естественные с идейной точки зрения для данного подхода, не только не являются метризуемыми, но и не удовлетворяют первой аксиоме счетности. Это касается, в частности, применяемого в пособии варианта $*$-слабой топологии и соответствующей индуцированной топологии. Мы, однако, получаем при этом компакт, что является важным для работоспособности схемы расширения. Итак,  мы имеем неметризуемый компакт, обладающий целым рядом полезных свойств. При исследовании вопроса, связанного с асимптотической нечувствительности исследуемой задачи при ослаблении части ограничений, потребовалась нульмерная топология $\tau_0(\ml),$ т.е. хаусдорфова топология, обладающая базой открыто-замкнутых множеств. Упомянутые обстоятельства, на наш взгляд, также играют полезную роль в изложении основных конструкций, поскольку при их рассмотрении заинтересованный читатель сможет повторить (или получить заново) весьма важные положения общей топологии, уже не связанные с метризацией и, тем не менее, удачно вписывающиеся в исходную содержательную задачу о достижимости. В этой связи отметим также более общие постановки абстрактных задач о достижимости с ограничениями асмптотического характера, рассматриваемые в \cite{46,47} (см., например, \cite[гл.\,3, 4]{46}) и знакомство с которыми может оказаться полезным для заинтересованного читателя.

 \vspace{5cm}

 \noindent Ченцов Александр Георгиевич, д.ф.-м.н., член-корреспондент РАН, главный научный сотрудник, отдел управляемых систем, Институт математики и механики им. Н.\,Н.~Красовского УрО РАН

  \vspace{0,5cm}

 \noindent chentsov@imm.uran.ru

  \vspace{2cm}

 \noindent Шапарь Юлия Викторовна, к.ф.-м.н., доцент, кафедра вычислительных методов и уравнений математической физики, Уральский федеральный университет

  \vspace{0,5cm}

 \noindent shaparuv@mail.ru

  \vspace{2cm}

\noindent г. Екатеринбург, 2016


\newpage

\begin{thebibliography}{40}
\bibitem{44}{\it Варга Дж.} Оптимальное управление дифференциальными и функциональными уравнениями.  М.: Наука. 1977, 624~с.
\bibitem{20}{\it Красовский~Н.\,Н., Субботин~А.\,И.} Позиционные дифференциальные игры. М.: Наука. 1974, 456~с.
\bibitem{21} {\it Понтрягин~Л.\,С., Болтянский~В.\,Г., Гамкрелидзе~Р.\,В., Мищенко~Е.\,Ф.} Математическая теория оптимальных процессов.  М.: Наука. 1961,  384 с.
\bibitem{22} {\it Гамкрелидзе~Р.\,В.} Основы оптимального управления. Тбилиси: Изд-во Тбилисского университета. 1977,  253~с.
 \bibitem{23} {\it Красовский~Н.\,Н.}  Теория управления движением.  М.: Наука. 1968, 475~с.
 \bibitem{24} { \it Завалищин~С.\,Т, Сесекин~А.\,Н.}  Импульсные процессы. Модели и приложения. М.: Наука. 1991, 256~с.
   \bibitem{25}{\it Chentsov~A.\,G.} Finitely additive measures and relaxations of extremal problems.  New York, London and Moscow:
 Plenum Publishing Corporation. 1996, 244~p.
\bibitem{46} {\it Chentsov~A.\,G.} Asymptotic attainability. Dordrecht-Boston-London: Kluwer Academic Publishers. 1997,  322~p.
\bibitem{47} {\it Chentsov~A.\,G.}  Finitely additive measures and relaxations of abstract control  problems.  Journal of Mathematical Sciences.  2006,  vol.133, \No~2, p.~1045--1206.
\bibitem{43} {\it Данфорд~H., Шварц~Дж.\,Т.}  Линейные операторы. Общая теория.  М.: Издательство иностранной литературы. 1962,  895~с.
\bibitem{26} {\it Бродская~Л.\,И., Ченцов~А.\,Г.} Некоторые примеры неустойчивых задач управления. Учебное пособие. Екатеринбург. Изд-во Уральского университета. 2014, 101~с.
\bibitem{41} {\it Chentsov~A.\,G., Morina~S.\,I.} Extensions and Relaxation. Dordrecht-Boston-London: Kluwer Academic Publishers. 2002,  408~p.
\bibitem{30} {\it Ченцов~А.\,Г.} Элементы конечно-аддитивной теории меры, I.  Екатеринбург: РИО УГТУ-УПИ. 2008, 388 с.
\bibitem{31}{\it Ченцов А.Г.} Множества, события, вероятность (основные структуры). Учебное пособие. УГТУ-УПИ. Екатеринбург. 2006, 199~с.
%\bibitem{31} {\it Ченцов~А.\,Г.}  Множества. События. Вероятность.
\bibitem{32} {\it Ченцов~А.\,Г.} Элементы теории множеств. Учебное пособие. Екатеринбург. УГТУ-УПИ. 2009, 55~с.
\bibitem{33} {\it Ченцов~А.\,Г., Шукшина~Н.\,В.} Теория множеств: простейшие конструкции.  Екатеринбург. УрФУ. 2012,  130~с.
\bibitem{34} {\it Куратовский~К., Мостовский~А.}  Теория множеств.  М.: Мир. 1970,  416~с.
\bibitem{35} {\it Кормен~Т., Лейзерсон~Ч., Ривест~Р., Штайн~K.}  Алгоритмы: построение и анализ. М.: ИД <<Вильямc>>. 2005, 1296~с.
\bibitem{36}{\it Ченцов~А.\,Г.}  Конечно-аддитивные меры и релаксации экстремальных задач. Екатеринбург: Наука. 1993, 229~с.
\bibitem{37}{\it Булинский~А.\,В., Ширяев~А.\,Н.} Теория случайных процессов. М.:~ ФИЗМАТЛИТ. 2005,  408~с.
\bibitem{38} {\it Ченцов~А.\,Г.} Элементы конечно-аддитивной теории меры, II.  Екатеринбург. 2010, 541~с.
\bibitem{39} {\it Энгелькинг~Р.} Общая топология.  М.: Мир. 1986, 751~с.
\bibitem{40} {\it Келли~Дж.\,Л.}  Общая топология.  М: Наука. 1981, 431~с.
\bibitem{45} {\it Бурбаки~Н.}  Общая топология. Основные структуры.  М.: Наука. 1968,  272~с.
 \bibitem{48} {\it Ченцов~А.\,Г.} Множества притяжения в абстрактных задачах о достижимости: эквивалентные представления и основные свойства. Известия вузов. Математика. 2013, \No~11, с.~33--50.
 \bibitem{49} {\it Ченцов~А.\,Г.} Фильтры и ультрафильтры в конструкциях множеств притяжения. Вестник Удмуртского университета. Математика. Механика. Компьютерные науки.  2011, вып. 1, с.~113--142.
 \bibitem{50} {\it Ченцов~А.\,Г.} Ультрафильтры в конструкциях множеств притяжения: задача соблюдения ограничений асимптотического характера. Дифференциальные уравнения.  2011, т. 47, \No~ 7, с.~1047--1064
 \bibitem{51} {\it Ченцов~А.\,Г., Бакланов~А.\,П.}  К вопросу о построении множества достижимости при огранчениях асимптотического характера.  Труды Института математики и механики. 2014, т. 20, \No~3, с.~309--323
 \bibitem{52} {\it Ченцов~А.\,Г., Бакланов~А.\,П.} Об одной задаче, связанной с асимптотической достижимостью в среднем. Доклады Академии Наук. 2014, т. 459, \No~6, с.~672--676
 \bibitem{53} {\it Ченцов~А.\,Г.}  Об одном примере представления пространства ультрафильтров алгебры множеств. Труды Института математики и механики.  2011, т. 27, \No~4, с.~293--311

 \end{thebibliography}
\end{document}